\documentclass[11pt]{amsart}

\usepackage{latexsym}
\usepackage{amscd}
\usepackage[all]{xy}
\usepackage[mathscr]{euscript}
\usepackage{dsfont}

\normalsize

\RequirePackage{xspace}
\newcommand{\thought}[1]{}
\renewcommand{\thought}[1]{ \textbf{[#1]}}

\usepackage{enumerate}        % Fancy enumerations;  used for Roman numbering.
\newenvironment{roenumerate}{\begin{enumerate}[\upshape (i)]}{\end{enumerate}}
\newcommand\nc {\newcommand}
\newcommand\rnc{\renewcommand}

\newcount\blopone
\newcount\xone
\newcount\xtwo
\newcount\ytwo

\newtheorem{theorem}{Theorem}[section]
\newtheorem{prop}[theorem]{Proposition}
\newtheorem{observation}[theorem]{Observation}
\newtheorem{com}[theorem]{Comment}
\newtheorem{elaboration}[theorem]{Elaboration}
\newtheorem{redu}[theorem]{Reduction}
\newtheorem{refinement}[theorem]{Refinement}
\newtheorem{summary}[theorem]{Summary}
\newtheorem{importnota}[theorem]{Important Notation}
\newtheorem{prblm}[theorem]{Problem}
\newtheorem{notation}[theorem]{Notation}
\newtheorem{defin}[theorem]{Definition}
\newtheorem{caution}[theorem]{Caution}
\newtheorem{remark}[theorem]{Remark}
\newtheorem{reminder}[theorem]{Reminder}
\newtheorem{illustration}[theorem]{Illustration}
\newtheorem{lemma}[theorem]{Lemma}
\newtheorem{convention}[theorem]{Convention}
\newtheorem{construction}[theorem]{Construction}
\newtheorem{corollary}[theorem]{Corollary}
\newtheorem{example}[theorem]{Example}
\newtheorem{conclusion}[theorem]{Conclusion}
\newtheorem{triviality}[theorem]{Triviality}
\newtheorem{proto}[theorem]{Prototype Quasifibration}
\newtheorem{cauex}[theorem]{Cautionary Example}
\newtheorem{hypo}[theorem]{Hypothesis}
\newtheorem{subth}{ }[theorem]
\newtheorem{case}{Case}[theorem]
\newtheorem{ssubth}{ }[subth]
\newtheorem{facts}[theorem]{Facts}

\nc\tri[1]{\begin{triviality}
\label{#1}}
\nc\fac[1]{\begin{facts}
\label{#1}
\begin{em}}
\nc\cas[1]{\begin{case}
\label{#1}
\begin{em}}
\nc\cvn[1]{\begin{convention}
\label{#1}
\begin{em}}
\nc\rfn[1]{\begin{refinement}
\label{#1}}
\nc\prt[1]{\begin{proto}
\label{#1}}
\nc\lem[1]{\begin{lemma}
\label{#1}}
\nc\pro[1]{\begin{prop}
\label{#1}}
\nc\thm[1]{\begin{theorem}
\label{#1}}
\nc\obs[1]{\begin{observation}
\label{#1}}
\nc\cor[1]{\begin{corollary}
\label{#1}}
\nc\dfn[1]{\begin{defin}
\label{#1}}
\nc\sthm[1]{\begin{subth}
\label{#1}}
\nc\exm[1]{\begin{example}
\label{#1}
\begin{em}}
\nc\plm[1]{\begin{prblm}
\label{#1}
\begin{em}}
\nc\rmk[1]{\begin{remark}
\label{#1}
\begin{em}}
\nc\rmd[1]{\begin{reminder}
\label{#1}
\begin{em}}
\nc\ntn[1]{\begin{notation}
\label{#1}
\begin{em}}
\nc\smr[1]{\begin{summary}
\label{#1}
\begin{em}}
\nc\cau[1]{\begin{caution}
\label{#1}
\begin{em}}
\nc\elb[1]{\begin{elaboration}
\label{#1}
\begin{em}}
\nc\hyp[1]{\begin{hypo}
\label{#1}}
\nc\imn[1]{\begin{importnota}
\label{#1}
\begin{em}}
\nc\rdn[1]{\begin{redu}
\label{#1}
\begin{em}}
\nc\cax[1]{\begin{cauex}
\label{#1}
\begin{em}}
\nc\cmt[1]{\begin{com}
\label{#1}
\begin{em}}
\nc\con[1]{\begin{construction}
\label{#1}
\begin{em}}
\nc\ill[1]{\begin{illustration}
\label{#1}
\begin{em}}
\nc\ssthm[1]{\begin{ssubth}
\label{#1}
\begin{em}}
\nc\cnc[1]{\begin{conclusion}
\label{#1}
\begin{em}}

\nc\elem{\end{lemma}}
\nc\erdn{\end{em}\end{redu}}
\nc\erfn{\end{refinement}}
\nc\eprt{\end{proto}}
\nc\ethm{\end{theorem}}
\nc\eobs{\end{observation}}
\nc\ecor{\end{corollary}}
\nc\edfn{\end{defin}}
\nc\esthm{\end{subth}}
\nc\epro{\end{prop}}
\nc\etri{\end{triviality}}
\nc\eexm{\end{em}
\end{example}}
\nc\ecvn{\end{em}
\end{convention}}
\nc\ecmt{\end{em}
\end{com}}
\nc\efac{\end{em}
\end{facts}}
\nc\ermk{\end{em}
\end{remark}}
\nc\ermd{\end{em}
\end{reminder}}
\nc\eill{\end{em}
\end{illustration}}
\nc\eplm{\end{em}
\end{prblm}}
\nc\ecas{\end{em}
\end{case}}
\nc\ecau{\end{em}
\end{caution}}
\nc\eelb{\end{em}
\end{elaboration}}
\nc\ecax{\end{em}
\end{cauex}}
\nc\eimn{\end{em}
\end{importnota}}
\nc\entn{\end{em}
\end{notation}}
\nc\econ{\end{em}
\end{construction}}
\nc\esmr{\end{em}
\end{summary}}
\nc\ehyp{
\end{hypo}}
\nc\ecnc{\end{em}
\end{conclusion}}
\nc\essthm{\end{em}
\end{ssubth}}

\nc\sst{\scriptstyle}
\newcommand{\comment}[1]{}
\newcommand{\ri}{\longrightarrow}

\newcommand{\zz}{{\mathbb Z}}
\nc\gm{{{\mathbb G}_m}}

\newcommand{\D}{{\mathbf D}}
\newcommand{\Dqc}{{\mathbf D_{\text{\bf qc}}}}
\newcommand{\qq}{{\mathbb Q}}

\newcommand{\pp}{{\mathbb P}}

\newcommand{\oo}{\otimes}

\newcommand{\ak}{{\mathbb A}}

\nc\op{^{\hbox{\rm\tiny op}}}
\nc\mth{^{\hbox{\rm\tiny th}}}

\nc\dcoh{\D_{\mathbf{coh}}^b}
\nc\script{\mathscr}
\nc\z{\zeta}
\nc\bc{{\mathbb{BC}}}
\nc\ct{{\script T}}
\nc\cf{{\script F}}
\nc\cl{{\script L}}
\nc\cv{{\script V}}
\nc\cq{{\script Q}}
\nc\cu{{\script U}}
\nc\ce{{\script E}}
\nc\cg{{\script G}}
\nc\ch{{\script H}}
\nc\cs{{\script S}}
\nc\car{{\script R}}
\nc\cd{{\script D}}
\nc\cc{{\script C}}
\nc\ck{{\script K}}
\nc\ca{{\script A}}
\nc\ci{{\script I}}
\nc\cj{{\script J}}
\nc\co{{\script O}}
\nc\cm{{\script M}}
\nc\cw{{\script W}}
\nc\cx{{\script X}}
\nc\cy{{\script Y}}
\nc\cz{{\script Z}}
\nc\bd{\begin{description}}
\nc\ed{\end{description}}
\nc\ctob{{\script C}at\big(\ci^{op},\ca\big)}
\nc\clim{{\ds\mathop{\rm lim}_{\ds\longleftarrow}}\,}
\nc\climi{\clim^{\!i}\,}
\nc\climn{\clim^{\!n}\,}
\nc\colim{{\ds\mathop{\rm colim}_{\ds\la}}}
\nc\oa{\overline{\ca}}
\nc\s{\sigma}
\nc\ta{\tau}
\nc\os{\overline\sigma}
\nc\ot{\overline\tau}
\nc\T{\Sigma}
\nc\Tm{\Sigma^{-1}}
\nc\de[1]{{\mathop{\rm deg(#1)}}}
\nc\Ad[1]{\mathop{\rm Ad}(#1)}
\nc\ad[1]{\mathop{\rm ad}(#1)}
\nc\kth{{\it K}--theory}
\nc\loc[1]{{\text{\rm Loc(#1)}}}
\nc\coloc[1]{{\text{\rm Coloc}(#1)}}

\nc\one{\mathds{1}}

\def\der #1 {D\left(#1\right)}
\nc\prf{\begin{proof}}
\nc\eprf{\end{proof}}
\nc\ds{\displaystyle}
\nc\Tor{\text{\rm Tor}}

\nc\cb{{\script B}}
\nc\ab{{\script A}b}

\nc\be{\begin{roenumerate}}
\nc\ee{\end{roenumerate}}

\nc\cat[1]{{\script C}at\Big({\big\{#1\big\}}\op\,\,,\,\,\ab\Big)}
\nc\csab{{\script C}at\big(\cs^{op},\ab\big)}
\nc\ctab{{\script C}at\Big({\{\ct^\alpha\}}^{op},\ab\Big)}
\nc\csex{{\script E}x\big(\cs^{op},\ab\big)}
\nc\ctex{{\script E}x\Big({\{\ct^\alpha\}}^{op},\ab\Big)}
\nc\sub{\qquad\subset\qquad}
\nc\ctr[1]{{\left.\ct\left(-,#1\right)\right|}_{\cs}}
\nc\ctrf[2]{{\left.\ct\left(#1,#2\right)\right|}_{\cs}}
\nc\Ctr[1]{{\left.\ct\left(-,#1\right)\right|}_{\ct^\alpha}}
\nc\Ctrf[2]{{\left.\ct\left(#1,#2\right)\right|}_{\ct^\alpha}}

\nc\la{\longrightarrow}
\nc\nin{\noindent}
\nc\cad[1]{\text{card}(#1)}
\nc\eq{\quad=\quad}
\nc\BA{\begin{array}{c}}
\nc\EA{\end{array}}
\nc\barr{
\[
\begin{array}{cccccccccccccccc}
}
\nc\earr{
\end{array}
\]
}
\nc\as[1]{{\langle S\rangle}^{#1}}
\nc\sh{\text{\it shift}}

\nc\yy[1]{{\left.\ct\left(-,#1\right)\right|}_{\ct^c}}
\nc\vrep[2]{{\left.\ct\left(#1,#2\right)\right|}_{\ct^\alpha}}
\nc\da{\downarrow}
\nc\Hom{{\mathop{\rm Hom}}}
\nc\RHom{{\mathop{\rm RHom}}}
\nc\HHom{{\script H}{\mathop{\rm om}}}
\nc\RHHom{{\script{RH}}{\mathop{\rm om}}}
\nc\End{{\mathop{\rm End}}}
\nc\Ext{{\mathop{\rm Ext}}}
\nc\mr\Modtc
\nc\PExt{{\mathop{\rm PExt}}}
\nc\stm{\text{\rm stmod}(kG)}
\nc\stM{\text{\rm StMod}(kG)}
\nc\e{\varepsilon}
\nc\p{\mathfrak{p}}

\nc\rs{\s^{-1}A}
\nc\br{{\{\s^{-1}A\}}}
\nc\ra\ri
\nc\y[1]{\mathbf{y}#1}
\nc\x[1]{\mathbf{z}#1}
\nc\mmod[1]{#1\text{--\rm mod}}
\nc\Mod[1]{#1\text{--\rm Mod}}
\nc\Md {\ensuremath{\mathop{\textup{Mod}}}}
\rnc\mod[1]{\ensuremath{\mathop{\textup{mod-}#1}}\xspace}
\nc\Modtc{\Mod{\ct^c}}
\nc\pgldim[1]{\mathop{\rm pgldim}\,#1}
\nc\tf{{\rm [TR5]}}
\nc\tfs{{\rm [TR5$^*$]}}
\nc\Fun{\text{\rm Funct}(F\op,\ab)}
\nc\sym{\text{\rm Sym}}
\nc\sgn{\text{\rm sgn}}
\nc\Pro{\text{\rm Prod}^{}_\alpha(F\op,\ab)}
\nc\Yt[1]{{\left.\Hom_\ct^{}\left(-,#1\right)\right|}_F^{}}
\nc\dl{\delta}
\nc\Proj[1]{#1\text{--\rm Proj}}
\nc\proj[1]{#1\text{--\rm proj}}
\nc\Flat[1]{\text{\rm Flat}\,#1}
\nc\Inj[1]{\text{\rm Inj\,}#1}
\nc\qc[1]{\text{\rm qc\,}#1}
\nc\ov{\overline}
\nc\wt{\widetilde}
\nc\wh{\widehat}
\nc\ph{\varphi}
\nc\tstr{{\it t}--structure}
\nc\spec[1]{{\text{\rm Spec}(#1)}}
\nc\Spec[1]{{\text{\rm Spec}\big(#1\big)}}

\newcommand{\Dqcpl}{{\mathbf D^+_{\mathrm{qc}}}}
\newcommand{\Dqcb}{{\mathbf D^b_{\mathrm{qc}}}}
\nc\EProd{\text{\rm EProd}}
\nc\ECoprod{\text{\rm ECoprod}}
\nc\Prod{\text{\rm Prod}}
\nc\Coprod{\text{\rm Coprod}}
\nc\ldimp{\text{\rm LDim}^{\prod}}
\nc\ldimc{\text{\rm LDim}^{\coprod}}

\nc\hoco{
\begin{picture}(40,10)
\put(20,0){\makebox(0,0)[b]{\text{\rm Hocolim}}}
\put(5,-2){\vector(1,0){30}}
\end{picture}\,}

\nc\holim{
\begin{picture}(40,10)
\put(20,0){\makebox(0,0)[b]{\text{\rm Holim}}}
\put(35,-2){\vector(-1,0){30}}
\end{picture}}

\nc\Cop{\text{\rm Coprod}}
\nc\seq{{\mathbb{S}_{\mathbf{e}}}}
\nc\nseq{{\mathbb{NS}_{\mathbf{e}}}}
\nc\lnseq{{\mathbb{LNS}_{\mathbf{e}}}}
\nc\fnseq{{\mathbb{FNS}_{\mathbf{e}}}}
\nc\rnseq{{\mathbb{RNS}_{\mathbf{e}}}}
\nc\gnseq{{\mathbf{G}\mathbb{NS}_{\mathbf{e}}}}
\nc\se{{S^{\text{\tt e}}}}
\nc\te{{T^{\text{\tt e}}}}
\nc\LL{{\text{\bf L}}}
\nc\R{\text{\bf R}}
\nc\id{\text{\rm id}}
\nc\supp{\text{\rm supp}}
\nc\Loc{\text{\rm Loc}}
\nc\Thick{\text{\rm Thick}}
\nc\Tri{\mathbb{T}\mathrm{ri}}
\nc\sq{\mathbf{Sq}}
\nc\nsq{\mathbb{N}\mathbf{Sq}}
\nc\rsq{\mathbb{R}\mathbf{Sq}}
\nc\hseq{{\mathbb{HS}_{\mathbf{e}}}}
\nc\vseq{{\mathbb{VS}_{\mathbf{e}}}}
\nc\fkr{\mathfrak{R}}
\nc\fks{\mathfrak{S}}
\nc\fkm{\mathfrak{M}}

\begin{document}

\author{Amnon Neeman}\thanks{The research was partly supported 
  by the Australian Research Council. Part of the work was carried out
  at the Max Planck Institute for Mathematics in Bonn, during the
 2016 Program on Higher Structures in Geometry and Physics.}
\address{Centre for Mathematics and its Applications \\
        Mathematical Sciences Institute\\
        Building 145\\
        The Australian National University\\
        Canberra, ACT 2601\\
        AUSTRALIA}
\email{Amnon.Neeman@anu.edu.au}

\title{An improvement on the base-change theorem and the functor $f^!$}

\begin{abstract}
In 2009 Lipman published his polished book~\cite{Lipman09}, 
the product of a decade's work, giving 
the definitive, state-of-the-art treatment of Grothendieck duality. In this article
we achieve a sharp improvement: we begin by giving a  new 
proof of the base change theorem, which can handle unbounded
complexes and work
in the generality of algebraic stacks
(subject to mild technical restrictions).
This means that our base-change theorem must be subtle and
delicate, the unbounded version
is right at the boundary of known counterexamples---counterexamples
(in the world
of schemes) that had led the experts to believe that major parts of
the theory could only be developed in the bounded-below derived
category.

Having proved our new base-change theorem,
we then use it to define the functor $f^!$ on 
the unbounded derived category and establish its functoriality
properties. In Section~\ref{S0} we will use this
to clarify the relation among all
the various constructions of Grothendieck duality.
One illustration of the power of the new methods is that
we can improve Lipman~\cite[Theorem~4.9.4]{Lipman09} to handle
complexes that are not necessarily bounded. There are also
applications to the theory developed by Avramov, Iyengar, Lipman
and Nayak on the connection between Grothendieck duality and 
Hochschild homology and cohomology but, to keep this paper from
becoming even longer, these are being relegated to separate
articles. See for example~\cite{NeemanTIFR}.
\end{abstract}



\maketitle

\tableofcontents

\setcounter{section}{-1}

\section{Executive summary}
\label{S-1}

Let $f:X\la Y$ be a separated morphism of noetherian schemes,
essentially of finite type. There is
a pushforward map $\R f_*:\Dqc(X)\la \Dqc(Y)$, which has a left
adjoint $\LL f^*:\Dqc(Y)\la\Dqc(X)$ 
and a right
adjoint $f^\times:\Dqc(Y)\la\Dqc(X)$. In duality theory it is
customary to consider yet another functor $f^!$; in all treatments
up to now it has been viewed as a functor
$f^!:\Dqcpl(Y)\la\Dqcpl(X)$. 
Here $\Dqcpl$ means the (cohomologically) bounded below derived category.
Duality theory studies the interplay among these functors.

Grothendieck duality is an old subject. There are basically two approaches 
to it: one can develop the theory
via residual and Cousin complexes as in Hartshorne~\cite{Hartshorne66} and 
Conrad~\cite{Conrad07}, or one can proceed more functorially, as outlined
by Deligne~\cite{Deligne66} and 
Verdier~\cite{Verdier68} and developed very fully in
Lipman~\cite{Lipman09}. 
The key to the second approach
is the base-change theorem of Verdier,  it is fundamental to
developing the properties of the functor $f^!$. There is also a proof 
of the base-change theorem due to
Hartshorne, but it comes at the end, after the theory has been 
substantially set up.

The first result in this paper is a new proof of the base-change theorem,
valid much more generally than the old proofs.
Unlike the old proofs it is based on the compact generation of the 
triangulated categories involved, see 
Lemmas~\ref{L4.5} and \ref{L4.7}. Let us give a version of
Lemma~\ref{L4.7} here---this isn't the most general statement we can
prove but gives the flavor:

\lem{L-1.1}
Suppose we are given a 2-cartesian square of quasi-compact, quasi-separated
algebraic stacks
\[
\CD
 W @>u>> X\\
@VfVV @VVgV \\
Y @>v>> Z
\endCD
\]
where $X$ and $Z$ have quasi-finite, separated diagonals. 
Suppose $v$ is flat, while the morphism $g$ is concentrated
as in \cite[Definition~2.4]{Hall-Rydh13} and 
quasi-proper.\footnote{If we assume the stacks $X$ and $Z$ noetherian, then
$g:X\la Z$  is
quasi-proper as long as $\R g_*$ takes bounded-above complexes
  with coherent cohomology
  to bounded-above complexes with coherent cohomology. This happens
  if $g$ is proper and representable,
  but there are also non-representable examples.}
Let $E$ be an object $E\in\Dqc(Z)$. 
Suppose further that one of the two hypotheses holds:
\be
\item
The map $f$
is of finite Tor-dimension, and the map $g$ is pseudo-coherent.\footnote{If the
stacks $X$ and $Z$ are noetherian and $g$ is of finite type then $g$ is pseudo-coherent.} 
\item
The object $E$ belongs to $\Dqcpl(Z)\subset\Dqc(Z)$.
\ee
Then the base-change map $\Phi:u^*g^\times E\la f^\times v^* E$
is an isomorphism.
\elem

The old theorem of Hartshorne and Verdier is
the special case of Lemma~\ref{L-1.1}(ii)
in the case where all the stacks are noetherian schemes,
and the modern proof given here is not hard--- see Section~\ref{S4}. 
The new Lemma~\ref{L-1.1}(i)
is more delicate; the 
proof involves an application of a generalized
Thomason's localization theorem whose original version, for schemes,
allows us to extend perfect complexes from open subsets of a 
scheme. The reader might wish to compare our proof to the existing ones
in the literature---the proof we give of Lemma~\ref{L-1.1}(ii)
is similar in spirit to the one in Lipman~\cite{Lipman09},
but the twist required
to prove the new Lemma~\ref{L-1.1}(i)
is a little subtle. And the new result is close to best possible, 
\cite[Example~6.5]{Neeman96} shows that unless we impose some
condition the base-change map
$\Phi$ need not be an isomorphism in the unbounded derived category.

We have a new proof of a strengthening of an old theorem.
We can say something for algebraic stacks, but for this
executive summary
let us focus on what the new result says for
noetherian schemes.
The treatments of Grothendieck duality
to date have used the old base-change theorem, with the result 
that many theorems were known
only for bounded-below complexes. 
Lemma~\ref{L-1.1}(i) will allow us to fix this; we will circumvent the
problems and define $f^!$ unconditionally, on the unbounded derived category.
The majority of the paper is devoted to developing the functor $f^!$ and
its functoriality properties. Given how much space it takes to 
set up the theory properly, in this article we confined ourselves
to only one application and the sketch of a second:
\be
\item
We extend Lipman~\cite[Theorem~4.9.4]{Lipman09}, one
of the major results in the book. The theorem has several ramifications, which
Lipman explores in his book but we don't have the space for here.
Nevertheless it is possible to 
reduce it to a simple statement we can include:
suppose $f:X\la Y$ is a proper morphism of noetherian schemes. 
Suppose the complex $E\in\Dqc(Y)$ is bounded above. Then $f^\times E$ does
not have to be bounded above. Let $u:U\la X$ be a flat morphism so that
$fu$ is of finite Tor-dimension. Is $u^*f^\times E$ bounded above?
Although 
Lipman~\cite{Lipman09} considers the question and its 
implications, the techniques available at the
time were able to prove this only when $E$ is bounded. 
We will extend this to any bounded-above complex.
\item
The theory developed here immediately applies to the results of 
Avramov, Iyengar, Lipman and 
Nayak~\cite{Avramov-Iyengar08,
  Avramov-Iyengar-Lipman-Nayak10}, see also
\cite{Iyengar-Lipman-Neeman13}; the relation they found, between Grothendieck
duality and Hochschild homology and cohomology, extends to the unbounded
derived category to give reduction isomorphisms valid unconditionally,
for unbounded complexes. To keep this article
from becoming even longer the exposition is being postponed to the
manuscript~\cite{NeemanTIFR}. 
\ee

In this article we do not
study formal schemes at all although there is a version
of the theory of
\cite{Avramov-Iyengar08,
Avramov-Iyengar-Lipman-Nayak10}
valid for formal schemes, see Shaul~\cite{Shaul13}.

\section{Introduction}
\label{S0}

Let us begin by recalling an old theorem of Nagata. Suppose
$X$ and $Y$ are noetherian schemes and $f:X\la Y$ is
a separated morphism of finite type. Then $f$ may be factored as
$X\stackrel u\la \ov X\stackrel p\la Y$ with $u$ an open immersion
and $p$ proper. Note that the open immersion $u$ is certainly flat,
and is a monomorphism in the category of schemes. And a
proper map of noetherian schemes is most definitely
of finite type and universally quasi-proper, see
the footnotes to Lemma~\ref{L-1.1}.

As above let $f:X\la Y$ be a separated morphism of finite type between
noetherian schemes, and let $X\stackrel u\la \ov X\stackrel p\la Y$
be a factorization as given by Nagata. 
The functor $f^!:\Dqcpl(Y)\la\Dqcpl(X)$ is traditionally defined 
to be the composite
$\Dqcpl(Y)\stackrel{p^\times}\la\Dqcpl(\ov X)\stackrel{u^*}\la\Dqcpl(X)$. 
It needs to be checked
that $f^!$ does not depend on the choice of Nagata compactification and is
pseudofunctorial. This has been done in the bounded-below derived
category
but not in the unbounded derived category. The key lemma in proving
that $f^!$ has good properties is the base-change theorem in 
Hartshorne~\cite[Corollary~3.4(a) on p.~383,
elaborated in (5) on p.~191]{Hartshorne66} or
Verdier~\cite[Theorem~2 on p.~394, proof pp.~400---407]{Verdier68}
which says that, if we have a cartesian square of noetherian schemes
\[
\CD
W @>u>> X\\
@VfVV @VVgV \\
Y @>v>> Z
\endCD
\]
with $g$ proper and $v$ flat, then the base-change map
$\Phi:u^*g^\times\la f^\times v^*$ is an isomorphism. The reader can
see the way base-change is used in the proof of the
pseudofunctoriality of $f^!$ in Lipman~\cite[proof of
Theorem~4.8.1]{Lipman09}.
Unfortunately the
base-change theorem is true as stated only in the bounded-below
derived category.
An (unbounded) counterexample may be found in \cite[Example~6.5]{Neeman96}.

In this paper we fix the problem---we prove statements that hold in
the unbounded derived category. 
In the Introduction we will not state our results in
maximal possible generality; instead we strive for clarity
and reasonably useful generality.
To formulate our results it will be helpful to introduce

\ntn{N0.1}
We will assume given a 2-subcategory $\seq$ of the 2-category
of noetherian algebraic stacks. The 2-morphisms are simple: any
2-morphism in the category of
algebraic stacks, between 1-morphisms in $\seq$, belongs to $\seq$.
The objects 
$X\in\seq$ will be assumed to be noetherian stacks with quasi-affine diagonals,
to admit finitely-presentable,
representable, separated and \'etale covers $\ov X\la X$ with $\ov X$
satisfying the resolution property, and to satisfy one
of the two conditions below
\be
\item
  either $X$ is a $\qq$--stack
\item
  or the diagonal of $X$ is quasi-finite.
  \setcounter{enumiii}{\value{enumi}}
  \ee
  All notherian Deligne-Mumford stacks satisfy the hypotheses---this includes
  all noetherian algebraic spaces, which includes noetherian schemes. 

Now for the 1-morphisms: we assume every
$f:X\la Y$ in $\seq$ to be separated, as well
as concentrated as in \cite[Definition~2.4]{Hall-Rydh13}. All morphisms
of noetherian algebraic spaces are concentrated. We furthermore
suppose
  \be
  \setcounter{enumi}{\value{enumiii}}
\item
  If $X$ is an object of $\seq$ and $u:U\la X$ is an open immersion then $u$
  belongs to $\seq$.
  \item
    Pullbacks in the 2-category of algebraic stacks, of diagrams
    \[
    \CD
    @. X \\
    @. @VVgV\\
    Y @>v>> Z
    \endCD
    \]
    in $\seq$, belong to $\seq$.
  \item
    Every morphism $f:X\la Y$ in $\seq$ admits in $\seq$ 
    a Nagata compactification.
    For us a Nagata compactification is a 2-isomorphism $pu\la f$,
    where $u,p$ are composable 1-morphisms
    $X\stackrel u\la \ov X\stackrel p\la Y$, where $u$
    is a dominant, flat monomorphism
    and $p$ of finite type and universally quasi-proper (for noetherian
    algebraic spaces the condition on $p$ means that it is proper).
    \ee
    \entn

\rmk{R0.10095}    
For the non-expert on stacks: if your major focus is the case of schemes or
algebraic spaces then all noetherian algebraic spaces satisfy the hypotheses on
the objects. The only hypothesis on the morphisms which might be tricky is the
existence of Nagata compactifications: it's known for morphisms of finite type
of algebraic spaces, and in the case of schemes one can even allow
morphisms essentially of finite type.

In the case of more general stacks most of the hypotheses are very mild.
There is the question of the existence of Nagata compactifications, it's still
open how generally they exist, see below for what is known. But an
important
restriction is that the 1-morphisms are assumed concentrated.

Let $G$ be a finite group, for example $\zz/p$ where $p$ is a prime number,
and let $k$ be a field of characteristic $p$. If $BG$ is the classifying stack
then the morphism $f:BG\la \spec k$ provides an
example of a non-concentrated morphism of Deligne-Mumford stacks.
In this case there are equivalences $\Dqc\big(\spec k\big)=\D(k)$
and $\Dqc(BG)=\D(kG)$, and under these equivalences the
functor $\R f_*$ identifies with $\RHom^{}_{kG}(k,-)$. It is known not to
respect coproducts, hence there is no right
adjoint $f^\times$. This makes it a challenge to extend Grothendieck
duality to such $f$'s. There are non-concentrated $f$'s out there and,
as presented below, the theory can't possibly work for them.

The class of concentrated morphisms was designed to exclude such
pathologies. A morphism $g:X\la Y$ is concentrated if $\R f_*$
is uniformly bounded above for all pullbacks $f$ of $g$.
That is: there is an integer $n$
so that, for all pullbacks $f$ of $g$,
we have $\R f_*\Dqc^{\leq 0}\subset\Dqc^{\leq n}$. This recent notion
was defined by Hall and Rydh~\cite[Definition~2.4]{Hall-Rydh13},
and their article goes on to show (among other things)
that for concentrated morphisms
many nice properties hold. To the extent that we are able to state
our results in the generality of algebraic stacks (rather
than algebraic spaces), this article hinges on
very recent work about algebraic stacks---we depend heavily on
the machinery developed in \cite{Hall-Rydh13}.
\ermk
    
\exm{E0.9991}
The following are examples of $\seq$'s as in Notation~\ref{N0.1}:
\be
\item
$\seq$ could be the category whose objects are noetherian
schemes, and whose morphisms are separated maps essentially of
finite type. The fact that any morphism
in this $\seq$ has a Nagata compactification
is due to Nayak~\cite[Theorem~3.6]{Nayak09}, although in the special case 
where $f$ is of finite type this is the old theorem of Nagata
discussed in the first paragraph of the Introduction. See
Nagata~\cite{Nagata62}
for the original proof or Conrad~\cite{Conrad07} for a more modern
treatment, following Deligne.
\item
  $\seq$ could be the category whose objects are noetherian
  algebraic spaces, and whose
  morphisms are separated maps of finite type. 
  The existence of Nagata compactifications for this $\seq$ is the
  main theorem of Conrad, Lieblich and
  Olsson~\cite{Conrad-Lieblich-Olsson12}, although the special case
  of algebraic spaces of finite type over an excellent scheme may
  essentially be
  found in Raoult's thesis (1974, unpublished), with sketches in his articles~\cite{Raoult71,Raoult74}.
\item
  We may take for $\seq$ the 2-category whose
  objects are noetherian Deligne-Mumford $\qq$--stacks, and
  whose 
  1-morphisms are those
  maps which are separated and of finite type. 
  The
  proof of the existence of Nagata compactifications for such
  maps may be
  found in Rydh's unpublished paper~\cite{Rydh09}, available
  on the author's web page.
\item
  $\seq$ could be the 2-category whose objects are
  noetherian
  stacks with quasi-finite and separated diagonal,
  and whose 1-morphisms are the maps which are separated, of
  finite type and representable. The
  proof of the existence of Nagata compactifications is again
  in Rydh's unpublished paper~\cite{Rydh09}. All representable maps
  are concentrated. The existence of an \'etale map $\ov X\la X$ with
  $\ov X$ satisfying the resolution property follows from
  Rydh~\cite[Theorem~7.2(iii)]{Rydh11}.
\ee
\eexm

\rmk{R0.90505}
There are other recent results giving Nagata compactifications for morphisms
of stacks. The first such theorem was Kresch~\cite[Theorem~5.3]{Kresch09},
but Rydh~\cite[Theorem~F]{Rydh09} is more general, we have already met the
$\seq$ it produces in Example~\ref{E0.9991}(iii).
There is also a result of Edidin, but for us
the problem is that in the factorization
$f=pu$ the $p$ is not in general concentrated. Edidin's
result is unpublished, but the proof may be found in Rydh's survey
article~\cite{Rydh09B}, which is available on the author's web page.

Rydh's methods allow him to compactify certain maps even without the
characteristic~zero hypothesis. However the current version of
Rydh~\cite{Rydh09} is undergoing revision and the correct, precise statement will
hopefully appear in a future manuscript.
\ermk

\rmk{R0.90506}
We have already discussed the importance of the hypothesis that the morphisms
in $\seq$ must be concentrated. Remember also that
they are assumed separated: I'm not sure how inescapable this is, but
if there is a theory that can handle non-separated morphisms then it must
somehow take care of pathologies like Example~\ref{EA1.1}.

Since morphisms in $\seq$ are assumed concentrated and separated
the following proposition becomes relevant. The author
thanks David Rydh for pointing this out and providing the proof. 
\ermk

\pro{P0.90507}
Any separated and concentrated morphism $f:X\la Y$,
of noetherian algebraic stacks,  has relatively tame
stabilizers. In the terminology of
Abramovich, Olsson and Vistoli~\cite{Abramovich-Olsson-Vistoli08}
any such $f$ is a relatively tame Artin stack.
\epro

\prf
The question is local in $Y$ so we can assume that $Y$ is affine.
Separated implies that the stabilizers of $X$ are proper.
For $X$ to be concentrated, it is necessary that $BG$
is concentrated for every stabilizer $G$, see
Hall and Rydh~\cite[Theorem~C, 1$\Longrightarrow$2]{Hall-Rydh15}.
This rules out non-finite $G$, as
$BA$ is never concentrated if $A$ is an abelian variety,
see \cite[Theorem~B or Proposition~1.5]{Hall-Rydh15}.
For finite $G$, we have that $BG$ is concentrated iff $G$ is linearly
reductive by~\cite[Theorem~B or Theorem~1.2]{Hall-Rydh15}.

Conversely: if all stabilizer groups of $X$ are finite and linearly reductive
then $X$ is concentrated; see~\cite[Theorem~2.1(2)]{Hall-Rydh15}.
\eprf

In this article we will define $f^!$ on the unbounded
derived category. For this we prove a technical
refinement of Lemma~\ref{L-1.1},
which will be our replacement for
the base-change theorem. Again we do not state the most
general version.

\thm{T0.1}
Suppose we are given a diagram, where the objects
satisfy the restrictions of Notation~\ref{N0.1}
\[
\CD
U @>u'>> W @>u>> X\\
@.  @VfVV @VVgV \\
@.  Y @>v>> Z
\endCD
\]
Assume that the square is 2-cartesian, the map $g$ is concentrated, of
finite type and universally
quasi-proper, the map $v$ is
flat, and the image of $u'$ is contained in the set on which $f$ is of
finite Tor-dimension. Let 
$\LL{u'}^*:\Dqc(W)\la\Dqc(U)$,
$u^*:\Dqc(X)\la\Dqc(W)$,
$v^*:\Dqc(Z)\la\Dqc(Y)$,
$f^\times:\Dqc(Y)\la\Dqc(W)$ and 
$g^\times:\Dqc(Z)\la\Dqc(X)$
be the usual pullback maps but in the unbounded derived category.
If $\Phi:u^*g^\times\la f^\times v^*$ is the
base-change map then $\LL{u'}^*\Phi:\LL{u'}^* u^*g^\times\la \LL{u'}^* f^\times
v^*$ is an isomorphism. Furthermore the isomorphic functors
$\LL{u'}^* u^*g^\times\cong \LL{u'}^* f^\times
v^*$
respect coproducts.
\ethm

Once we have this theorem we can make deductions.
The next theorem will give a list of formal
consequences. It follows that any $\seq$, satisfying the conditions
of Notation~\ref{N0.1},
also satisfies all the conclusions of Theorem~\ref{T0.5}.
The short
summary
is that we will construct on $\seq$ a structure, some 2-functors and natural
transformations
among them, see Theorem~\ref{T0.5}~(i)--(vii).
To fix ideas we
make the blanket assumption that all our 2-functors and natural transformations
strictly respect identities. 
After introducing the players, the 2-functors and natural transformations,
comes a long
list of compatibility properties which our structure satisfies. 
After stating the theorem, most of the remainder of the introduction will use
the compatibility properties in Theorem~\ref{T0.5} to show that the
structure
is unique up to canonical isomorphism. The end of the introduction
will mention a couple of applications.

\thm{T0.5}
Let the 2-category $\seq$ be as above, and let $\Tri$ be the
2-category of triangulated categories. There are three contravariant
2-functors $\seq\la\Tri$ which we will denote
$(-)^*$, $(-)^\times$ and $(-)^!$. The formulas are
\be
\item
  On objects: for any object $X\in\seq$ we have $X^*=X^\times=X^!=\Dqc(X)$. We
  emphasize: on objects the three 2-functors are identical and unsurprising,
  all three send an object $X$ to $\Dqc(X)$. It is what they do to 1-morphisms
  and 2-morphisms that distinguishes them.
\item
  On 1-morphisms: for a 1-morphism $f:X\la Y$ in $\seq$ our
  $f^*:\Dqc(Y)\la\Dqc(X)$
  is the usual derived pullback $\LL f^*$, our $f^\times:\Dqc(Y)\la\Dqc(X)$
  is the usual $f^\times$
  [i.e.~the right adjoint of the right adjoint of $f^*$],
  and our $f^!:\Dqc(Y)\la\Dqc(X)$ is a new functor
  (at least new in this generality).
\setcounter{enumiii}{\value{enumi}}
\ee
The 2-functors $(-)^*$ 
  and $(-)^\times$ are pseudofunctorial, meaning that for composable
  morphisms $X\stackrel f\la Y\stackrel g\la Z$ in $\seq$ we have isomorphisms
  $\tau(f,g):(gf)^*\la f^*g^*$ and $\delta(f,g):(gf)^\times\la f^\times g^\times$,
  but the 2-functor $(-)^!$ is only oplax, we just have a natural map
  $\rho(f,g):(gf)^!\la f^!g^!$, it needn't be an isomorphism in general.
  
These are the 2-functors in the theory, now for
the natural transformations. We begin with the straightforward ones:
\be
\setcounter{enumi}{\value{enumiii}}
\item
  There is an oplax natural transformation $\psi:(-)^\times\la(-)^!$.
\item
  The 2-functor $(-)^*$ is a premonoid, there is a pseudonatural transformation
  $\mu:(-)^*\times(-)^*\la(-)^*$
  satisfying an obvious associativity property.
\setcounter{enumiii}{\value{enumi}}
\ee
\begin{em}
  The formalism of 2-functors which are premonoids is treated
  extensively and in glorious generality in the category-theory literature.
  The reader might wish to look (for example)
  at 
  Chikhladze, Lack and Street~\cite{Chikhladze-Lack-Street10}; the modern
  terminology for our premonoids is \emph{monoidales.}
  In this article we give a minimal discussion in
  \S\ref{S805} and \S\ref{S703}, with only the results we absolutely
  need. For the Introduction, suffice it to say that the natural
  transformation $\mu$ takes the object $X\in\seq$ to the usual
  tensor product functor $\mu_X^{}:\Dqc(X)\times\Dqc(X)\la\Dqc(X)$. To
  say that this is part of a 
  pseudonatural transformation is just a concise way of packaging the
  usual information that, for each 1-morphism $f\in\seq$, the functor $f^*$ is
  strong monoidal and this strong monoidal structure is compatible with
  composition.
\end{em}
\be
\setcounter{enumi}{\value{enumiii}}
\item
  The functors $(-)^\times$ and $(-)^!$ are oplax modules over $(-)^*$.
  That is
  we have oplax natural transformations $\chi:(-)^*\times(-)^\times\la(-)^\times$
  and $\s:(-)^*\times(-)^!\la(-)^!$ satisfying the obvious associativity
  property. 
\setcounter{enumiii}{\value{enumi}}
\ee
\begin{em}
  Once again, there is considerable category-theoretic literature on
  functors acted on by monoidales (our premonoids), see
  for example Day and Street~\cite{Day-Street97} or
  Lack~\cite{Lack00}.
  The usual setting is infinitely more general than ours,
  we will come back to this in Remark~\ref{R805.99999}. In
  the literature the current name
  for our modules is \emph{actegories.} In this
  article our treatment is minimal and very restrictive,
  narrowly tailored to the application we have in mind---see \S\ref{S805}
  for our formalism.
Part of our definition, of a module over a 2-functor which is
  a premonoid, insists that on objects the natural transformations coincide.
  That is, if $X\in\seq$ is an object, then $\chi(X):X^*\times X^\times\la X^\times$ and $\s(X):X^*\times X^!\la X^!$ are also just the tensor-product
  functor $\mu_X^{}:\Dqc(X)\times\Dqc(X)\la\Dqc(X)$. More detail is
  provided in \S\ref{S805}, for the Introduction let this suffice.
  
Back to the Theorem: there are some more 
natural transformations, but to define them we need to introduce the
2-categories $\hseq$ and $\vseq$. The 1-morphisms in $\seq$ are
the objects of $\hseq$, and the flat 1-morphisms in $\seq$
are the objects of $\vseq$. The 1-morphisms, in both $\hseq$ and $\vseq$,
are the 2-cartesian squares in $\seq$
\[\xymatrix{
  W\ar[r]^u\ar[d]_f\ar@{}[dr]|-{(\diamondsuit)} & X\ar[d]^g \\
  Y\ar[r]^v   & Z
}\]
with flat horizontal maps.
The difference is that in $\hseq$ we view $(\diamondsuit)$ as a morphism
$f\la g$, while in $\vseq$ it is a morphism $u\la v$. On both $\vseq$ and
$\hseq$ there are projection 2-functors $p_1^{}$ and $p_2^{}$ to $\seq$,
taking an object $f:X\la Y$ in either $\hseq$ or $\vseq$ to
$p_1^{}(f)=X$ and $p_2^{}(f)=Y$.
Consider the 2-functors
\[\xymatrix@C+50pt{
  \hseq \ar@<0.5ex>[r]^-{p_1^{}} \ar@<-0.5ex>[r]_-{p_2^{}}&\seq\ar[r]^-{(-)^*} & \Tri
}\]
and 
\[\xymatrix@C+30pt@R=15pt{
  \vseq \ar@<0.5ex>[r]^-{p_1^{}} \ar@<-0.5ex>[r]_-{p_2^{}}   & \seq
  \ar@/^2pc/[rr]^-{(-)^*}
  \ar[rr]^-{(-)^\times}
  \ar@/_2pc/[rr]_-{(-)^!} & &\Tri
}\]
Now we are ready for the remaining
natural transformations.
\end{em}
\be
\setcounter{enumi}{\value{enumiii}}
\item
  On the category $\hseq$: there are three lax natural transformations
  $\tau,\Phi,\theta:(-)^*\circ p_2^{}\la (-)^*\circ p_1^{}$. Of them
  $\tau$ is pseudonatural.
\item
  On the category $\vseq$: there is a pseudonatural transformation
  $\tau:(-)^*\circ p_2^{}\la (-)^*\circ p_1^{}$,
  and oplax natural transformations
  $\Phi:(-)^\times\circ p_2^{}\la (-)^\times\circ p_1^{}$ and 
 $\theta:(-)^!\circ p_2^{}\la (-)^!\circ p_1^{}$.
  \setcounter{enumiii}{\value{enumi}}
\ee
\begin{em}
These are the 2-functors and natural transformations, now it's time for
the relations among them. Perhaps we should begin with (vi) and (vii), to explain
what might look like ambiguous notation:
on the face of it $\tau$, $\theta$ and $\Phi$ have
different meanings, depending on whether the input category is $\hseq$ or
$\vseq$.

On objects: in the category $\vseq$ the natural transformations
$\tau$, $\Phi$ and $\theta$ all take an object $u:W\la X$ to
$\tau(u)=\Phi(u)=\theta(u)=u^*$.
In the category
$\hseq$, the object $f:W\la Y$ maps under $\tau$ to $\tau(f)=f^*$,
under $\Phi$ to $\Phi(f)=f^\times$, while
$\theta$ takes it to $\theta(f)=f^!$.

The 1-morphisms,  in both $\hseq$ and  $\vseq$, are 2-cartesian squares
$(\diamondsuit)$ as
above. And the naturality requires us to provide 2-morphisms. If we view
$(\diamondsuit)$ as a morphism $(\diamondsuit):f\la g$ in $\hseq$,
the natural transformation $\Phi$ must provide 
a 2-morphism comparing $\Phi(f)\big[p_2^{}(\diamondsuit)\big]^*$ with
$\big[p_1^{}(\diamondsuit)\big]^*\Phi(g)$ [in a 1-category they would just be
  equal]. The direction we choose is lax, the 2-morphism goes in the
direction $\Phi(\diamondsuit):\big[p_1^{}(\diamondsuit)\big]^*\Phi(g)\la
\Phi(f)\big[p_2^{}(\diamondsuit)\big]^*$, and is nothing
other than the usual base-change map
$\Phi(\diamondsuit):u^*g^\times\la f^\times v^*$.
Similarly the natural transformation $\tau$ must give us a map
$\tau(\diamondsuit):\big[p_1^{}(\diamondsuit)\big]^*\tau(g)\la
\tau(f)\big[p_2^{}(\diamondsuit)\big]^*$;
we choose the canonical
isomorphism $\tau(\diamondsuit)=(\diamondsuit)^*:u^*g^*\la f^*v^*$.
The fact that $\tau(\diamondsuit)$  is an isomorphism qualifies
$\tau$ to be pseudonatural, in the case of $\Phi$ there is a direction.

Of course $(\diamondsuit)$ can also be viewed as a 1-morphism
$(\diamondsuit):u\la v$ in the 2-category $\vseq$. The natural transformation
$\Phi$ takes the objects $u,v\in\vseq$ to $\Phi(u)=u^*$ and $\Phi(v)=v^*$,
and the naturality requires us to provide a comparison 2-morphism
between $\Phi(u)\big[p_2^{}(\diamondsuit)\big]^\times$ and
$\big[p_1^{}(\diamondsuit)\big]^\times\Phi(v)$. This time the direction
is oplax: the map we choose has the direction
$\Phi(\diamondsuit):\Phi(u)\big[p_2^{}(\diamondsuit)\big]^\times\la
\big[p_1^{}(\diamondsuit)\big]^\times\Phi(v)$,
and the reason
for what seems confusing notation is that it is the identical map
to the horizontal version, namely $\Phi(\diamondsuit):u^*g^\times\la f^\times v^*$.
Similarly $\tau(\diamondsuit):\tau(u)\big[p_2^{}(\diamondsuit)\big]^*\la
\big[p_1^{}(\diamondsuit)\big]^*\tau(v)$ is just
the isomorphism $\tau(\diamondsuit)=(\diamondsuit)^*:u^*g^*\la f^*v^*$.
The assertion that we have natural transformations really comes down
to the standard fact that the base-change map is compatible with
the concatenation of 2-cartesian squares, both horizontally and vertically,
as is the map $(\diamondsuit)^*$.

What is really being asserted is that there is a parallel which
works for
$(-)^!$. We assert the existence of a 2-morphism
$\theta(\diamondsuit):u^*g^!\la f^!v^*$, which serves as a lax natural
transformation in the case of $\hseq$ and an oplax natural transformation
for $\vseq$, that is respects the concatenation of 2-cartesian squares,
both vertical and horizontal.

We have explained at length the natural transformations
$\tau$, $\Phi$ and $\theta$, and now we come to the other natural
transformations and their interplay.
\end{em}
\sthm{ST0.5.-1}
On objects, the oplax natural transformation $\psi:(-)^\times\la(-)^!$
must provide, for every object $X\in\seq$, a 1-morphism
$\psi(X):X^\times\la X^!$. This map
is $\psi(X)=\id:\Dqc(X)\la\Dqc(X)$.
\esthm
\sthm{ST0.5.1}
The natural transformation $\psi:(-)^\times\la(-)^!$ is a homomorphism
of $(-)^*$--modules. That is the square of natural transformations
\[\xymatrix@C+40pt{
  (-)^*\times(-)^\times\ar[r]^-{\chi}
  \ar[d]_-{\id\times\psi}& (-)^\times\ar[d]^-{\psi} \\
  (-)^*\times(-)^! \ar[r]^-{\s} & (-)^!
}\]
strictly commutes.
\esthm
\sthm{ST0.5.3}
On the 2-category $\vseq$ the natural transformations $\Phi$ and $\theta$
are strictly compatible with $\psi$. That is the following square
of natural transformations strictly commutes
\[\xymatrix{
(-)^\times\circ p_2^{}\ar[d]_-{\psi p_2^{}} \ar[r]^-{\Phi} & (-)^\times\circ p_1^{}\ar[d]^-{\psi p_1^{}} \\
(-)^!\circ p_2^{} \ar[r]^-{\theta} & (-)^!\circ p_1^{} 
}\]
\esthm
\sthm{ST0.5.3.1}
On the category $\hseq$ the compatibility of $\Phi$, $\theta$
is less strict. If we consider the diagram of 2-functors
on $\hseq$
\[\xymatrix{
(-)^*\circ p_2^{}\ar[d]_-{\id} \ar[r]^-{\Phi} & (-)^*\circ p_1^{}\ar[d]^-{\id} \\
(-)^*\circ p_2^{} \ar[r]^-{\theta} & (-)^*\circ p_1^{} 
}\]
then the object $f\in\hseq$ maps under $\Phi$ to $f^\times$ and under $\theta$
to $f^!$. The natural transformation $\psi:(-)^\times\la(-)^!$
provides us with a
2-morphism $\psi(f):f^\times\la f^!$ (see Remark~\ref{R0.953} for details),
and the assignment taking the
object $f\in\hseq$ to the 2-morphism $\psi(f):f^\times\la f^!$ extends to
a modification of natural transformations, that is a morphism
$\Phi\la\theta$.
\esthm
\sthm{ST0.5.3.2}
On the category $\vseq$ the squares
\[\xymatrix@C+30pt{
[(-)^*\circ p_2^{}]\times [(-)^\times\circ p_2^{}]\ar[d]_-{\tau\times\Phi} \ar[r]^-{\chi} & (-)^\times\circ p_2^{}\ar[d]^-{\Phi} \\
[(-)^*\circ p_1^{}]\times [(-)^\times\circ p_1^{}] \ar[r]^-{\chi} & (-)^\times\circ p_1^{}  \\
}\]
and
\[\xymatrix@C+30pt{
[(-)^*\circ p_2^{}]\times [(-)^!\circ p_2^{}]\ar[d]_-{\tau\times\theta} \ar[r]^-{\s} & (-)^!\circ p_2^{}\ar[d]^-{\theta} \\
[(-)^*\circ p_1^{}]\times [(-)^!\circ p_1^{}] \ar[r]^-{\s} & (-)^!\circ p_1^{}  \\
}\]
both 2-commute. More precisely:
in both cases the composites of the shape $\llcorner$
take the object $u:W\la X$ in the 2-category $\vseq$ to
$\mu_W^{}(u^*\times u^*)$, while the composites
of the shape $\urcorner$ take $u$ to $u^*\mu_X^{}$.
The pseudonatural transformation $\mu:(-)^*\times(-)^*\la(-)^*$
gives a 2-isomorphism
$\mu(W,X):\mu_W^{}(u^*\times u^*)\la u^*\mu_X^{}$, and we
assert that this extends to modifications making
both squares 2-commute.

We leave to the reader the analogous statement about $\hseq$.
\esthm

\begin{em}
We have had a string of results telling us that pairs of composites
of our natural transformations agree, or maybe only agree up to
modification. The remaining results are conditions under which
these maps induce isomorphisms.
\end{em}
\sthm{S0.5.3.7}
Let $X\stackrel f\la Y\stackrel g\la Z$
be a pair of composable 1-morphisms
in $\seq$. The map  $\rho(f,g):{(gf)}^!\la f^!g^!$, which
is part of the structure of the oplax 2-functor $(-)^!$, 
is an isomorphism if one of the conditions below holds:
\begin{enumerate}
\item
$f$ is of finite Tor-dimension.
\item
  $g$ is of finite type and universally quasi-proper.
\item
  The composite $gf$ is of finite type and universally quasi-proper. 
\item
We restrict to the subcategory $\Dqcpl(Z)\subset\Dqc(Z)$.
\end{enumerate}
\esthm
\sthm{ST0.5.4}
Let $f:X\la Y$ be a 1-morphism in $\seq$.
The map $\psi(f):f^\times\la f^!$ is an isomorphism if $f$ is
of finite type and
universally quasi-proper.
\esthm
\sthm{ST0.5.7}
Let $(\diamondsuit)$ be the 2-cartesian square pictured between (v) and (vi)
above.
The base-change map $\theta(\diamondsuit):u^*g^!\la f^!v^*$ is an isomorphism
if one of the following holds:
\begin{enumerate}
\item
  $f$ is of finite Tor-dimension. More generally if $u':U\la W$ is
  a map whose image lies in
the subset of $\,W$ on which $f$ is of finite Tor-dimension, then 
${u'}^*\theta(\diamondsuit):{u'}^*u^*g^!\la {u'}^*f^!v^*$
is an isomorphism.
\item
We restrict to the subcategory $\Dqcpl(Z)\subset\Dqc(Z)$.
\end{enumerate}
\esthm
\sthm{ST0.5.17}
Let $f:X\la Y$ be a 1-morphism in $\seq$ 
and let $E,F$ be objects in $\Dqc(Y)$.
Then $\chi(f):\mu_X^{}(f^*\times f^\times)\la f^\times\mu_Y^{}$
and $\s(f):\mu_X^{}(f^*\times f^!)\la f^!\mu_Y^{}$
are both natural transformation of functors $\Dqc(Y)\times\Dqc(Y)\la\Dqc(X)$,
and we can evaluate them at the pair $(E,F)$. Let
us write $\chi(f,E,F):f^*E\oo f^\times F\la f^\times(E\oo F)$
and $\s(f,E,F):f^*E\oo f^! F\la f^!(E\oo F)$ for the resulting
1-morphisms in the category $\Dqc(X)$. Then
\begin{enumerate}
\item
Both $\chi(f,E,F)$ and $\sigma(f,E,F)$ are isomorphisms if $E$ is a 
perfect complex.
\item
  $\sigma(f):\mu_X^{}(f^*\times f^!)\la f^!\mu_Y^{}$ is an isomorphism
  as long as $f$ is of finite Tor-dimension. That is: if $f$ is of finite
  Tor-dimension then $\s(f)$ evaluates to an isomorphism on every pair
  of objects $(E,F)\in\Dqc(Y)$.

  More generally: if $u:W\la X$ is a flat morphism so that $fu$ is
  of finite Tor-dimension, then
  $u^*\sigma(f):u^*\mu_X^{}(f^*\times f^!)\la u^*f^!\mu_Y^{}$
  is an isomorphism.
\end{enumerate}
\esthm
\sthm{ST0.5.999}
Let the notation be as in \ref{ST0.5.17}. The morphism
$\chi(f,E,F):f^*E\oo f^\times F\la f^\times (E\oo F)$ corresponds, under
the adjunction $\R f_*\dashv f^\times$, to a morphism
$\xi(f,E,F):\R f_*(f^*E\oo f^\times F)\la E\oo F$. The map $\xi(f,E,F)$
is just the composite
\[\xymatrix@C+30pt{
  \R f_*(f^*E\oo f^\times F) \ar[r]^-{\cong} & E\oo \R f_*f^*F\ar[r]^-{\id\oo\e} &
  E\oo F\ ,
}\]
where the isomorphism is by the projection formula and the
map $\e:\R f_*f^\times\la\id $ is the counit of adjunction. 
\esthm
\ethm

It might help a little if we work out, explicitly, what some of the
2-categorical 
formalism says. We begin with

\rmk{R0.959}
We made the blanket assumption that our 2-functors
strictly respect identities. Hence $\id^*=\id^\times=\id^!=\id$ and the
2-morphisms
\[\xymatrix@R-15pt{
  f^*\id^* &&(\id\circ f)^* \ar@{=}[r]\ar[ll]_-{\tau(f,\id)} &  f^*\ar@{=}[r]& (f\circ\id)^*\ar[rr]^-{\tau(\id,f)}&&\id^* f^*\\
  f^\times\id^\times &&(\id\circ f)^\times \ar@{=}[r]\ar[ll]_-{\delta(f,\id)} &  f^\times\ar@{=}[r]& (f\circ\id)^\times\ar[rr]^-{\delta(\id,f)}&&\id^\times f^\times\\
   f^!\id^! &&(\id\circ f)^!\ar@{=}[r]\ar[ll]_-{\rho(f,\id)} &  f^!\ar@{=}[r]& (f\circ\id)^!\ar[rr]^-{\rho(\id,f)}&&\id^! f^!
}\]
are all identities. We also made the blanket assumption that our
natural transformations strictly respect identities. The 
2-cartesian squares
\[\xymatrix{
W \ar[r]^u \ar[d]_\id\ar@{}[dr]|-{(\clubsuit)} & X \ar[d]^\id & &W \ar[r]^\id \ar[d]_f\ar@{}[dr]|-{(\diamondsuit)} & W \ar[d]^f  \\
W \ar[r]^u & X & & Y\ar[r]^\id & Y
}\]
can be viewed as identity morphisms $(\clubsuit)\in\vseq$ and
$(\diamondsuit)\in\hseq$, and we learn that
\[\xymatrix@C+15pt@R-15pt{
  u^*\id^*\ar[r]^-{\tau(\clubsuit)} & \id^* u^* &
  u^*\id^\times\ar[r]^-{\Phi(\clubsuit)} & \id^\times u^* &
  u^*\id^!\ar[r]^-{\theta(\clubsuit)} & \id^! u^*\\
   \id^*f^*\ar[r]^-{\tau(\diamondsuit)} & f^* \id^* &
  \id^*f^\times\ar[r]^-{\Phi(\diamondsuit)} & f^\times \id^* &
   \id^*f^!\ar[r]^-{\theta(\diamondsuit)} & f^! \id^*
}\]
are all identities. Similarly for the natural transformation $\psi$: the
2-morphism $\psi(\id):\id^\times\la\id^!$ is the identity.
The fact that the natural transformations
$\mu:(-)^*\times(-)^*\la(-)^*$, $\chi:(-)^*\times(-)^\times\la(-)^\times$
and $\s:(-)^*\times(-)^!\la(-)^!$ respect identities comes down to saying
that, for any object $X\in\seq$ and
for the identity map $\id:X\la X$, the 2-morphisms
\[
\xymatrix@C-5pt{
 \mu_X^{}(\id^*\times\id^*)\ar[r]^-{\mu} & \id^*\mu_X^{} &
 \mu_X^{}(\id^*\times\id^\times)\ar[r]^-{\chi} & \id^\times\mu_X^{} &
  \mu_X^{}(\id^*\times\id^!)\ar[r]^-{\s} & \id^!\mu_X^{} 
}
\]
are all identities.
\ermk

\rmk{R0.2007}
In Theorem~\ref{T0.5}(iii) we learned about the
existence of an oplax natural transformation
$\psi:(-)^\times\la(-)^!$, while in \ref{ST0.5.-1} we were told that, for
$X$ an object of $\seq$, the map $\psi(X):X^\times\la X^!$ is the identity
$\id:\Dqc(X)\la\Dqc(X)$. What this means is the following.

Let $f:X\la Y$ be a 1-morphism in $\seq$. The under the 2-functors
$(-)^\times$ and $(-)^!$ it is mapped (respectively) to
$f^\times:\Dqc(Y)\la\Dqc(X)$ and $f^!:\Dqc(Y)\la\Dqc(X)$. The natural
transformation must provide us with a comparison map between
$f^\times=\psi(X)f^\times$ and $f^!=f^!\psi(Y)$, and the direction is
oplax, the map $\psi(f)$ goes $\psi(X)f^\times\la f^!\psi(Y)$,
or more simply we have a map $\psi(f):f^\times\la f^!$. The assertion
that this is a natural transformation includes the statement that
$\psi(f)$ must be compatible with composition in the obvious way.

Because this particular natural transformation is such that $\psi(X)=\id_X^{}$
for every object $X$, we may also view it as a lax natural transformation
$(-)^!\la(-)^\times$. If we do this then the map
goes in the direction $f^\times\psi(Y)\la \psi(X)f^!$, that is the lax direction.
\ermk

\rmk{R0.953}
Let us evaluate the strictly commutative diagram of \ref{ST0.5.3} at the
1-morphism $(\diamondsuit):u\la v$ in the 2-category $\vseq$.
The 2-functors
\[(-)^\times\circ p_2^{},\quad
(-)^\times\circ p_1^{},\quad
(-)^!\circ p_2^{},\quad
(-)^!\circ p_1^{} 
\]
take $(\diamondsuit)$ respectively to $g^\times$, $f^\times$, $g^!$ and
$f^!$. The natural transformation $\Phi$ and $\theta$ both take the
object $u$ to the 1-morphism $u^*$ and the object
$v$ to the 1-morphism $v^*$. And the natural transformation $\psi$
is the identity on objects. In other words: on objects both composite
natural transformations
are simple, they take $u$ to $u^*$ and $v$ to $v^*$. 

The evaluation of the diagram on $(\diamondsuit)$ therefore becomes
\[
\xymatrix@C+5pt@R+5pt{
u^*g^\times \ar[d]_{u^*\psi(g)}\ar[r]^{\Phi(\diamondsuit)} & f^\times v^*\ar[d]^{\psi(f)v^*}\\
u^*g^! \ar[r]^{\theta(\diamondsuit)} & f^! v^*
}
\]
and \ref{ST0.5.3} asserts that this commutes. The reader can check
the $\hseq$ version: 
if we evaluate the 2-commutative diagram of \ref{ST0.5.3.1} at the
1-morphism $(\diamondsuit)\in\hseq$ we end up with the
commutative square above. In other words
\ref{ST0.5.3} and  \ref{ST0.5.3.1}, evaluated at the 1-morphism
$(\diamondsuit)$ in either $\vseq$ or $\hseq$, give rise to identical
commutative squares.
\ermk

\rmk{R0.1025}
Let us evaluate the second 2-commutative diagram
of \ref{ST0.5.3.2} at $(\diamondsuit)$. It yields the commutative
diagram
\[
\xymatrix@C+40pt{
  \mu_X^{}(u^*\times u^*)(g^*\times g^!) \ar[r]^-{\mu(W,X)}\ar[d]_{\tau(\diamondsuit)\times\theta(\diamondsuit)} &
  u^*\mu_X^{}(g^*\times g^!)\ar[r]^-{\s(X,Z)} & u^*g^!\mu_Z^{}\ar[d]^-{\theta(\diamondsuit)}\\
   \mu_X^{}(f^*\times f^!)(v^*\times v^*) \ar[r]^-{\s(W,Y)} &
  f^!\mu_Y^{}(v^*\times v^*)\ar[r]^-{\mu(Y,Z)} & f^!v^*\mu_Z^{}
}\]
where the maps labeled $\mu(W,X)$ and $\mu(Y,Z)$ come from
the modification, which we must apply both at the source and at the
target. The vertices of the diagram are all 
functors $\Dqc(Z)\times\Dqc(Z)\la \Dqc(W)$ and the arrows are all natural
transformations. Hence the top left and bottom right give
two functors
$\xymatrix{\Dqc(Z)\times\Dqc(Z)\ar@<0.5ex>[r]\ar@<-0.5ex>[r] & \Dqc(W),}$
and the two composites give natural transformations between them.
If we evaluate at an object $(E,F)\in\Dqc(Z)\times\Dqc(Z)$ we
obtain two morphisms in $\Dqc(W)$ which the theorem asserts must
be equal. Concretely this comes down to the commutative diagram
\[
\xymatrix@C+30pt{
  u^*g^*E\oo u^*g^!F \ar[d]_{\tau\times\theta}\ar[r]^-{\mu} & u^*(g^*E\oo g^!F) \ar[r]^-{\s}& u^*g^!(E\oo F)\ar[d]^\theta \\
  f^*v^*E\oo f^!v^*F \ar[r]^-{\s} & f^!(v^*E\oo v^*F) \ar[r]^-{\mu}& f^!v^*(E\oo F) 
}
\]
in the category $\Dqc(W)$.
\ermk

In the next few remarks we will show how to
 use the properties listed in Theorem~\ref{T0.5}
to conclude that the functor $f^!$, the map $\psi:f^\times\la f^!$, the map
$\rho(f,g):{(gf)}^!\la f^! g^!$, the base-change map $\theta(\diamondsuit)$ for
a 2-cartesian square $(\diamondsuit)$ and the map 
$\sigma:f^*E\oo f^!F\la f^!(E\oo F)$ are all determined up to unique
isomorphism.

\rmk{R0.7}
If $f:X\la Y$ is a flat monomorphism we can consider the 
2-cartesian square
\[
\xymatrix@C+10pt@R+10pt{
X \ar[r]^{\id}\ar[d]_{\id} & X\ar[d]^{f}\\
X \ar[r]^{f} & Y\ar@{}[ul]|{(\diamondsuit)}
}
\]
Since $\id:X\la X$ is definitely of finite
Tor-dimension \ref{ST0.5.7}(1) applies
and $\theta(\diamondsuit):\id^*f^!\la \id^! f^*$ is an isomorphism.
It gives us a canonical isomorphism
$\theta(\diamondsuit):f^!\la f^*$.

Still considering the simple 2-cartesian square $(\diamondsuit)$ above,
the general commutative diagram of Remark~\ref{R0.953} specializes to
\[
\xymatrix@C+10pt@R+10pt{
\id^* f^\times \ar[r]^{\Phi(\diamondsuit)}\ar[d]_{\id^*\psi(f)} &
\id^\times f^*\ar[d]^{\psi(\id)f^*}\\
\id^* f^! \ar[r]^{_\theta(\diamondsuit)} & \id^!f^*
}
\]
From the previous paragraph we know that $\theta(\diamondsuit):f^!\la f^*$ is 
an isomorphism,
and now we learn that the map $\psi(f):f^\times\la f^!$ can be
computed as the composite 
\[
\CD
f^\times @>\Phi(\diamondsuit)>> f^* @>{\theta(\diamondsuit)}^{-1}>> f^!.
\endCD
\] 
\ermk

\rmk{R0.8}
Suppose we are given two composable 1-morphisms 
$X\stackrel f\la Y\stackrel g\la Z$ in $\seq$, with $f$ a flat
monomorphism
and $g$  of finite type and universally quasi-proper.
Because $(-)^!$ is an oplax 2-functor we
have a map $\rho(f,g):{(gf)}^!\la f^!g^!$,
which is an isomorphism by \ref{S0.5.3.7}~(1) or (2).
By Remark~\ref{R0.7} the map
$ \theta(\diamondsuit):f^!\la  f^*$ is an isomorphism,
while \ref{ST0.5.4} tells is that $\psi(g):g^\times\la g^!$ is an isomorphism. 
Combining these we obtain a canonical isomorphism
\begin{equation}
\label{eq0.1}
\xymatrix@C+30pt{
{(gf)}^!\ar[r]^{\rho(f,g)} &
f^!g^!\ar[r]^{\theta(\diamondsuit)\psi{(g)}^{-1}} & 
f^*g^\times.
}\end{equation}
Moreover, because $\psi:(-)^\times\la (-)^!$ is a natural transformation
it respects composition, and
we have 
a commutative square  
\[\xymatrix@C+15pt{
{(gf)}^\times \ar[r]^{\delta(f,g)}\ar[d]_{\psi(gf)} &f^\times g^\times \ar[d]^{\psi(f)\psi(g)}\\
{(gf)}^! \ar[r]_{\rho(f,g)} & f^!g^!
}\]
If we identify ${(gf)}^\times$ with $f^\times g^\times$ via the
canonical isomorphism $\delta(f,g)$, and identify ${(gf)}^!$ with
$f^*g^\times$ via the isomorphism (\ref{eq0.1}) above, then the map
$\psi(gf):{(gf)}^\times\la {(gf)}^!$ reduces to
$\Phi(\diamondsuit)g^\times: f^\times g^\times\la f^* g^\times$.

Since every morphism in $\seq$ is 2-isomorphic to a
composite $gf$, with $f$ a
flat monomorphism (even dominant) and
$g$ of finite type and universally
quasi-proper, these computations tell us that Theorem~\ref{T0.5}
gives us no choice on how to define $h^!$ and
$\psi(h):h^\times\la h^!$
for any 1-morphism $h\in\seq$.
\ermk

\rmk{R0.9}
Suppose we are given the 2-cartesian square $(\diamondsuit)$ below
\[
\xymatrix@C+10pt@R+10pt{
W \ar[r]^{u}\ar[d]_{f} & X\ar[d]^{g}\\
Y \ar[r]^{v} & Z\ar@{}[ul]|{(\diamondsuit)}
}
\]
Assume
that $u$ and $v$ are flat, and that $f$ and $g$ are of finite
type and universally quasi-proper.
By Remark~\ref{R0.953} the square
\[
\xymatrix@C+5pt@R+5pt{
u^*g^\times \ar[d]_{u^*\psi(g)}\ar[r]^{\Phi(\diamondsuit)} & f^\times v^*\ar[d]^{\psi(f)v^*}\\
u^*g^! \ar[r]^{\theta(\diamondsuit)} & f^! v^*
}
\]
commutes, but by \ref{ST0.5.4} the vertical maps are isomorphisms. Hence
for 2-cartesian squares $(\diamondsuit)$, with
flat horizontal maps
and vertical maps that are of finite
type and universally quasi-proper, the base-change 
map $\theta(\diamondsuit):u^*g^!\la f^!v^*$
can be canonically identified with the ordinary base-change map 
$\Phi(\diamondsuit):u^*g^\times
\la f^\times v^*$.

Now suppose we have a 2-cartesian square
\[
\xymatrix@C+5pt@R+5pt{
W \ar[r]^{u}\ar[d]_{f} & X\ar[d]^{g}\\
Y \ar[r]^{v} & Z\ar@{}[ul]|{(\heartsuit)}
}
\]
where $u,v$ are flat as usual, but where
$f$ and $g$ are flat monomorphisms. We may consider the 
larger diagrams
\[
\xymatrix@C+5pt@R+5pt{
W \ar[r]^{\id}\ar[d]_{\id}  & W \ar[r]^{u}\ar[d]_{f} & X\ar[d]^{g} 
& &
W \ar[r]^{u}\ar[d]_{\id}  & X\ar[r]^{\id}\ar[d]_{\id} & X\ar[d]^{g} 
\\
W \ar[r]^{f} &Y\ar@{}[ul]|{(\spadesuit)} \ar[r]^{v} &
Z\ar@{}[ul]|{(\heartsuit)} & & 
W \ar[r]^{u} &X\ar@{}[ul]|{(\clubsuit)} \ar[r]^{g} &
Z\ar@{}[ul]|{(\#)}
}
\]
The squares are still all 2-cartesian, and
the diagrams concatenate to isomorphic 
squares
\[
\xymatrix@C+20pt@R+5pt{
W \ar[r]^{u}\ar[d]_{\id} & X\ar[d]^{g}\\
X \ar[r]^{vf\cong gu} & Z\ar@{}[ul]|{(\bullet)}
}
\]
In other words: in the category $\hseq$ we have a 2-isomorphism
$(\heartsuit)\circ(\spadesuit)\cong (\bullet)\cong
(\#)\circ(\clubsuit)$.
The natural transformation $\theta$ respects composition, allowing
us to compute $\theta(\bullet)$ as two different composites.
This yields a commutative diagram
\[
\xymatrix{
u^*g^!\ar[drr]|{\theta(\bullet)} \ar[d]_{u^*\theta(\#)}
\ar[r]^{\theta(\heartsuit)} & f^! v^* \ar[r]^{\theta(\spadesuit)v^*} &
f^* v^*\ar[d]^{\wr} \\
u^*g^*
\ar[rr]_-{\sim} & & {(gu)}^*={(vf)}^*
}
\]
Reading the perimeter of this diagram we see that, up to identifying
via
the isomorphisms 
$\theta(\spadesuit):f^!\la f^*$ and $\theta(\#):g^!\la g^*$, the map
$\theta(\heartsuit)$ is just the canonical isomorphism $u^*g^*\la
f^*v^*$.

We have computed the base-change map $\theta(\diamondsuit)$ in the
special cases where $f,g$ are either flat monomorphisms or
of finite type and universally quasi-proper. Now let $(\diamondsuit)$
be arbitrary and factor $g:X\la Z$ up to
2-isomorphism as $X\stackrel{t'}\la S\stackrel{p'}\la
Z$, with $t'$ a dominant, flat monomorphism and $p'$ of finite type
and universally quasi-proper. Pull back to form the
2-cartesian squares
\[\xymatrix{
  W\ar[r]^u\ar[d]_t\ar@{}[dr]|-{(\dagger)} & X\ar[d]^{t'} \\
  R\ar[r]^w\ar[d]_{p}\ar@{}[dr]|-{(\dagger\dagger)} &S\ar[d]^{p'}\\
  Y\ar[r]^v   & Z
}\]
We have just factored the 1-morphism $(\diamondsuit)$, in the 2-category
$\vseq$, as $(\diamondsuit)\cong(\dagger\dagger)\circ(\dagger)$.
Since $\theta$ is a natural transformation of functors on
$\vseq$ it respects composition, and we may compute
$\theta(\diamondsuit)$ in terms of $\theta(\dagger\dagger)$
and $\theta(\dagger)$. This allows us to
compute 
$\theta(\diamondsuit)$ for all $(\diamondsuit)$.
\ermk

\rmk{R0.11}
Suppose we have a 2-cartesian square in $\seq$ 
\[
\xymatrix@C+10pt@R+2pt{
W \ar[r]^{u}\ar[d]_{f} & X\ar[d]^{g}\\
Y \ar[r]^{v} & Z\ar@{}[ul]|{(\spadesuit)}
}
\]
where $u,v$ are flat monomorphisms. In the following diagrams
all the squares are 2-cartesian
\[
\xymatrix@C+10pt@R+2pt{
W \ar[r]^{\id}\ar[d]_{f} & W\ar[d]^{f} & & 
    W \ar[r]^{\id}\ar[d]_{\id} & W\ar[d]^{u} \\
Y \ar[r]^{\id}\ar[d]_{\id} & Y\ar[d]^{v}\ar@{}[ul]|{(\clubsuit)} & &
    W\ar[r]^{u}\ar[d]_{f} & X\ar[d]^{g}\ar@{}[ul]|{(\heartsuit)}\\
Y \ar[r]^{v} & Z\ar@{}[ul]|{(\diamondsuit)} & & 
  Y \ar[r]^{v} & Z\ar@{}[ul]|{(\spadesuit)}
}
\]
The concatenation of these two diagrams are 2-isomorphic, namely
\[
\xymatrix@C+10pt@R+2pt{
W \ar[r]^{\id}\ar[d]_{f} & W\ar[d]^{vf=gu}\\
Y \ar[r]^{v} & Z\ar@{}[ul]|{(\#)}
}
\]
and because $\theta$ is a natural transformation
on the 2-category $\vseq$, we now compute $\theta(\#)$ in two different ways.
The reader can check that the following diagram commutes
\[
\xymatrix@C+20pt@R+10pt{
{(gu)}^!={(vf)}^!\ar[rrrrd]|{\theta(\#)} \ar[rr]^-{\rho(u,g)}\ar[d]_{\rho(f,v)} & & u^!g^! 
      \ar[rr]^{\theta(\heartsuit)g^!} & & u^*g^!\ar[d]^{\theta(\spadesuit)}\\
f^!v^!\ar[rr]_{\theta(\clubsuit)v^!=\id}& &  f^!v^! \ar[rr]_{f^!\theta(\diamondsuit)} & &f^! v^*
}\]
In this diagram $\rho(u,g)$ is an isomorphism because 
$u$ is of finite Tor-dimension, and $\theta(\diamondsuit)$ and 
$\theta(\heartsuit)$ are the isomorphisms of Remark~\ref{R0.7}.
Up to the isomorphisms in the diagram we have that $\rho(f,v)$ agrees
with $\theta(\spadesuit)$.

If $f$ and $g$ are proper then Remark~\ref{R0.9} informs us that,
up to canonical isomorphism, $\theta(\spadesuit)$ agrees with the usual
base-change map $\Phi(\spadesuit)$. In \cite[Example~6.5]{Neeman96}
we see that $\Phi(\spadesuit)$ does not have to be an isomorphism; in
the example the stacks are all noetherian, affine
schemes, the maps $f$ and $g$ are finite and $u,v$ are
dominant open immersions.
Hence $\rho(f,v)$ need not be an isomorphism in general,
and placing further restrictions on the stacks 
is unlikely to help.
\ermk

\rmk{R0.994}
If $f:X\la Y$ is a morphism in $\seq$ then the subset of $X$ on which
$f$ is of finite Tor-dimension is open. I am not aware of a place in
the
literature where this is stated explicitly, even for schemes---but
the reader can find a proof in
Corollary~\ref{C4.497.-111}.
\ermk

\rmk{R0.13}
Suppose we are given two composable 1-morphisms 
$X\stackrel f\la Y\stackrel g\la Z$ in $\seq$. Let $u:U\la X$ be the inclusion
of the open subset on which $f$ is of finite Tor-dimension. Then 
$u$ and $fu$ are of finite Tor-dimension. Because $(-)^!$ is
an oplax 2-functor the square below commutes 
\[\xymatrix@C+15pt{
{(gfu)}^! \ar[r]^{\rho(u,gf)}\ar[d]_{\rho(fu,g)} & u^!{(gf)}^! \ar[d]^{u^!\rho(f,g)}\\
{(fu)}^!g^! \ar[r]_{\rho(u,f)g^!} & u^!f^!g^!
}\]
In this square $\rho(fu,g)$, $\rho(u,gf)$ and $\rho(u,f)$ must be isomorphisms
by \ref{ST0.5.3}, and hence $u^!\rho(f,g)=u^*\rho(f,g)$ is an 
isomorphism.

More generally, if $w:W\la X$ is a map in $\seq$ whose image is contained
in the set on which $f$ is of finite Tor-dimension then $w$ must factor
through $u$, and $\LL w^*\rho(f,g)$ will also be an isomorphism.
\ermk

\rmk{R0.17}
Suppose $X\stackrel f\la Y\stackrel g\la Z$  are composable 1-morphisms in 
$\seq$. Suppose further that we are given a 2-commutative diagram 
in $\seq$
\[
\CD
@. X @>u'>> X'  @>u>> \ov X\\
@. @.   @Vp'VV    @VVpV\\
(\dagger)\qquad\qquad @. @.    Y @>v>> \ov Y\\
 @. @. @. @VVqV \\
@. @.  @. Z
\endCD
\]
where the horizontal maps are all dominant, flat monomorphisms,
the vertical 
maps are all of finite type and universally quasi-proper,
the square is 2-cartesian and $p'u'\cong f$, $qv\cong g$.
Because $(-)^!$ is an oplax
2-functor we obtain a commutative diagram
\[\xymatrix@C-20pt{
{(qvp'u')}^! \ar[rr]^-{\rho(u',qvp')} \ar[rd]_{\rho(p'u',qv)}
      & & {u'}^!{(qvp')}^! \ar[rr]^-{{u'}^!\rho(vp',q)}&&
{u'}^!{(vp')}^!q^!  \ar[ld]^{{u'}^!\rho(v,p')q^!}\\
 & {(p'u')}^!{(qv)}^! \ar[rr]_{\rho(u',p')\rho(v,q)} &&
   {u'}^!{p'}^!v^!q^! &
}\]
By \ref{S0.5.3.7}~(1) and (2) the horizontal maps are 
isomorphisms and, up to these isomorphisms, the map
$\rho(f,g)=\rho(p'u',qv)$ is computed by ${u'}^!\rho(p',v)q^!$.
But in the last paragraph of Remark~\ref{R0.11} we noted that
$\rho(p',v)$ can be computed as the composite
\[\xymatrix@C+40pt{
{(vp')}^!={(pu)}^! 
\ar[r]_-{\sim} \ar@/^2pc/[rrr]^{\rho(p',v)} & u^*p^\times \ar[r]_-{\Phi}
& v^*{p'}^\times \ar[r]_{\sim} & v^! {p'}^!
}\]
 where 
$\Phi:u^*p^\times\la {p'}^\times{u'}^*$ is the base-change
map 
of the cartesian square
in $(\dagger)$. This means that $\rho(f,g)$ can be expressed as the
composite
\[\xymatrix@C+40pt{
{(gf)}^! 
\ar[r]_-{\sim} \ar@/^2pc/[rrr]^{\rho(f,g)} & {u'}^*u^*p^\times
q^\times \ar[r]_-{{u'}^*\Phi q^\times}
& ({u'}^*{p'}^\times )( v^* q^\times)\ar[r]_-{\sim} & f^! g^!
}\]
Thus  as long as we have enough diagrams 
$(\dagger)$ we can compute the map $\rho(f,g)$ for every 
composable
$f,g\in\seq$.
In Lemma~\ref{L6.2058} we produce enough diagrams
$(\dagger)$. 
\ermk

\rmk{R0.19}
From \ref{ST0.5.999} we know the map
$\chi$ explicitly 
for every $f$.
If we evaluate the strictly commutative square
of \ref{ST0.5.1} at a 1-morphism $f:X\la Y$ in $\seq$ we obtain
the commutative square
\[\xymatrix{
\mu_X^{}(f^*\times f^\times) \ar[r]^-{\chi}\ar[d]_-{\mu_X^{}(\id\times\psi(f))} & f^\times \mu_Y^{} \ar[d]^{\psi(f)\mu_Y^{}}\\
\mu_X^{}(f^*\times f^!) \ar[r]^-{\s} & f^! \mu_Y^{}
}\]
If $f$ is of finite type and universally quasi-proper the
vertical maps are isomorphisms by \ref{ST0.5.4}, hence we 
know $\s(f):\mu_X^{}(f^*\times f^!)\la f^! \mu_Y^{}$ whenever $f$
is of finite type and universally quasi-proper.

Now suppose $f$ is a flat monomorphism.  
We have the 2-cartesian square
\[
\xymatrix@C+10pt@R+2pt{
X \ar[r]^{\id}\ar[d]_{\id} & X\ar[d]^{f}\\
X \ar[r]^{f} & Y\ar@{}[ul]|{(\diamondsuit)}
}
\] 
which we may view as a 1-morphism $(\diamondsuit)\in\vseq$,
and we evaluate the second 2-commutative diagram
of \ref{ST0.5.3.2} at $(\diamondsuit)$; the reader might wish to
refer to Remark~\ref{R0.1025}, where this was made a little more
explicit. Anyway: we obtain the commutative diagram
\[
\xymatrix@C+40pt{
  \mu_X^{}(\id^*\times \id^*)(f^*\times f^!) \ar[r]^-{\mu(X,X)}\ar[d]_{\tau(\diamondsuit)\times\theta(\diamondsuit)} &
  \id^*\mu_X^{}(f^*\times f^!)\ar[r]^-{\s(X,Y)} & \id^*f^!\mu_Y^{}\ar[d]^-{\theta(\diamondsuit)}\\
   \mu_X^{}(\id^*\times \id^!)(f^*\times f^*) \ar[r]^-{\s(X,X)} &
  \id^!\mu_X^{}(f^*\times f^*)\ar[r]^-{\mu(X,Y)} & \id^!f^*\mu_Y^{}
}\]
Recall that $\theta(\diamondsuit):f^!\la f^*$ is the isomorphism
of Remark~\ref{R0.7}, and we discover that, up to this isomorphism, the map
$\s(f):\mu_X^{}(f^*\times f^!)\la f^!\mu_Y^{}$ agrees with
$\mu(f):\mu_X^{}(f^*\times f^*)\la f^*\mu_Y^{}$.

This means that we know $\s(f)$ if $f$ is either a flat monomorphism
or of finite type and universally quasi-proper.
But every morphism in $\seq$
is 2-isomorphic to a 
composite $gf$, with $g$ of finite type and
universally quasi-proper and $f$ a (dominant) flat monomorphism.
Because $\s$ is a natural transformation of 2-functors it respects
composition, and if we write out explicitly what this means we discover
that we can compute $\s(gf)$ in terms of
$\s(f)$ and $\s(g)$.
There is a formula for $\s(gf)$ which we leave to the reader.
\ermk

This concludes the series of remarks in which we showed that the 
compatibilities of \ref{ST0.5.-1} through \ref{ST0.5.999}
force upon us a canonically unique choice for
each of the data specified in  Theorem~\ref{T0.5}~(i)--(vii). The
difficulty,
which will occupy us for most of the paper, is to show
that the recipe that is forced on us works. This means the 
following.

In Remark~\ref{R0.8} we learned that if we
choose a Nagata compactification for a 1-morphism $f\in\seq$,
meaning we write the map $f$ as
$f\cong pu$
with $u$ a dominant, flat monomorphism and $p$ of finite type and
universally quasi-proper, then 
the composite \ref{eq0.1} is an isomorphism which we will write
$\ph(f,pu):f^!\la u^*p^\times$. This means that, if $f\cong pu\cong qv$ are two
Nagata compactifications of $f$,
then we have isomorphisms
\[\xymatrix@C+40pt{
u^*p^\times \ar[r]_{{\ph(f,pu)}^{-1}} \ar@/^2pc/[rr]^{\ph(pu,qv)} & f^! \ar[r]_{\ph(f,qv)} & v^*q^\times
}\]
where $\ph(pu,qv)$ is defined to be the composite. It is automatic
from this definition that the $\ph(pu,qv)$ must compose correctly:
given
three factorizations $f\cong pu\cong qv\cong rw$ we must have
$\ph(qv,rw)\ph(pu,qv)=\ph(pu,rw)$.

What we will do in this paper is introduce
a category of $\nseq(X,Y)$, whose objects are morphisms $f:X\la Y$
together with Nagata compactifications $pu\la f$. There is a
functor $F:\nseq(X,Y)\la\seq(X,Y)$ which forgets the Nagata compactification,
and the main point is that this functor is a groupoid completion. It has the
property that any functor from $\nseq(X,Y)$ to a groupoid factors
uniquely through $F$, and any natural transformation
of functors that factor through $F$ must factor through $F$.
In order to define $(-)^!$ on the
1-category $\seq(X,Y)$ it suffices to produce a functor
$\nseq(X,Y)\la \Tri\big[\Dqc(Y),\Dqc(X)\big]$ which
takes all morphisms in $\nseq(X,Y)$ to isomorphisms
in $\Tri$. And thereafter all our constructions will
be made on categories like $\nseq(X,Y)$, and will be shown to factor
through $F$ by the universal property.

This article is foundational, and one of the points is that the
results could be used as a starting point for developing the theory
of the functor $(-)^!$. Therefore the body of the article, as well as
the entire Introduction up to this point, goes to some
lengths to avoid using any results in the
literature on Grothendieck duality. In what little remains
of the Introduction this will change: we are about to mention
an application and compare to other results in the literature, and
for this there is no virtue in being self-contained.

\rmk{R0.33}
We mention in passing one application.
Suppose $g:Y\la Z$ is a 1-morphism in $\seq$ of finite type and universally
quasi-proper,
$f:X\la Y$ is flat morphism of stacks,
and $gf:X\la Z$ is of finite Tor-dimension. Then
$f^*g^\times :\Dqc(Z)\la\Dqc(X)$ is a bounded functor. This means
there exist integers $A$ and $B$ so that 
\begin{eqnarray*}
\{\ch^i(E)=0\text{ for }i\geq 0\} &\quad\Longrightarrow\quad&
\{\ch^i(f^*g^\times E)=0\text{ for }i\geq A\},\\
\{\ch^i(E)=0\text{ for }i\leq 0\} &\quad\Longrightarrow\quad&
\{\ch^i(f^*g^\times E)=0\text{ for }i\leq B\}
\end{eqnarray*}
where $\ch^i(E)$ means the $i$th cohomology sheaf of the complex
$E$. In 
Hartshorne's old terminology of~\cite{Hartshorne66} bounded functors
would be called ``way out
left'' and ``way out right''.

Now $f^*$ is bounded and $g^\times$ is always bounded below, meaning the second
implication is standard. We have to prove $f^*g^\times$ bounded above. Note 
that the image of $f$ must be contained in the open subset $U\subset Y$ on which
$g$ is of finite Tor-dimension, and hence $f$ factors through the open
immersion $U\hookrightarrow Y$. We may
therefore confine ourselves to the case where $f$ is an open
immersion, in particular $f$ is a morphism in $\seq$. 
But then the map $\rho(f,g):{(gf)}^!\la f^!g^!\cong f^*g^\times$ is an
isomorphism (for example by \ref{S0.5.3.7}(2)),
and we are reduced to proving that $h^!$ is bounded
above for $h:X\la Z$ a 1-morphism in $\seq$ of 
finite type and of finite Tor-dimension. Choose an affine scheme
$Y$\and a faithfully flat map $v:Y\la Z$ and form the pullback
\[\xymatrix{
  W \ar[r]^u\ar[d]_i & X\ar[d]^-h\\
  Y \ar[r]^v & Z
}\]
The map $h$ is of finite Tor-dimension, hence so is its flat base-change
$i:W\la Y$. By \ref{ST0.5.7}(1) the map $\theta:u^*h^!\la i^!v^*$
is an isomorphism. As $u$ is faithfully flat it suffices to
prove that $u^*h^!\cong i^!v^*$ is  bounded above,
and since $v$ is faithfully flat
$v^*$ is bounded, and it suffices to prove that $i^!$ is bounded above.
We are reduced to the case where $Z$ is an affine scheme.

Now choose an affine scheme $W$ and a smooth, surjective map $\ell:W\la X$.
Because $\ell$ is flat \ref{S0.5.3.7}(1) tells us that
$\rho(\ell,h):(h\ell)^!\la\ell^!h^!$ is an isomorphism. We wish to
show that $h^!$ is bounded, but $\ell^!(-)=\omega\oo\ell^*(-)$ where
$\omega$ is some shift of the relative canonical bundle. Therefore
$h^!$ will be bounded if and only if $\ell^!h^!\cong(h\ell)^!$
is---we are reduced
to the case where $h:X\la Z$ is a morphism of affine schemes, still
of finite Tor-dimension and of finite type.

The map $h$ is quasiprojective, meaning it
factors as $X\stackrel j\la \pp_Z^n\stackrel\pi\la Z$ with 
$j$ a locally closed immersion. 
Factor $j$ as $X\stackrel i\la V\stackrel v\la\pp_Z^n$
with $v$ an open immersion and $i$ a closed immersion.
The maps $v$ and $\pi$ are both smooth, and as $h=\pi vi$ is
of finite Tor-dimension it follows that $i$ is
a proper map of finite Tor-dimension.
Because all three maps are of finite Tor-dimension we have that
$f^!=i^!v^!\pi^!$,
and $i^!=i^\times$ and $\pi^!=\pi^\times$
are bounded above by 
\cite[Corollary~4.3.1]{Lipman-Neeman07} while $v^!=v^*$ is obviously
bounded. This completes the proof of our claim.

Thus our techniques allow us to extend to unbounded derived
categories of stacks
the boundedness results of Lipman~\cite[Theorem~4.9.4]{Lipman09}.
The reader is also
referred to Nayak~\cite[Theorem~5.9]{Nayak09}, who
improved on Lipman's results by relaxing Lipman's hypothesis that the maps
must be of finite type. 
Lipman's book~\cite{Lipman09} studies the boundedness of the functor $f^!$ 
extensively. The boundedness of $f^!$ can be seen to be equivalent
to \ref{ST0.5.17}(2), and the
reader might wish to look at Lipman~\cite[Exercise~4.7.6(f)]{Lipman09} for
further equivalent statements. Also interesting is that the boundedness of
$f^!$, for $f$ of finite Tor-dimension,
can be used to prove Theorem~\ref{T0.1}.
As it stands this would be circular since we used 
 Theorem~\ref{T0.1} to prove the boundedness of $f^!$; but it is possible
to give a direct proof of the boundedness, and hence a second proof
of Theorem~\ref{T0.1}. We hope to include this in a subsequent paper.

It should be noted that the direct proof I have of the boundedness of $f^!$ 
still hinges on Lemma~\ref{L4.5}; somehow the key to everything is 
Thomason's localization theorem.
\ermk

\rmk{R0.21}
If we restrict to the 2-subcategory $\mathbb{S}_{\mathbf{f}}\subset\seq$ where
the morphisms are the maps of finite Tor-dimension then the theory simplifies,
and $(-)^!$ becomes a pseudofunctor. On maps of finite Tor-dimension 
$\rho(f,g)$ is an isomorphism. 

If we restrict further to the 1-subcategory
$S_{\mathbf{f}}\subset\mathbb{S}_{\mathbf{f}} $ where the
objects are schemes, one can also proceed as in
Alonso, Jerem\'{\i}as and
Lipman~\cite[Section~5.7]{Alonso-Jeremias-Lipman11}: they
defined $f^!E$ to be $f^*E\oo f^!\co_Y^{}$, which they could do
without
the results of this article because $\co_Y^{}$ belongs to $\Dqcpl(Y)$
and the classical theory (for schemes) defined $f^!$ on
$\Dqcpl(Y)$. By \ref{ST0.5.17} 
the map $\sigma(f,E,\co_Y^{}):f^*E\oo f^!\co_Y^{}\la f^!E$ 
is an isomorphism for $f\in \mathbb{S}_{\mathbf{f}} $, and hence
their functor agrees on $S_{\mathbf{f}} $ 
with the one in this article. Their methods
allowed them to produce an isomorphism ${(gf)}^!\la f^!g^!$ when 
$f,g$ are of finite Tor-dimension, and to
prove associativity as well as a base-change theorem.

The recent work of Avramov, Iyengar, Lipman and 
Nayak~\cite{Avramov-Iyengar08,
Avramov-Iyengar-Lipman-Nayak10} makes it especially interesting to
understand base-change for morphisms of finite Tor-dimension. In the
articles \cite{Avramov-Iyengar08,
Avramov-Iyengar-Lipman-Nayak10} 
the authors find new formulas for $f^!$, which do not involve
compactifications,
but are only valid for morphisms of finite Tor-dimension. It was in trying
to understand this recent work that I came to the theorems in this article.
Since this article is already very long, the discussion of the
relation with the theory of
\cite{Avramov-Iyengar08,
  Avramov-Iyengar-Lipman-Nayak10} is left to~\cite{NeemanTIFR}.
\ermk

\medskip
\nin
{\bf Acknowledgements.}\ \  The author would like to thank Jack Hall,
Joe Lipman and David Rydh, for extensive
comments and helpful suggestions that led
to major improvements on earlier versions of this manuscript.

\section{Notation}
\label{S20}

All our algebraic stacks
 will be assumed quasi-compact and quasi-separated, Furthermore: for most of the
article we'll assume our stacks noetherian---the reason for this restriction 
will be discussed in Example~\ref{E4.7.53}, as well as in 
Remarks~\ref{R4.7.9997} and \ref{R4.10.10.10.3}.

In this entire article we work only in derived categories of the form
$\Dqc(X)$ for algebraic stacks $X$. Therefore all functors 
we consider will be assumed
derived: if $f:X\la Y$ is a morphism of stacks we will 
always write
$f^*$ for $\LL f^*$, $f_*$ will stand for $\R f_*$, the tensor $\oo$ will mean
derived tensor $\oo^\LL_{}$, and $\HHom$ will be the
right adjoint of $\oo$, the internal derived Hom of $\Dqc(X)$. We distinguish 
the internal $\HHom$ from $\Hom$, which takes its values in complexes of
abelian groups. In more traditional notation our $\HHom$ would be written
$\mathbf{q}\RHHom$\footnote{The $\mathbf q$ stands for the
  derived quasi-coherator, the right adjoint to the
  inclusion $\Dqc(X)\la\D\big(\text{Mod(X)}\big)$. As we have already said: the classical construction, using the quasi-coherator, is irrelevant to us. We only want the formal properties of $\HHom$, and can construct it without ever leaving $\Dqc(X)$, using Brown representability to give the existence
  of a right adjoint
  to the tensor product.}
and our $\Hom$ would be the usual $\RHom$.
 The reader will note
that our convention differs 
from much of the classical literature, we never leave the category
$\Dqc(X)$---one
facet is that our $\HHom$ is rarely 
local. On page~\pageref{p25.200}
we present an argument showing that, if $P$ is a perfect complex and $R$
is arbitrary, then $\HHom(P,R)$ is isomorphic to $R\oo P^\vee$ for some
$P^\vee$. From the nature of $P^\vee$ it then follows
that $\HHom(P,R)$ has to be local.
The only other useful theorem I know,
providing a situation in which our $\HHom(P,R)$ can be proved local, is when
$P$ is pseudo-coherent and $R$ is bounded below---this is an old result of
Illusie~\cite[Proposition~3.7]{Illusie71A}, we will meet it again in 
the proof of
Lemma~\ref{L4.-4}.

Suppose $f:X\la Y$ is a concentrated morphism 
of stacks, as in \cite[Definition~2.4]{Hall-Rydh13}.
Then
\cite[Theorem~2.6(3)]{Hall-Rydh13} tells us that
$f_*:\Dqc(X)
\la\Dqc(Y)$
respects
coproducts, while from \cite[Theorem~A.3]{Hall-Neeman-Rydh14}
we know that $\Dqc(X)$ is a well generated triangulated category---hence 
$\Dqc(X)$ satisfies Brown representability by 
\cite[Proposition~8.4.2]{Neeman99} and $f_*$ has a 
right adjoint by \cite[Theorem~8.4.4]{Neeman99}.
The right adjoint of $f_*:\Dqc(X)
\la\Dqc(Y)$ will henceforth be denoted
$f^\times:\Dqc(Y)
\la\Dqc(X)$. 

Still assuming that the map $f$ is concentrated: from
\cite[Definition~2.4]{Hall-Rydh13} it follows that there exists an integer 
$\ell$ with $f_*\Dqc(X)^{\leq n}\subset \Dqc(Y)^{\leq n+\ell}$. 
Hence the adjoint $f^\times$ takes $\big( \Dqc(Y)^{\leq n+\ell}\big)^\perp$
to $\big( \Dqc(X)^{\leq n}\big)^\perp$; in other words
$f^\times\Dqc(Y)^{> n+\ell}\subset \Dqc(X)^{> n}$.
Consequently $f^\times\Dqcpl(Y)\subset\Dqcpl(X)$.

Finally: the category-theoretic conventions are as in MacLane's 
book~\cite{MacLane71}. If $F:\ca\la\cb$, $G,G':\cb\la\cc$ and $H:\cc\la\cd$ are
functors, and $\ph:G\la G'$ is a natural transformation, then
$H\ph F$ is the obvious natural transformation $HGF\la HG'F$. That is,
the natural transformation $H\ph F$ sends an object $a\in\ca$ to the morphism
$H\ph_{Fa}^{}:HGFa\la HG'Fa$, where $Fa\in\cb$ is the image of the object $a\in\ca$ under $F:\ca\la\cb$, the morphism $\ph_{Fa}^{}:GFa\la G'Fa$ is where the 
object $Fa\in\cb$ goes under
the natural transformation $\ph:G\la G'$, and $H\ph_{Fa}^{}$ is the image of this
morphism under $H$.

\section{A review of Thomason's localization theorem}
\label{SThomason}

The idea that working in derived categories can be facilitated by looking 
at objects satisfying finiteness properties is old---it goes back to 
SGA6, where there are several articles by Illusie introducing perfect 
complexes and variants. Illusie works in great generality, dealing with 
derived categories of sheaves of $\co_X^{}$--modules on a ringed
topos $X$, and its various subcategories. In this paper we consider the
special case of an algebraic stack $X$, and the derived category $\Dqc(X)$
of sheaves of $\co_X^{}$--modules with quasi-coherent 
cohomology. The perfect complexes are those which, when you pull them back via
a faithfully flat map $\spec R\la X$, become isomorphic to bounded
complexes of vector bundles on $\spec R$.

The modern approach to Illusie's ideas is from
the homotopy-theoretic perspective---we briefly remind the reader.
Perfect complexes on quasicompact, quasi-separated schemes have two
key properties: they are both compact objects and
strongly dualizable objects
in $\Dqc(X)$. We recall

\dfn{DThom.1}
Let $\ct$ be a triangulated category with coproducts. An object $P\in\ct$ is
\emph{compact} if $\Hom(P,-)$ commutes with coproducts.
\edfn

\dfn{DThom.3}
Let $\ct$ be a symmetric monoidal category---this means $\ct$ has a symmetric
tensor product and this tensor product has a unit $\one$.
 An object $P\in\ct$ is \emph{strongly dualizable} 
if there exists an object $P^\vee\in\ct$ and morphisms $\e:P^\vee\oo P\la\one$
and $\eta:\one\la P\oo P^\vee$ so that the two composites
\[
\CD
P @>\eta\oo\id>> P\oo P^\vee\oo P @>\id\oo\e>> P, 
  @.\qquad\qquad @. P^\vee @>\id\oo\eta>> P^\vee\oo P\oo P^\vee @>\e\oo\id>> P^\vee
\endCD
\]
are identities. 
\edfn

\rmk{RThom.1}
We follow the conventions of homotopy theory, where the objects
of Definition~\ref{DThom.3}  are called
strongly dualizable. There are parts of the literature which refer to the 
same 
objects as \emph{rigid.}
\ermk

If $X$ is a quasi-compact, quasi-separated
stack the perfect complexes in $\Dqc(X)$
are precisely the strongly dualizable objects.
 From
\cite[Lemma~4.4(1)]{Hall-Rydh13} we know that
the compact objects in $\Dqc(X)$ 
are all perfect. For many stacks the inclusion is an equality: 
more precisely this happens exactly for
concentrated stacks. A stack is $X$ is concentrated if the
map $X\la\spec\zz$ is concentrated. All 
tame Deligne-Mumford stacks are concentrated, hence
for those stacks the compact objects and
(strongly dualizable objects)= (perfect complexes) are the
same. Note that our results \emph{do not} assume
that the stacks are concentrated, we do not want to assume tameness.
The morphisms between them will be assumed concentrated, but
not the objects. Hence in our arguments care will be
taken to deal with the case when the compact objects are properly
contained in the perfect complexes.

If $\ct$ is a symmetric monoidal category and $P,Q,R\in\ct$ are objects, with
$P$ strongly dualizable, we have natural maps
\[
\Hom(Q\oo P, R) \stackrel\alpha\la \Hom(Q,R\oo P^\vee),\qquad
\Hom(Q,R\oo P^\vee) \stackrel\beta\la \Hom(Q\oo P, R)
\]
where $\alpha$ takes $\ph:Q\oo P\la R$ to the composite
$Q\stackrel{\id\oo\eta}\la Q\oo P\oo P^\vee \stackrel{\ph\oo\id}\la
R\oo P^\vee$, 
while $\beta$ takes the map $\rho:Q\la R\oo P^\vee$ to the composite
$Q\oo P\stackrel{\rho\oo\id}\la R\oo P^\vee \oo P\stackrel{\id\oo\e}\la
R$.
The reader will easily check that $\beta\alpha=\id$ and $\alpha\beta=\id$,
and since these inverse isomorphisms are natural in $Q$ and $R$ they
give us a canonical isomorphism $\HHom(P,R)\cong R\oo P^\vee$, meaning
that 
$R\oo P^\vee$ satisfies the universal property defining the internal Hom-object
$\HHom(P,R)$.
\label{p25.200}
The special case $R=\one$ tells us that $P^\vee$ is
unique up to canonical isomorphism.
If $f:X\la Y$ is a morphism of stacks then $f^*:\Dqc(Y)\la\Dqc(X)$ is
a strong monoidal functor, meaning it respects the tensor product.
It easily follows that $f^*$ takes
strongly dualizable objects to strongly
dualizable objects, more precisely
$f^*(P^\vee)$ is canonically isomorphic to $(f^*P)^\vee$. And 
the argument above shows that $f^*\HHom(P,R)$ is canonically 
isomorphic to $\HHom(f^*P,f^*R)$. We will often use this.

We recall in this preliminary section
some useful facts, well-known 
for the $f^*$ that arise from morphisms of schemes. The argument below
notes that the statements are formal.

\ntn{NThom.73.63}
Let $f^*:\ca\la\cb$ be an oplax monoidal functor of 
symmetric monoidal 
categories, which means that we are given the data of
\be
\item
 A natural map $\ell:f^*\one_\ca^{}\la\one_\cb^{}$.
\item
For every pair of objects
$A,A'\in\ca$, a natural map 
$\mu=\mu(A,A'):f^*(A\oo A')\la f^*A\oo f^*A'$.
\setcounter{enumiii}{\value{enumi}}
\ee
Furthermore
these maps satisfy the obvious 
compatibilities with the symmetric 
monoidal structures of $\ca$ and $\cb$. 
If $\ell$ and $\mu(A,A')$ are all isomorphisms then
 $f^*$ is called \emph{strong monoidal.}

Next we recall some natural transformations.
\be
\setcounter{enumi}{\value{enumiii}}
\item
  Assume $f^*$ has a right adjoint $f_*$. [If $f:X\la Y$ is
    a morphism of stacks then $f^*:\Dqc(Y)\la\Dqc(X)$ is an
    example.]
The \emph{projection natural transformation $p$} 
is a natural 
transformation between two functors $\ca\times\cb\la \ca$.
It takes a pair of objects
$A\in\ca$ and $B\in\cb$ to a map $p=p(f,A,B):A\oo f_*B \la f_*(f^*A\oo B)$.
The map $p$ corresponds under the adjunction $f^*\dashv f_*$ to the
composite
\[
\CD
f^*(A\oo f_*B) @>\mu>> f^*A\oo f^*f_*B @>\id\oo\e'>> f^*A\oo B \ ,
\endCD
\]
where $\e':f^*f_*\la\id $ is the counit of the adjunction $f^*\dashv f_*$.
\item
Assume that $f^*$ is such that the projection
natural transformation $p$ in (iii) is an
isomorphism---in this case we say that the \emph{projection formula}
holds. Suppose further that $f_*$ has a right adjoint $f^\times$.
[An $f^*$ coming
  from a concentrated morphism of algebraic stacks $f:X\la Y$ is an example.
  We saw the existence on $f^\times$ in \S\ref{S20}, and the projection formula
  holds
by \cite[Corollary~4.12]{Hall-Rydh13}\footnote{Another large
class of $f^*$
for which the projection formula holds may be found in
\protect{\cite[Proposition~2.15]{Balmer-Dellambrogio-Sanders16}}.
There is overlap with \protect{\cite[Corollary~4.12]{Hall-Rydh13}}, but 
to apply \protect{\cite[Proposition~2.15]{Balmer-Dellambrogio-Sanders16}}
to morphisms of stacks $f:X\la Y$ 
it isn't enough to assume $f$ concentrated,
we must also suppose $\Dqc(Y)$ compactly generated.
}.]
Suppose we are
given a pair of objects $E,F\in\ca$;
we define $\xi(f,E,F)$ to be the composite
\[
\xymatrix@C+30pt{
f_*^{}\left[f^*E\oo f^\times F\right] \ar[r]_-{p(f,E,f^\times F)^{-1}}
\ar@/^2pc/[rr]^{\xi(f,E,F)} &
E\oo f_*^{} f^\times F \ar[r]_-{1\oo \e'' } & E\oo F
}
\]
where the first map is the inverse of
the projection formula map $p$ of (iii),
while
$\e'':f_*^{}f^\times\la \id$ is the counit of the adjunction 
$f_*^{}\dashv f^\times$.
By the adjunction $f_*\dashv f^\times$
the map
$\xi(f,E,F)$ corresponds to a map
\[
\xymatrix@C+50pt{
f^*E\oo f^\times F 
\ar[r]^-{\chi(f,E,F)}& f^\times (E\oo F).
}
\]
\ee
The 
maps $\xi(f,E,F)$ and  
$\chi(f,E,F)$ are natural in  $E$ and $F$, and there is a sense in
which they are also natural in $f$.
In \S\ref{S805} we will meet a fancy formulation of this
naturality.
\entn

\pro{PThom.73.66}
Let the situation be as in Notation~\ref{NThom.73.63}(iv), in particular
the projection formula holds for $f^*$. Assume
further that $f^*$ is strong monoidal (not only lax monoidal). 
If $E$ is strongly dualizable then the map 
$\chi:f^*E\oo f^\times F\la f^\times(E\oo F)$ is an isomorphism.
\epro

\prf
Because $E$ is strongly dualizable there exists an object
$E^\vee\in\ca$ with morphisms $\e:E^\vee\oo E\la\one_\ca^{}$ and
$\eta:\one_\ca^{}\la E\oo E^\vee$ as in Definition~\ref{DThom.3}, and because
$f^*$ is strong monoidal the morphisms $f^*\e:f^*E^\vee\oo f^*E\la\one_\cb^{}$ and
$f^*\eta:\one_\cb^{}\la f^*E\oo f^*E^\vee$ also satisfy the hypotheses of 
Definition~\ref{DThom.3}.
Let $B\in\cb$ be arbitrary, and observe the following isomorphisms
\begin{eqnarray*}
\Hom_\cb^{}(B\,\,,\,\,f^*E\oo f^\times F) &\cong & \Hom_\cb^{}(f^*E^\vee\oo B\,\,,\,\,f^\times F) \\
&\cong& \Hom_\ca\big(f_*(f^*E^\vee\oo B)\,\,,\,\,F\big) \\
&\cong& \Hom_\ca(E^\vee\oo f_*B\,\,,\,\,F)\\
&\cong& \Hom_\ca(f_*B\,\,,\,\,E\oo F)\\
&\cong& \Hom_\cb\big(B\,\,,\,\,f^\times(E\oo F)\big)\ .
\end{eqnarray*}
The first and fourth isomorphisms come because $E$ and $f^*E$ are strongly
dualizable, see the inverse
isomorphisms $\alpha$ and $\beta$ of the paragraph before
Notation~\ref{NThom.73.63}.
The second and fifth isomorphisms are by the adjunction $f_*\dashv
f^\times$, and the third isomorphism is by the projection formula.

Of course it remains to show that this composite of isomorphisms agrees
with $\Hom(B,-)$ applied to the morphism $\chi:f^*E\oo f^\times F\la
f^\times (E\oo F)$. Since the maps are all explicitly described in
Notation~\ref{NThom.73.63}
this amounts to a straightforward diagram chase. Unfortunately I
couldn't make this large diagram fit on the page, we leave it to
the reader.
\eprf

One of the natural questions is whether it is possible to extend a perfect
complex from an open substack of a stack to the whole stack. An old theorem
of Serre says that any coherent sheaf on an open subset of a noetherian scheme
$X$ extends to all of $X$, 
and perfect complexes can be thought of as
an analog of coherent sheaves. The question might be natural enough, but the 
answer is No---there exist noetherian schemes $X$, open subsets $U\subset X$,
and perfect complexes on $U$ that do not extend. There is a $K$--theoretic
obstruction: a perfect complex $P$ on $X$ defines a class $[P]\in K_0(X)$,
and the class in $K_0(U)$ of the restriction of $P|_U^{}$ is the image
of  $[P]\in K_0(X)$ under the restriction map $K_0(X)\la K_0(U)$. If $X$
is singular it is possible for the map $K_0(X)\la K_0(U)$
not to be surjective, in which
case there will exist perfect complexes on $U$ that cannot possibly extend.

This much was known early on. Until Thomason and
Trobaugh~\cite{ThomTro} no one realized that this is
the only obstruction. In fact \cite{ThomTro} proved more: 
if $X$ is a quasi-compact, semi-separated
scheme, $U\subset X$ is a quasi-compact open subset with
$j:U\la X$ the immersion, $P\in\Dqc(U)$ 
is a perfect complex with $[P]\in K_0(U)$ lying
in the image of the map $K_0(X)\la K_0(U)$, and if $f:P\la j^*G$
is any morphism in $\Dqc(U)$, then there is a perfect complex $P'\in\Dqc(X)$
and a morphism $f':P'\la G$ so that $f$ is isomorphic to $j^*f'$.

The proof in \cite{ThomTro} is algebro-geometric, very much in the 
spirit of SGA6. In \cite{Neeman92A} there is a short, sweet proof of a much
more general theorem, using homotopy-theoretic techniques. 
The general version of Thomason's localization theorem, and its proofs,
can be formalized as facts about compact objects in triangulated
categories. To minimize the number of
abstract definitions introduced, let us not state the 
most general version
here---we 
mention only that this general theorem applies to many quasi-compact, 
quasi-separated stacks, and this underpins our argument.

Before we leave the
 realm of abstraction there is one more
concept we will use.

\dfn{DThom.5}
Let $\ct$ be a triangulated category. We say that $\ct$ is 
\emph{compactly generated} if it has
coproducts, and there is a set $\cg$ of compact
objects in $\ct$ satisfying the following equivalent conditions:
\be
\item
If $\Hom(G,X)=0$ for all $G\in\cg$ then $X=0$.
\item
Any localizing subcategory of $\ct$ containing $\cg$ must be all of $\ct$.
[A full subcategory $\cs\subset\ct$ is \emph{localizing} if it contains
zero and is closed under coproducts and triangles].
\ee
\edfn
 
\nin
The equivalence of (i) and (ii) is not supposed to be obvious, but it
is a special case of \cite[Proposition~8.4.1]{Neeman99}. A 
set $\cg$ satisfying the equivalent
conditions of Definition~\ref{DThom.5} is called a \emph{set of 
compact generators.}
The categories $\Dqc(X)$ are known to be compactly generated for
all quasi-compact, quasi-separated algebraic spaces, and also for a large
class of stacks. As we will see in Sections~\ref{S4} and \ref{S901},
what makes perfect complexes so useful
is the fact that
they  contain a set of compact
generators.

\section{A review of formal nonsense concerning base-change}
\label{S27}

Base-change maps will play a key role in this article, so we include here
a little reminder of the formal theory.
We give a fairly abstract treatment: let $\ca$
be a 2-subcategory of the category of categories which contains
all natural transformations between its 1-morphisms---for example $\ca$
could be the 2-category
 whose objects are triangulated
categories,
the 1-morphisms are triangulated functors, and the 2-morphisms are any
natural transformations.
\dfn{D27.1}
Suppose we are given in $\ca$ a 2-commutative square
\[
\xymatrix@C+15pt@R+15pt{
  \cw  & \cx\ar[l]_-{u^*} \\
  \cy\ar@{=>}[ur]^-{\tau} \ar[u]^-{f^*} & \cz \ar[l]^-{v^*} \ar[u]_-{g^*}
}\]
This means that we are given in $\ca$ the 1-morphisms $u^*:\cx\la\cw$,
$g^*:\cz\la\cx$, $f^*:\cy\la\cw$
and $v^*:\cz\la\cy$, as well as the natural transformation
$\tau:f^*v^*\la u^*g^*$. Suppose also that $f^*$ and $g^*$ have
in $\ca$ right adjoints
$f_*$ and $g_*$ respectively.
The induced \emph{base change map} is the 2-commutative square
\[
\xymatrix@C+15pt@R+15pt{
  \cw\ar[d]_-{f_*}  & \cx\ar[l]_-{u^*}\ar[d]^-{g_*} \\
  \cy  & \cz \ar[l]^-{v^*} \ar@{=>}[ul]_-{\beta}
}\]
where $\beta:v^*g_*\la f_*u^*$ is defined to be the map corresponding
under the adjunction $f^*\dashv f_*$ to the composite
\[
\CD
f^*v^*g_* @>\tau g_*>> u^*g^*g_* @>u^*\e>> u^*\ ,
\endCD
\]
where $\e:g^*g_*\la\id$ is the counit of the adjunction $g^*\dashv g_*$.
\edfn

\rmk{R27.2}
Recall that the adjunction $f^*\dashv f_*$ takes a map $\beta:x\la f_* y$
to the composite $f^*x\stackrel{f^*\beta}\la f^*f_* y\stackrel{\e'}\la y$,
with $\e':f^*f_*\la\id$ the counit of the adjunction.
Applying this to the $\beta$ of Definition~\ref{D27.1}, we deduce the
commutativity
of the square
\[
\xymatrix{
  f^*v^*g_* \ar[r]^-{f^*\beta}\ar[d]_-{\tau g_*} & f^*f_*u^*\ar[d]^-{\e'u^*}\\
  u^*g^*g_* \ar[r]^-{u^*\e} & u^*
}\]
\ermk

\dfn{D27.3}
Let the notation be as in Definition~\ref{D27.1}, but assume further:
\be
\item
  The map $\beta:v^*g_*\la f_*u^*$ is an isomorphism.
\item
  The functors $f_*$ and $g_*$ both have right adjoints, denoted
  $f^\times$ and $g^\times$ respectively.
  \ee
Then the base-change map
  $\Phi:u^*g^\times\la f^\times v^*$ is obtained by applying
Definition~\ref{D27.1} to the 2-commutative square
\[
\xymatrix@C+15pt@R+15pt{
  \cy   & \cz \ar[l]_-{v^*}\\
  \cw\ar[u]^-{f_*}\ar@{=>}[ur]|-{\beta^{-1}} & \cx\ar[l]^-{u^*}\ar[u]_-{g_*} 
}\]
Concretely: the map $\Phi:u^*g^\times\la f^\times v^*$ corresponds
under the adjunction $f_*\dashv f^\times$ to the composite
\[
\CD
f_*u^*g^\times @>\beta^{-1}g^\times>> v^*g_*g^\times @>v^*\wt\e>> v^*\ ,
\endCD
\]
where $\wt\e:g_*g^\times\la\id$ is the counit of the adjunction
$g_*\dashv g^\times$.
\edfn

\rmk{R27.4}
Applying Remark~\ref{R27.2} to the base-change map of Definition~\ref{D27.3}
we obtain the commutativity of the square
\[
\xymatrix{
  f_*u^*g^\times \ar[r]^-{f_*\Phi}\ar[d]_-{\beta^{-1} g^\times} & f_*f^\times v^*\ar[d]^-{\widehat\e v^*}\\
  v^*g_*g^\times \ar[r]^-{v^*\wt\e} & v^*
}\]
where $\widehat\e:f_*f^\times\la\id$ is the counit of the adjunction
$f_*\dashv f^\times$.
\ermk

\exm{E27.5}
Suppose we are given a 2-commutative square of
quasicompact, quasiseparated algebraic stacks
\[
\xymatrix{
  W \ar[r]^-{u}\ar[d]_-{f}\ar@{}[dr]|{\diamondsuit}  & X\ar[d]^-{g} \\
  Y\ar[r]_-{v} & Z
}\]
Applying $(-)^*$ we obtain a 2-commutative square
\[
\xymatrix@C+15pt@R+15pt{
  \Dqc(W)  & \Dqc(X)\ar[l]_-{u^*} \\
  \Dqc(Y)\ar@{=>}[ur]^-{\tau} \ar[u]^-{f^*} & \Dqc(Z)\ar[l]^-{v^*} \ar[u]_-{g^*}
}\]
with
$\tau:f^*v^*\la u^*g^*$
an isomorphism. The maps $f^*$ and $g^*$ have
right adjoints $f_*$ and $g_*$, and Definition~\ref{D27.1} gives
the base-change
square
\[
\xymatrix@C+15pt@R+15pt{
  \Dqc(W)\ar[d]_-{f_*}  & \Dqc(X)\ar[l]_-{u^*}\ar[d]^-{g_*} \\
  \Dqc(Y)  & \Dqc(Z) \ar[l]^-{v^*} \ar@{=>}[ul]_-{\beta}
}\]
This map $\beta$ need not always be an isomorphism, but 
\cite[Theorem~2.6(4)]{Hall-Rydh13} tells us that, as long as
\be
\item
  The square $\diamondsuit$ is 2-cartesian,
\item
  The morphism $v$ is flat while $g$ is concentrated,
\ee
then $\beta$ is indeed an isomorphism. Since $g$ is concentrated
\cite[Lemma~2.5(1)]{Hall-Rydh13} guarantees that so is its pullback
$f$, and 
in \S\ref{S20}
we saw that when $f$, $g$ are concentrated then $f_*$, $g_*$ have right
adjoints $f^\times$, $g^\times$. We are in the situation
of Definition~\ref{D27.3}, and have the base-change
2-commutative square
\[
\xymatrix@C+15pt@R+15pt{
  \Dqc(W)  & \Dqc(X)\ar[l]_-{u^*}\ar@{=>}[dl]_-{\Phi} \\
  \Dqc(Y)\ar[u]^-{f^\times}  & \Dqc(Z)\ar[u]_-{g^\times} \ar[l]^-{v^*} 
}\]
\eexm

Back to gorgeous generality. In Notation~\ref{NThom.73.63} we learned
that, if $f$ is an oplax monoidal functor, then one can construct
a projection natural transformation $p:A\oo f_*B\la f_*(f^*A\oo B)$. If
$p$ is an isomorphism one can go further and construct a morphism
$\chi:f^*E\oo f^\times F\la f^\times(E\oo F)$. Next we recall how these
are compatible with the base-change maps above; the next Lemma
is virtually the same as \cite[Proposition~3.7.3]{Lipman09}.

\lem{L27.7}
Suppose we are given in $\ca$ a 2-commutative square
\[
\xymatrix@C+15pt@R+15pt{
  \cw  & \cx\ar[l]_-{u^*} \\
  \cy\ar@{=>}[ur]^-{\tau} \ar[u]^-{f^*} & \cz \ar[l]^-{v^*} \ar[u]_-{g^*}
}\]
Suppose the categories $\cw$, $\cx$, $\cy$ and $\cz$ are all symmetric
monoidal, and the functors $f^*$, $g^*$, $u^*$ and $v^*$
are all oplax monoidal. Suppose
further that the map $\tau:f^*v^*\la u^*g^*$ is compatible with the oplax
monoidal structure maps $\mu$, in the sense that the following
hexagon commutes
\[
\xymatrix@C+20pt{
  f^*v^*(A\oo B) \ar[d]_{\tau}\ar[r]^-{f^*\mu_v^{}} & f^*(v^*A\oo v^*B) \ar[r]^-{\mu_f^{}} & f^*v^*A\oo f^*v^*B\ar[d]^{\tau\oo\tau}\\
u^*g^*(A\oo B)  \ar[r]^-{u^*\mu_g^{}} & u^*(g^*A\oo g^*B) \ar[r]^-{\mu_u^{}} & u^*g^*A\oo u^*g^*B
}\]
Finally, assume also that $f^*$ and $g^*$ have right adjoints $f_*$ and $g_*$
in $\ca$.

Then, for any objects $A\in\cz$ and $B\in\cx$, the following diagram
commutes
\[\xymatrix@C+5pt{
  v^*(A\oo g_*B)   \ar[r]^-{\mu_v^{}} \ar[dd]^-{v^*p(g,A,B)} & v^*A\oo v^*g_*B\ar[r]^-{\id\oo\beta} & v^*A\oo f_*u^*B \ar[dr]|-{p(f,v^*A,u^*B)}\\
  & & & f_*(f^*v^*A\oo u^*B) \ar[dl]^-{f_*(\tau\oo\id)}\\
  v^*g_*(g^*A\oo B) \ar[r]^-{\beta} & f_*u^*(g^*A\oo B) \ar[r]^-{f_*\mu_u^{}} &
  f_*(u^*g^*A\oo u^*B)
  }\]
\elem

\prf
We have to prove the equality of two composites of the form $x\la f_*y$,
and the adjunction $f^*\dashv f_*$ reduces us to proving the equality
of the corresponding composites $f^*x\la y$. Concretely, we need to
prove the commutativity of
\[\xymatrix@C+5pt{
  f^*v^*(A\oo g_*B)   \ar[r]^-{f^*\mu_v^{}} \ar[d]_-{f^*v^*p(g,A,B)} & f^*(v^*A\oo v^*g_*B)\ar[r]^-{f^*(\id\oo\beta)} & f^*(v^*A\oo f_*u^*B) \ar[r]^-{\mu_f^{}}\ar[dr]&
 f^*v^*A\oo f^*f_*u^*B\ar[d]^-{\id\oo\e'} \\
f^*v^*g_*(g^*A\oo B)\ar[d]_-{\tau}\ar[dr]  & & & f^*v^*A\oo u^*B \ar[dl]^-{\tau\oo\id}\\
 u^*g^*g_*(g^*A\oo B)  \ar[r]^-{u^*\e} & u^*(g^*A\oo B) \ar[r]^-{\mu_u^{}} &
  u^*g^*A\oo u^*B
}\]
where the commutative triangles are simply places where we have replaced the
diagonal 
morphism $x\la f_*y$ by the corresponding $f^*x\la y$, as in the definitions.
The naturality of $\tau$ and $\mu_f^{}$ gives the commutativity of the
squares in the diagram
\[\xymatrix@R+15pt@C+5pt{
  f^*v^*(A\oo g_*B) \ar[dr]^-{\tau}  \ar[r]^-{f^*\mu_v^{}} \ar[d]|-{f^*v^*p(g,A,B)} & f^*(v^*A\oo v^*g_*B)\ar[r]^-{f^*(\id\oo\beta)}\ar[dr]|-{\mu_f^{}} & f^*(v^*A\oo f_*u^*B) \ar[r]^-{\mu_f^{}}&
 f^*v^*A\oo f^*f_*u^*B\ar[d]^-{\id\oo\e'} \\
 f^*v^*g_*(g^*A\oo B)\ar[d]_-{\tau}  & u^*g^*(A\oo g_*B)
 \ar[dl]|-{u^*g^*p(g,A,B)}&  f^*v^*A\oo f^*v^*g_*B\ar[ur]|-{\id\oo f^*\beta}& f^*v^*A\oo u^*B \ar[dl]^-{\tau\oo\id}\\
 u^*g^*g_*(g^*A\oo B)  \ar[r]^-{u^*\e} & u^*(g^*A\oo B) \ar[r]^-{\mu_u^{}} &
  u^*g^*A\oo u^*B
}\]
Hence what remains is to prove the commutativity of the perimeter in the
diagram
\[\xymatrix@R+10pt@C+5pt{
  f^*v^*(A\oo g_*B) \ar[d]_-{\tau}  \ar[r]^-{f^*\mu_v^{}}  & f^*(v^*A\oo v^*g_*B)\ar[r]^-{\mu_f^{}} &f^*v^*A\oo f^*v^*g_*B\ar[r]^-{\id\oo f^*\beta}\ar[d]_-{\id\oo\tau}& 
 f^*v^*A\oo f^*f_*u^*B\ar[d]^-{\id\oo\e'} \\
  u^*g^*(A\oo g_*B)
  \ar[d]|-{u^*g^*p(g,A,B)}\ar[r]^-{u^*\mu_g^{}}&
  u^*(g^*A\oo g^*g_*B)\ar[d]^-{u^*(\id\oo\e)}& f^*v^*A\oo u^*g^*g_*B\ar[r]^-{\id\oo u^*\e} \ar@{}[ur]|-{(\clubsuit)}& f^*v^*A\oo u^*B \ar[dl]^-{\tau\oo\id}\\
 u^*g^*g_*(g^*A\oo B)  \ar@{}[ur]|-{(\spadesuit)}\ar[r]^-{u^*\e} & u^*(g^*A\oo B) \ar[r]^-{\mu_u^{}} &
  u^*g^*A\oo u^*B
}\]
The square $(\clubsuit)$ is the tensor product of $f^*v^*A$ with the commutative
square of Remark~\ref{R27.2}, hence commutes. The square $(\spadesuit)$
commutes because it is $u^*$ of the square spelling
out the definition of the morphism $g^*(A\oo g_*B)\la g^*A\oo B$
corresponding to the map $p_g^{}:A\oo g_*B\la g_*(g^*A\oo B)$; see
Notation~\ref{NThom.73.63}(iii).

It remains to show the commutativity of the perimeter in the diagram
\[\xymatrix@C+5pt{
  f^*v^*(A\oo g_*B) \ar[d]_-{\tau}  \ar[r]^-{f^*\mu_v^{}}  & f^*(v^*A\oo v^*g_*B)\ar[r]^-{\mu_f^{}} &f^*v^*A\oo f^*v^*g_*B\ar[r]^-{\id\oo\tau}&
  f^*v^*A\oo u^*g^*g_*B\ar[d]^-{\id\oo u^*\e} \ar[dl]|-{\tau\oo\id}\ar@{}[ddl]|-{(\diamondsuit)}
  \\
  u^*g^*(A\oo g_*B)
  \ar[r]^-{u^*\mu_g^{}}&
  u^*(g^*A\oo g^*g_*B)\ar[d]_-{u^*(\id\oo\e)}\ar[r]^-{\mu_u^{}}&  u^*g^*A\oo u^*g^*g_*B\ar[d]|-{\id\oo u^*\e}& f^*v^*A\oo u^*B \ar[dl]^-{\tau\oo\id}\\
 & u^*(g^*A\oo B)\ar@{}[ur]|-{(\heartsuit)} \ar[r]^-{\mu_u^{}} &
  u^*g^*A\oo u^*B
}\]
The square $(\diamondsuit)$ obviously commutes, the square $(\heartsuit)$
commutes by the naturality of $\mu_u^{}$, and the heptagon commutes
because $\tau$ is compatible with the oplax monoidal
structures of $v^*f^*$ and $u^*g^*$,
see the hypotheses of the Lemma.
\eprf

We follow with another relatively standard fact, see
for example \cite[Exercise~4.7.3.4(c)]{Lipman09}.

\lem{L27.9}
Let the notation be as in Lemma~\ref{L27.7}, but assume in addition
that
\be
\item
  The map $\tau:f^*v^*\la u^*g^*$ is an isomorphism.
\item
  The functors $f_*$ and $g_*$ have right adjoints $f^\times$ and $g^\times$.
\item
  The base-change map $\beta:v^*g_*\la f_*u^*$ is invertible, and hence
  the base-change map $\Phi:u^*g^\times\la f^\times v^*$ is given
  as in Definition~\ref{D27.3}.
\item
  The projection formula holds for $f^*$ and $g^*$. 
\ee
Then, for any pair of objects $E,F\in\cz$, the following
square commutes
\[
\xymatrix@C+30pt{
  u^*g^*E \oo u^*g^\times F \ar[d]_-{\tau^{-1}\oo\Phi} &
    u^*(g^*E\oo g^\times F)\ar[l]_-{\mu_u^{}}
    \ar[r]^-{u^*\chi(g,E,F)}&  u^*g^\times(E\oo F)\ar[d]^-{\Phi} \\
    f^*v^*E\oo f^\times v^* F\ar[r]^-{\chi(f,v^*E,v^*F)} &
    f^\times(v^*E\oo v^*F) &
    f^\times v^*(E\oo F)\ar[l]_-{f^\times\mu_v^{}}
}\]
\elem

\prf
We need to show the equality of two maps $x\la f^\times y$, and it suffices
to prove the equality of the corresponding maps $f_*x\la y$. Concretely this
means we must show the commutativity
of the perimeter of the diagram
\[
\xymatrix@C+30pt{
  f_*(u^*g^*E \oo u^*g^\times F) \ar[dd]_-{f_*(\tau^{-1}\oo\Phi)} &
    f_*u^*(g^*E\oo g^\times F)\ar[l]_-{f_*\mu_u^{}}
    \ar[r]^-{f_*u^*\chi(g,E,F)}\ar[d]_-{\beta^{-1}}&  f_*u^*g^\times(E\oo F)\ar[d]^-{\beta^{-1}g^\times} \\
    & v^*g_*(g^*E\oo g^\times F)\ar[r]^-{v^*g_*\chi(g,E,F)}\ar[dr]|-{v^*\xi(g,E,F)}& v^*g_*g^\times(E\oo F)\ar[d]^-{v^*\wt\e}\ar@{}[ul]|-{(\heartsuit)} \\
    f_*(f^*v^*E\oo f^\times v^* F)\ar[r]_-{\xi(f,v^*E,v^*F)} &
    v^*E\oo v^*F &
    v^*(E\oo F)\ar[l]^-{\mu_v^{}}
}\]
In passing from the diagram of the Lemma to the perimeter
of the diagram above
the map $\chi(f,v^*E,v^*F):
f^*v^*E\oo f^\times v^*F\la f^\times v^*(E\oo F)$ was replaced by
the corresponding morphism
$\xi(f,v^*E,v^*F):f_*(f^*v^*E\oo f^\times v^*F\la
v^*(E\oo F)$,
see Notation~\ref{NThom.73.63}(iv), while the vertical map on the right
$\Phi:u^*g^\times (E\oo F)\la f^\times u^*(E\oo F)$ was replaced by
the composite $f_*u^*g^\times (E\oo F)\la v^*(E\oo F)$
as in Definition~\ref{D27.3}. Now the square $(\heartsuit)$ commutes by
the naturality of $\beta$, and the triangle commutes just because
$\xi(g,E,F):g_*(g^*E\oo g^\times F)\la E\oo F$ is the map
corresponding to $\chi(g,E,F):g^*E\oo g^\times F\la g^\times(E\oo F)$
under the adjunction $g_*\dashv g^\times$. It suffices therefore to prove the
commutativity of the remaining hexagon. If we expand out the maps
$\xi(g,E,F)$ and $\xi(f,v^*E,v^*F)$ as in Notation~\ref{NThom.73.63}(iv)
we are left with proving the commutativity of the perimeter of
\[
\xymatrix@C+30pt{
  f_*(u^*g^*E \oo u^*g^\times F) \ar[d]_-{f_*(\tau^{-1}\oo\Phi)} &
    f_*u^*(g^*E\oo g^\times F)\ar[l]_-{f_*\mu_u^{}}
    \ar[r]^-{\beta^{-1}}& v^*g_*(g^*E\oo g^\times F)\ar[d]^-{v^*p(g,E,g^\times F)^{-1}} \\
    f_*(f^*v^*E\oo f^\times v^* F)\ar[d]_-{p(f,v^*E,f^\times v^*F)^{-1}}  &
    v^*E\oo v^*g_*g^\times F\ar[d]_-{\id\oo v^*\wt\e}&\ar[l]_-{\mu_v^{}} v^*(E\oo g_*g^\times F)\ar[d]^-{v^*(\id\oo\wt\e)}\ar@{}[dl]|-{(\diamondsuit)} \\
    v^*E\oo f_*f^\times v^* F\ar[r]_-{\id\oo\widehat\e} &
    v^*E\oo v^*F &
    v^*(E\oo F)\ar[l]^-{\mu_v^{}}
}\]
The square $(\diamondsuit)$ commutes by the naturality of $\mu_v^{}$, hence
we are reduced to proving the commutativity of the perimeter of
\[
\xymatrix@C+30pt{
  f_*(u^*g^*E \oo u^*g^\times F) \ar[d]_-{f_*(\tau^{-1}\oo\id)} &
    f_*u^*(g^*E\oo g^\times F)\ar[l]_-{f_*\mu_u^{}}
    \ar[r]^-{\beta^{-1}}& v^*g_*(g^*E\oo g^\times F)\ar[dd]^-{v^*p(g,E,g^\times F)^{-1}} \\
 f_*(f^*v^*E \oo u^*g^\times F)\ar[d]_-{f_*(\id\oo\Phi)} \ar[dr]^-{p(f,v^*E,u^*g^\times F)^{-1}}  &  & \\
    f_*(f^*v^*E\oo f^\times v^* F)\ar[d]_-{p(f,v^*E,f^\times v^*F)^{-1}}  \ar@{}[r]|-{(\spadesuit)}&
  v^*E\oo f_*u^*g^\times F\ar[ld]_-{\id\oo f_*\Phi}\ar[dr]^-{\id\oo\beta^{-1}g^\times}\ar@{}[d]|-{(\clubsuit)} &\ar[d]_-{\mu_v^{}} v^*(E\oo g_*g^\times F) \\
    v^*E\oo f_*f^\times v^* F\ar[r]_-{\id\oo\widehat\e} &
    v^*E\oo v^*F &
    v^*E\oo v^*g_*g^\times F\ar[l]^-{\id\oo v^*\wt\e}
}\]
The square $(\spadesuit)$ commutes by the naturality of $p$, while
$(\clubsuit)$ is the tensor product of $v^*E$ with the commutative square
of Remark~\ref{R27.4}. Finally the remaining diagram is the commutative
heptagon of Lemma~\ref{L27.7}, for the pair of objects
$E\in\cz$ and $g^\times F\in\cx$.
\eprf

\rmk{R27.11}
The example we care about is Example~\ref{E27.5}, where we start with a
2-commutative square of algebraic stacks
\[
\xymatrix{
  W \ar[r]^-{u}\ar[d]_-{f}\ar@{}[dr]|{\diamondsuit}  & X\ar[d]^-{g} \\
  Y\ar[r]_-{v} & Z
}\]
and apply $\Dqc(-)$ to obtain a 2-commutative square
\[
\xymatrix@C+15pt@R+15pt{
  \Dqc(W)  & \Dqc(X)\ar[l]_-{u^*} \\
  \Dqc(Y)\ar@{=>}[ur]^-{\tau} \ar[u]^-{f^*} & \Dqc(Z)\ar[l]^-{v^*} \ar[u]_-{g^*}
}\]
The hypotheses of Lemma~\ref{L27.7} are always satisfied.
The hypotheses of Lemma~\ref{L27.9} hold as long as
the square $(\diamondsuit)$ is 2-cartesian, $g$ is concentrated and $v$ is flat.
\ermk

\section{An improved flat base-change theorem}
\label{S4}

This section is devoted to the proof of
Theorem~\ref{T0.1}=Theorem~\ref{T4.13} and proceeds by a series of lemmas.
If $\alpha:V\la Z$ is an open immersion of 
noetherian schemes, the first lemma is easy and can be 
proved in many ways. We want the non-noetherian version, and to state
it we begin with

\rmd{R4.-7}
Let $Z$ be an algebraic stack. 
 An object $D\in\Dqc(Z)$ is called 
\emph{pseudo-coherent} if there exists a flat cover $f:\spec R\la Z$
so that $f^*D$ is isomorphic in $\Dqc\big(\spec R\big)\cong \D(R)$ to a
bounded-above complex of finitely-generated, projective $R$--modules.

One small historical 
comment: if the reader compares the above
with Illusie's definition it's obvious
that a pseudo-coherent complex in our sense is also 
pseudo-coherent in Illusie's. For the converse we need to show that, if a 
complex in $\Dqc\big(\spec R\big)$ is pseudo-coherent in Illusie's sense,
then it is isomorphic in $\Dqc\big(\spec R\big)$ to a 
bounded-above complex of finitely-generated, projective $R$--modules.
For this we need the equivalence $\D(R)\cong\Dqc\big(\spec R\big)$ of
\cite[Theorem~5.1]{Bokstedt-Neeman93}, which wasn't known at the time Illusie 
wrote his expos\'es in SGA6. The point is that finite generation of an 
 $R$--module is local in the flat topology, and hence the top cohomology
module of a complex in $\Dqc\big(\spec R\big)$, 
which is pseudo-coherent in Illusie's sense, must be finitely generated.
But now an easy induction produces a resolution for the complex
as in our definition.
\ermd

\lem{L4.-4}
Let $Z$ be an algebraic stack and let  $\alpha:V\la Z$ be a
flat morphism. Suppose 
$D\in\Dqcb(Z)$ is pseudo-coherent, and assume $\alpha^*D=0$. 
Then, for any $E\in\Dqc(Z)$, we have 
$\alpha^*\HHom_Z^{}(D,E)=0$.
\elem

\prf
The hypotheses tell us that $D$ is pseudo-coherent and belongs to
$\Dqcpl(Z)$; hence we may apply  \cite[Proposition~3.7]{Illusie71A}
with $F=G=D$ to deduce
 that the natural map
$\alpha^*\HHom(D,D)\la\HHom(\alpha^*D, \alpha^*D)$
is an isomorphism. Since $\alpha^*D=0$ we conclude that
$\alpha^*\HHom(D,D)=0$.

But now $\HHom(D,D) $ is a monoid in the monoidal category $\Dqc(Z)$
and $\HHom(D,E)$ is a module over it---there are canonical maps
\[
\CD
\co_Z^{}\oo \HHom(D,E) @>i\oo\id>> \HHom(D,D)\otimes \HHom(D,E) @>\mu>> \HHom(D,E)
\endCD
\] 
which compose to an isomorphism. Applying $\alpha^*$ we have that the
composite
\[
\CD
\co_V^{}\oo \alpha^*\HHom(D,E) @>(\alpha^*i)\oo\id>> \alpha^*\HHom(D,D)\otimes \alpha^*\HHom(D,E) @>\alpha^*\mu>> \alpha^*\HHom(D,E)
\endCD
\] 
is also an isomorphism. Hence $\alpha^*\HHom(D,E)$ is a direct summand
of the middle term, which vanishes since $\alpha^*\HHom(D,D)=0$.
\eprf

The next easy Lemma (in the case of schemes) 
may be found in \cite[1.6.1]{Iyengar-Lipman-Neeman13};
we include the proof for the convenience of the reader.

\lem{L4.97}
Let $f:X\la Y$ be a concentrated morphism of quasi-compact, quasi-separated
algebraic stacks. Suppose $f_*:\Dqc(X)\la \Dqc(Y)$ is the derived
pushforward map, and $f^\times:\Dqc(Y)\la\Dqc(X)$ its right
adjoint. Then the natural map is an isomorphism
\[
\CD
f_*\HHom(E,f^\times F) @>>>\HHom(f_*E, F) .
\endCD
\]
\elem

\prf
Let $C\in\Dqc(Y)$ be arbitrary. The lemma follows from 
the sequence of  isomorphisms
\begin{eqnarray*}
 \Hom_Y^{}\big(C\,\,,\,\, f_*\HHom_X^{}(E,f^\times F)\big)&\cong &
 \Hom_X^{}\big({} f^*C\,\,,\,\, \HHom_X^{}(E,f^\times F)\big)\\
&\cong &
 \Hom_X^{}\big({} f^*C\oo_X^{} E\,\,,\,\, f^\times F\big)\\
&\cong &
\Hom_Y^{}\big(f_*[{} f^*C\oo_X^{} E]\,\,,\,\, F\big)\\
&\cong &
\Hom_Y^{}(C\oo_Y^{} f_*E\,\,,\,\, F)\\
&\cong &
\Hom_Y^{}\big(C\,\,,\,\,\HHom_Y^{}(f_*E, F)\big)
\end{eqnarray*}
The 
first isomorphism is by the adjunction
${} f^*\dashv f_*$,
the second and fifth are by the adjunction of $\oo$ and
$\HHom$,
the third is by
the
adjunction $f_*\dashv f^\times$, and the 
fourth is by the projection formula (which holds 
for concentrated morphisms by \cite[Corollary~4.12]{Hall-Rydh13}).
\eprf

\rmk{R4.2.55}
There exist algebraic stacks $X$ for which the derived category $\Dqc(X)$ 
is not compactly generated---see \cite[Theorem~1.1]{Hall-Neeman-Rydh14}
for examples.
In this article we view the phenomenon as pathological, our results will mostly
focus on the good stacks $X$ for which $\Dqc(X)$ is compactly
generated. There are many interesting classes
of such stacks and much work has gone into studying them,
for the best theorems to date the reader is referred to Hall and 
Rydh~\cite[Introduction]{Hall-Rydh13}.

Because in this section we will make some effort to state our results in
the maximal generality in which the proofs hold, it will help to make a
convention.
\ermk

\cvn{C4.2.553}
We will say that a stack $X$ \emph{satisfies Thomason's condition} if
$\Dqc(X)$ is compactly generated, and moreover for any quasi-compact open set
$U\subset X$ the subcategory $\D_{\mathbf{qc},X-U}^{}(X)$ is generated by 
the compact objects in $\Dqc(X)$ that happen to lie in the subcategory.
We remind the reader: 
a complex belongs to $\D_{\mathbf{qc},X-U}^{}(X)\subset\Dqc(X)$ 
if its restriction to $U$ is
acyclic.
Thomason~\cite{ThomTro} proved that every quasi-compact, semi-separated
scheme satisfies Thomason's condition. The current best theorems
for algebraic stacks
tell us that a  stack $X$ satisfies Thomason's condition if either
$X$ is a $\qq$--stack of $s$--global type
in the notation of \cite[Section~2]{Rydh15}, or $X$ is quasi-compact and its
diagonal is separated and quasi-finite--- see Hall and 
Rydh~\cite[Theorems~A, B and 4.10(2)]{Hall-Rydh13}. The following
observation will be useful.
\ecvn

\obs{O4.2.1999}
If $X$ and $X'$ are quasi-compact, quasi-separated stacks,
if $X$ satisfies Thomason's condition and if $i:X'\la X$ is a
locally closed immersion, then $X'$ satisfies Thomason's condition.
\eobs

\prf
The map $i$ is a locally closed immersion of quasi-compact,
quasi-separated stacks, hence it is quasi-affine. From
\cite[Lemma~8.2]{Hall-Rydh13} we have that the objects $i^*C$, with
$C\in\Dqc(X)$ compact, are compact generators of $\Dqc(X')$.

If $U'\subset X'$ is a quasi-compact open subset then there exists a 
quasi-compact open subset $U\subset X$ with $U'=U\cap X'$. As 
$X$ satisfies Thomason's condition the objects
$C\in\D_{\mathbf{qc},\,X-U}^{}(X)$, with compact image in $\Dqc(X)$, 
generate
$\D_{\mathbf{qc},\,X-U}^{}(X)$. But then the objects $i^*C$ generate
$\D_{\mathbf{qc},\,X'-U'}^{}(X')$
and are compact in $\Dqc(X')$.
\eprf

We will find useful the following little fact.

\lem{L4.705.332}
Let $g:X\la Z$ be a morphism of stacks, and assume
$\Dqc(Z)$ is compactly generated. Then the following
\[
S=\{C\oo g^*\wt C\mid C\text{ compact in }\Dqc(X),\,\,\wt C\text{ compact in }\Dqc(Z)\}
\]
is a class of compact objects in $\Dqc(X)$. Moreover the thick subcategory
generated by $S$---that is the smallest subcategory containing $S$ and closed
under triangles and direct summands---is the
category $\Dqc(X)^c$ of all compact objects in $\Dqc(X)$.
\elem

\prf
If $\wt C$ is a compact object in $\Dqc(Z)$ then
\cite[Lemma~4.4(1)]{Hall-Rydh13}
tells us that $\wt C$ is a perfect complex, hence so is $g^*\wt C$.
If $C$ is a compact object in $\Dqc(X)$ we learn from
\cite[Lemma~4.4(2)]{Hall-Rydh13} that $C\oo g^*\wt C$ is also a compact object in
$\Dqc(X)$.

Next let $\ct\subset\Dqc(X)$ be the localizing subcategory generated by
the class of compacts $S$. The objects of $S$ lie in $\ct$ and are compact
in the larger $\Dqc(X)$, hence are certainly compact in $\ct$.
Since $S$ generates $\ct$ we
may apply \cite[Lemma~4.4.5]{Neeman99} to conclude
that $\text{Thick}(S)=\ct^c$,
where $\text{Thick}(S)$ stands for the thick subcategory generated by $S$.
Because $S\subset\Dqc(X)^c$ and $\Dqc(X)^c$ is thick we conclude that
$\ct^c=\text{Thick}(S)\subset\Dqc(X)^c$. We need to prove the
inclusion $\Dqc(X)^c\subset\ct^c$.

Choose therefore any compact $C\in\Dqc(X)$; the subcategory
$\cs(C)\subset\Dqc(Z)$ defined by
\[
\cs(C)=\{E\in\Dqc(Z)\mid C\oo g^*E\in\ct\}
\]
is obviously a localizing subcategory of $\Dqc(Z)$ containing all the compact
objects $\wt C\in\Dqc(Z)$. But $\Dqc(Z)$ is compactly generated,
hence $\cs(C)=\Dqc(Z)$, and in particular $\co_Z^{}\in\cs(C)$. Therefore
$C=C\oo g^*\co_Z^{}\in\ct$. Thus $\ct$ contains all the compact objects
of $\Dqc(X)$: in symbols $\Dqc(X)^c\subset\ct$. But it's clear that, for an
object in $\ct$, being compact in the larger $\Dqc(X)$ implies
compactness in $\ct$. Hence $\Dqc(X)^c\subset\ct^c$.
\eprf

\ntn{N4.1}
For most of this section we will suppose we are given a 2-cartesian square of
stacks
\[
\CD
W @>u>> X\\
@VfVV @VVgV \\
Y @>v>> Z
\endCD
\]
We always assume the stacks quasi-compact and quasi-separated, and the 
morphism $g$ will be
concentrated as in \cite[Definition~2.4]{Hall-Rydh13}. 
Note that any noetherian stack is automatically 
quasi-compact and quasi-separated, and any representable
morphism is concentrated. As all morphisms of quasi-compact,
quasi-separated algebraic spaces are representable they are concentrated.
\entn

\lem{L4.3}
With the notation as in~\ref{N4.1}, assume the stack $Z$ has
quasi-affine diagonal and that $\Dqc(X)$
is compactly generated. Let $S$ be a set of compact generators in $\Dqc(X)$.
Then 
\be
\item
A morphism $\ph:E\la E'$ in $\Dqc(W)$ is an isomorphism if and only
if, for all compact objects $C\in S$, the functor
$f_*\HHom(u^*C,-)$ takes $\ph$ to an isomorphism.
\item
An object $E\in\Dqc(W)$ vanishes if and only if, for every compact
object
$C\in S$, the functor $f_*\HHom(u^*C,-)$ takes $E$ to zero.
\ee
\elem

\prf
It suffices to prove (ii), as (i) follows by applying (ii) to the
mapping cone of the morphism $\ph:E\la E'$.

Suppose therefore that $E$ is an object of $\Dqc(W)$ and that
$f_*\HHom_W^{}(u^*C,E)=0$ for all $C\in S$. Let $Y'$ 
be an affine scheme and $v':Y'\la Y$
a faithfully flat map; we may extend the diagram of~\ref{N4.1} by pulling back
along $v'$, obtaining
\[
\CD
W' @>u'>> W @>u>> X\\
@Vf'VV @VfVV @VVgV \\
Y' @>v'>> Y @>v>> Z
\endCD
\]
We are assuming that $g$ is concentrated,
and \cite[Lemma~2.5(1)]{Hall-Rydh13} tells us that so are
its pullbacks $f'$ and $f$.
We are given that $f_*\HHom_W^{}(u^*C,E)=0$ for all 
$C\in S$,
hence also ${v'}^*f_*\HHom_W^{}(u^*C,E)=0$ for all $C\in S$.
By \cite[Theorem~2.6(4)]{Hall-Rydh13} the base-change map is an isomorphism
${v'}^*f_*\cong f'_*{u'}^*$, hence $f'_*{u'}^*\HHom_W^{}(u^*C,E)=0$.
Now \cite[Lemma~4.4(1)]{Hall-Rydh13} tells us that the compact object
$C\in\Dqc(X)$ is a perfect complex, therefore so is $u^*C\in\Dqc(W)$, and the
natural map ${u'}^*\HHom_W^{}(u^*C,E)\la \HHom_{W'}^{}({u'}^*u^*C,{u'}^*E)$
is an isomorphism; we conclude that $f'_*\HHom_{W'}^{}({u'}^*u^*C,{u'}^*E)=0$.
Taking global sections we have that $\Hom_{W'}^{}({u'}^*u^*C,{u'}^*E)
=\Hom_X^{}(C,u_*u'_*{u'}^*E)=0$ for every $C\in S$. But the
objects of $S$ generate, hence $u_*u'_*{u'}^*E=0$.

But $Y'$ is assumed affine and $Z$ has quasi-affine diagonal, hence the map
$vv':Y'\la Z$ is quasi-affine. Its pullback $uu':W'\la X$ is also quasi-affine,
and by \cite[Corollary~2.8]{Hall-Rydh13} the
functor $u_*u'_*:\Dqc(W')\la\Dqc(X)$ is conservative. The vanishing
of $u_*u'_*{u'}^*E$ implies the vanishing of ${u'}^*E$. The map $v'$ is
faithfully flat and hence so is its pullback $u'$, therefore
the vanishing of ${u'}^*E$ implies the vanishing of $E$.
\eprf

\cor{C4.7993}
With the notation as in~\ref{N4.1}, assume the stack $Z$ has
quasi-affine diagonal, and the
categories $\Dqc(X)$ and $\Dqc(Y)$ are both
compactly generated. Then $\Dqc(W)$ is generated by the compact
objects $u^*C\oo f^*D$, where $C\in\Dqc(X)$ and $D\in\Dqc(Y)$ are
compact.
\ecor

\nin
The case of the corollary
where the diagram is of
quasi-compact, quasi-separated schemes was proved in
Bondal and Van den Bergh~\cite[Lemma~3.4.1]{BondalvandenBergh04},
while the case where $Z$ has affine diagonal may be
found in Ben-Zvi, Francis and 
Nadler~\cite[Proposition~3.24]{BenZvi-Francis-Nadler10}.

\prf
In Notation~\ref{N4.1} we assume that
$g$ is concentrated, hence so is its pullback $f$; by
\cite[Theorem~2.6(3)]{Hall-Rydh13} we know that $f_*$ preserves coproducts,
while from the proof of \cite[Theorem~5.1]{Neeman96} it follows that $f^*$ takes
compacts to compacts. Hence $f^*D$ is compact for any compact $D\in\Dqc(Y)$.
The object $C\in\Dqc(X)$ is compact, hence perfect by
\cite[Lemma~4.4(1)]{Hall-Rydh13}. Therefore 
$u^*C$ is perfect in $\Dqc(W)$. 
Now \cite[Lemma~4.4(2)]{Hall-Rydh13} allows us to 
conclude that $u^*C\oo f^*D$ is compact in $\Dqc(W)$.

It remains to show that the objects $u^*C\oo f^*D$ generate. Suppose that $E$
is an object of $\Dqc(W)$ so that, for all compact 
$C\in\Dqc(X)$ and $D\in\Dqc(Y)$, we have $\Hom_W^{}(u^*C\oo f^*D,E)=0$; we
need to show $E=0$. To this end consider the isomorphisms
\begin{eqnarray*}
\Hom_W^{}(u^*C\oo f^*D,E) &\cong & 
   \Hom_W^{}\big(f^*D,\HHom_W^{}(u^*C,E)\big) \\
&\cong & 
   \Hom_Y^{}\big(D,f_*\HHom_W^{}(u^*C,E)\big)
\end{eqnarray*}
The vanishing of this for all compact $D\in\Dqc(Y)$ tells us that
$f_*\HHom_W^{}(u^*C,E)=0$ for all compact $C\in\Dqc(X)$, and Lemma~\ref{L4.3}(ii)
now gives that $E=0$.
\eprf

\rmk{R4.1054}
As a curiosity we note that it isn't clear whether $W$ satisfies Thomason's 
condition whenever $X$ and $Y$ do. Even the simplest cases are
unclear, for example when
the maps are all representable, \'etale and separated.
\ermk

We are building up to
Lemma~\ref{L4.5}; in the generality in which we will state it
we need 

\rmd{R4.-1099}
Let $g:X\la Z$ be a 
concentrated morphism of quasi-compact, quasi-separated stacks.
The morphism $g$ is called \emph{quasi-proper} if
 $g_*$ takes 
 pseudo-coherent objects of $\Dqc(X)$
 to pseudo-coherent objects of $\Dqc(Z)$; see 
Reminder~\ref{R4.-7} for the definition of pseudo-coherence of
objects.
The morphism $g$ is \emph{universally quasi-proper} if all pullbacks of $g$
are quasi-proper.
The morphism $g$ is called \emph{pseudo-coherent} if there exists
a 2-commutative square
\[
\CD
\spec S @>u''>> X \\
@Vg'VV   @VVgV \\
\spec R @>u'>> Z
\endCD
\]
with $u'$ smooth and surjective, with the map
$\spec S\la X\times_Z^{}\spec R$ also smooth
and surjective, and where $g'$ 
admits
a factorization $\spec S\stackrel i\la\ak^n_R\la \spec R$ 
with $i$ a closed immersion and so 
that  $i_*$
takes the structure sheaf $\co_{\spec S}^{}$ to a pseudo-coherent object
in $\Dqc(\ak^n_R)\cong\D\big(R[x_1^{},x_2^{},\ldots, x_n^{}]\big)$.
\ermd

\rmk{R4.-20099}
In Reminder~\ref{R4.-1099} we do not assume $g$ to be representable.
Pseudo-coherence of morphisms is preserved by flat base-change, with
quasi-properness this is less clear. 
All finite-type morphisms of noetherian stacks
are pseudo-coherent, and 
all proper maps of noetherian stacks are quasi-proper. For non-noetherian
stacks, we learn from
Kiehl~\cite[p.~315, Theorem~2.2]{Kiehl72} that every proper
pseudo-coherent map is quasi-proper.
But there are also non-representable examples of quasi-proper maps.
Let $X$ be a scheme of finite type over $S$ and $G$ a linearly
reductive group over $S$ acting on $X$. If the GIT quotient $X/\!\!/G$ exists
then the map $[X/G]\la X/\!\!/G$, from the stack $[X/G]$ to the GIT
quotient, is always quasi-proper but rarely representable.
\ermk

\lem{L4.497.-5}
Let $X$ be a quasi-compact,
quasi-separated stack and $C\in\Dqcb(X)$ a pseudo-coherent object. Then the 
locus where $C$ is perfect is open. More precisely:
let
\[
V_n=\left\{
x\in X\left| \begin{array}{c}\text{\rm there exists a flat map }f:\spec R\la X\text{\ \rm with }x\in f\big(\spec R\big)\\
\text{\rm and with $f^*C$ having Tor-amplitude contained in $[-n,\infty)$}
\end{array}\right.\right\}
\] 
Then for $n\gg0$ the sets $V_n$ are open. Therefore so is 
$V=\cup_{n\gg0}^{} V_n$, which is the set of points at which $C$ is perfect.
\elem

\prf
The question is local in the flat topology on $X$, hence we may assume $X$
to be an affine scheme $X=\spec R$. Being pseudo-coherent, the object 
$C
\in\Dqc(X)\cong\D(R)$ 
is isomorphic to a 
bounded above complex of finitely generated
projective $R$--modules---choose and fix such 
an isomorph $C$. As we are assuming that $C$ has bounded cohomology we may
(after shifting) suppose
that $H^i(C)=0$ for all $i<0$.
We assert that $V_n$ is open
for all $n\geq0$. 

Let $\p\in\spec R$ be a prime ideal 
of $R$ belonging to $V_n$ with $n\geq0$. In the definition
of $V_n$ we may take $f$ to be the map $\spec{R_\p}\la\spec R$,
where  $R_\p$ 
is the localization of $R$ at $\p$---from
the fact that $\p\in V_n$
it follows that $f^*C=R_\p\oo_R^{}C$
has Tor-amplitude contained in $[-n,\infty)$. Hence for
any $R_\p$ module $N$ we have $H^i(C\oo N)=0$ for all
$i<-n\leq0$. Thus the brutal truncation $C_{\leq -n}^{}$ of the 
complex $C$, given below 
\[
\CD
\cdots@>>> C^{-n-2} @>>> C^{-n-1} @>\partial>> C^{-n} @>>> 0
\endCD
\]
is acyclic except in degree $-n$, as is the complex $C_{\leq -n}^{}\oo N$.
The complex $C_{\leq -n}^{}$
is a projective resolution over $R$ of the finitely-presented
module $M=\mathrm{Coker}(\partial)$, and the fact that,
for every $R_\p$--module $N$, the complex $C_{\leq -n}^{}\oo N$
is acyclic in degrees $\neq -n$ says that $R_\p\oo M$ is flat as
an $R_\p$ module, hence projective. 
But then $M$ must be projective and finitely 
generated in a Zariski-open neighborhood of $\p$, and on this
neighborhood the complex $C$ is quasi-isomorphic to the perfect 
complex
\[
\CD
\cdots@>>> 0 @>>> M @>\partial>> C^{-n+1} @>>> C^{-n+2} @>>> \cdots
\endCD
\]
\eprf

\cor{C4.497.-111}
Let $g:X\la Z$ be a pseudo-coherent morphism of quasi-compact,
quasi-separated stacks. Then the subset $U\subset X$ of points where
$g$ is of finite Tor-dimension is open. More precisely: for any $n$
let $U_n\subset X$ be defined as the set of all points $x\in X$ at which
the Tor-amplitude of $g$ is contained in $[-n,\infty)$. Then $U_n$ is open for
$n\gg0$, and hence so is the union $U=\cup_{n\gg0}^{}U_n$.
\ecor

\prf
The pseudo-coherence of $g$ means that there exists a 2-commutative square
\[
\CD
\spec S @>u''>> X \\
@Vg'VV   @VVgV \\
\spec R @>u'>> Z
\endCD
\]
as in Reminder~\ref{R4.-1099}. Since $u'$ and $u''$ are both faithfully flat,
it suffices to prove the assertion for $g'$. Explicitly: we need 
to show the openness for large $n$
of the set $U_n\subset\spec S$, of all prime ideals $\p\subset S$
where the localization $S_\p$ has Tor-amplitude over $R$
 contained in the interval $[-n,\infty)$.
But $g'$ 
admits
a factorization $\spec S\stackrel i\la\ak^n_R\stackrel j\la \spec R$
with $j$ flat, and it suffices to prove the assertion for 
$i$. We are given that $i$ is a closed immersion and  
$i_*\co_{\spec S}^{}\in\Dqcpl(\ak^n_R)$  is pseudo-coherent,
and we may apply Lemma~\ref{L4.497.-5} to deduce that, for $n\gg0$, the set 
$V_n\subset\ak^n_R$ 
on which $i_*\co_{\spec S}^{}$ has Tor-amplitude contained in $[-n,\infty)$
is open in $\ak^n_R$. Hence $U_n=V_n\cap\spec S$ is open in $\spec S$.
\eprf

\lem{L4.497.-100}
Let $g:X\la Z$ be a pseudo-coherent morphism of quasi-compact,
quasi-separated stacks. Suppose $z\in Z$ is a point so that the 
map $g$ is of finite Tor-dimension at every
$x\in X$ with $z$ in the closure of $g(x)$.
Then there exist an integer $\ell$
and a 2-cartesian square
\[
\CD
W @>u>> X\\
@VfVV @VVgV \\
Y @>v>> Z
\endCD
\]
with $Y$ an affine scheme, 
with $v$ flat, with $z\in v(Y)$, and so that
$f^*\Dqc(Y)^{\geq n}\subset\Dqc(W)^{\geq n-\ell}$.
\elem

\prf
We assume $g$ pseudo-coherent, hence there exists a 2-commutative square
\[
\CD
\spec S @>u''>> X \\
@Vg'VV   @VVgV \\
\spec R @>u'>> Z
\endCD
\] 
as in Reminder~\ref{R4.-1099}. The map $u'$ is surjective (and smooth), hence
we may choose a point $\p\in\spec R$ with $u'(\p)=z$. Let $R_\p$ be
the localization of $R$ at the prime ideal $\p$, let $v:Y\la Z$ be the
composite $\spec{R_\p}\la\spec R\stackrel{u'}\la Z$, and form the
 pullback square
\[
\CD
W @>u>> X\\
@VfVV @VVgV \\
Y @>v>> Z
\endCD
\]
We will prove that this square satisfies the requirements of the Lemma.

By construction $z$ belong to the image of the flat map $v:Y\la Z$,
and
$Y=\spec{R_\p}$ is an affine scheme. It remains to show that, 
for some integer $\ell$, we have
$f^*\Dqc(Y)^{\geq 0}\subset \Dqc(W)^{\geq -\ell}$.
We have a smooth and surjective
map $\rho:\spec{S_\p}\la W$, which is the pullback along the map 
$\spec{R_\p}\la
\spec R$
 of the smooth and surjective map $\spec S\la X\times_Z^{}\spec R$,
and it clearly suffices
to show that
$\rho^*f^*\Dqc(Y)^{\geq 0}\subset \Dqc\big(\spec{S_\p}\big)^{\geq -\ell}$.

By hypothesis $g$ is of finite Tor-dimension at every
point $x\in X$ such that the closure of $g(x)$ contains $z$. 
Hence $f$ is of finite Tor-dimension at every point of 
$W$. 
The map  $\rho$ is flat
by construction, and so $f\rho$ is of finite Tor-dimension at
every point of $\spec{S_\p}$.
But the map $f\rho$ factors as 
$\spec{S_\p}\stackrel i\la
\ak^n_{R_\p}\stackrel j\la \spec{R_\p}$, 
and as $j$ is flat $j^*\Dqc(Y)= j^*\Dqc\big(\spec{R_\p}\big)^{\geq 0}\subset
\Dqc(\ak^n_{R_\p})^{\geq 0}$. Thus it suffices
to prove that 
$i^*\Dqc(\ak^n_{R_\p})^{\geq 0}\subset 
\Dqc\big(\spec{S_\p}\big)^{\geq -\ell}$ 
for some integer $\ell$. The map $i$ is a closed immersion hence affine, and it
therefore suffices to prove that 
$i_*i^*\Dqc(\ak^n_{R_\p})^{\geq 0}\subset \Dqc(\ak^n_{R_\p})^{\geq -\ell}$.

The projection formula tells us that
 $i_*i^*(-)\cong(-)\oo i_*\co_{\spec{S_\p}}^{}$, and
we know that $i_*\co_{\spec{S_\p}}^{}$ is pseudo-coherent. The fact that 
$ji$ and therefore $i$ are of finite Tor-dimension at every point 
$s\in\spec{S_\p}$ tells us first
that $i_*\co_{\spec{S_\p}}^{}$
is perfect at every point of the 
closed subset $\spec{S_\p}\subset\ak^n_{R_\p}$, while on the complement
of $\spec{S_\p}\subset\ak^n_{R_\p}$
the complex $i_*\co_{\spec{S_\p}}^{}$
vanishes---thus it is perfect at every
point of $\ak^n_{R_\p}$. Lemma~\ref{L4.497.-5}, applied to
the pseudo-coherent object $i_*\co_{\spec{S_\p}}^{}$
on the stack $\ak^n_{R_\p}$, expresses $\ak^n_{R_\p}$ as a union
$\ak^n_{R_\p}=\cup_{\ell\gg0}^{}V_\ell$ with $V_\ell$ increasing and 
open for large $\ell$. 
As $\ak^n_{R_\p}$ is quasi-compact
there exists an integer $\ell$ with $V_\ell=\ak^n_{R_\p}$,
that is the Tor-amplitude of $i_*\co_{\spec{S_\p}}^{}$
is contained in the interval $[-\ell,\infty)$. The result 
follows.
\eprf

\lem{L4.497}
Let $g:X\la Z$ be a concentrated, pseudo-coherent  morphism of quasi-separated, quasi-compact
stacks, let
$C\in\Dqc(X)$ be a perfect complex, and assume $g_*C$ is
pseudo-coherent. [If $g:X\la Z$ is quasi-proper, the last assumption
is a consequence of $C$ being perfect.] Let $V\subset Z$ be the open
set obtained by applying Lemma~\ref{L4.497.-5} to the
pseudo-coherent $g_*C\in\Dqc(Z)$---that is, $V$
is the set of points at which $g_*C$ is perfect.
Then $V$ contains the set of points $z\in Z$ such that $g$ is of
finite Tor-dimension at every point $x\in X$ with $z$ in
the closure
of $g(x)$.
\elem

\prf
Let $z\in Z$ be a point such that $g$ is of finite Tor-dimension
at  any $x\in X$ with $z$ in the
closure of $g(x)$. By Lemma~\ref{L4.497.-100} there exists an integer 
$\ell$ and a 2-cartesian square
\[
\CD
W @>u>> X\\
@VfVV @VVgV \\
Y @>v>> Z
\endCD
\]
with $Y=\spec R$ an affine scheme, 
with $v$ flat, with $z\in v(Y)$, and so that 
$f^*\Dqc(Y)^{\geq
  n}\subset\Dqc(W)^{\geq n-\ell}$.
As $C$ is perfect the Tor-amplitude of $C$ is in
some bounded interval $[-m,m]$, and as $u$ is flat the Tor-amplitude
of $u^*C$ is also in the interval $[-m,m]$.  Therefore
$u^*C\oo f^*\Dqc(Y)^{\geq
  n} \subset 
\Dqc(W)^{\geq n-\ell-m}$, and hence
\begin{eqnarray*}
f_*u^*C\oo \Dqc(Y)^{\geq
  n} &=& 
f_*\big[u^*C\oo f^*\Dqc(Y)^{\geq
  n}\big] \\
&\subset &
f_*\Dqc(W)^{\geq n-\ell-m} \\
&\subset &
\Dqc(Y)^{\geq n-\ell-m}\ .
\end{eqnarray*}
Here the morphism $f$, being the fullback of the concentrated morphism $g$,  is
concentrated by \cite[Lemma~2.5(1)]{Hall-Rydh13}, and
the equality in the sequence of inclusions above is the projection formula for the
concentrated morphism $f$, which holds by \cite[Corollary~4.12]{Hall-Rydh13}.
The inclusion  
$f_*u^*C\oo \Dqc(Y)^{\geq
  n} \subset 
\Dqc(Y)^{\geq n-\ell-m}$ tells us that
$f_*u^*C$ is an object of $\Dqc(Y)=\Dqc\big(\spec R\big)\cong\D(R)$ with Tor-amplitude contained in the interval
$[-\ell-m,\infty)$.
On the other hand 
\cite[Theorem~2.6(4)]{Hall-Rydh13}, applied to the cartesian square
above in which the $v$ is flat and $g$ is concentrated,
tells us that $v^*g_*\cong f_*u^*$, and we deduce that 
$v^*g_*C\cong f_*u^*C$. As $g_*C$ is assumed pseudo-coherent so is
$v^*g_*C$, and being of  
Tor-amplitude contained in $[-\ell-m,\infty)$ it must be
perfect. Hence $z$ is contained in the open set $V$ of points at which
$g_*C$ is perfect.
\eprf

\cor{C4.5.-100}
Let $g:X\la Z$ be a concentrated, pseudo-coherent  morphism of 
quasi-separated, quasi-compact
stacks, and suppose $g$ is of finite Tor-dimension. We have
\be
\item
If $C\in\Dqc(X)$ is perfect and $g_*C\in\Dqc(Z)$ is pseudo-coherent, then
$g_*C$ is perfect.
\item
If $g$ is quasi-proper then $g_*$ takes perfect complexes
to perfect complexes.
\item
  Assume that $g$ is quasi-proper and that $\Dqc(Z)$ is
 compactly generated. Then
  \begin{enumerate}
  \item
    $g_*$ takes compact objects in $\Dqc(X)$ to compact objects in $\Dqc(Z)$.
  \item
 If $\Dqc(X)$ is compactly generated then $g^\times$ respects coproducts, and
 the
 natural transformation $\chi(E,F):g^*E\oo g^\times F\la g^\times(E\oo F)$, which was defined in
 Notation~\ref{NThom.73.63}(iv), 
 is an isomorphism for every $E,F\in\Dqc(Z)$.
  \end{enumerate}
\ee
\ecor

\prf
In proving (i) apply Lemma~\ref{L4.497} noting that, because $g$
is of finite Tor-dimension, every
 $z\in Z$ satisfies the technical condition of the Lemma. Hence $V=Z$,
 that is the complex $g_*C\in\Dqcpl(Z)$ is perfect at every point in $Z$. In the
 notation of Lemma~\ref{L4.497.-5} we have $Z=V=\cup_{n\gg0}^{}V_n$
 with $V_n$ open, and the quasi-compactness of $Z$ gives that $Z=V_n$
 for some $n$. Thus $g_*C$ is a pseudo-coherent complex of
 Tor-amplitude contained in $[-n,\infty)$, hence $g_*C$ is
 perfect. 

Part (ii) is immediate from (i): if $g$ is quasi-proper then $g_*$
takes any perfect complex $C\in\Dqc(X)$ to a pseudo-coherent $g_*C$, which
must be perfect by (i).

To prove (iii)(a) we use Lemma~\ref{L4.705.332}: the class
$S=\{C\oo g^*\wt C\mid C\in\Dqc(X)^c,\,\,\wt C\in\Dqc(Z)^c\}$
is a class of compact objects in $\Dqc(X)$ and 
$\Dqc(X)^c$ is its thick closure.
If we let $\cl$ be the
full subcategory of objects $L\in\Dqc(X)$ with $g_*L$ compact, then
$\cl$ is obviously a thick subcategory. We wish to show that
$\Dqc(X)^c\subset\cl$, and it suffices to show that $S\subset\cl$.

But now we are in business: objects in $S$ are of the
form $C\oo g^*\wt C$, with $C\in\Dqc(X)$ and $\wt C\in\Dqc(Z)$
both compact. The projection formula holds for the
concentrated morphism $g:X\la Z$ by \cite[Corollary~4.12]{Hall-Rydh13},
hence $g_*(C\oo g^*\wt C)\cong g_*C\oo \wt C$.
By \cite[Lemma~4.4(1)]{Hall-Rydh13}
the compact object $C\in\Dqc(X)$ is perfect, and by (ii) above so is
$g_*C\in\Dqc(Z)$. But now \cite[Lemma~4.4(2)]{Hall-Rydh13} tells us
that $g_*C\oo \wt C$ is compact.

This proved (iii)(a). The first part (iii)(b) is a formal consequence: 
we know that $\Dqc(X)$ is compactly generated while (iii)(a) tells us
that $g_*:\Dqc(X)\la\Dqc(Y)$
takes compacts to compacts, and \cite[Theorem~5.1]{Neeman96}
allows us to conclude that $g^\times$ respects coproducts.

Finally note that $\chi$ is a natural transformation
between functors both of which respect coproducts. Given
an object $F\in\Dqc(Z)$, the full subcategory $\car\subset\Dqc(Z)$
of all objects $E$
for which map $\chi(E,F)$ is an isomorphism is localizing.
On the other hand 
Proposition~\ref{PThom.73.66} tells us that every compact object belongs to
$\car$, and as $\Dqc(Z)$ is assumed compactly generated we
have $\car=\Dqc(Z)$.
\eprf

\lem{L4.5}
Let 
$v:Y\la
Z$ 
be a flat morphism of quasi-compact, quasi-separated stacks, and let 
$C',E$ be objects in $\Dqc(Z)$. Assume $C'=C''\oo\wt C$,
with $C''$ pseudo-coherent and $\wt C$ compact.
Assume further that at least one of the two conditions below 
is satisfied:
\be
\item
$Z$ 
satisfies Thomason's condition, and
the image of 
$v:Y\la Z$ 
is contained in the subset $V=\cup_{n\gg0}^{}V_n\subset Z$
given by Lemma~\ref{L4.497.-5} applied
to the pseudo-coherent object $C''\in\Dqc(Z)$.
\item
$E$ belongs to $\Dqcpl(Z)\subset\Dqc(Z)$.
\ee
Then the map
\[
\CD
v^*\HHom_Z^{}(C',E) @>\ph>> \HHom_Y^{}(v^*C',v^*E) 
\endCD
\]
is an isomorphism.
\elem

\prf
If $E$ is bounded below then 
$\ph:v^*\HHom_Z^{}(C',E) \la \HHom_Y^{}(v^*C',v^*E) $
is an isomorphism
by \cite[Proposition~3.7]{Illusie71A}. 
The difficulty is in case (i) of the Lemma.

Suppose therefore that we are in case (i) of the Lemma; the fact 
that the image of $v$ is contained in $V$ means
that $Y=\cup_{n\gg0}^{}v^{-1}V_n$.  Lemma~\ref{L4.497.-5} tells us that
 $V_n$ is open for $n\gg0$, hence so is $v^{-1}V_n$, and they
obviously increase with $n$. As $Y$ is quasi-compact there exists 
an integer $n\gg0$ with $Y=v^{-1}V_n$. Now $V_n$ is the union of its
quasi-compact open subsets $V_{n,\alpha}$, hence $Y$ is the union
of $v^{-1}V_{n,\alpha}$, and being quasi-compact $Y$ must be a finite
union of $v^{-1}V_{n,\alpha}$. Since the set of quasi-compact open substacks
of $V_n$ is closed under finite unions, we may choose a quasi-compact
open subset $U\subset V_n$ with $v^{-1}U=Y$. 
The map  
$v:Y\la Z$ factors as
$Y\stackrel\beta\la U\stackrel\alpha\la Z$ with $\alpha$ an open
immersion and $\alpha^*C''$ perfect. We wish to show that the 
composite
\[
\CD
\beta^*\alpha^*\HHom_Z^{}(C',E) @>>> 
\beta^*\HHom_Y^{}(\alpha^*C',\alpha^*E) 
@>\cong>> \HHom_Y^{}(\beta^*\alpha^*C',\beta^*\alpha^*E) 
\endCD
\]
is an isomorphism, and the second map is an isomorphism because $\alpha^*C'
\cong\alpha^*C''\oo\alpha^*\wt C$
is perfect. Hence it suffices to prove
that the map $\alpha^*\HHom_Z^{}(C',E)\la \HHom_Y^{}(\alpha^*C',\alpha^*E) $
is an isomorphism---we are reduced to the case where $v:Y\la Z$ 
is an open immersion with $v^*C''$ perfect and $Y$ quasi-compact.
The 
quasi-compact open immersion
$v:Y\la Z$ is representable, therefore 
a concentrated morphism by 
\cite[Lemma~2.5(3)]{Hall-Rydh13}, and
\cite[Theorem~2.6(3)]{Hall-Rydh13} gives that $v_*$ respects
coproducts. From~\cite[Theorem~5.1]{Neeman96} it
follows that $v^*$ takes compact objects
to compact objects---we are given that $\wt C\in\Dqc(Z)$ is compact,
hence $v^*\wt C$ is compact, and  \cite[Lemma~4.4(2)]{Hall-Rydh13}
allows us to conclude that $v^*C'\cong v^*C''\oo v^*\wt C$ is a compact
object in $\Dqc(Y)$.

The object $C'\oplus\T
C'\in\Dqc(Z)$ restricts to a compact object in $\Dqc(Y)$, whose 
image in $K_0(Y)$ vanishes. We are assuming that $Z$ satisfies
Thomason's condition. Thomason's
localization theorem therefore guarantees the existence of a compact object
$G\in\Dqc(Z)$
and a morphism $\alpha:G\la C'\oplus\T C'$, so that the restriction to
the quasi-compact open subset $Y\subset Z$ is a quasi-isomorphism---for
the case of semi-separated schemes see
Thomason-Trobaugh~\cite{ThomTro}, in the generality required here see
\cite{Neeman92A}, whose 
relevant results are summarized in \cite[statements 2.1.4 and 2.1.5]{Neeman96}.

Now complete $\alpha$ to a triangle
\[
\CD
D @>>> G @>\alpha >> C'\oplus\T C' @>>> \T D
\endCD
\]
The complex $D\in\Dqcb(Z)$ is pseudo-coherent and $v^*D=0$. 
Lemma~\ref{L4.-4} tells
us that, for every object $E\in\Dqc(Z)$, we have $v^*\HHom(D,E)=0$.
 Look at the morphism
of
triangles
\[
\CD
v^*\HHom(D,E) @<<< v^*\HHom(G,E) @<\alpha << v^*\HHom(C'\oplus\T
C'\,,\,E) \\
@V\beta VV @V\gamma VV @V\delta VV \\
\HHom(v^* D,v^* E) @<<< \HHom(v^* G,v^* E) @<\alpha << \HHom(v^* 
C'\oplus v^*\T
C'\,,\,v^* E) 
\endCD
\]
The object $G$ is compact, hence perfect by \cite[Lemma~4.4(1)]{Hall-Rydh13}.
This guarantees that $\gamma$ is an isomorphism, and $\beta$
is an isomorphism since the domain and codomain both vanish. From
the morphism of triangles we learn that $\delta$ is an isomorphism.
\eprf

\lem{L4.7}
As in Notation~\ref{N4.1}
suppose we are given a 2-cartesian square of quasi-compact, quasi-separated
stacks
\[
\CD
 W @>u>> X\\
@VfVV @VVgV \\
Y @>v>> Z
\endCD
\]
Assume that $\Dqc(X)$ and $\Dqc(Z)$ are compactly generated, that
$Z$ has quasi-affine diagonal, that $v$ is flat and that $g$ is
quasi-proper and concentrated.
Assume further that at least one of the two conditions below 
is satisfied:
\be
\item
The stack $Z$ satisfies Thomason's
condition, the map $g$ is pseudo-coherent,  and 
the map $f$ has finite Tor-dimension. 
\item
$E$ belongs to $\Dqcpl(Z)\subset\Dqc(Z)$.
\setcounter{enumiii}{\value{enumi}}
\ee
Then the base-change map $\Phi:u^*g^\times E\la f^\times v^* E$
is an isomorphism.
\elem

\prf
We wish to show that the base-change map $\Phi:u^*g^\times\la f^\times v^*$ is an 
isomorphism in $\Dqc(W)$, and by Lemmas~\ref{L4.705.332} and \ref{L4.3}
it suffices
to prove that the functor 
$f_*\HHom(u^*C,-\big)$ takes
$\Phi$ to an isomorphism for every object $C=C''\oo g^*\wt C$,
with $C''\in\Dqc(X)$ and $\wt C\in\Dqc(Z)$ both compact.
Fix therefore a compact object $C''\in\Dqc(X)$ and
a compact object $\wt C\in\Dqc(Z)$, and put $C =C''\oo g^*\wt C$.
I assert:
\be
\setcounter{enumi}{\value{enumiii}}
\item
The natural map
$v^*\HHom_Z^{}(g_*C,E)\la \HHom_Y^{}(v^*g_*C,v^*E)$ is an isomorphism.
\ee
We will return to the proof of (iii) later, let us first see that our Lemma follows.

To prove the Lemma consider the commutative diagram below (for a careful check
of the commutativity
see Lipman~\cite[Lemma~4.6.5]{Lipman09})
\[
\xymatrix@C+60pt{
f_*u^*\HHom_X^{}(C,g^\times E)\ar[dd]_{\beta^{-1}} \ar[r]^{(1)} & f_*\HHom_W^{}(u^*C,u^*g^\times E)\ar[d]\\
& \HHom_Y^{}(f_*u^*C,f_*u^*g^\times E)\ar[d]^{\HHom(\beta,\beta^{-1})}\\
v^*g_*\HHom_X^{}(C,g^\times E)\ar[r]\ar[d]_{(2)}& \HHom_Y^{}(v^*g_*C,v^*g_*g^\times E)\ar[d]^{\HHom(\id,v^*\e)}\\
v^*\HHom_Z^{}(g_*C, E)\ar[r]^{(3)} &  \HHom_Y^{}(v^*g_*C,v^* E)\ar[d]^{\HHom(\beta^{-1},\id)}\\
 \HHom_Y^{}(f_*u^*C,f_*f^\times v^* E)\ar[r]^-{\HHom(\id,\e')} &  \HHom_Y^{}(f_*u^*C,v^* E)\\
 & f_*\HHom_W^{}(u^*C,f^\times v^* E) \ar[ul]\ar[u]_{(4)}
}\]
In this diagram $\beta:v^*g_*\la f_*v^*$ is the
base-change map of Definition~\ref{D27.1}, the map $\e:g_*g^\times\la\id$ is
the counit of the adjunction $g_*\dashv g^\times$,  
and $\e':f_*f^\times\la\id$ is
the counit of the adjunction $f_*\dashv f^\times$.
The commutative triangle at the bottom is there merely to remind us of
the definition of the isomorphism $(4)$ of Lemma~\ref{L4.97}
\[
f_*\HHom_W^{}(A, f^\times B)\quad\cong\quad\HHom_Y^{}(f_*A,B).
\]
The map $(1)$ is an isomorphism because $C$ is perfect, the map $(2)$ is
the isomorphism 
$ g_*\HHom_X^{}(A, g^\times B)\la\HHom_Z^{}(g_*A,B)$ of Lemma~\ref{L4.97}, and 
$(3)$ is an isomorphism by
(iii) above, to be proved later. 
The vertical arrows on the right of the diagram compute 
$f_*\HHom_W^{}(u^*C,-)$ applied to the composite 
defining $\Phi:u^*g^\times E\la
f^\times v^*E$, essentially by the
definition of $\Phi$. And the commutative diagram shows how to factor this
map through a series of isomorphisms.

It remains to prove assertion (iii) above. The object $C''$ is assumed compact
in $\Dqc(X)$, and 
\cite[Lemma~4.4(1)]{Hall-Rydh13} tells us that $C''$ is perfect.
Since $g$ is quasi-proper the object $g_*C''\in\Dqc(Z)$ is pseudo-coherent.
This makes $g_*C=g_*(C''\oo g^*\wt C)\cong g_*C''\oo\wt C$ the tensor
product of a pseudocoherent object and a compact object.
We're in the situation of Lemma~\ref{L4.5}. 
If $E\in\Dqcpl(Z)$ then
(iii) follows immediately from Lemma~\ref{L4.5}(ii).
We have to prove that, under the hypotheses of case (i) of
our current Lemma, we may apply Lemma~\ref{L4.5}(i) to the
object $C'=g_*C=g_*C''\oo\wt C\in\Dqc(Z)$ and the map $v:Y\la Z$.
We know that $C''$ is compact in
$\Dqc(X)$, hence perfect.  Therefore $u^*C''\in\Dqc(W)$
is perfect.
We also know that
$g_*C''$ is pseudo-coherent, thus so is $v^*g_*C''\cong f_*u^*C''$.
We apply Corollary~\ref{C4.5.-100}(i) to the morphism $f:W\la Y$
and the perfect complex $u^*C''\in\Dqc(X)$---by hypothesis $f$
is of finite Tor-dimension and, being
the pullback of $g$, it is concentrated and pseudo-coherent. The 
conclusion is that $v^*g_*C''\cong f_*u^*C''$ is perfect, and hence $v$
takes $Y$ to the set of points $V\subset Z$ where $g_*C''$ is perfect.
Lemma~\ref{L4.5}(i) applies.
\eprf

\rmk{R4.105}
For noetherian schemes
case (ii) of Lemma~\ref{L4.7} is classical and due to
Hartshorne~\cite{Hartshorne66} and
Verdier~\cite{Verdier68}. In fact one of the two classical approaches 
to Grothendieck
duality used Lemma~\ref{L4.7}(ii) as the cornerstone on which
the theory is built---in the category $\Dqcpl$, of course. In
\cite[Example~6.5]{Neeman96} we observed that the base-change map need
not be an isomorphism on unbounded objects. The main new tool of this
article is  Lemma~\ref{L4.7}(i), and the rest of
the article shows how to deduce from this little lemma
that everything works fine in the larger category $\Dqc$, 
as long as one is just a little careful with composites.

There is a modern treatment of Lemma~\ref{L4.7}(ii)
in Lipman~\cite[Section 4.6]{Lipman09}, quite similar to the one
above. What is
really new and different in this paper is that the crucial
Lemma~\ref{L4.5} 
has a version (i), leading to 
Lemma~\ref{L4.7}(i). The reader will note that 
the proof of Lemma~\ref{L4.5}(i) is 
rather subtler than that of the easy, classical Lemma~\ref{L4.5}(ii).
\ermk

\section{Geometric manipulations to obtain a sharper
  version of the base-change Lemma~\protect{\ref{L4.7}}}
\label{S901}

In this section we go through some geometric gymnastics to obtain
a refinement of Lemma~\ref{L4.7}---we will still be appealing to
the topological techniques involving compact generation,
but in combination with the construction
of auxiliary stacks and diagrams of morphisms among them.
Before we start let us note
  
\rmk{R4.7.9997}
In Example~\ref{E4.7.53} we will see
 a cartesian square of affine schemes,
satisfying all the 
hypotheses
of Lemma~\ref{L4.7} except for the pseudo-coherence of $g$.
But the map $u^*g^\times\la f^\times
v^*$ is not an isomorphism, not
even when restricted to $\Dqcpl(Z)$. 
The pseudo-coherence seems essential.

When the stacks are noetherian all finite-type maps are pseudo-coherent and
all proper maps are quasi-proper,
making the hypotheses of Lemma~\ref{L4.7} easy to fulfill.
In the non-noetherian case quasi-proper, pseudo-coherent maps are hard
to come by---there is a theorem of
Kiehl~\cite[p.~315, Theorem~2.2]{Kiehl72} saying that every proper,
pseudo-coherent map is quasi-proper but, while 
proper maps are plentiful,
much of the geometric acrobatics we're about to witness
will not
preserve the pseudo-coherence of morphisms.

This means that very soon we'll start assuming all our stacks
noetherian. But to better pinpoint where the noetherian hypothesis seems 
essential, the next few lemmas remain in gorgeous generality. 
\ermk

\lem{L4.8}
Suppose we are given a diagram of 2-cartesian squares
of quasi-compact, quasi-separated stacks
\[
\CD
U @>u'>> W @>u>> X\\
@VeVV @VfVV @VVgV \\
V @>v'>> Y @>v>> Z
\endCD
\]
Assume $Y$ and $Z$ have quasi-affine diagonals and 
satisfy Thomason's condition, 
the category $\Dqc(X)$ is
compactly generated, 
the morphisms $f$ and $g$ are quasi-proper, the morphism $g$ 
 is concentrated and pseudo-coherent, $v'$ and $v$ are flat and $e$
is
of finite Tor-dimension. 
If $\Phi:u^*g^\times\la f^\times v^*$ is the base-change map of
the square
\[
\CD
W @>u>> X\\
@VfVV @VVgV \\
Y @>v>> Z
\endCD
\]
then ${u'}^*\Phi: {u'}^*u^*g^\times\la {u'}^*f^\times v^*$ is an
isomorphism.
\elem

\prf
Because $\Dqc(X)$ and $\Dqc(Y)$ 
are compactly generated and
$Z$ has quasi-affine diagonal, Corollary~\ref{C4.7993} tells us that $\Dqc(W)$
is compactly generated. Pseudo-coherence and concentratedness are
preserved by 
flat base-change
and $f$, being the pullback of the concentrated and pseudo-coherent
$g$ by the flat 
map $v$, must also be
concentrated and pseudo-coherent.
Consider the two 2-cartesian squares
\[
\CD
U @>u'>> W  @. @. @.     U @>uu'>>  X\\
@VeVV @VfVV  @.\qquad \qquad\text{and}\qquad\qquad @.    @VeVV @VVgV \\
V @>v'>> Y  @. @. @.   V @>vv'>> Z
\endCD
\]
Since $e$ is of finite Tor-dimension Lemma~\ref{L4.7} applies to
both squares, and tells us that
the base-change maps for these two squares are isomorphisms. That is
we have isomorphisms $\Phi_1^{}:{u'}^*f^\times\la e^\times {v'}^*$ and
$\Phi_2^{}:{u'}^*u^*g^\times\la e^\times {v'}^*v^*$. But the isomorphism
$\Phi_2^{}$ can be written as the composite
\[
\CD
{u'}^*u^*g^\times @>{u'}^*\Phi>> {u'}^*f^\times v^* @>\Phi_1^{}v^*>> 
{e}^\times{v'}^*v^*
\endCD
\]
where we know that $\Phi_1^{}$ is an isomorphism. Hence so is 
${u'}^*\Phi={\big[\Phi_1^{}v^*\big]}^{-1}\Phi_2^{}$.
\eprf

The next few lemmas are similar; it might help to set up common
notation for all of them. We begin with a reminder of some standard notation.

\rmd{R4.911}
Recall that, if $\cs'$ is a finitely presented quasicoherent sheaf of
$Z$, then $\pp(\cs')$ stands for the stack
$\text{Proj}\big(\text{Sym}(\cs')\big)$ where $\text{Sym}(\cs')$ is the
symmetric
algebra of $\cs'$ over $\co_X^{}$. The stack $Z$ is said to \emph{have
the resolution property} if every finitely presented sheaf $\cs'$ over $Z$
admits an epimorphism $\cv'\la\cs'$, with $\cv'$ a finite-rank vector
bundle on $Z$. If $g:X\la Z$ has a factorization $X\stackrel
j\la\pp(\cs')\la Z$, with $j$ a closed immersion and $\cs'$ a
finitely-presented
quasicoherent sheaf on $Z$, and if $Z$ has the resolution property,
then we can factor $g$ as $X\stackrel
j\la\pp(\cs')\stackrel{j'}\la\pp(\cv')\la Z$ with $j'$ also a closed immersion.
\ermd

\ntn{N4.9.99312}
For the next three results the setup will be as follows. As in
Notation~\ref{N4.1} we will
suppose given a 2-cartesian square of quasi-compact, 
quasi-separated stacks
\[
\CD
W @>u>> X\\
@VfVV @VVgV \\
Y @>v>> Z
\endCD
\]
We'll assume that $v$ is flat and $g$ is
concentrated and pseudo-coherent. We will also assume given a separated,
finitely-presentable, 
representable, \'etale
morphism
$\theta_Z^{}:\ov Z\la Z$ (which may be
the identity $\id:Z\la Z$). Assume that
$\Dqc(X)$ is compactly generated and
$Y$ and $Z$ have quasi-affine diagonals. Suppose further that, for all
separated,
finitely-presentable, 
representable \'etale maps $\alpha:Y'\la Y$ and
$\beta:Z'\la Z$ and for every finitely presented quasicoherent
sheaf
$\cs$ on either $Y'$ or 
$Z'$, the stack $\pp(\cs)$ satisfies  Thomason's condition.
Assume also that $g$ factors as 
$X\stackrel j\la\pp(\cs')\stackrel
q\la Z$, 
where $\cs'$ is a finitely presented 
quasicoherent sheaf on $Z$ and $j$ is a closed
immersion. By Corollary~\ref{C4.497.-111}
the set $U\subset W$ on which $f$
is of finite Tor-dimension is open---let  $u':U\la W$ be the open immersion.
 Let 
$\Phi:u^*g^\times\la f^\times v^*$ be the base-change map.

The next three results all have the following shape: under some added
hypotheses on the data above, which become less restrictive as we
progress, some slight variants of some of (i), (ii) and possibly (iii)
below hold.
\be
\item
 The natural transformation ${u'}^*\Phi \theta_{Z*}^{}:
{u'}^*u^*g^\times \theta_{Z*}^{}\la {u'}^*f^\times v^*\theta_{Z*}^{}$ is an isomorphism.
\item
The isomorphic functors ${u'}^*u^*g^\times \theta_{Z*}^{}\cong {u'}^*f^\times v^*\theta_{Z*}^{}$
respect coproducts.
\item
If we further assume that all the stacks in the diagram are noetherian, then
for any object $z\in\Dqc(Z)$ we have that
${u'}^*u^*g^\times z=0$ 
if and only if 
${u'}^*u^*g^*z=0$. 
\ee
\entn

\lem{L4.9}
With the conventions of Notation~\ref{N4.9.99312}, assume further that
$\theta_Z^{}:\ov Z\la Z$ is the identity map $\id:Z\la Z$, and that the
quasicoherent sheaf $\cs'=\cv'$ is a finite-rank vector bundle on $Z$. Then the conclusions
of
Notation~\ref{N4.9.99312} (i), (ii) and (iii) are true.
\elem

\prf
The map $g$ is assumed to factor as $g=qj$, where 
$j:X\la \pp(\cv')$ is a closed immersion and $q:\pp(\cv')\la Z$ is 
the projection. 
Let $\cv$ be the pullback of the vector bundle $\cv'$ by the map
$v:Y\la Z$, then
our 2-cartesian square factors as a concatenation
of two 2-cartesian squares
\[
\CD
W @>u>> X\\
@ViVV @VVjV \\
\pp(\cv) @>w>> \pp(\cv')\\
@VpVV @VVqV \\
Y @>v>> Z
\endCD
\]
The vertical morphisms in this diagram are all concentrated.
The stacks $Z$, $\pp(\cv)$ and  $\pp(\cv')$ have quasi-affine diagonals
and satisfy Thomason's condition. The maps $p_*$ and $q_*$ obviously
take pseudo-coherent complexes to pseudo-coherent complexes, that is
$p$
and $q$ are quasi-proper.
 The map $q$ is smooth while
$g=qj$ is pseudo-coherent by hypothesis, hence $j$ is
pseudo-coherent. Pseudo-coherence is stable by flat base change,
%(see~\cite[Cor~1.10 and Cor~4.7.2]{Illusie71A}), 
hence
$i$ is also pseudo-coherent. Since $i$ and $j$ are proper and
pseudo-coherent, Kiehl~\cite[p.~315, Theorem~2.2]{Kiehl72} tells us they 
are quasi-proper. 

Since $p$ is 
flat it is certainly of finite Tor-dimension, and Lemma~\ref{L4.7}
tells us that the base-change map $\Phi_1^{}:w^*q^\times\la p^\times v^*$ is
an isomorphism.

Since $p$ is smooth, the open
subset $U\subset W$ on which
$f=pi$ is of finite Tor-dimension is equal to the
subset on which $i$ is of finite Tor-dimension;
see~\cite[page 246, 3.6]{Illusie71C}.
Choose an open set $V\subset\pp(\cv)$ with $U=V\cap W$. 
If $u':U\la W$ is the open immersion, we wish to
show that ${u'}^*\Phi:
{u'}^*u^*g^\times\la {u'}^*f^\times v^*$ is an isomorphism. The question
is local in $U$, and replacing $V$ by a 
quasi-compact open subset we may assume $V$ is quasi-compact and
open in $\pp(\cv) $. Now consider the diagram with 2-cartesian squares
\[
\CD
U @>u'>> W @>u>> X\\
@Vi'VV @ViVV @VVjV \\
V@>w'>> \pp(\cv) @>w>> \pp(\cv')
\endCD
\]
then the map $i':U\la V$ is of finite Tor-dimension. By Lemma~\ref{L4.8}
we have that, if $\Phi_2^{}:u^*j^\times\la i^\times w^*$ is the base 
change map of the square
\[
\CD
W @>u>> X\\
@ViVV @VVjV \\
 \pp(\cv) @>w>> \pp(\cv')
\endCD
\]
then ${u'}^*\Phi_2^{}:{u'}^*u^*j^\times\la {u'}^*i^\times w^*$ is an
isomorphism.
Let us summarize what we have so far: we began with a diagram of 
2-cartesian squares
\[
\CD
W @>u>> X\\
@ViVV @VVjV \\
\pp(\cv) @>w>> \pp(\cv')\\
@VpVV @VVqV \\
Y @>v>> Z
\endCD
\]
We have studied the base-change maps of the squares
\[
\CD
W @>u>> X      @. @. @.       \pp(\cv) @>w>> \pp(\cv')      \\
@ViVV @VVjV  @.\qquad \qquad\text{and}\qquad\qquad @.          @VpVV @VVqV \\
\pp(\cv) @>w>> \pp(\cv') @. @. @.     Y @>v>> Z
\endCD
\]
and proved that 
\be
\item
$\Phi_1^{}:w^*q^\times\la p^\times v^*$ is an isomorphism.
\item
The functor ${u'}^*$ takes the base-change map
$\Phi_2^{}:u^*j^\times\la i^\times w^*$ to an isomorphism.
\ee
If $\Phi:u^*g^\times=u^*j^\times q^\times\la i^\times p^\times v^*=f^\times v^*$ is the 
base-change map of the 2-cartesian square
\[
\CD
W @>u>> X\\
@VfVV @VVgV \\
Y @>v>> Z
\endCD
\]
then we can write it as the composite
\[
\CD
u^*j^\times q^\times @>\Phi_2^{}q^\times>>
i^\times w^* q^\times @>i^\times\Phi_1^{}>> 
i^\times p^\times v^*
\endCD
\]
where $i^\times\Phi_1^{}$ is an isomorphism, while $\Phi_2^{}q^\times$ becomes 
an isomorphism if we compose it with ${u'}^*$. Therefore  ${u'}^*\Phi$
is an isomorphism, and we have proved (i).

Next we prove (ii), that is we show that the isomorphic
functors
${u'}^*u^*g^\times\cong {u'}^*f^\times v^*$ respect coproducts. 
The question is still local in $U$, hence we may continue
assuming that, in the diagrams above, $U$ and $V$ are quasi-compact.
Write ${u'}^*f^\times v^*$ as ${u'}^*i^\times p^\times v^*$. The base-change map
${u'}^*i^\times\la {i'}^\times{w'}^*$ is an isomorphism by 
Lemma~\ref{L4.7}(i), hence
the functor ${u'}^*f^\times v^*$ is isomorphic to ${i'}^\times{w'}^*p^\times v^*$.
The functors ${w'}^*$ and $v^*$ preserve comproducts because
they have right adjoints. It therefore suffices to show that
${i'}^\times$ and $p^\times$ preserve coproducts, and
we will prove this by applying
Corollary~\ref {C4.5.-100}(iii)(b).
We know that $i'$ and $p$ are pseudocoherent and proper---hence quasi-proper
by Kiehl's theorem. They are also concentrated and of finite
Tor-dimension, hence to apply Corollary~\ref{C4.5.-100}(iii)(b)
we only need to check that the categories
$\Dqc(U)$, $\Dqc(V)$, $\Dqc\big(\pp(\cv)\big)$ and $\Dqc(Y)$ are
compactly generated. For $\Dqc(Y)$ and $\Dqc\big(\pp(\cv)\big)$ this is
part of the hypotheses of Notation~\ref{N4.9.99312}. Note also that
by the hypotheses of Notation~\ref{N4.9.99312} we know that $Z$ has
quasi-affine diagonal and $\Dqc(Y)$ and $\Dqc(X)$ are compactly generated,
and Corollary~\ref{C4.7993} guarantees that $\Dqc(W)$ is compactly
generated. But the maps $u':U\la W$ and $v':V\la\pp(\cv)$
are  open immersions and \cite[Lemma~8.2]{Hall-Rydh13} (see also
Observation~\ref{O4.2.1999}) tell us that
the compact generation of
$\Dqc(W)$ and $\Dqc\big(\pp(\cv)\big)$ implies the
compact generation of $\Dqc(U)$ and $\Dqc(V)$.

It remains to prove (iii), in which all the 
stacks are assumed noetherian.
With the notation as above choose an affine scheme
$S'$ of finite type and faithfully flat over $V$,
and form the diagram of 2-cartesian
squares of finite-type maps of noetherian stacks
\[
\CD
S @>\sigma>> U @>u'>> W @>u>> X \\
@V\gamma VV  @Vi'VV @VVi V @VVjV\\
S' @>\theta>> V @>w'>> \pp(\cv) @>w>> \pp(\cv')
\endCD
\]
By construction the horizontal maps are flat,  
the vertical maps 
are closed immersions, and $i'$ and therefore
$\gamma$ are of finite Tor-dimension.
By the faithful flatness of $\sigma$ we have that, for any
object $\mathfrak{u}\in\Dqc(U)$, the object $\s^*\mathfrak{u}$ vanishes if
and only if $\mathfrak{u}$ does. Hence given $z\in\Dqc(Z)$, we have
${u'}^*f^\times v^*z=0$
if and only if $\sigma^*{u'}^*f^\times v^*z\cong
\sigma^*{u'}^*i^\times p^\times v^*z=0$.
But Lemma~\ref{L4.7}, applied to the cartesian square
\[
\CD
S  @>u'\s>> W \\
@V\gamma VV @VVi V\\
S' @>w'\theta>>  \pp(\cv)
\endCD
\] 
gives an isomorphism 
$\sigma^*{u'}^*i^\times\cong \gamma^\times\theta^*{w'}^*$,
and it follows that ${u'}^*f^\times v^*z=0$
 if and only if 
$\gamma^\times\theta^*{w'}^*p^\times v^*z=0$.

But now $\gamma$ is a closed immersion of affine noetherian schemes, and 
the functors $\gamma^*$ and $\gamma^\times$ are very concrete. Let $S'=\spec {R'}$
and $S=\spec R$ and identify $\Dqc(S)\cong D(R)$ and $\Dqc(S')\cong\D(R')$,
then the functor $\gamma^*$ identifies with $R\oo_{R'}^{}(-)$ while
$\gamma^\times$ is the functor $\Hom_{R'}^{}(R,-)$.
Support theory, more concretely 
\cite[Proposition~A.3(1)]{Iyengar-Lipman-Neeman13} applied to
the object $E=R\in\D(R')$ with Zariski-closed support and to the object
$F=\theta^*{w'}^*p^\times v^*z$
 in 
$\D(R')$, tells us that
$\gamma^\times\theta^*{w'}^*p^\times v^*z=0$
if and only if 
$\gamma^*\theta^*{w'}^*p^\times v^*z=
\sigma^*{u'}^*i^*p^\times v^*z=0$,
 and the faithful flatness of
$\sigma$ tells us that this happens if and only if 
${u'}^*i^*p^\times v^*z=0$. Summarizing what
we have so far: ${u'}^*f^\times v^*z=0$
if and only if ${u'}^*i^*p^\times v^*z=0$.

It remains to show that 
${u'}^*i^*p^\times v^*z=0$ if and only if 
${u'}^*i^*p^* v^*z={u'}^*f^*v^*z=0$. But the map 
$iu':U\la \pp(\cv)$ is a locally closed immersion, hence
$i_*u'_*$ is conservative---it reflects vanishing.
Hence it remains to show that $i_*u'_*{u'}^*i^*p^\times v^*z=0$
if and only if 
$i_*u'_*{u'}^*i^*p^*v^*z=0$. 
But the projection formula---which holds for the concentrated 
morphism $iu'$ by 
\cite[Corollary~4.12]{Hall-Rydh13}---gives an isomorphism 
$i_*u'_*{u'}^*i^*(-)\cong
(-)\oo i_*u'_*{u'}^*i^*\co_{\pp(\cv)}^{}$. 
Hence it will certainly
suffice to show that, for any objects $F\in\Dqc\big(\pp(\cv)\big)$
and $y\in\Dqc(Y)$, we have
$F\oo p^*y=0$ if and only if $F\oo p^\times y=0$.
Because the map $p$ is flat and proper Lemma~\ref{L4.7} says that 
$p^\times$ commutes with flat base change, and hence the assertion can
be checked flat-locally. Replacing $Y$ by an affine faithfully flat cover
$\spec R$
on which the vector bundle $\cv$ trivializes, what remains to be 
proved is
\be
\setcounter{enumi}{\value{enumiii}}
\item
Let $R$ be a noetherian ring.
Let $p:\pp^n_R\la\spec R$ be the natural projection, let $y$ be
an object in $\Dqc\big(\spec R\big)\cong\D(R)$, and let
$F$ be an object in $\Dqc(\pp^n_R)$. Then $F\oo p^*y=0$ 
if and only if $F\oo p^\times y=0$.
\setcounter{enumiii}{\value{enumi}}
\ee
Of course at this point we could really stop---in the special case
of $p:\pp^n_R\la\spec R$ the functor $p^\times$ is well-understood
and we can appeal to computations in the literature. But in this article
we're partly interested in how much of the theory can be developed on very
formal, category-theoretic grounds, without any computation;
in the manuscript \cite{NeemanTIFR} we will explain how the 
formal theory can be used to render the traditional computations very easy. Let
us therefore proceed to prove (iii) using zero concrete computations.

Our particular $p$ is flat and proper, hence Corollary~\ref{C4.5.-100}(iii)
tells us that $p_*$ 
respects compacts and $p^\times$
respects coproducts. By Corollary~\ref{C4.5.-100}(iii)(b), with $F=R=\co_R^{}$,
 it follows
that there
is a natural isomorphism
$p^\times(-)\cong p^*(-)\oo p^\times R$.
Therefore proving (iii) is equivalent to proving that
$F\oo p^*y\oo p^\times R=0$ if and only if
$F\oo p^*y=0$. Support theory, more precisely 
\cite[Proposition~A.3(1)]{Iyengar-Lipman-Neeman13}, tells us that
it suffices to prove that $\text{supp}(p^\times R)=\pp_R^n$.

Given a point $x\in\pp^n_R$ with residue field $k(x)$, we
wish to show that $k(x)\oo p^\times R\neq 0$.
Let $\p\in\spec R$ be $p(x)$; it is a prime ideal
$\p\subset R$. Let $k=k(\p)$ be its residue field.
Consider the cartesian square
\[
\CD
\pp^n_{k} @>\gamma>> \pp^n_R \\
@Vp'VV      @VVp V \\
\spec{k} @>\gamma'>> \spec R
\endCD
\]
The map $p$ is flat, hence Tor-independent base-change gives an
isomorphism ${\gamma'}^*p_*\la {p'}_*\gamma^*$. Taking right
adjoints
gives an isomorphism $\gamma_*{p'}^\times\la p^\times
\gamma'_*$. 
Applying this to the object 
$k\in\D(k)\cong\Dqc\big(\spec
k\big)$ 
we have an isomorphism 
$\gamma_*{p'}^\times k\la
p^\times k
\cong p^*k\oo_{\pp_R^n}^{}p^\times R$.
 We wish to show that
$k(x)\oo_{\pp_R^n}^{}p^\times R\neq0$, 
and because 
$k(x)
\oo_{\pp_R^n}^{}p^*k$ 
is a (nonzero) direct sum of suspensions of
$k(x)$,
it suffices to prove the non-vanishing of 
\[
\big[k(x)
\oo_{\pp_R^n}^{}p^*k\big] \oo_{\pp_R^n}^{}p^\times
R\quad\cong\quad
k(x)
\oo_{\pp_R^n}^{}\big[p^*k \oo_{\pp_R^n}^{}p^\times
R\big]\quad\cong\quad
k(x)
\oo_{\pp_R^n}^{}[\gamma_*{p'}^\times k]
\]
The projection formula tells us that $k(x)
\oo_{\pp_R^n}^{}[\gamma_*{p'}^\times k]\cong
\gamma_*\big[\gamma^*k(x)\oo_{\pp_k^n}^{} {p'}^\times k]$. The morphism
$\gamma$ is affine,
hence $\gamma_*$ is conservative---it suffices to prove the
non-vanishing of $\gamma^*k(x)\oo_{\pp_k^n}^{}  {p'}^\times k$. But
$\gamma^*k(x)$ is a non-vanishing direct sum of suspensions of the residue field
$k(\gamma^{-1}x)$, and it suffices to prove
$k(\gamma^{-1}x)\oo_{\pp_k^n}^{}  {p'}^\times k$
does not vanish. In other words we are reduced to the case where $R=k$
is a field.

Assume therefore that $R=k$ is a field, therefore $\pp^n_k$ is smooth 
and projective over
the field $k$. The subcategory of compact
objects
in $\Dqc(\pp^n_k)$ is equal to $\dcoh(\pp^n_k)$, and by
\cite[Theorem~1.1]{BondalvandenBergh04}
every finite-type homological functor on $\dcoh(\pp^n_k)$ is representable. In particular
there is an object $E\in\dcoh(\pp^n_k)$ and an isomorphism of functors
$\Hom_{\pp_k^n}^{}(-,E)\la
\Hom_k^{}\big(p_*(-),k\big)\cong\Hom_{\pp_k^n}^{}(-,p^\times
k)$
on the category $\dcoh(\pp^n_k)$. That is, even though we don't (yet) know
whether $p^\times k$ belongs to $\dcoh(\pp^n_k)$, when we restrict
the functor $\Hom_{\pp_k^n}^{}(-,p^\times
k)$ to the subcategory $\dcoh(\pp^n_k)$ we obtain a representable
functor, represented by $E$. The identity
in $\Hom_{\pp_k^n}^{}(E,E)$ maps to an element of $\Hom_{\pp_k^n}^{}(E,p^\times
k)$, that is a morphism $\rho:E\la p^\times k$, and Yoneda's lemma tells us that $\rho$
induces the isomorphism of functors on $\dcoh(\pp^n_k)$. Hence for
every
object $F\in\dcoh(\pp^n_k)$, the functor $\Hom(F,-)$ takes $\rho$ to an
isomorphism. But now the category of all $F'\in\Dqc(\pp^n_k)$ such that $\Hom(F',-)$ takes $\rho$ to an
isomorphism is a localizing subcategory of $\Dqc(\pp^n_k)$ containing
all the compacts, hence it is all of $\Dqc(\pp^n_k)$ . It follows that
$\rho$ is an isomorphism, and $p^\times k$ is compact.

This means that $\text{supp}(p^\times k)$ is a closed subset of
$\pp^n_k$, and Hilbert's Nulstellensatz says the closed points
are dense. It suffices to show that every closed point in
$\pp^n_k$ belongs to $\text{supp}(p^\times k)$. But for a closed
point $x$ the support of $k(x)$ is closed, and
\cite[Proposition~A.3(1)]{Iyengar-Lipman-Neeman13}
tells us that $k(x)\oo p^\times k=0$ if and only if
$\Hom\big(k(x),p^\times k\big)=0$. But 
$\Hom_{\pp_k^n}^{}\big(k(x),p^\times k\big)=\Hom_{k}^{}\big(p_*k(x),k\big)\neq0$.
\eprf

\rmk{R4.973885}
In the 
proof of 
Lemma~\ref{L4.9} the noetherian hypothesis enters when we use 
support theory. Support theory is known to work beautifully for
noetherian schemes, while the obvious generalizations fail
miserably in the non-noetherian case.
But the noetherian hypothesis will become even more crucial later.
\ermk

\lem{L4.93115}
With the conventions of Notation~\ref{N4.9.99312}, assume further that
the stack $\ov Z$ in the given \'etale map 
$\theta_Z^{}:\ov Z\la Z$ satisfies the resolution property. Then the conclusions
of
Notation~\ref{N4.9.99312} (i) and (ii) are true.
\elem

\prf
Pull back along the map ${\theta_Z^{}}:\ov Z\la Z$ to obtain the diagram
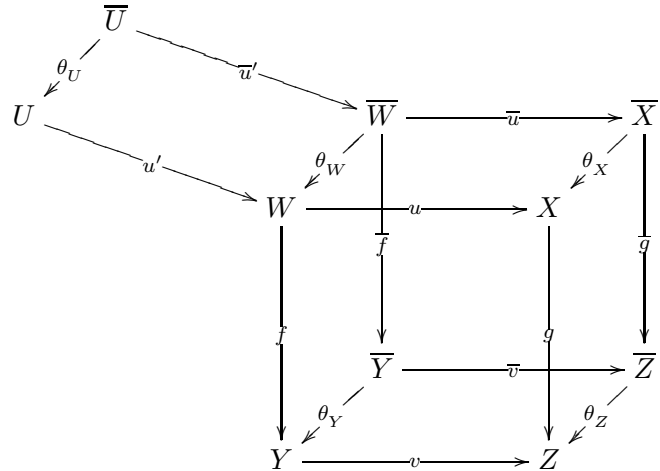
\[
\xymatrix@R-5pt@C-5pt{ 
 & \ov U\ar[rrrd]|{\ov u'}\ar[dl]|{{\theta_U^{}}} &&& &&& \\
U\ar[rrrd]|{u'} &&& &\ov W \ar[rrr]|{\ov u}\ar[ddd]|{\ov
  f}\ar[dl]|{{\theta_W^{}}} && & \ov X\ar[ddd]|{\ov g}\ar[dl]|{{\theta_X^{}}}\\
&&& W \ar[rrr]|{u}\ar[ddd]|{f} &&&  X\ar[ddd]|{g}& \\
\\
& &&& \ov Y \ar[dl]|{{\theta_Y^{}}} \ar[rrr]|{\ov v} &&& \ov Z \ar[dl]|{{\theta_Z^{}}}\\
&&& Y \ar[rrr]|{v} &&& Z &
}
\]
We are assuming that $\Dqc(X)$ is compactly generated, that $\ov Z$
satisfies Thomason's condition (and in particular $\Dqc(\ov Z)$ is
compactly generated), while $Z$ is assumed to have quasi-affine
diagonal. From Corollary~\ref{C4.7993} we deduce that 
$\Dqc(\ov
X)$ 
is compactly generated. 
The map ${\theta_Z^{}}$ is assumed representable, hence so are its
pullbacks ${\theta_U^{}}$, ${\theta_W^{}}$, ${\theta_X^{}}$ and ${\theta_Y^{}}$. Being
representable these maps are also concentrated. The map $\ov g$ is the pullback
of the pseudo-coherent map $g$ which factors as $X\stackrel{j}\la\pp(\cs')\la Z$,
hence $\ov g$ is pseudo-coherent and factors as 
$\ov X\stackrel{\ov j}\la\pp(\ov\cs')\la \ov Z$. Since
$\ov Z$ has the resolution property Reminder~\ref{R4.911} permits us to
factor $\ov g$ further as
$\ov X\stackrel{\ov j}\la\pp(\ov\cs')\stackrel{\ov{j}}\la\pp(\cv')\la
\ov Z$,
with $\cv'$ a finite-rank vector bundle on $\ov Z$. Lemma~\ref{L4.9} now
applies to the diagram
\[
\CD
\ov U @>\ov u'>>\ov W @>\ov u>> \ov X\\
@.  @V\ov fVV @VV\ov gV \\
 @.  \ov Y @>\ov v>> \ov Z
\endCD
\]
and we deduce that 
\be
\item
The map ${\ov u'}^*\ov\Phi:{\ov u'}^*\ov u^*\ov g^\times
\la {\ov u'}^*\ov f^\times\ov v^*$
is an isomorphism.
\item
The isomorphic functors 
${\ov u'}^*\ov u^*\ov g^\times
\cong {\ov u'}^*\ov f^\times\ov v^*$
respect coproducts.
\item
In the case where the stacks are all noetherian,
if $z\in\Dqc(\ov Z)$
is an object then 
${\ov u'}^*\ov u^*\ov g^\times z
\cong {\ov u'}^*\ov f^\times\ov v^*z=0$ 
if and only if ${\ov u'}^*\ov u^*\ov g^*z
\cong {\ov u'}^*\ov f^*\ov v^*z=0$.
\setcounter{enumiii}{\value{enumi}}
\ee
We will not use (iii) in the proof of the current Lemma, but will refer
back to it in 
the proof of the next Proposition. In fact we will be referring back
to the diagrams of this proof. 

Back to the proof of the current Lemma:
flat base-change for the
concentrated maps ${\theta_W^{}}$, ${\theta_Y^{}}$ and ${\theta_Z^{}}$ gives isomorphisms
\be
\setcounter{enumi}{\value{enumiii}}
\item
${u'}^*{\theta_{W*}^{}}\cong{\theta_{U*}^{}}{\ov u'}^*$,\ \ \ 
$u^*{\theta_{X*}^{}}\cong{\theta_{W*}^{}}\ov u^*$,\ \ \ 
$v^*{\theta_{Z*}^{}}\cong{\theta_{Y*}^{}}\ov v^*$.
\setcounter{enumiii}{\value{enumi}}
\ee
while flat base-change for the concentrated morphisms $g$ and $f$
gives isomorphisms
\be
\setcounter{enumi}{\value{enumiii}}
\item
$\theta_Z^*g_*\cong\ov g_*\theta_X^*$,\ \ \ $\theta_Y^*f_*\cong\ov
f_*\theta_W^*$, 
and taking right adjoints 
${\theta_{X*}^{}}\ov
g^\times\cong g^\times{\theta_{Z*}^{}}$,\ \ \ 
${\theta_{W*}^{}}\ov
f^\times\cong f^\times{\theta_{Y*}^{}}$.
\setcounter{enumiii}{\value{enumi}}
\ee
Combining (iv) and (v) and a diagram chase we have that
${u'}^*\Phi{\theta_{Z*}^{}}:{u'}^*u^*g^\times{\theta_{Z*}^{}}\la {u'}^*f^\times
v^*{\theta_{Z*}^{}}$ is isomorphic to ${\theta_{U*}^{}}{\ov u'}^*\ov\Phi:
{\theta_{U*}^{}}{\ov u'}^*\ov u^*\ov g^\times
\la {\theta_{U*}^{}} {\ov u'}^*\ov f^\times\ov v^*$. By (i) we know
that ${\ov u'}^*\ov\Phi$ is an isomorphism, hence so is
${\theta_{U*}^{}}{\ov u'}^*\ov\Phi\cong
{u'}^*\Phi{\theta_{Z*}^{}}$. From (ii) we know that the isomorphic
functors 
${\ov u'}^*\ov u^*\ov g^\times
\cong {\ov u'}^*\ov f^\times\ov v^*$ 
respect coproducts, as does the
functor ${\theta_{U*}^{}}$ for the concentrated morphism
$\theta_U$. Hence the isomorphic functors ${\theta_{U*}^{}}{\ov u'}^*\ov u^*\ov g^\times
\cong {\theta_{U*}^{}} {\ov u'}^*\ov f^\times\ov v^*$ respect
coproducts,
but they are also isomorphic to
${u'}^*u^*g^\times{\theta_{Z*}^{}}\cong {u'}^*f^\times
v^*{\theta_{Z*}^{}}$. 
\eprf

\pro{P4.9112233}
Let the conventions be as in Notation~\ref{N4.9.99312}, but now assume 
that $\theta_Z^{}:\ov Z\la Z$ is an \'etale \emph{cover} 
with $\ov Z$ satisfying 
the resolution property.
Then part (iii) of 
Notation~\ref{N4.9.99312} is true, as are the following strengthenings of 
(i) and (ii) of Notation~\ref{N4.9.99312}:
\be
\item
 The natural transformation ${u'}^*\Phi :
{u'}^*u^*g^\times \la {u'}^*f^\times v^*$ is an isomorphism.
\item
The isomorphic functors ${u'}^*u^*g^\times\cong {u'}^*f^\times v^*$
respect coproducts.
\setcounter{enumiii}{\value{enumi}}
\ee
\epro

\prf
By Lemma~\ref{L4.93115} we have that
${u'}^*\Phi\theta_{Z*}^{}$ is an isomorphism 
and the isomorphic functors ${u'}^*u^*g^\times\theta_{Z*}^{}\la {u'}^*f^\times
v^*\theta_{Z*}^{}$ 
respect coproducts. If $P\in\Dqc(\ov Z)$ is any object, then 
$P\oo\theta_{Z}^{*}(-)$ is a functor $\Dqc(Z)\la\Dqc(\ov Z)$ respecting
coproducts, and hence
\be
\setcounter{enumi}{\value{enumiii}}
\item
The natural transformation 
${u'}^*\Phi\theta_{Z*}^{}\big[P\oo\theta_{Z}^{*}(-)\big]$
is an isomorphism.
\item
The isomorphic functors 
${u'}^*u^*g^\times\theta_{Z*}^{}\big[P\oo\theta_{Z}^{*}(-)\big]\,\,\cong\,\, {u'}^*f^\times
v^*\theta_{Z*}^{}\big[P\oo\theta_{Z}^{*}(-)\big]$ respect coproducts.
\setcounter{enumiii}{\value{enumi}}
\ee
The projection formula for the concentrated morphism 
$\theta_Z^{}:\ov Z\la
Z$
 gives an isomorphism 
$\theta_{Z*}^{}\big[P\oo\theta_{Z}^{*}(-)\big]\,\,\cong\,\,(\theta_{Z*}^{}P)\oo(-)$,
and (iii) and (iv) rewrite as
\be
\setcounter{enumi}{\value{enumiii}}
\item
The natural transformation ${u'}^*\Phi\big[(\theta_{Z*}^{}P)\oo(-)\big]$
is an isomorphism.
\item
The isomorphic functors 
${u'}^*u^*g^\times\big[(\theta_{Z*}^{}P)\oo(-)\big]\,\,\cong\,\,{u'}^*f^\times
v^*\big[(\theta_{Z*}^{}P)\oo(-)\big]$ 
respect coproducts.
\setcounter{enumiii}{\value{enumi}}
\ee
Let $\ch$ be the full subcategory of $\Dqc(Z)$ of all objects $E \in\Dqc(Z)$ so
that
\be
\setcounter{enumi}{\value{enumiii}}
\item
The natural transformation ${u'}^*\Phi(E\oo -)$
is an isomorphism.
\item
The isomorphic functors 
${u'}^*u^*g^\times(E\oo -)\,\,\cong\,\,{u'}^*f^\times
v^*(E\oo -)$ respect coproducts.
\setcounter{enumiii}{\value{enumi}}
\ee
By (v) and (vi) we know that, 
for any object $P\in\Dqc(\ov Z)$, the
object $\theta_{Z*}^{}P$ belongs to $\ch$.
Clearly $\ch$ is a thick subcategory of $\Dqc(Z)$. 
By \cite[Lemma~A.1]{Hall-Rydh13} we have that every 
compact object in
$\Dqc(Z)$ lies in $\ch$.

Let $E$ be a compact object of $\Dqc(Z)$ and let $F\in\Dqc(Z)$
be arbitrary. By Lemma~\ref{L27.9} we have a commutative diagram
\[
\xymatrix@C+30pt{
  u^*g^*E \oo u^*g^\times F \ar[d]_-{\tau^{-1}\oo\Phi} &
    u^*(g^*E\oo g^\times F)\ar[l]_-{\mu_u^{}}
    \ar[r]^-{u^*\chi(g,E,F)}&  u^*g^\times(E\oo F)\ar[d]^-{\Phi} \\
    f^*v^*E\oo f^\times v^* F\ar[r]^-{\chi(f,v^*E,v^*F)} &
    f^\times(v^*E\oo v^*F) &
    f^\times v^*(E\oo F)\ar[l]_-{f^\times\mu_v^{}}
}\]
The horizontal maps are isomorphisms; for $\mu_u^{}$ and $f^\times\mu_v^{}$
this is obvious, and for $\chi(g,E,F)$ and $\chi(f,v^*E,v^*F)$
we note that $E$ is compact, hence perfect, hence $v^*E$ is also
perfect, and Proposition~\ref{PThom.73.66} applies. Applying
the functor ${u'}^*$ we obtain the commutativity
of the right-hand square in
\[
\xymatrix@C+30pt{
  {u'}^*u^*g^*E \oo {u'}^*u^*g^\times F \ar[d]_-{{u'}^*\tau^{-1}\oo{u'}^*\Phi} \ar[r]^-\cong &{u'}^*(u^*g^*E \oo u^*g^\times F) \ar[d]_-{{u'}^*(\tau^{-1}\oo\Phi)} \ar[r]^-\cong &  {u'}^*u^*g^\times(E\oo F)\ar[d]^-{{u'}^*\Phi} \\
    {u'}^*f^*v^*E\oo {u'}^*f^\times v^* F\ar[r]^-\cong &
    {u'}^*(f^*v^*E\oo f^\times v^* F)\ar[r]^-\cong &
    {u'}^*f^\times v^*(E\oo F)
}\]
and the left-hand square commutes by the naturality of the oplax
structure map $\mu_{{u'}^*}^{}:{u'}^*(A\oo B)\la{u'}^*A\oo{u'}^*B$.
The fact that $\ch$ contains all the compact objects $E$ says that the
vertical map on the right is an isomorphism and
the isomorphic functors respect coproducts. The commutativity says
that the same is true on the left.
Therefore define the subcategory $\cl\subset\Dqc(Z)$ to contain all
the objects $E\in\Dqc(Z)$ so that
\be
\setcounter{enumi}{\value{enumiii}}
\item
The natural transformation ${u'}^*u^*g^*(E)\oo{u'}^*\Phi(-)$
is an isomorphism.
\item
The isomorphic functors 
${u'}^*u^*g^*(E)\oo{u'}^*u^*g^\times(-)\,\,\cong\,\,{u'}^*u^*g^*(E)\oo{u'}^*f^\times
v^*(-)$ respect coproducts.
\setcounter{enumiii}{\value{enumi}}
\ee
The above argument showed that $\Dqc(Z)^c$ is contained in $\cl$, while
$\cl$ is obviously localizing. As $\Dqc(Z)$ is compactly generated
we conclude that $\cl=\Dqc(Z)$, in particular $\co_Z^{}\in\cl$.
Therefore (i) and (ii) follow.

It remains to prove that, under the hypotheses
of Proposition~\ref{P4.9112233}, 
part (iii) of Notation~\ref{N4.9.99312} is true.
Let $\cl_1^{},\cl_2^{}\subset\Dqc(Z)$ be the full 
subcategories annihilated (respectively)
by the functors ${u'}^*u^*g^\times$ and ${u'}^*u^*g^*$; we need to prove that
$\cl_1^{}=\cl_2^{}$. With the notation as in the proof
of Lemma~\ref{L4.93115},  part (iii) of the proof
of Lemma~\ref{L4.93115} tells us that $\cm_1^{}=\cm_2^{}$, 
where $\cm_1^{}\subset\Dqc(\ov Z)$
is full subcategory annihilated by ${\ov u'}^*\ov u^*\ov g^\times$ 
while $\cm_1^{}\subset\Dqc(\ov Z)$ is the full subcategory annihilated
by  ${\ov u'}^*\ov u^*\ov g^*$. 
The idea will be to relate $\cl_1^{}$ to $\cm_1^{}$ and $\cl_2^{}$ to $\cm_2^{}$. 

Recall the diagrams in the proof of Lemma~\ref{L4.93115}:
because the map $\theta_U^{}:\ov U\la U$ is finite-type, 
representable, separated
 and \'etale it
is quasi-affine, hence $\theta_{U*}^{}$ is conservative. Thus the kernel 
$\cm_1^{}$
(respectively $\cm_2^{})$
of ${\ov u'}^*\ov u^*\ov g^\times$ 
(respectively ${\ov u'}^*\ov u^*\ov g^*$) is equal to the kernel
of $\theta_{U*}^{}{\ov u'}^*\ov u^*\ov g^\times$ 
(respectively $\theta_{U*}^{}{\ov u'}^*\ov u^*\ov g^*$). In the proof of
Lemma~\ref{L4.93115} we saw that 
$\theta_{U*}^{}{\ov u'}^*\ov u^*\ov g^\times
\cong{u'}^*u^*g^\times\theta_{Z*}^{}$,
 while the isomorphism
$\theta_{U*}^{}{\ov u'}^*\ov u^*\ov g^*
\cong{u'}^*u^*g^*\theta_{Z*}^{}$ 
is by Tor-independent base change
($\theta_Z^{}$ is \'etale, in particular flat).
This gives
\be
\setcounter{enumi}{\value{enumiii}}
\item
For $i\in\{1,2\}$ we have 
$\cm_i=\theta_{Z*}^{-1}\cl_i$, and hence
$\theta_{Z*}^{}\cm_i=\theta_{Z*}^{}\theta_{Z*}^{-1}\cl_i\subset\cl_i$. 
\setcounter{enumiii}{\value{enumi}}
\ee

Notation~\ref{NThom.73.63}(iv) gives us maps natural in
$E,F\in\Dqc(Z)$ 
of the form
$\chi(g,E,F):g^*E\oo g^\times F\la g^\times(E\oo F)$,
and these induce maps 
${u'}^*u^*\chi(g,E,F):{u'}^*u^*g^*E\oo {u'}^*u^*g^\times F\la 
{u'}^*u^*g^\times(E\oo F)$.
Let $\cs\subset\Dqc(Z)$ be
the full subcategory
\[
\cs=\{E\in\Dqc(Z)\mid {u'}^*u^*\chi(g,E,F)\text{ is an isomorphism for
every }F\in\Dqc(Z)\}.
\]
Proposition~\ref{PThom.73.66} tells us that the map $\chi(g,E,F)$
is an isomorphism whenever $E$ is compact, hence
${u'}^*u^*\chi(g,E,F)$ is also an isomorphism when $E$ is compact---in
other words $\Dqc(Z)^c$ is contained in $\cs$.
On the other hand 
by Proposition~\ref{P4.9112233}(ii) we know that the functor 
${u'}^* u^*g^\times$ 
respects coproducts, and this makes ${u'}^*u^*\chi(g,E,F)$ a natural
transformation between functors both of which respect coproducts (in
both the variable $E$ and the variable $F$). It follows that $\cs$ is
a localizing subcategory, and since $\Dqc(Z)^c\subset\cs$ and 
$\Dqc(Z)$ 
is compactly generated we have that $\cs=\Dqc(Z)$. That is
${u'}^*u^*\chi(g,E,F)$ is an isomorphism for all $E,F\in\Dqc(Z)$. 
If $F\in\cl_1^{}$,
that is if ${u'}^*u^*g^\times F=0$, the isomorphism
${u'}^*u^*g^\times(E\oo F)\cong {u'}^*u^*g^*E\oo
{u'}^*u^*g^\times F$ tells us that ${u'}^*u^*g^\times(E\oo F)=0$,
that is $E\oo F\in\cl_1^{}$. Thus $\cl_1^{}$ is not only localizing,
it is also a tensor ideal. Because ${u'}^*u^*g^*$ is
a strong monoidal functor its kernel $\cl_2^{}$ is also a localizing
tensor ideal. 

For $i\in\{1,2\}$ we have the inclusions
\[
\theta_{Z*}^{}\big[\Dqc(\ov Z)\oo\theta_{Z}^{*}\cl_i\big]\quad=\quad
\big[\theta_{Z*}^{}\Dqc(\ov Z)\big]\oo\cl_i\quad\subset\quad\cl_i\ ,
\]
where the equality is by the projection formula and the inclusion
because $\cl_i$ is a tensor ideal. From (xi) it follows that
$\Dqc(\ov Z)\oo\theta_{Z}^{*}\cl_i\subset\theta_{Z*}^{-1}\cl_i=\cm_i$. But then
\be
\setcounter{enumi}{\value{enumiii}}
\item
$\quad
\big[\theta_{Z*}^{}\Dqc(\ov Z)\big]\oo\cl_i\quad=\quad
\theta_{Z*}^{}\big[\Dqc(\ov Z)\oo\theta_{Z}^{*}\cl_i\big]\quad\subset\quad
\theta_{Z*}^{}\cm_i$\ .
\setcounter{enumiii}{\value{enumi}}
\ee
 Combining (xi) and  (xii) we have inclusions
\[
\big[\theta_{Z*}^{}\Dqc(\ov Z)\big]\oo\cl_i\quad\subset\quad
\theta_{Z*}^{}\cm_i\quad\subset\quad
\cl_i\ ,
\]
and if we let $\text{\rm Loc}(-)$ be the 
operation taking a subcategory of $\Dqc(Z)$
to the localizing subcategory it generates, then
\[
\text{\rm Loc}\big[\theta_{Z*}^{}\Dqc(\ov Z)\big]\oo\cl_i
\quad\subset\quad\text{Loc}(\theta_{Z*}^{}\cm_i)\quad\subset\quad
\cl_i\ .
\]
By \cite[Lemma~A.1]{Hall-Rydh13} we have that 
$\co_Z^{}\in\text{\rm Loc}\big[\theta_{Z*}^{}\Dqc(\ov Z)\big]$, and
so $\text{\rm Loc}(\theta_{Z*}^{}\cm_i)=\cl_i$.
The equality $\cm_1^{}=\cm_2^{}$ now gives $\cl_1^{}=\cl_2^{}$.
\eprf

In order to use Proposition~\ref{P4.9112233} we need a result that produces
for us maps $g:X\la Z$ which factor as $X\stackrel j\la\pp(\cs')\la
Z$. 
The following little Lemma is what we will use.

\lem{L4.13.0.7}
Let $\alpha:V\la X$ be a finite-type morphism of quasi-compact,
quasi-separated stacks. Assume $V$ is 
a quasi-affine scheme and $X$ has quasi-affine diagonal. Then
the map $\alpha:V\la X$ factors as 
$V\stackrel\beta\la\pp(\cs)\stackrel\gamma\la X$, with $\beta$
a locally closed immersion and $\cs$ a  finitely presented
quasicoherent sheaf on $X$.
\elem

\prf
Since $V$ is quasi-affine and  $X$ has
quasi-affine diagonal the map $\alpha$ must be
quasi-affine. The sheaf $\alpha_*\co_V^{}$ is a quasicoherent
sheaf on $X$, and can be expressed 
as a filtered direct limit of finitely presented quasicoherent
sheaves: that is  
$\alpha_*\co_V^{}=\colim\,\cs_\lambda$ with each $\cs_\lambda$
finitely presented. For each  $\lambda$ the map
$\alpha:V\la X$ factors as 
$V\stackrel{i_\lambda^{}}\la
\ak(\cs_\lambda)\la X$, 
where $\ak(\cs_\lambda) $ is Spec of the
symmetric algebra on the sheaf $\cs_\lambda$. I assert 
\be
\item
We may choose $\lambda$ so
that $i_\lambda^{}$ is a locally closed immersion. 
\ee
Assuming (i) we let
$\cs=\cs_\lambda\oplus\one$,
where $\one$ is the trivial line
bundle on $X$.
Then the following factorization does the trick:
$V\stackrel{i_\lambda^{}}\la \ak(\cs_\lambda)\subset\pp(\cs)\la X$
with $\ak(\cs_\lambda)\subset\pp(\cs)$ the standard open immersion.

It remains to prove (i). 
Fortunately the question is local in the flat topology. Choose an
affine
scheme $S$ of finite type over $X$ and a faithfully flat map $\s:S\la X$,
and form the pullback square
\[
\CD
V' @>>> V \\
@V\alpha'VV @VV\alpha V \\
S @>\s>> X
\endCD
\]
Because $\alpha$ is quasi-affine and of finite type so is the
pullback $\alpha'$. 
We have that the map $\alpha'$ is a finite-type morphism from a
quasi-compact open subset of an affine scheme to an affine scheme, and
that
$\alpha'_*\co_{V'}^{}$ is the filtered colimit of sheaves
$\s^*\cs_\lambda$.
It is an easy exercise to show
that we may choose an $\cs_\lambda$ so that, in the
factorization $V'\stackrel{i'_\lambda}\la \ak(\s^*\cs_\lambda)\la S$,
the map $i'_\lambda$ is a locally closed immersion. 
\eprf

Until now we have given 
almost all of our arguments in non-noetherian, 
general forms. The next theorem is the
 first point at which I have no
idea how to proceed without the noetherian hypothesis.

\thm{T4.13}
As in Notation~\ref{N4.1}
suppose we are given a 2-cartesian square  of \emph{noetherian} stacks
\[
\CD
W @>u>> X\\
@VfVV @VVgV \\
Y @>v>> Z
\endCD
\]
Assume that $g$ is of finite type, separated, concentrated
and universally quasi-proper,  
and that $v$ is flat. Assume that $W$, $X$, $Y$
and $Z$ 
have quasi-affine diagonals. Furthermore assume that there exist
finite type, separated, representable \'etale covers 
$\ov W\la W$, $\ov X\la X$
and $\ov Z\la Z$, with $\ov W$, $\ov X$ and $\ov Z$ 
having the resolution property. 
Suppose that $\pp(\cs)$ 
satisfies Thomason's condition for any 
coherent sheaf $\cs$
on any of the stacks $\ov W$, $\ov X$, $\ov Y$ or $\ov Z$ 
which admit
morphisms to $W$, $X$, $Y$
or $Z$ which are finite-type, 
separated, representable and \'etale. Let 
$\ov U\subset W$ be the open subset of all the points at which $f$ is of finite
Tor-dimension.
Consider the base-change map 
$\Phi:u^*g^\times\la f^\times v^*$; 
if $\rho:U\la W$ is a morphism whose image lies in the open 
subset $\ov U\subset W$, then $\rho^*\Phi:
\rho^*u^*g^\times\la \rho^*f^\times v^*$ is an isomorphism, and
the isomorphic functors $\rho^*u^*g^\times\cong\rho^*f^\times v^*$
respect coproducts.
\ethm

\prf
The openness of the subset $\ov U\subset W$ is by Corollary~\ref{C4.497.-111},
which applies because any finite-type morphism $f:W\la Y$ of noetherian
stacks is pseudo-coherent.
The assertion of the Theorem is local in the flat topology on $W$: it suffices to
produce a single flat morphism $\rho:U\la W$, whose image
is all of $\ov U$,
so that $\rho^*\Phi$ is an isomorphism and
the isomorphic functors $\rho^*u^*g^\times\cong\rho^*f^\times v^*$
respect coproducts. We will now do the construction.

Choose a smooth surjection $\alpha:V\la X$ with $V$
an affine scheme, and apply Lemma~\ref{L4.13.0.7} to the
maps $\alpha:V\la X$ and $g\alpha:V\la Z$.
We obtain commutative squares of finite-type maps
\[
\CD
V @>\alpha>> X @. @.  @. V @>\alpha>> X\\
@V\beta VV    @| @. \qquad\qquad@. @V\beta'VV @VVgV\\
\pp(\cs') @>\gamma>> X @. @.  @.\pp(\ct') @>\gamma'>> Z
\endCD
\]
with $\cs'$ and $\ct'$  
coherent sheaves
on $Z$ and $X$ respectively, and with 
$\beta$ and $\beta'$ locally closed immersions.
The map $g$ is concentrated by hypothesis, while $\alpha$, $\beta$,
$\beta'$, $\gamma$ and
$\gamma'$ are representable, hence concentrated. 
Form the diagram with  2-cartesian squares
\[
\xymatrix{
\pp(\cs')\times_X^{}\pp(g^*\ct') \ar[rr] \ar[d]^\pi&& \pp(g^*\ct') \ar[r] \ar[d]& \pp(\ct')\ar[d] \\
\pp(\cs')  \ar[rr]& &X \ar[r]^{g}  &Z
}
\]
The two horizontal maps in the bottom row are separated and
universally quasi-proper---the morphism
$g$ by hypothesis, the morphism $\pp(\cs')\la X$ by
construction. Hence their pullbacks, the horizontal maps on the top
row, are also separated and universally quasi-proper. We have an obvious map
$V\la \pp(\cs')\times_X^{}\pp(g^*\ct') $, and we know that the two
composites
\[
\xymatrix{
V\ar[rr] &&  \pp(\cs')\times_X^{}\pp(g^*\ct') \ar[rr] \ar[d]^\pi&& \pp(\ct') \\
 && \pp(\cs') & &
}
\]
are locally closed immersions. On the other hand the map $\pi$ is
separated and representable, hence the morphism 
$V\la
\pp(\cs')\times_X^{}\pp(g^*\ct') $
 is also a locally closed
immersion. Let $X''$ be the stack-theoretic closure of $V$ in
$\pp(\cs')\times_X^{}\pp(g^*\ct') $ and let $X'$ be the
stack-theoretic closure of $V$ in
$\pp(\ct')$.

This gives a 2-commutative square of concentrated morphisms
\[
\CD
X'' @>g'>> X'\\
  @Vq'VV @VVqV \\
 X @>g>> Z
\endCD
\]
Here $X'$ and $X''$ are two compactifications of $V$, while the map
$V\la X$ is faithfully flat.
The maps $g$ and $g'$ are separated and universally 
quasi-proper, while the morphisms
$q$ and $q'$ have factorizations $X'\la  \pp(\ct)\la Z$ and
$X'' \la\pp(\cs)\la X$, with 
$X'\la  \pp(\ct)$ and $X'' \la\pp(\cs)$
closed immersions. If
$V'$ is the pullback in the square 
\[
\CD
V' @>>> V \\
@VVV  @VVV\\
W @>u>> X
\endCD
\]
then we obtain diagrams where all the squares are 2-cartesian
\[
\xymatrix@C+10pt@R+10pt{
V' \ar[r]^i &  W''\ar@{}[dr]|{(1)}\ar[d]_{p'} \ar[r]^{u''} & 
   X''\ar[d]^{q'} & &    
   V' \ar@{}[dr]|{(2)} \ar[r]^i \ar@{=}[d]&
        W'' \ar@{}[dr]|{(3)}\ar[r]^{u''}\ar[d]_{f'}&
      X''\ar[d]^{g'}\\
 &  W \ar@{}[dr]|{(0)}\ar[d]_f\ar[r]^u & X \ar[d]^g    &\text{and} &   
           V'  \ar[r]^{f'i}&  W'\ar@{}[dr]|{(4)}\ar[d]_p \ar[r]^{u'} 
        & X'\ar[d]^q \\
& Y \ar[r]^v & Z    & & &     Y \ar[r]^v & Z
}
\]
The horizontal maps are all flat, in fact the maps $i$ and $f'i$
are even open immersions. The vertical maps are all quasi-proper maps
of noetherian stacks.
All the vertical morphisms are concentrated by construction.
The stacks $W$, $X$, $Y$ and $Z$ are assumed to have quasi-affine diagonals,
and the stacks $W'$ and $X'$ are closed substacks of 
$\pp(\cs)$ and $\pp(\cs')$,
with $\cs$ a coherent sheaf 
 on $Y$ and $\cs'$ a coherent sheaf on $Z$. Hence
$W'$ and $X'$ also have quasi-affine diagonals. For any 
coherent sheaf 
 $\cs$ on any $\ov V$, admitting 
a map to $W$, $X$, $Y$ or $Z$ which is finite-type, 
representable, separated
and \'etale, we assumed that  $\pp(\cs)$
satisfies Thomason's condition. The stacks $W'$ and $X'$ 
 were constructed as closed substacks of 
$\pp(\cs)$
and $\pp(\cs')$, with $\cs$ a coherent sheaf
 on $Y$ and $\cs'$ a coherent sheaf
 on $Z$, and Observation~\ref{O4.2.1999} guarantees that
$W'$ and $X'$ satisfy Thomason's condition.
Finally $X''$ is a closed 
substack of some $\pp(\cs)$, where $\cs$ is a 
coherent sheaf on 
$X$---hence $\Dqc(X'')$ is compactly
generated.

Let 
$\Phi_3^{}:{u''}^*{g'}^\times\la
{f'}^\times{u'}^*$ 
be the base-change map of the square
labeled $(3)$ in the diagram. From the diagram
\[
\xymatrix@C+10pt@R+10pt{ 
   V'  \ar[r]^i \ar@{}[dr]|{(2)}\ar@{=}[d]&  W'' \ar@{}[dr]|{(3)}\ar[r]^{u''}\ar[d]_{f'}&
      X''\ar[d]^{g'}\\  
           V'  \ar[r]^{f'i}&  W'\ar[r]^{u'} 
        & X' \\
}
\]
and Lemma~\ref{L4.8} we learn that 
$i^*\Phi_3^{}:i^*{u''}^*{g'}^\times\la
i^*{f'}^\times{u'}^*$ is an isomorphism.

Since the map $q$ has a factorization as $X'\la\pp(\cs')\la Z$,
and $Z$ admits a finite-type, 
representable, separated and \'etale morphism 
$\ov Z\la Z$ with $\ov Z$ satisfying the resolution
property, we can
apply Proposition~\ref{P4.9112233} to the square labeled 
$(4)$ in the diagram.
Let $\Phi_4:{u'}^*q^\times\la p^\times v^*$ be its base-change map;
we learn that any flat morphism $\e:U\la W'$, whose image is
contained in the open set of points $x\in W'$ on which $p$ has 
finite Tor-dimension, yields an isomorphism
$\e^*\Phi_4^{}:\e^*{u'}^*q^\times\la \e^*p^\times v^*$,
and the isomorphic functors
$\e^*{u'}^*q^\times\cong\e^*p^\times v^*$
respect coproducts. We are given a flat
map $f'i:V'\la W'$, let $U=(f'i)^{-1}F$ be the inverse image
of the open substack $F\subset W'$ on which the map $p$
is of finite Tor-dimension, and let $\s:U\la V'$ be the open
immersion. We may apply the above to $\e=f'i\s:U\la W'$, and deduce
that $\s^*i^*{f'}^*\Phi_4$ is an isomorphism and
$\s^*i^*{f'}^*p^\times v^*$ respects coproducts. 
 Consider the square labeled $(2)$; 
together with Lemma~\ref{L4.7} this square tells  us that the base-change 
map 
$\Phi_2^{}:i^*{f'}^\times\la
\id^\times{(f'i)}^*$ is an isomorphism, and thus $i^*{f'}^*$ is isomorphic
to $i^*{f'}^\times$. Combining with the above we have that
$\s^*i^*{f'}^\times\Phi_4^{}\cong \s^*i^*{f'}^*\Phi_4$ is an isomorphism,
and the functor $\s^*i^*{f'}^\times p^\times v^*\cong\s^*i^*{f'}^*p^\times v^*$
respects coproducts.
Hence $\s^*i^*$ takes each of the 
two composable morphisms
\[
\CD
{u''}^*{g'}^\times q^\times @>\Phi_3^{}q^\times>> 
{f'}^\times{u'}^*q^\times @>{f'}^\times\Phi_4^{}>> {f'}^\times p^\times v^*
\endCD
\]
to an  isomorphism. But the composite is equal to the composite
\[
\CD
{u''}^*{q'}^\times g^\times @>\Phi_1^{}g^\times>> 
{p'}^\times{u}^*g^\times @>{p'}^\times\Phi>> {p'}^\times f^\times v^*
\endCD
\]
where $\Phi$ is the base-change map of the square labeled $(0)$ and
$\Phi_1$ is the base-change map of the square labeled $(1)$. 
In the case of the square labeled $(1)$ the map $q':X''\la X$ factors
as $X''\la\pp(\cs)\la X$, where $\cs$ is a coherent
sheaf on 
$X$, where the map $X''\la\pp(\cs)$ is a closed immersion,
and where $X$ admits a finite-type,
separated, representable and \'etale
cover $\ov X\la X$ with $\ov X$ satisfying the resolution 
property. Hence we may
apply Proposition~\ref{P4.9112233}. We have an open immersion $i:V'\la W''$, and
the composite $V'\stackrel i\la W''\stackrel {p'}\la W$ is the pullback of
the faithfully flat map $V\la X$. Therefore $i$ maps $V'$ into a substack of
$W''$ where $p'$ is even flat---certainly of finite Tor-dimension.
Hence $i^*\Phi_1^{}$ is an isomorphism and the isomorphic
functors $i^*{p'}^\times u^*\cong i^*{u''}^*q^\times$
respect coproducts. Since $i^*\Phi_1^{}$ is an isomorphism
so is $\s^*i^*\Phi_1^{}g^\times$.
Summarizing what we have so far:
\be
\item
$\s^*i^*$ takes ${p'}^\times\Phi$ to an isomorphism, or more
simply that $\s^*i^*{p'}^\times$ takes $\Phi$ to an isomorphism.
\item
 The 
functor $\s^*i^*{p'}^\times f^\times v^*\cong\s^*i^*{p'}^\times u^* g^\times $
respects coproducts because it is isomorphic to the functor
$\s^*i^*{f'}^\times p^\times v^*$ which was proved to respect
coproducts above. Concretely: if
$\{\cs_\lambda^{},\,\lambda\in\Lambda\}$ is a set of objects
in the category $\Dqc(Z)$, the natural map
\[
\CD
 \coprod_{\lambda\in\Lambda}\s^*i^*{p'}^\times u^* g^\times\cs_\lambda^{} @>>>
\s^*i^*{p'}^\times u^*g^\times\left(\coprod_{\lambda\in\Lambda}\cs_\lambda^{}\right)
\endCD
\]
is an isomorphism
\item
The functor $\s^*i^*{p'}^\times u^*$  respects coproducts
since $\s^*$ obviously does and $i^*{p'}^\times u^*$ was proved to
respect coproducts above.
This certainly implies that if
$\{\cs_\lambda^{},\,\lambda\in\Lambda\}$ is a set of objects
in the category $\Dqc(Z)$, the natural map
\[
\CD
 \coprod_{\lambda\in\Lambda}\s^*i^*{p'}^\times u^* g^\times\cs_\lambda^{} @>>>
\s^*i^*{p'}^\times u^*\left(\coprod_{\lambda\in\Lambda}g^\times\cs_\lambda^{}\right)
\endCD
\]
is an isomorphism.
\setcounter{enumiii}{\value{enumi}}
\ee
Let $\ph$ be the natural map
\[
\CD
 \ph\ :\  u^* \left(\coprod_{\lambda\in\Lambda}g^\times\cs_\lambda^{} \right)@>>>
u^*g^\times\left(\coprod_{\lambda\in\Lambda}\cs_\lambda^{}\right)
\endCD
\]
and consider the composite
 \[
\CD
 \coprod_{\lambda\in\Lambda}\s^*i^*{p'}^\times u^* g^\times\cs_\lambda^{} @>\alpha>>
 \s^*i^*{p'}^\times u^*\left(\coprod_{\lambda\in\Lambda}g^\times\cs_\lambda^{}\right)
 @>\s^*i^*{p'}^\times\ph>>
 \s^*i^*{p'}^\times u^*g^\times\left(\coprod_{\lambda\in\Lambda}\cs_\lambda^{}\right)
\endCD
\]
In (iii) we saw that $\alpha$ is an isomorphism while in (ii) we learned
that so is the composite $(\s^*i^*{p'}^\times\ph)\circ\alpha$. Hence
\be
\setcounter{enumi}{\value{enumiii}}
\item
 The functor $\s^*i^*{p'}^\times$
takes the natural map
\[
\CD
  u^* \left(\coprod_{\lambda\in\Lambda}g^\times\cs_\lambda^{} \right)@>\ph>>
u^*g^\times\left(\coprod_{\lambda\in\Lambda}\cs_\lambda^{}\right)
\endCD
\]
to an isomorphism.
\setcounter{enumiii}{\value{enumi}}
\ee

Next we apply Proposition~\ref{P4.9112233}(iii) to the (trivial) 2-cartesian square
\[
\CD
W''  @>\id >>  W'' \\
@Vp'VV @VVp'V \\
W @>\id>> W
\endCD
\]
and the map $i\sigma:U\la W''$. The map $i\s$ is an open 
immersion (hence flat), the map $p'$ is of finite
Tor-dimension on the image of $i\s$ and admits
a factorization $W''\stackrel\alpha\la \pp(\cs)\stackrel\beta\la W$
with $\alpha$ a closed immersion, and $W$ has a finite-type,
representable, separated and \'etale cover $\ov W\la W$ with
$\ov W$ satisfying the resolution property. The
hypotheses of Proposition~\ref{P4.9112233} hold and the stacks are noetherian.
We have 
\be
\setcounter{enumi}{\value{enumiii}}
\item
If $\phi:E\la E'$ is a morphism
in $\Dqc(W)$, by applying
assertion (iii) of Notation~\ref{N4.9.99312}
to the mapping cone of $\phi$ 
we deduce that $\sigma^*i^*{p'}^*\phi$ is an isomorphism
if and only if $\sigma^*i^*{p'}^\times \phi$ is. 
\setcounter{enumiii}{\value{enumi}}
\ee
In (i) we saw that $\sigma^*i^*{p'}^\times \Phi$ is an isomorphism for
the base-change map $\Phi:u^*g^\times\la f^\times v^*$, while in
(iv) we saw that $\sigma^*i^*{p'}^\times \ph$ is an isomorphism when
$\ph$ is the natural map 
\[
\CD
  u^* \left(\coprod_{\lambda\in\Lambda}g^\times\cs_\lambda^{} \right)@>\ph>>
u^*g^\times\left(\coprod_{\lambda\in\Lambda}\cs_\lambda^{}\right)
\endCD
\] 
By (v) we now conclude that $\sigma^*i^*{p'}^*\Phi$ and
$\sigma^*i^*{p'}^*\ph$ are also isomorphisms. Thus we have that the
map $\rho=p'i\sigma:U\la W$ satisfies 
\be
\setcounter{enumi}{\value{enumiii}}
\item
The map $\rho^*\Phi$  is an isomorphism.
\item
$\rho^*\ph$ is an isomorphism, and because $\rho^*$ and $u^*$  respect
coproducts it follows that,
for any set of objects $\{\cs_\lambda^{},\,\lambda\in\Lambda\}$
in the category $\Dqc(Z)$, the natural map
\[
\CD
  \coprod_{\lambda\in\Lambda}\rho^*u^*g^\times\cs_\lambda^{} @>>>
\rho^*u^*g^\times\left(\coprod_{\lambda\in\Lambda}\cs_\lambda^{}\right)
\endCD
\] 
is an isomorphism. In other words $\rho^*u^*g^\times$ respects coproducts.
\setcounter{enumiii}{\value{enumi}}
\ee
It remains to show that
the Theorem follows. By our construction $F$ was the open subset of $W'$ where
$p$ is of finite Tor-dimension, the map $f'i$ is an open immersion $V'\la W'$,
and $U=F\cap V'$ is just the open subset of $V'$ at which the composite
$pf'i=fp'i$ is of finite Tor-dimension.
Let $\ov U\subset W$ be the open set at which $f:W\la Y$ is of
finite Tor-dimension; by construction the map $p'i:V'\la W$ is faithfully flat,
and hence $U=(p'i)^{-1}\ov U$ and $p'iU=\ov U$. The image of the flat
map $\rho:U\la W$ is precisely $\ov U\subset W$, and the Theorem now follows
from (vi), (vii) and the first paragraph of the proof.
\eprf

\rmk{R4.10.10.10.3}
In the proof
of Theorem~\ref{T4.13} 
there are two points at which we appealed to the 
noetherian hypothesis:
it entered when we cited Proposition~\ref{P4.9112233}(iii)  in the proof
of assertion~(i) above, the proof of Proposition~\ref{P4.9112233}(iii)
used support theory---in
Remark~\ref{R4.973885} we already
mentioned the failure of support theory
in the non-noetherian situation. But we
also relied on the noetherian hypothesis when we used compactifications. And it 
is the compactifications where the hypothesis 
seems crucial.

All our lemmas 
proving that base-change is an isomorphism assumed that the vertical
maps are quasi-proper, and the geometric way to achieve this is 
with proper, pseudo-coherent maps. But when we start with a 
pseudo-coherent, proper map
$g:X\la Z$, and then produce $q:X'\la Z$, $q':X''\la X$ and $g':X''\la X'$
by taking the closure of some $V$ in a suitable $\pp(\cs)$, then the
maps constructed will be proper but there is no 
reason to expect them
to be pseudo-coherent. 
Taking compactifications cannot, as far as I know, be done in a way
that achieves pseudo-coherence. And this problem will plague
the rest of the article.

And as Example~\ref{E4.7.53} and Remark~\ref{R4.7.9997} show,
base-change fails to be an isomorphism without something
like quasi-properness.
\ermk

\section{Elementary properties of Nagata compactifications}
\label{SNag}

In Notation~\ref{N0.1} we defined \emph{Nagata compactifications} in
the 2-category of algebraic stacks: they are pairs of composable
1-morphisms
$X\stackrel u\la \ov X\stackrel p\la Y$, with $p$ of finite type and
universally quasi-proper
and $u$ a dominant, flat monomorphism. It might be useful to
assemble together in this section the elementary properties of the maps
$u$ and $p$ in a Nagata compactification. We begin with

\lem{LNag.-1}
Let $f:X\la Y$ be a 1-morphism of noetherian stacks, and
assume $f$ is concentrated, of finite type and quasi-proper.
Then $f$ is closed.
\elem

\prf
Because $f$ is of finite type and the stacks are noetherian, the image
$f(Z)$ of
any closed subset $Z\subset X$ is constructible. It suffices therefore
to show that $f(Z)$ is specialization-closed; in other words we need to
show that, if $z\in Z$ is a point, then $f(Z)$ contains the closure
of $f(z)$. Choose therefore a field $k$ and a morphism $g:\spec k\la X$
whose image is $z$, and in an abuse of notation we will write
$k$ for the structure sheaf of $\spec k$, and $g_*k$ will be its
(derived) pushforward to $\Dqc(X)$.

Now $g_*k$ is supported on $Z$, hence it can be expressed as the
homotopy colimit of complexes $C_\lambda\in\dcoh(Z)$. Because $f$ is concentrated
it respects coproducts and hence homotopy colimits, and $f_*g_*k$
can be expressed as the homotopy colimit of the $f_*C_\lambda$. Since
$f_*g_*k$ is nonzero at the point $f(z)\in Y$ there must be
a $C_\lambda$ with $f_*C_\lambda$ nonzero at $f(z)$. Now because
$f$ is quasi-proper $f_*C_\lambda$ belongs to $\dcoh(Y)$ and its support is
closed. Thus the support of $f_*C_\lambda$ is a closed set in $Y$ containing
$f(z)$, and the proof will be finished if we show that the support
of $f_*C_\lambda$ is contained in the image of $Z$.

But this follows from the projection formula, which holds for the concentrated
morphism $f$: if $y\in Y$ is a point outside $f(Z)$ and $k(y)$ is any
object of $\Dqc(Y)$ supported at $y$, then
$f_*C_\lambda\oo k(y)\cong f_*\big(C_\lambda\oo f^*k(y)\big)=0$,
where the vanishing is because $C_\lambda$ is
supported on $Z$ and
$f^*k(y)$ is supported on the disjoint $f^{-1}y$.
\eprf

\lem{LNag.1}
Let $X\stackrel p\la Y\stackrel q\la Z$ be two composable 1-morphisms
of noetherian stacks with
$q$ of
finite type and separated, and with $qp$ of finite type and
universally quasi-proper.
Then $p$ is of finite type and universally
quasi-proper. Furthermore: if $qp$ is concentrated and/or separated
then so is $p$.
\elem

\prf
We are assuming $q$ separated so the diagonal map $Y\la Y\times_Z^{}Y$ is
proper and representable,
as is its pullback $X\la X\times_Z^{}Y$. Since $Y$ is
of finite type over $Z$ the stack $X\times_Z^{}Y$ is of finite type over
the noetherian stack $X$. Thus $X\la X\times_Z^{}Y$ is a proper,
representable 
morphism of noetherian algebraic stacks, hence concentrated,
separated, of finite
type and universally quasi-proper.

But the map $qp:X\la Z$ is assumed of
finite type and univerally quasi-proper, and we allow ourselves to
consider the option that it might also be concentrated and/or
separated. Hence its
pullback $X\times_Z^{}Y\la Y$ also has these properties.
Therefore the composite
$X\la X\times_Z^{}Y\la Y$ is of finite type  and universally quasi-proper,
and if $qp$ is concentrated and/or separated then so is the composite.
\eprf

To avoid possible confusion we recall

\dfn{DNag.3}
Let $\cc$ be a 2-category. A \emph{monomorphism} $f:Y\la Z$ is a 1-morphism
so that the (strictly) commutative square
\[
\CD
Y @>\id>> Y \\
@V\id VV  @VVfV \\
Y @>f>> Z
\endCD
\]
is 2-cartesian.
\edfn

\lem{LNag.4}
A surjective flat monomorphism is an isomorphism.
\elem

\prf
Let $f:Y\la Z$ be a surjective, flat monomorphism. The surjectivity and
flatness imply that $f$ is faithfully flat, while the fact that $f$
is a monomorphism tells us that the square
  \[
  \CD
  Y @>\id>>  Y\\
  @V\id VV       @VVfV \\
  Y  @>f>>   Z
  \endCD
  \]
is 2-cartesian. 
But then the pullback of $f$ by the faithfully flat map
  $f$ is an isomorphism, and faithfully flat descent
 assures us that $f$ must be an isomorphism.
\eprf

\lem{LNag.5}
Let $X\stackrel f\la Y\stackrel g\la Z$ be composable 1-morphisms in
a 2-category $\cc$, with $g$ a monomorphism.
Then the square
\[
\CD
X @>f>>  Y \\
@V\id VV    @VVgV \\
X @>gf>> Z
\endCD
\]
is 2-cartesian.
\elem

\prf
In the diagram
\[
\CD
X @>f>> Y @>\id>>  Y \\
@V\id VV @V\id VV    @VVgV \\
X @>f>> Y @>g>> Z
\endCD
\]
the square on the right is 2-cartesian because $g$ is a monomorphism,
while the square on the left is trivially 2-cartesian. Hence the
concatenated square is 2-cartesian.
\eprf

\rmk{RNag.7}
Let $\cc$ be a 2-category. Suppose $P$ and $Q$
are properties
of 1-morphisms in $\cc$ so that
\be
\item
  The identity 1-morphisms have property $P$ and property $Q$,
  any 1-morphism
  isomorphic to a morphism having property $P$ (respectively $Q$)
  has property $P$ (respectively $Q$), and composites of morphisms having
  property $P$ (respectively $Q$) have property $P$ (respectively $Q$).
\item
  In any 2-cartesian square
  \[
  \xymatrix{
    W \ar[r]^u \ar[d]_f & X\ar[d]^g\\
    Y \ar[r]^v & Z
  }\]
  if $v$ has property $P$ and $g$ has property $Q$ then $u$ has
  property $P$ and $f$ has property $Q$.
  \setcounter{enumiii}{\value{enumi}}
  \ee
If $X\stackrel f\la Y\stackrel g\la Z$ are composable
morphisms such that
\be
\setcounter{enumi}{\value{enumiii}}
\item
  $g$ is a monomorphism,
\item
 and $gf$ has property $P$ and $g$ has property $Q$.
 \ee
 then $f$ has property $P$, by (ii) applied to the 2-cartesian
 square of Lemma~\ref{LNag.5}.
  \ermk

 \exm{ENag.9}
 We apply Remark~\ref{RNag.7} to the 2-category of algebraic stacks, with
 $P$ the property of being a flat monomorphism and $Q$ the class of all
 morphisms.
 If $X\stackrel f\la Y\stackrel g\la Z$ are composable morphisms of algebraic
 stacks, with $gf$ a flat monomorphism and $g$ a monomorphism, then $f$
 is a flat monomorphism.
 \eexm

 \exm{ENag.10}
 Next we apply Remark~\ref{RNag.7} but with $P$ being the class of all
 dominant, flat monomorphisms and $Q$ the class of all
 flat morphisms. Note that the pullback of a dominant map by a flat map
 is dominant, so the hypotheses of Remark~\ref{RNag.7} hold.

 From Remark~\ref{RNag.7} we learn that if
 $X\stackrel f\la Y\stackrel g\la Z$ is
 a pair of composable 1-morphisms with $gf$ a dominant, flat monomorphism
 and $g$
 a flat monomorphism, then $f$ must be a dominant, flat monomorphism. And
 of course $g$ must also be dominant: the stack-theoretic closure
 of the image of $g$ contains the stack-theoretic closure of the
 image of $gf$, but since $gf$ is dominant this is all of $Z$.
 \eexm

 We end the section with
 
\cor{CNag.22}
Suppose we are given a 2-commutative square of noetherian algebraic stacks 
\[\xymatrix@C+15pt@R+2pt{
W\ar[r]^{u''}\ar[d]_f \ar@{}[dr]|-{(\dagger)}& X\ar[d]^g \\
Y\ar[r]^v & Z
}\]
with $g$ 
separated and of finite type, with $f$ concentrated, universally
quasi-proper
and of finite type, with $v$ a flat monomorphism and with $u''$
a dominant, flat monomorphism.
Then the square $(\dagger)$ is 2-cartesian.
\ecor

\prf
Form the pullback square
\[\xymatrix@C+15pt@R+2pt{
P\ar[r]^u\ar[d]_{p}\ar@{}[dr]|-{(\bullet)} & X\ar[d]^g \\
Y\ar[r]^v & Z
}\]
and let $u':W\la P$ be the map given by the morphism from
the 2-commutative $(\dagger)$ to the 2-cartesian $(\bullet)$.
We obtain a 2-commutative diagram
\[\xymatrix@C+10pt{
W \ar[r]^{u'}\ar[dr]_f & P\ar[r]^u\ar[d]^{p} & X\ar[d]^g \\
&Y\ar[r]^v & Z
}\]
We are given that $v$ is a flat
monomorphism, and its pullback $u$ must also be. But
$u''\cong uu'$ is a dominant,
flat monomorphism by hypothesis,
and Example~\ref{ENag.10},
applied to the composable maps $X\stackrel{u'}\la X'\stackrel u\la X''$,
now
informs us that both $u$ and $u'$ are dominant, flat monomorphisms.

On the other hand $p$, being the pullback of $g$, is separated
and of finite type, and
$f\cong pu'$ is concentrated, of finite type
and universally quasi-proper. Lemma~\ref{LNag.1} informs us
that $u'$ must be concentrated,
of finite type and universally quasi-proper. From Lemma~\ref{LNag.-1}
we deduce that $u'$ is closed: in particular $u'(W)\subset P$ is a
closed subset. Since it is dominant it must be all of $P$. This makes
$u'$ a surjective, flat monomorphism, and Lemma~\ref{LNag.4} guarantees
that $u'$ is an isomorphism.
\eprf

\section{The categories $\nseq(X,Z)$ and $\nseq(X,Y,Z)$}
\label{S6}

Let $\seq$ be a 2-subcategory of the 2-category of noetherian algebraic
stacks, satisfying the hypotheses of Notation~\ref{N0.1}. Part of the
structure of 2-categories gives us, for every pair of objects
$X,Z\in\seq$, a category of 1-morphisms $X\la Z$. As is customary
we denote this category $\seq(X,Z)$. For our
purposes it will be better to work with slightly different categories.

\dfn{D6.2051}
Let $\seq$ be a 2-category as in  Notation~\ref{N0.1}, and let $X,Z$
be any two objects of $\seq$. We define a category $\fnseq(X,Z)$ as
follows
\be
\item
  The objects of $\fnseq(X,Z)$ are 2-commutative triangles in $\seq$
  \[
  \xymatrix@C+30pt@R-20pt{
    & Y\ar[dr]^-{p} & \\
    X\ar[rr]_-{f}\ar[ur]^-{u}  & & Z 
    }\]
Note that the isomorphism $pu\la f$ giving the 2-commutativity
  is part of the data of an object in $\fnseq(X,Z)$.
\item
  A morphism in  $\fnseq(X,Z)$, from object
   \[
  \xymatrix@C-5pt@R-30pt{
    && Y\ar[ddrr]^-{p} &&&& & &Y'\ar[ddrr]^-{p'}& & \\
    &&&&&\text{to object}&&&&&\\
    X\ar[rrrr]_-{f}\ar[uurr]^-{u}& &  & & Z&  & X\ar[rrrr]_-{f'}\ar[uurr]^-{u'} && & & Z
  }\]
  is an equivalence class of data in $\seq$. First let us give the data,
  and then define the equivalence relation. The following 
information suffices to determine a morphism in $\fnseq$:
  \begin{enumerate}
  \item
    A 2-morphism $f\la f'$ in $\seq$;
  \item
    A 2-commutative diagram
\[\xymatrix@C+30pt@R-20pt{
       & Y\ar^p[dr]\ar_\alpha[dd] & \\
X\ar^u[ur]\ar_{u'}[dr] & &  Z \\
       & Y'\ar_{p'}[ur] & 
}\]
in $\seq$, where again the 2-isomorphisms $\alpha u\la u'$ and
$p\la p'\alpha$ are part of the data;
\item
These should be compatible in that the square
\[\CD
pu  @>>> p'u'\\
@VVV  @VVV \\
f @>>> f'
\endCD
\]
must commute in the category $\seq(X,Z)$.
  \end{enumerate}
  But we want to identify these data up to equivalence: we allow ourselves
  to replace $\alpha$ by an isomorph $\alpha'$. 
  More fully, what this means is that given a 2-commutative diagram
  \[\xymatrix@C+30pt@R-20pt{
       & Y\ar^p[dr]\ar@/_/[dd]_\alpha\ar@/^/[dd]^{\alpha'} & \\
X\ar^u[ur]\ar_{u'}[dr] & &  Z \\
       & Y'\ar_{p'}[ur] & 
}\] 
  that is 2-isomorphisms $\ph:\alpha u\la u'$, $\psi:\alpha'\la\alpha$
  and $\rho:p\la p'\alpha'$, then we are allowed to replace
   \[\xymatrix@C+30pt@R-20pt{
       & Y\ar^p[dr]\ar[dd]|\alpha & &       & Y\ar^p[dr]\ar[dd]|{\alpha'} & \\
X\ar^u[ur]\ar_{u'}[dr]\ar@{}[r]|-{\ph} & \ar@{}[r]|-{\psi\rho}&  Z \ar@{}[r]|-{\ds\text{by}} & X\ar^u[ur]\ar_{u'}[dr]\ar@{}[r]|-{\ph\psi} &\ar@{}[r]|-{\rho} &  Z \\
       & Y'\ar_{p'}[ur] &  &         & Y'\ar_{p'}[ur] & 
}\] 
   To avoid any possible confusion: we are free to replace $\alpha$ by $\alpha'$, as long
   as the natural isomorphisms are changed from the pair
   \[
\ph:\alpha u\la u',\,\,
\psi\rho:p\la p'\alpha\qquad\text{\rm to the pair}\quad
\ph\psi:\alpha' u\la u',\,\,
\rho:p\la p'\alpha'.
   \]
 \item
  Composition is obvious, once we note that the equivalence relation is
  preserved by composition.
\item
  We record the following full subcategories of $\fnseq(X,Z)$: 
  the objects of $\lnseq(X,Z)$, $\rnseq(X,Z)$ and $\nseq(X,Z)$ 
  are the 2-commutative diagrams
   \[
  \xymatrix@C+30pt@R-20pt{
    & Y\ar[dr]^-{p} & \\
    X\ar[rr]_-{f}\ar[ur]^-{u}  & & Z 
    }\]
  where:
  \begin{enumerate}
  \item
    In  $\lnseq(X,Z)$ we assume that $u$ is a flat morphism.
  \item
    In $\rnseq(X,Z)$ we suppose that $p$ is of finite type and
    universally quasi-proper.
  \item $\nseq(X,Z)$ is the full subcategory of $\lnseq(X,Z)\cap\rnseq(X,Z)$
    where $u$ is assumed not only to be flat, but also to be a dominant
    monomorphism.
  \end{enumerate}
  \ee
\edfn

\rmk{R6.2052}
The variant we care about most is $\nseq(X,Z)$.
To recast the definition in simple words:
an object of $\nseq(X,Z)$ is a morphism
$f:X\la Z$ in $\seq$ together with a Nagata compactification, and the
morphisms in $\nseq(X,Z)$ are the obvious maps of the data. The one quirk is
that we mod out by the possible 2-categorical structure,
identifying morphisms up to equivalence. The lax versions $\lnseq(X,Z)$
and $\rnseq(X,Z)$ will come up in proofs, hence we recorded them.

There is a forgetful functor $F:\fnseq(X,Y)\la\seq(X,Y)$ which
remembers only $f:X\la Z$. This map is not an equivalence;
after all the category $\seq(X,Y)$ is a groupoid while
$\fnseq(X,Y)$ is not, and neither is
any of the three subcategories
we introduced. What we will show in the next few lemmas is that the
variant $\nseq(X,Y)$ admits a formal calculus of right fractions,
allowing us to easily prove that
the map from the groupoid completion of $\nseq(X,Z)$ to
$\seq(X,Y)$ is an equivalence.
\ermk

\rmk{R6.2052.3}
Suppose we choose a representative of a morphism in $\rnseq(X,Z)$;
part of the data this specifies, according to
Definition~\ref{D6.2051}(ii), is the 2-commutative diagram
\[\xymatrix@C+30pt@R-20pt{
       & Y\ar^p[dr]\ar[dd]_\alpha & \\
X\ar^u[ur]\ar_{u'}[dr] & &  Z \\
       & Y'\ar_{p'}[ur] & 
}\]
We assert that the morphism $\alpha:Y\la Y'$ must be of finite
type and universally quasi-proper.

The proof is by Lemma~\ref{LNag.1}: since the morphism
lies in $\rnseq(X,Z)$ both $p'$ and $p\cong\alpha p'$ are of finite
type and universally quasi-proper. All morphisms of $\seq$ are assumed
separated, hence  Lemma~\ref{LNag.1} applies and tells us that $\alpha$
is of finite type and universally quasi-proper.
\ermk

\rmk{R6.2052.5}
If we assume that all our stacks are schemes then life becomes much
simpler: given two objects in $\nseq(X,Z)$, i.e.~two diagrams
   \[
  \xymatrix@C-5pt@R-30pt{
    && Y\ar[ddrr]^-{p} &&&& & &Y'\ar[ddrr]^-{p'}& & \\
    &&&&&\text{and}&&&&&\\
    X\ar[rrrr]_-{f}\ar[uurr]^-{u}& &  & & Z&  & X\ar[rrrr]_-{f'}\ar[uurr]^-{u'} && & & Z
  }\]
  there is at most one commutative diagram
\[\xymatrix@C+30pt@R-20pt{
       & Y\ar^p[dr]\ar[dd]_\alpha & \\
X\ar^u[ur]\ar_{u'}[dr] & &  Z \\
       & Y'\ar_{p'}[ur] & 
}\] 
After all the map $u$ is assumed dominant and $p'$ is a separated morphism,
hence the equality $u'=\alpha u$ uniquely specifies $\alpha$ (if it exists).
If all the objects of $\seq$ are schemes then $\nseq(X,Z)$ is a
partially ordered set.
\ermk

\lem{L6.2053}
Let $F:\rnseq(X,Y)\la\seq(X,Y)$ be as in Remark~\ref{R6.2052}. Given
two objects $b,c\in\rnseq(X,Y)$ and a morphism $\ph:F(b)\la F(c)$ in $\seq$,
there
is an object
$a\in\nseq(X,Z)$  and a
pair of morphisms $b\stackrel\alpha\longleftarrow a\stackrel\beta\la c$
in $\rnseq(X,Y)$ so that $\ph=F(\beta)F(\alpha)^{-1}$.
\elem

\prf
The objects $b$ and $c$ are a pair of factorizations
$pu\la f$ and $p'u'\la f'$, and the morphism $\ph:F(b)\la F(c)$ comes
down to a map $f\la f'$ in $\seq$. These assemble in $\seq$
to a 2-commutative square
\[
\CD
X @>u>> Y \\
@Vu'VV @VVpV \\
Y' @>p'>> Z
\endCD
\]
Now form the 2-cartesian square
\[
\CD
Y\times_Z^{}Y' @>q'>> Y \\
@VqVV @VVpV \\
Y' @>p'>> Z
\endCD
\]
and consider the induced map $X\la Y\times_Z^{} Y'$. It is a morphism
in $\seq$, hence we may choose a Nagata compactification $X\stackrel{u''}\la
\ov Y
\stackrel r\la Y\times_Z^{}Y'$.

Because $p$ is of finite type and universally quasi-proper so is its
pullback $q$, and because $p'$ and $r$ are also of finite type and universally
quasi-proper so is the composite $p'qr:\ov Y\la Z$. The map $u'':X\la \ov Y$
is a dominant, flat monomorphism by construction. Hence
$X\stackrel{u''}\la
\ov Y\stackrel{p'qr}\la Z$
is a Nagata compactification of the composite; we let this be our object
$a\in\nseq(X,Y)$. The maps $b\stackrel\alpha\longleftarrow a\stackrel\beta\la c$
are now obvious.
\eprf

\lem{L6.2054}
Let $\ph,\psi$ be two morphisms $\xymatrix{b\ar@<0.5ex>[r]\ar@<-0.5ex>[r]& c}$
  in $\rnseq(X,Z)$ with $F(\ph)=F(\psi)$.
Then there exists an object $a\in\nseq(X,Z)$ and  
a morphism $\rho:a\la b$ with $\ph\rho=\psi\rho$.
\elem

\prf
If we choose representatives
in the equivalence classes giving $\ph$ and $\psi$ we obtain two 2-commutitive
diagrams in $\seq$
\[\xymatrix@C+30pt@R-20pt{
       & Y\ar^p[dr]\ar[dd]^\alpha & &        & Y\ar^p[dr]\ar[dd]^{\alpha'} & \\
X\ar^u[ur]\ar_{u'}[dr] & &  Z  & X\ar^u[ur]\ar_{u'}[dr] & &  Z \\
       & Y'\ar_{p'}[ur] & &    & Y'\ar_{p'}[ur] & 
}\]
The assertion that $F(\ph)=F(\psi)$ says that in $\seq$ the two 2-morphisms
$pu\la p'u'$ agree. Let $\Delta$ be the diagonal map
$Y'\la Y'\times_Z^{}Y'$, and let $s:Y\la Y'\times_Z^{}Y'$ be the map
corresponding to the pair of 1-morphisms 
$\xymatrix{Y\ar@<0.5ex>[r]^\alpha \ar@<-0.5ex>[r]_{\alpha'} & Y'}$,
coupled with the two
2-morphisms $p'\alpha\la p\longleftarrow p'\alpha'$.
Now consider the 2-commutative square
\[
\xymatrix@C+10pt{
X \ar[r]^-u \ar[d]_-{u'} & Y\ar[d]^s \\
Y' \ar[r]^-\Delta & Y'\times_Z^{}Y'
}
\]
Form the pullback
\[
\xymatrix@C+10pt{
P \ar[r]^-q \ar[d]_-{q'} & Y\ar[d]^s \\
Y' \ar[r]^-\Delta & Y'\times_Z^{}Y'
}
\]
and choose a 1-morphism $X\la P$ giving a map
from the 2-commutative to the 2-cartesian square. Choose a Nagata compactification
$X\stackrel w \la Y''\stackrel r\la P$. The map $p':Y'\la Z$ is separated
and of finite type, hence
\be
\item
  The pullback $Y'\times_Z^{}Y'\la Y'$ is of finite type, and since
  $Y'$ is noetherian so is $Y'\times_Z^{}Y'$.
\item
  The diagonal map
  $\Delta:Y'\la Y'\times_Z^{}Y'$ is a proper map
  of noetherian stacks, hence of finite type and universally
  quasi-proper.
  \ee
  Therefore the pullback $q:P\la Y$ is also of finite type and
  universally quasi-proper, as are $r:Y''\la P$ and $p:Y\la Z$. Therefore
  the composite $pqr:Y''\la Z$ is of finite type and universally quasi-proper.
  This makes the 2-commutative diagram
\[\xymatrix@C+30pt@R-20pt{
       & Y''\ar^-{pqr}[dr]\ar[dd]_{qr} & \\
X\ar^w[ur]\ar_{u}[dr] & &  Z \\
       & Y\ar_{p}[ur] & 
}\]
a morphism $\rho:a\la b$ in $\rnseq(X,Z)$, and
$a$ belongs to the subcategory $\nseq(X,Z)$.
We
leave it to the reader to check that $\ph\rho=\psi\rho$.

In this check the reader should note that it isn't true that
$\alpha qr=
\alpha'qr$.
If $\pi_1^{}$, $\pi_2^{}$ are the two projections $Y'\times_Z^{}Y'\la Y'$,
all we 
have  are isomorphisms
$\alpha qr \cong \pi_1^{}sqr\cong\pi_1^{}\Delta q'r\cong q'r
\cong\pi_2^{}\Delta q'r\cong \pi_2^{}sqr\cong\alpha' qr$, and that these
2-isomorphisms are compatible with all the given 
2-morphisms of the data. But 2-isomorphisms are not equalities: before
passing to the equivalence relation it isn't true that $\ph\rho=\psi\rho$.
\eprf

We proved Lemmas~\ref{L6.2053} and \ref{L6.2054} in the generality of
$\rnseq(X,Z)$, but it says something
interesting even if we concentrate on the case 
where all the objects lie in $\nseq(X,Z)$. For future use we 
state this abstractly.

\rmk{R6.2054.5}
We suppose given a 
functor $F:\ca\la\cb$ satisfying
\be
\item
  The category $\cb$ is a groupoid, meaning all the morphisms in $\cb$
  are invertible.
\item
  Given a morphism $\ph:F(b)\la F(c)$ in $\cb$, there is a pair of
  morphisms $b\stackrel\alpha\longleftarrow a\stackrel\beta\la c$
  in the category $\ca$ with $\ph=F(\beta)F(\alpha)^{-1}$.
\item
  Let $\ph,\psi$ be two morphisms in $\ca$ with $F(\ph)=F(\psi)$.
Then there exists
a morphism $\rho\in\ca$ with $\ph\rho=\psi\rho$
\ee
Note that the functor $F:\nseq(X,Z)\la\seq(X,Z)$ satisfies (i), (ii) and (iii);
(i) is obvious since the 2-category of stacks is a groupoid,
(ii) by Lemma~\ref{L6.2053},
and (iii) by Lemma~\ref{L6.2054}. The following is a formal
consequence:
\ermk

\lem{L6.2055}
Let
$F:\ca\la\cb$ be a functor satisfying conditions
(i), (ii) and (iii) of Remark~\ref{R6.2054.5}.
Let $\mathbf{G}\ca$ be the groupoid completion of $\ca$, that
is the category where we formally invert all the morphisms.
Then the category $\ca$ satisfies the right Ore condition, and
in
the canonical factorization of the functor $F$  as
$\ca\stackrel\pi\la\mathbf{G}\ca\stackrel{F'}\la\cb$,
the functor $F'$ is fully faithful.
\elem

\prf
We should remind ourselves what it means for
the category $\ca$ to satisfy the right Ore condition. The following two
conditions must be satisfied:
\be
\item
  Any pair of
  morphisms $c\la d\longleftarrow b$ in $\ca$
  can be completed to a commutative
square
\[
\xymatrix{
  a\ar[r]\ar[d] & b\ar[d] \\
  c\ar[r] & d
  }
\]
\item
  Any pair of morphisms $\xymatrix{b\ar@<0.5ex>[r]\ar@<-0.5ex>[r] & c}$
   in $\ca$ admitting a coequalizer
   also admits an equalizer. That is, if there exists a morphism
   $c\la d$ so that the composites
   $\xymatrix{b\ar@<0.5ex>[r]\ar@<-0.5ex>[r] & c\ar[r] & d}$
   are equal, then there also
   exists a morphism $a\la b$ so that
the composites $\xymatrix{a\ar[r] &b\ar@<0.5ex>[r]\ar@<-0.5ex>[r] & c}$
are equal.
\ee
Let us begin by proving that these two conditions are satisfied. 

We start with (ii): 
suppose we are given in $\ca(X,Y)$ morphisms 
$\xymatrix{b\ar@<0.5ex>[r]\ar@<-0.5ex>[r] & c\ar[r] & d}$
so that the two composites $b\la d$ are equal. 
Applying the functor $F$ we have that the two composites
$\xymatrix{F(b)\ar@<0.5ex>[r]\ar@<-0.5ex>[r] & F(c)\ar[r] & F(d)}$
are also equal. But the category $\cb$ is a groupoid
by Remark~\ref{R6.2054.5}(i),
hence the morphism $F(c)\la F(d)$ is an isomorphism,
and the two morphisms
$\xymatrix{F(b)\ar@<0.5ex>[r]\ar@<-0.5ex>[r] & F(c)}$ must be equal.
Remark~\ref{R6.2054.5}(iii) now informs us that there exists in $\ca$
a morphism
$a\la b$ so that the composites
$\xymatrix{a\ar[r] &b\ar@<0.5ex>[r]\ar@<-0.5ex>[r] & c}$ are equal.

Next we prove (i): suppose that we are given in $\ca$ morphisms
$c\la d\longleftarrow b$. Applying the functor $F$ we obtain morphisms
$F(c)\la F(d)\longleftarrow F(b)$, but $\cb$ is a groupoid so
the map $F(d)\longleftarrow F(b)$ is invertible. The composite
$F(c)\la F(d)\la F(b)$
is a morphism $\ph:F(c)\la F(b)$, and by Remark~\ref{R6.2054.5}(ii)
it may be written as
$\ph=F(\beta)F(\alpha)^{-1}$
for some morphisms $c\stackrel\alpha\longleftarrow \wt a\stackrel\beta\la b$
in the category $\ca$. We have produced in $\ca$ morphisms
\[
\xymatrix{
  \wt a\ar[r]^\beta\ar[d]_\alpha & b\ar[d] \\
  c\ar[r] & d
  }
\]
and what we know is that the functor $F$ takes this to a commutative square.
Consider the two composites
$\xymatrix{\wt a\ar@<0.5ex>[r]\ar@<-0.5ex>[r] & d}$; we don't know that
they agree in $\ca$, but the functor $F$ takes them to
equal maps. Applying Remark~\ref{R6.2054.5}(iii) there exists a morphism
$a\la \wt a$ so that the two composites
$\xymatrix{a\ar[r] &\wt a\ar@<0.5ex>[r]\ar@<-0.5ex>[r] & d}$
are equal, and hence we obtain in $\ca$ a commutative square
\[
\xymatrix{
  a\ar[r]\ar[d] & b\ar[d] \\
  c\ar[r] & d
  }
\]
This completes the proof of (i).

Let $\mathbf{G}\ca$ be
the groupoid completion of $\ca$. Because $\ca$
satisfies the right Ore conditions the calculus of right fractions is very
simple: any morphism $\ph:b\la c$ in $\mathbf{G}\ca$ can be represented as
$\beta\alpha^{-1}$ for morphisms
$b\stackrel\alpha\longleftarrow a\stackrel\beta\la c$
in the category $\ca$, any pair of morphisms 
$\xymatrix{b\ar@<0.5ex>[r]^\ph\ar@<-0.5ex>[r]_{\ph'} & c}$ in $\mathbf{G}\ca$
may be put on a ``common denominator'', meaning there are morphisms
$\xymatrix{b& \ar[l]_\alpha a\ar@<0.5ex>[r]^\beta\ar@<-0.5ex>[r]_{\beta'} & c}$
in $\ca$ with $\ph=\beta\alpha^{-1}$ and $\ph'=\beta'\alpha^{-1}$.
Furthermore if we represent two morphisms
$\xymatrix{b\ar@<0.5ex>[r]\ar@<-0.5ex>[r]& c}$ in $\mathbf{G}\ca$ as
$\ph=\beta\alpha^{-1}$ and $\ph'=\beta'\alpha^{-1}$, then $\ph=\ph'$ if and
only if there exists in $\ca$
a morphism $\gamma$ with
$\beta\gamma=\beta'\gamma$.

The functor $F$ is a functor from $\ca$ to the groupoid $\cb$,
hence factors uniquely as $\ca\stackrel\pi\la \mathbf{G}\ca\stackrel{F'}\la
\cb$.
Remark~\ref{R6.2054.5}(ii) shows that the functor
$F'$ is full, while
Remark~\ref{R6.2054.5}(iii) proves it faithful.
\eprf

\pro{P6.2055.5}
In the special case where the functor $F:\ca\la\cb$ of Lemma~\ref{L6.2055}
is the functor $F:\nseq(X,Z)\la\seq(X,Z)$ of Remark~\ref{R6.2052}, the
factorization $\nseq(X,Z)\stackrel\pi\la\gnseq(X,Z)\stackrel{F'}\la\seq(X,Z)$
is such that $F'$ is an equivalence.
\epro

\prf
By Lemma~\ref{L6.2055} we know that $F'$ is fully faithful. But in the
category $\seq$ we assume that any morphism $f:X\la Z$ has a Nagata
compactification, and the functors $F$ and  $F'$ are therefore essentially
surjective.
\eprf

For every pair of objects $X,Z\in\seq$ we have found a category $\gnseq(X,Z)$
equivalent to $\seq(X,Z)$. We now want to do a similar construction
for triples of objects.

\dfn{D6.2056}
Let $\seq$ be a 2-category as in Notation~\ref{N0.1}, and let $X,Y,Z$
be three objects in $\seq$. We define the category $\nseq(X,Y,Z)$ as
follows:
\be
\item
  The objects in $\nseq(X,Y,Z)$ are
  2-commutative diagrams in $\seq$
  \[
  \xymatrix{
    X\ar[dr]_-f \ar[r]^-{u'}  &X' \ar[r]^{u} \ar[d]^{p'} &  X''\ar[d]^p\\
 &   Y\ar[dr]_-g \ar[r]^-v & Y'\ar[d]^q\\
 & &  Z
}
\]
where the horizontal maps are all dominant, flat monomorphisms and
the vertical maps are of finite type and universally quasi-proper.
Corollary~\ref{CNag.22} tells us that the one square in the
diagram is automatically 2-cartesian.
\item
  The morphisms in $\nseq(X,Y,Z)$, from the object
 \[
  \xymatrix{
    X\ar[dr]_-f \ar[r]^-{u'}  &X' \ar[r]^{u} \ar[d]^{p'} &  X''\ar[d]^p & &
    X\ar[dr]_-{\ov f} \ar[r]^-{u'}  &\ov X' \ar[r]^{\ov u'} \ar[d]^{\ov p'} &  \ov X''\ar[d]^{\ov p}\\
 &   Y\ar[dr]_-g \ar[r]^-v & Y'\ar[d]^q &\text{to} &&   Y\ar[dr]_-{\ov g} \ar[r]^-{\ov v} &\ov Y'\ar[d]^{\ov q} \\
 & &  Z & & & & Z
}
  \]
    are equivalence classes of 2-commutative diagrams in $\seq$
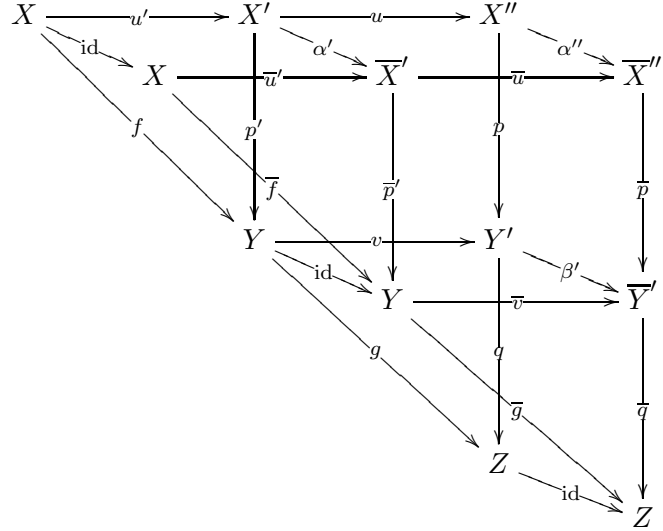
\[
\xymatrix@R-18pt@C-18pt{
X\ar[dddddrrrrr]|-f \ar[drrr]|{\id}\ar[rrrrr]|{u'} &&&&&
 X' \ar[drrr]|{\alpha'}\ar[ddddd]|{p'}  
          \ar[rrrrr]|{u} &&&&& X''
     \ar[drrr]|{\alpha''}\ar[ddddd]|{p}&&&\\
&&& X\ar[dddddrrrrr]|-{\ov f} \ar[rrrrr]|{\ov u'} &&&&& \ov X'\ar[ddddd]|{\ov p'}  
          \ar[rrrrr]|{\ov u} &&&&& \ov X''\ar[ddddd]|{\ov p}\\
&&&&&&&&&&&&&\\
&&&&&&&&&&&&&\\
&&&&&&&&&&&&&\\
&&&&&
Y\ar[dddddrrrrr]|-{g} \ar[drrr]|{\id} \ar[rrrrr]|{v} &&&&& Y' \ar[drrr]|{\beta'}\ar[ddddd]|{q}
&&&\\
& &&&&&&&  Y\ar[dddddrrrrr]|-{\ov g} \ar[rrrrr]|{\ov v} &&&&& \ov Y'\ar[ddddd]|{\ov q}\\
&&&&&&&&&&&&&\\
&&&&&&&&&&&&&\\
&&&&&&&&&&&&&\\
&&&&&&&&&&Z \ar[drrr]|{\id}&&&\\
&&&&&&&&&&&&& Z
}
\]
\item
  Two such diagrams
  connecting a given pair of objects give, among other things,
  three pairs of 1-morphisms
$\xymatrix{X'\ar@<0.5ex>[r]^{\alpha'}\ar@<-0.5ex>[r]_{\ov \alpha'} & \ov X'}$,
$\xymatrix{X''\ar@<0.5ex>[r]^{\alpha''}\ar@<-0.5ex>[r]_{\ov \alpha''} & \ov X''}$
  and $\xymatrix{Y'\ar@<0.5ex>[r]^{\beta'}\ar@<-0.5ex>[r]_{\ov \beta'} & \ov Y'}$.
  These diagrams are declared to be equivalent if there are 2-isomorphisms
  $\alpha'\cong\ov\alpha'$, $\alpha''\cong\ov\alpha''$ and
  $\beta'\cong\ov\beta'$ compatible with all the other 2-morphisms in
  the diagrams.
\item
  For future use we record the minor variant $\lnseq(X,Y,Z)$:
  the objects are 2-commutative
  diagrams as 
  in $\nseq(X,Y,Z)$, the square is 2-cartesian, but the vertical maps
  are unrestricted and the only restriction on the horizontal maps is
  that they are assumed flat.
\ee
\edfn

\rmk{R6.2057}
There is a functor $F:\nseq(X,Y,Z)\la\seq(X,Y)\times\seq(Y,Z)$; it takes
an object of $\nseq(X,Y,Z)$ to the pair of composable 1-morphisms
$X\stackrel f\la Y\stackrel g\la Z$. An object of $\nseq(X,Y,Z)$
should be thought of as the composable 1-morphisms together with
a bunch of compatible Nagata compactifications, and the morphisms
are equivalence classes of maps respecting the data.

The category $\nseq(X,Y,Z)$ also has three obvious functors
\[
\pi_{12}^{}:\nseq(X,Y,Z)\la\nseq(X,Y),\quad
\pi_{23}^{}:\nseq(X,Y,Z)\la\nseq(Y,Z)
\]
\[\text{and}\quad
\pi_{13}^{}:\nseq(X,Y,Z)\la\nseq(X,Z)
\]
each of which forgets some of the data. For example $\pi_{23}^{}$
takes the object
  \[
  \xymatrix{
    X\ar[dr]_-f \ar[r]^-{u'}  &X' \ar[r]^{u} \ar[d]^{p'} &  X''\ar[d]^p\\
 &   Y\ar[dr]_-g \ar[r]^-v & Y'\ar[d]^q\\
 & &  Z
}
  \]
  of the category $\nseq(X,Y,Z)$ to the object
 \[
  \xymatrix{
Y\ar[dr]_-g \ar[r]^-v & Y'\ar[d]^q\\
 &  Z
}
  \]
  of the category $\nseq(Y,Z)$.
\ermk
  
\rmk{R6.2057.5}
In Definition~\ref{D6.2056}(ii) we saw that a representative of a morphism
$\ph:a\la b$ in the category $\nseq(X,Y,Z)$, where the objects $a$ and
$b$ are the diagrams
\[
  \xymatrix{
    X\ar[dr]_-f \ar[r]^-{u'}  &X' \ar[r]^{u} \ar[d]^{p'} &  X''\ar[d]^p & &
    X\ar[dr]_-{\ov f} \ar[r]^-{u'}  &\ov X' \ar[r]^{\ov u'} \ar[d]^{\ov p'} &  \ov X''\ar[d]^{\ov p}\\
 &   Y\ar[dr]_-g \ar[r]^-v & Y'\ar[d]^q &\text{and} &&   Y\ar[dr]_-{\ov g} \ar[r]^-{\ov v} &\ov Y'\ar[d]^{\ov q} \\
 & &  Z & & & & Z
}
  \]
is a 2-commutative diagram where the ``new'' 1-morphisms---the ones
which do not
form part of the data of the objects---are the triple
$\alpha':X'\la\ov X'$, $\alpha'':X''\la\ov X''$ and $\beta':Y'\la \ov Y'$.
We assert that the 1-morphisms $\alpha'$, $\alpha''$ and $\beta'$ are
all of finite type and universally quasi-proper.

The proof is by Remark~\ref{R6.2052.3}: if we apply the functors of
Remark~\ref{R6.2057} we obtain three maps
$\pi_{}^{}(\ph):\pi_{}^{}(a)
\la \pi_{12}^{}(b)$,
$\pi_{12}^{}(\ph):\pi_{12}^{}(a)
\la \pi_{23}^{}(b)$
and
$\pi_{13}^{}(\ph):\pi_{13}^{}(a)
\la \pi_{13}^{}(b)$,
and each of them is represented by a diagram where the ``new'' 1-morphism
is (respectively) $\alpha'$, $\alpha''$ and $\beta'$.
Remark~\ref{R6.2052.3} guarantees that these 1-morphisms are of finite
type and universally quasi-proper.
\ermk

\lem{L6.2058}
The functor $F:\nseq(X,Y,Z)\la\seq(X,Y)\times\seq(Y,Z)$
is essentially surjective. 
\elem

\prf
Suppose we are given an object in $\seq(X,Y)\times\seq(Y,Z)$, that is a
pair of composable morphisms $X\stackrel f\la Y\stackrel g\la Z$. Choose
Nagata compactification for $f$ and $g$. So far this constructs for us
an object $a\in\nseq(X,Y)$ and an object $b\in\nseq(Y,Z)$,
which we assemble together in the 2-commutative diagram
  \[
  \xymatrix{
X \ar[dr]_-f \ar[r]^{u'} & X'\ar[d]^{p'} \\
& Y\ar[dr]_-g \ar[r]^-v & Y'\ar[d]^q\\
&  &  Z
}
  \]
Now the composite $X'\stackrel {p'}\la Y\stackrel v\la Y'$ is a morphism
in the category $\seq$, and has a Nagata compactification
$X'\stackrel {u}\la X''\stackrel p\la Y'$. This means that the composites
are 2-isomorphic, $u$ is a dominant, flat monomorphism and $p$ is of
finite type and universally quasi-proper. This completes for us the diagram
 \[
  \xymatrix{
X \ar[dr]_-f \ar[r]^{u'} & X'\ar[d]^{p'}\ar[r]^u & X''\ar[d]^p \\
& Y\ar[dr]_-g \ar[r]^-v & Y'\ar[d]^q\\
&  &  Z
}
  \]
which is an object of $\nseq(X,Y,Z)$ lifting  $X\stackrel {f}\la Y\stackrel g\la Z$.
\eprf

\lem{L6.2059}
Let $F:\nseq(X,Y,Z)\la\seq(X,Y)\times\seq(Y,Z)$ be the functor of
Remark~\ref{R6.2057}, and suppose we are
given a morphism $\ph:F(b)\la F(c)$ in $\seq(X,Y)\times\seq(Y,Z)$. There
exists a pair of morphisms $b\stackrel\sigma\longleftarrow a\stackrel\theta\la c$
in the category $\nseq(X,Y,Z)$, with $\ph=F(\theta)F(\sigma)^{-1}$.
\elem

\prf
We are given two objects $b,c\in\nseq(X,Y,Z)$, that is diagrams
in $\seq$ satisfying all the requirements of Definition~\ref{D6.2056}(i)
 \[
  \xymatrix{
    X\ar[dr]_-f \ar[r]^-{u'}  &X' \ar[r]^{u} \ar[d]^{p'} &  X''\ar[d]^p & &
    X\ar[dr]_-{\ov f} \ar[r]^-{\ov u'}  &\ov X' \ar[r]^{\ov u} \ar[d]^{\ov p'} &  \ov X''\ar[d]^{\ov p}\\
 &   Y\ar[dr]_-g \ar[r]^-v & Y'\ar[d]^q &\text{and} &&   Y\ar[dr]_-{\ov g} \ar[r]^-{\ov v} &\ov Y'\ar[d]^{\ov q} \\
 & &  Z & & & & Z
}
  \]
We are also given a morphism $\ph:F(b)\la F(c)$, meaning a
pair of 2-morphisms $\ph_f^{}:f\la \ov f$ and $\ph_g^{}:g\la \ov g$
in $\seq$. Apply the functors $\pi_{12}^{}$ and $\pi_{23}^{}$ of
Remark~\ref{R6.2057}, and then Lemma~\ref{L6.2053},
and we construct the diagram $A$ below 
\[
\xymatrix{
  X\ar[dr]_-{\wh f} \ar[r]^{\wh u'} & \wh X'\ar[d]^{\wh p'} \\
 &       Y\ar[dr]_-{\wh g} \ar[r]^-{\wh v} &\wh Y'\ar[d]^{\wh q} \\
 & &  Z
}
\]
together with maps
\[
\begin{array}{cl}
\sigma':\pi_{12}^{}(A)\la\pi_{12}^{}(b),\quad &\theta':\pi_{12}^{}(A)\la\pi_{12}^{}(c),\\
\sigma'':\pi_{23}^{}(A)\la\pi_{23}^{}(b),\quad&\theta'':\pi_{23}^{}(A)\la\pi_{23}^{}(c)
\end{array}
\]
so that
\be
\item
The horizontal maps in $A$ are dominant, flat monomorphisms and the vertical
maps are of finite type and universally quasi-proper.
\item
  $\ph_f^{}=F(\theta')F(\sigma')^{-1}$ and
  $\ph_{g}^{}=F(\theta'')F(\sigma'')^{-1}$.
\ee
What we are about to prove is a refinement of Lemma~\ref{L6.2059}. Since we
will need the refinement in the proof of Lemma~\ref{L6.2060} we state it
formally here.

\sthm{L6.2059.5}
Suppose we have chosen an object $A$ and maps
$\sigma'$, $\sigma''$, $\theta'$ and $\theta''$ satisfying (i) and (ii) above.
The object $a$ in the statement of Lemma~\ref{L6.2059} can be chosen so
that $\pi_{12}^{}(a)=\pi_{12}^{}(A)$ and $\pi_{23}^{}(a)=\pi_{23}^{}(A)$, and
the morphisms $b\stackrel\sigma\longleftarrow a\stackrel\theta\la c$
in the statement of Lemma~\ref{L6.2059} can be chosen so
that
\[ \pi_{12}^{}(\sigma)=\sigma',\quad
\pi_{23}^{}(\sigma)=\sigma'',\quad
 \pi_{12}^{}(\theta)=\theta',\quad
\pi_{23}^{}(\theta)=\theta''.
\]
\esthm

It remains to prove \ref{L6.2059.5}
and hence also Lemma~\ref{L6.2059}.
Choose representatives of the maps $\sigma'$, $\sigma''$, $\theta'$ and
$\theta''$. Concretely this means we choose maps $\alpha'_\sigma:\wh X'\la X'$,
$\alpha'_\theta:\wh X'\la\ov X'$, $\beta'_\sigma:\wh Y'\la Y'$
and $\beta'_\theta:\wh Y'\la\ov Y'$ rendering some triangles 2-commutative.
By Remark~\ref{R6.2052.3} the 1-morphisms $\alpha'_\sigma$, $\alpha'_\theta$,
$\beta'_\sigma$
and $\beta'_\theta$ are of finite type and universally quasi-proper. 

Choose a Nagata compactification for the composite
$\wh X'\stackrel{\wh p'}\la Y\stackrel {\wh v}\la \wh Y'$, giving the
2-commutative triangle 
   \[
  \xymatrix@C-2pt@R-27pt{
    && P\ar[ddrr]^-{r} &&\\
    &&&&\\
    \wh X'\ar[rrrr]_-{\wh v\wh p'}\ar[uurr]^-{w}& &  & & \wh Y'
  }\]
  which we view as an object $D\in\nseq(\wh X',\wh Y')$.
We know that $\beta'_\sigma$ and $\beta'_\theta$ are of finite type and
universally quasi-proper, hence
\[
  \xymatrix@C-5pt@R-30pt{
    && P\ar[ddrr]^-{\beta'_\sigma r} &&&& & &P\ar[ddrr]^-{\beta'_\theta r}& & \\
    &&&&&\text{and}&&&&&\\
    \wh X'\ar[rrrr]_-{\beta'_\sigma\wh v\wh p'}\ar[uurr]^-{w}& &  & & Y'&  & \wh X'\ar[rrrr]_-{\beta'_\theta \wh v\wh p'}\ar[uurr]^-{w} && & & \ov Y'
  }\]
can be viewed as objects $E\in\nseq(\wh X',Y')$ and $G\in\nseq(\wh X',\ov Y')$.
But we also have the diagrams
\[
  \xymatrix@C-5pt@R-30pt{
    && X''\ar[ddrr]^-{p} &&&& & &\ov X''\ar[ddrr]^-{\ov p}& & \\
    &&&&&\text{and}&&&&&\\
    \wh X'\ar[rrrr]_-{vp'\alpha'_\sigma}\ar[uurr]^-{u\alpha'_\sigma}& &  & & Y'&  & \wh X'\ar[rrrr]_-{\ov v\ov p'\alpha'_\theta }\ar[uurr]^-{\ov u\alpha'_\theta} && & & \ov Y'
  }\]
which define objects $E'\in\rnseq(\wh X',Y')$ and
$G'\in\rnseq(\wh X',\ov Y')$. Furthemore we are given 2-morphisms
$\lambda:\beta'_\sigma \wh v\wh p'\la vp'\alpha'_\sigma$ and
$\mu:\beta'_\theta \wh v\wh p'\la\ov v\ov p'\alpha'_\theta$,
in other words morphisms $\lambda:F(E)\la F(E')$ and $\mu:F(G)\la F(G')$.
Lemma~\ref{L6.2053} applies, and we produce objects
$H\in\nseq(\wh X',Y')$ and $\ov H\in\nseq(\wh X',\ov Y')$,
together with maps $E\stackrel{\gamma}\longleftarrow H\stackrel\delta\la E'$
and $G\stackrel{\gamma'}\longleftarrow \ov H\stackrel{\delta'}\la G'$
with $\lambda=F(\delta)F(\gamma)^{-1}$ and
$\mu=F(\delta')F(\gamma')^{-1}$. Let $H$ and $\ov H$ be the diagrams
\[
  \xymatrix@C-5pt@R-30pt{
    && Q\ar[ddrr]^-{t} &&&& & &\ov Q\ar[ddrr]^-{\ov t}& & \\
    &&&&&\text{and}&&&&&\\
    \wh X'\ar[rrrr]_-{s}\ar[uurr]^-{w'}& &  & & Y'&  & \wh X'\ar[rrrr]_-{\ov s}\ar[uurr]^-{\ov w'} && & & \ov Y'
  }\]
The morphisms $\gamma:H\la E$ and $\gamma':\ov H\la G$ have representatives
giving the 2-commutative triangles
\[
\xymatrix@C+30pt@R-20pt{
&  Q\ar[dr]^-{t}\ar[dd]_-{\zeta} &  & &
  \ov Q\ar[dr]^-{\ov t}\ar[dd]_-{\xi} &  \\
\wh X'\ar[ru]^-{w'}\ar[dr]_{w}&  & Y'& \wh X'\ar[ur]^-{\ov w'}\ar[dr]_-{w} && \ov Y'\\
 & P\ar[ur]_-{\beta'_\sigma r} &
& &P\ar[ur]_-{\beta'_\theta r} &
}\]
Replacing $H$ and $\ov H$ by isomorphs we may assume that
the 2-isomorphisms $t\la \beta'_\sigma r\zeta$
and
$\ov t\la \beta'_\theta r\xi$ are identities.
Remark~\ref{R6.2052.3} tells us that  $\zeta$ and $\xi$
are of finite type and universally quasi-proper, 
therefore the following are
morphisms $D'\stackrel{i}\la D\stackrel j\longleftarrow \ov D'$
in $\nseq(\wh X',\wh Y')$
\[
\xymatrix@C+30pt@R-20pt{
&  Q\ar[dr]^-{r\zeta}\ar[dd]_-{\zeta} &  & &
  \ov Q\ar[dr]^-{r\xi}\ar[dd]_-{\xi} &  \\
\wh X'\ar[ru]^-{w'}\ar[dr]_{w}&  & \wh Y'& \wh X'\ar[ur]^-{\ov w'}\ar[dr]_-{w} && \wh Y'\\
 & P\ar[ur]_-{ r} &
& &P\ar[ur]_-{r} &
}\]
and the morphisms $\gamma:H\la E$ and $\gamma':\ov H\la G$ satisfy
$\gamma=\beta'_\sigma i$ and $\gamma'=\beta'_\theta j$.

But Lemmas~\ref{L6.2053}, \ref{L6.2054} and \ref{L6.2055}
combine to tell us that the morphisms of $\nseq(\wh X',\wh Y')$ satisfy the right Ore
condition, hence the morphisms $D'\stackrel{i}\la D\stackrel j\longleftarrow \ov D'$
can be completed in $\nseq(\wh X',\wh Y')$ to a commutative square
\[\xymatrix{
  \wh D \ar[r]^-{k}\ar[d]_{\ell} & D'\ar[d]^-{i} \\
  \ov D' \ar[r]^-{j} & D
}\]  
Let the morphisms $k:\wh D\la D'$ and $\ell:\wh D\la \ov D'$ be
represented by the 2-commutative diagrams
\[
\xymatrix@C+30pt@R-20pt{
&  \wh X''\ar[dr]^-{\wh p}\ar[dd]_-{\tau} &  & &
  \wh X''\ar[dr]^-{\wh p}\ar[dd]_-{\ov\tau} &  \\
\wh X'\ar[ru]^-{\wh u}\ar[dr]_{w'}&  & \wh Y'& \wh X'\ar[ur]^-{\wh u}\ar[dr]_-{\ov w'} && \wh Y'\\
 & Q\ar[ur]_-{ r\zeta} &
& &\ov Q\ar[ur]_-{r\xi} &
}\]
and let the morphisms $\delta:H\la E'$ and $\delta':\ov H\la G'$ be represented
by the 2-commutative diagrams
\[
\xymatrix@C+30pt@R-20pt{
&  Q\ar[dr]^-{\beta'_\sigma r\zeta}\ar[dd]_-{\chi} &  & &
  \ov Q\ar[dr]^-{\beta'_\sigma r\xi}\ar[dd]_-{\ov\chi} &  \\
\wh X'\ar[ru]^-{w'}\ar[dr]_{u}&  & Y'& \wh X'\ar[ur]^-{\ov w'}\ar[dr]_-{\wh u} && \ov Y'\\
 & X''\ar[ur]_-{ p} &
& &\ov X''\ar[ur]_-{\ov p} &
}\]
This produces for us the 2-commutative diagram
\[
\xymatrix{
  X\ar[dr]_-{\wh f} \ar[r]^{\wh u'} & \wh X'\ar[d]^{\wh p'}\ar[r]^{\wh u} &\wh X''\ar[d]^{\wh p} \\
 &       Y\ar[dr]_-{\wh g} \ar[r]^-{\wh v} &\wh Y'\ar[d]^{\wh q} \\
 & &  Z
}
\]
which we view as an object $a\in\nseq(X,Y,Z)$, as well as 
1-morphisms $\alpha'_\sigma:\wh X'\la X'$,
$\alpha'_\theta:\wh X'\la\ov X'$, $\beta'_\sigma:\wh Y'\la Y'$
and $\beta'_\theta:\wh Y'\la\ov Y'$, $\chi\tau:\wh X''\la X''$
and $\ov\chi\ov \tau:\wh X''\la \ov X''$, together with a great deal
of 2-commutativity that assembles to give the morphisms
$b\stackrel\sigma\longleftarrow a\stackrel\theta\la c$ as in \ref{L6.2059.5}.
We leave the details to the reader.
\eprf

  \lem{L6.2060}
Let $F:\nseq(X,Y,Z)\la\seq(X,Y)\times\seq(Y,Z)$ be the functor of
Remark~\ref{R6.2057}, and 
let $\ph,\psi$ be two morphisms in $\nseq(X,Y,Z)$ with $F(\ph)=F(\psi)$.
Then there exists
a morphism $\rho\in\nseq(X,Y,Z)$ with $\ph\rho=\psi\rho$.
\elem

\prf
We are given two morphisms
  $\xymatrix@C+15pt{b\ar@<0.5ex>[r]^{\ph} \ar@<-0.5ex>[r]_{\psi} & c }$
  in the category $\nseq(X,Y,Z)$, satisfying $F(\ph)=F(\psi)$. 
Applying the functors 
$\pi_{12}^{}$ and $\pi_{23}^{}$ of Remark~\ref{R6.2057}
we obtain the pairs of  morphisms
\[
\xymatrix@C+15pt{\pi_{12}^{}(b)\ar@<0.5ex>[r]^{\pi_{12}^{}(\ph)} \ar@<-0.5ex>[r]_{\pi_{12}^{}(\psi)} & \pi_{12}^{}(c)  & &
  \pi_{23}^{}(b)\ar@<0.5ex>[r]^{\pi_{23}^{}(\ph)} \ar@<-0.5ex>[r]_{\pi_{23}^{}(\psi)} & \pi_{23}^{}(c)
}
\]
in the categories $\nseq(X,Y)$ and $\nseq(Y,Z)$,
with $F\big(\pi_{12}^{}(\ph)\big)=F\big(\pi_{12}^{}(\psi)\big)$
and $F\big(\pi_{23}^{}(\ph)\big)=F\big(\pi_{23}^{}(\psi)\big)$.
  By Lemma~\ref{L6.2054} there exist morphisms $\rho^{}_{12}\in\nseq(X,Y)$
and $\rho^{}_{23}\in\nseq(Y,Z)$
so that the two composites
\[
\xymatrix@C-5pt{a_{12}^{}\ar[r]^-{\rho^{}_{12}} &\pi_{12}^{}(b)\ar@<0.5ex>[rr]^{\pi_{12}^{}(\ph)} \ar@<-0.5ex>[rr]_{\pi_{12}^{}(\psi)} && \pi_{12}^{}(c)
  & &
 a_{23}^{}\ar[r]^-{\rho^{}_{23}} &\pi_{23}^{}(b)\ar@<0.5ex>[rr]^{\pi_{23}^{}(\ph)} \ar@<-0.5ex>[rr]_{\pi_{23}^{}(\psi)} && \pi_{23}^{}(c)}
 \]
 are equal.
We assemble $a_{12}^{}$ and $a_{23}^{}$ together into the diagram $A$ below
\[
\xymatrix{
  X\ar[dr]_-{\wh f} \ar[r]^{\wh u'} & \wh X'\ar[d]^{\wh p'} \\
 &       Y\ar[dr]_-{\wh g} \ar[r]^-{\wh v} &\wh Y'\ar[d]^{\wh q} \\
 & &  Z
}
\]
where $a_{12}^{}=\pi_{12}^{}(A)$ and $a_{23}^{}=\pi_{23}^{}(A)$, and the morphisms
$\rho_{12}^{}$ and $\rho_{23}^{}$ can be rewritten as
 maps
\[
\rho_{12}^{}:\pi_{12}^{}(A)\la\pi_{12}^{}(b),\quad 
\rho_{23}^{}:\pi_{23}^{}(A)\la\pi_{23}^{}(b).
\]
Now we apply \ref{L6.2059.5} to the identity map
$\id:F(b)\la F(b)$, to the diagram $A$ above, and with
$\sigma'=\theta'=\rho_{12}^{}$ and $\sigma''=\theta''=\rho_{23}^{}$.
From \ref{L6.2059.5} we learn that there exists an object
$a\in\nseq(X,Y,Z)$ with $\pi_{12}^{}(a)=\pi_{12}^{}(A)$ and $\pi_{23}^{}(a)=\pi_{23}^{}(A)$,
and a morphism $\rho:a\la b$ with
$\pi_{12}^{}(\rho)=\rho_{12}^{}$ and 
$\pi_{23}^{}(\rho)=\rho_{23}^{}$. Replacing the pair of morphisms
\[
\xymatrix@C+15pt{
  b\ar@<0.5ex>[r]^{\ph} \ar@<-0.5ex>[r]_{\psi} & c 
  &\text{by the composites} &
  a\ar[r]^{\rho} &b\ar@<0.5ex>[r]^{\ph} \ar@<-0.5ex>[r]_{\psi} & c 
}\]
we are reduced to proving the Lemma in the special case
where $\pi_{12}^{}(\ph)=\pi_{12}^{}(\psi)$ and $\pi_{23}^{}(\ph)=\pi_{23}^{}(\psi)$.
To explain it all more concretely: the objects
  $b,c\in\nseq(X,Y,Z)$ are given by 2-commutative diagrams
in $\seq$ satisfying all the requirements of Definition~\ref{D6.2056}(i)
 \[
  \xymatrix{
    X\ar[dr]_-f \ar[r]^-{u'}  &X' \ar[r]^{u} \ar[d]^{p'} &  X''\ar[d]^p & &
    X\ar[dr]_-{\ov f} \ar[r]^-{\ov u'}  &\ov X' \ar[r]^{\ov u} \ar[d]^{\ov p'} &  \ov X''\ar[d]^{\ov p}\\
 &   Y\ar[dr]_-g \ar[r]^-v & Y'\ar[d]^q &\text{and} &&   Y\ar[dr]_-{\ov g} \ar[r]^-{\ov v} &\ov Y'\ar[d]^{\ov q} \\
 & &  Z & & & & Z
}
  \]
  The morphisms $\xymatrix{b\ar@<0.5ex>[r]^\ph \ar@<-0.5ex>[r]_\psi & c}$
  are represented by pairs $\xymatrix{X'\ar@<0.5ex>[r]^{\alpha'_\ph}\ar@<-0.5ex>[r]_{\alpha'_\psi} & \ov X'}$,
$\xymatrix{X''\ar@<0.5ex>[r]^{\alpha''_\ph}\ar@<-0.5ex>[r]_{\alpha''_\psi} & \ov X''}$
  and $\xymatrix{Y'\ar@<0.5ex>[r]^{\beta'_\ph}\ar@<-0.5ex>[r]_{\beta'_\psi} & \ov Y'}$,
  together with a bunch of 2-isomorphisms giving the 2-commutativity
  in Definition~\ref{D6.2056}(ii). The reduction we have proved so far is
  that
  we may assume $\pi_{12}^{}(\ph)=\pi_{12}^{}(\psi)$ and
  $\pi_{23}^{}(\ph)=\pi_{23}^{}(\psi)$, meaning there are
  2-commutative diagrams
 \[\xymatrix@C+30pt@R-20pt{
   & X'\ar^{p'}[dr]\ar@/_/[dd]_{\alpha'_\ph}\ar@/^/[dd]^{\alpha'_\psi} & &
          & Y'\ar^q[dr]\ar@/_/[dd]_{\beta'_\ph}\ar@/^/[dd]^{\beta'_\psi} & \\
   X\ar^{u'}[ur]\ar_{\ov u'}[dr] & &  Y &
   Y\ar^{v}[ur]\ar_{\ov v}[dr] & &  Z\\
   & \ov X'\ar_{\ov p'}[ur] & &
   & \ov Y'\ar_{\ov q}[ur] &
}\] 
where the curved arrows represent the part of the data of $\ph$ and and $\psi$
that pertains to this portion of the diagram. But the equivalence relation
allows us to replace $\alpha'_\ph$ and $\beta'_\ph$ by the isomorphic
$\alpha'_\psi$ and $\beta'_\psi$, as long as we're careful to perturb
the rest of the data to match. Doing this we may
choose representatives of $\ph$ and $\psi$ so that on the $\pi_{12}^{}$ and
$\pi_{23}^{}$ portions they are identical: we may assume that our representatives
for $\ph$ and $\psi$ give the same data
 \[\xymatrix@C+30pt@R-20pt{
   & X'\ar^{p'}[dr]\ar[dd]_{\alpha'} & &
          & Y'\ar^q[dr]\ar[dd]_{\beta'} & \\
   X\ar^{u'}[ur]\ar_{\ov u'}[dr] & &  Y &
   Y\ar^{v}[ur]\ar_{\ov v}[dr] & &  Z\\
   & \ov X'\ar_{\ov p'}[ur] & &
   & \ov Y'\ar_{\ov q}[ur] &
}\] 
and Remark~\ref{R6.2052.3} guarantees that $\alpha'$ and $\beta'$ are
of finite type and universally quasi-proper.
 
 The above gives a single 2-morphism $\beta'vp'\la\ov v\ov p'\alpha'$, but the
 morphisms $\ph$ and $\psi$ give two liftings of it
  \[\xymatrix@C+30pt@R-20pt{
   & X''\ar^{\beta'p}[dr]\ar@/_/[dd]_{\alpha''_\ph}\ar@/^/[dd]^{\alpha''_\psi} &  \\
   X'\ar^{u}[ur]\ar_{\ov u\alpha'}[dr] & &  \ov Y'\\
   & \ov X''\ar_{\ov p}[ur] &
  }\]
which are not necessarily equivalent, in other words what is pictured above
is  two 2-commutative diagrams, one with $\alpha''_\ph$ and one with $\alpha''_\psi$,
which might not compress into a single 2-commutative diagram. Anyway
these define two morphisms
in $\rnseq(X',\ov Y')$ with the same image under $F$.
Lemma~\ref{L6.2054} assures us that there exists
an equalizer;
concretely the above can be extended to a diagram
\[\xymatrix@C+50pt@R-5pt{
   & \wh X''\ar^{\wh s}[dr]\ar[d]_{\alpha''_\rho} &  \\
   X'\ar^{\wh u}[ur]\ar|-{u}[r]\ar_{\ov u\alpha'}[dr] & X''\ar|-{\beta'p}[r]\ar@/_/[d]_{\alpha''_\ph}\ar@/^/[d]^{\alpha''_\psi}  &  \ov Y'\\
   &\ov X''\ar_{\ov p}[ur] &
}\]
which becomes 2-commutative if we delete the middle row.
Now Lemma~\ref{L6.2054} allows us to choose
$X'\stackrel{\wh u}\la X''\stackrel{\wh s}\la \ov Y'$
to be in $\nseq(X',\ov Y')$, meaning we may assume
that $\wh u$ is a dominant, flat monomorphism and,
by Remark~\ref{R6.2052.3}, we may also assume $\alpha''_\rho$
of
finite type and univerally quasi-proper.
We construct the 2-commutative
diagram $a$ below
 \[
  \xymatrix{
    X\ar[dr]_-f \ar[r]^-{u'}  &X' \ar[r]^{\wh u} \ar[d]^{p'} &  \wh X''\ar[d]^{p\alpha''_\rho} \\
 &   Y\ar[dr]_-g \ar[r]^-v & Y'\ar[d]^q \\
 & &  Z 
}
  \]
The horizontal maps are dominant, flat monomorphisms
and the vertical maps are of finite type and universally quasi-proper,
hence $a$ is an object of $\nseq(X,Y,Z)$.
And the maps $\id :X'\la X'$, $\alpha''_\rho:\wh X''\la X''$ and $\id:Y'\la Y'$
assemble to give a morphism $\rho:a\la b$, and we remind the reader that $a$ and $b$
are the 2-commutative diagrams
 \[
  \xymatrix{
    X\ar[dr]_-f \ar[r]^-{u'}  &X' \ar[r]^{\wh u} \ar[d]^{p'} &  \wh X''\ar[d]^{p\alpha''_\rho} &&    X\ar[dr]_-f \ar[r]^-{u'}  &X' \ar[r]^{u} \ar[d]^{p'} &  X''\ar[d]^{p} \\
 &   Y\ar[dr]_-g \ar[r]^-v & Y'\ar[d]^q & &  &   Y\ar[dr]_-g \ar[r]^-v & Y'\ar[d]^q \\
 & &  Z & &  & &  Z 
}
  \]
  We leave it to the reader to check that the composites 
$\xymatrix@C+15pt{a\ar[r]^\rho &b\ar@<0.5ex>[r]^{\ph} \ar@<-0.5ex>[r]_{\psi} & c }$
  are equal in $\nseq(X,Y,Z)$, meaning the representatives we
  constructed are equivalent.
  \eprf

\pro{P6.2061}
In Definition~\ref{D6.2056} we constructed a category $\nseq(X,Y,Z)$ and
in Remark~\ref{R6.2057} we noted that there is an obvious functor
$F:\nseq(X,Y,Z)\la\seq(X,Y)\times\seq(Y,Z)$. Since $\seq(X,Y)\times\seq(Y,Z)$
is a groupoid the functor $F$ factors canonically as
$\nseq(X,Y,Z)\stackrel{\pi}\la\gnseq(X,Y,Z)\stackrel{F'}\la\seq(X,Y)\times
\seq(Y,Z)$,
where $\gnseq(X,Y,Z)$ is the groupoid completion of $\nseq(X,Y,Z)$.

In this canonical factorization the functor $F'$ is an equivalence.
\epro

\prf
Lemma~\ref{L6.2055}, coupled with Lemmas~\ref{L6.2059} and \ref{L6.2060},
tell us that the functor $F'$ is fully faithful, while Lemma~\ref{L6.2058}
establishes the essential surjectivity.
\eprf

\dfn{D6.2062}
We will also need to refer to the four-object version, but fortunately
we don't care about its categorical structure, only the class of objects.
We therefore define $\nseq(W,X,Y,Z)\subset\lnseq(W,X,Y,Z)$ to be the classes
whose elements
are 2-commutative diagrams 
  \[
  \xymatrix{
    W\ar[dr]_-e \ar[r]^-{t''}  &W' \ar[r]^{t'} \ar[d]^{o''} &  W''\ar[d]^{o'}\ar[r]^t & W'''\ar[d]^o\\
&    X\ar[dr]_-f \ar[r]^-{u'}  &X' \ar[r]^{u} \ar[d]^{p'} &  X''\ar[d]^p\\
& &   Y\ar[dr]_-g \ar[r]^-v & Y'\ar[d]^q\\
& & &  Z
}
  \]
In both $\nseq(W,X,Y,Z)$ and $\lnseq(W,X,Y,Z)$ the squares
are assumed 2-cartesian and the horizontal maps are assumed flat.
For $\lnseq(W,X,Y,Z)$ these are the only restrictions, while in
$\nseq(W,X,Y,Z)$ we further restrict 
the horizontal maps to be dominant monomorphisms, 
and the vertical maps to be of finite type and
universally quasi-proper.
\edfn

\rmk{R6.2062.5}
As in Remark~\ref{R6.2057} we will want to consider some forgetful
maps out of $\nseq(W,X,Y,Z)$. For example the map $\pi_{124}^{}$ forgets
the column and row containing  $Y$ in the diagram of Definition~\ref{D6.2062}.
That is it leaves us with the diagram
  \[
  \xymatrix{
    W\ar[dr]_-e \ar[r]^-{t''}  &W' \ar[r]^{tt'} \ar[d]^{o''} & W'''\ar[d]^o\\
&    X\ar[dr]_-{gf} \ar[r]^-{uu'}  &  X''\ar[d]^{qp}\\
& &  Z
}
\]
Mercifully the only fact we will need about $\nseq(W,X,Y,Z)$ is that
it has plenty of objects.
\ermk

\lem{L6.2063}
The obvious map $F:\nseq(W,X,Y,Z)\la\seq(W,X)\times\seq(X,Y)\times\seq(Y,Z)$
is essentially surjective.
\elem

\prf
Suppose we are given an object $W\stackrel e\la X\stackrel f\la Y\stackrel g\la Z$
in $\seq(W,X)\times\seq(X,Y)\times\seq(Y,Z)$. 
Choose Nagata compactifications for $e$, $f$ and $g$, that is construct a
2-commutative diagram
\[
\xymatrix{
  W\ar[r]^{t''} \ar[dr]_e & W'\ar[d]^{o''}  \\
&    X\ar[dr]_-f \ar[r]^-{u'}  &X'  \ar[d]^{p'} \\
 & &   Y\ar[dr]_-g \ar[r]^-v & Y'\ar[d]^q\\
 && &  Z
  }
\]
where the horizontal maps are dominant, flat monomorphisms and the vertical maps
are of finite type and universally quasi-proper. Now
apply Lemma~\ref{L6.2058} to the composable morphisms
$W'\stackrel{u'o''}\la X'\stackrel{vp'}\la Y$
to complete to a 2-commutative diagram
\[
  \xymatrix{
    W\ar[dr]_-e \ar[r]^-{t''}  &W'\ar[d]^{o''}\ar[r]^{t'} & W''\ar[r]^{t}\ar[d]^{o'} & W'''\ar[d]^o\\
&    X\ar[dr]_-f \ar[r]^-{u'}  &X' \ar[r]^{u} \ar[d]^{p'} &  X''\ar[d]^p\\
& &   Y\ar[dr]_-g \ar[r]^-v & Y'\ar[d]^q\\
& & &  Z
}
  \]
  where the horizontal maps are dominant, flat monomorphisms and the vertical maps are
  of finite type and universally quasi-proper. The fact that the squares are
  2-cartesian is by Corollary~\ref{CNag.22}.
  \eprf

\section{The 2-functor $(-)^!$}
\label{S200}

In the last section we met a great many auxiliary categories, and the
reader might wonder what the point is. Now we begin to explain, starting
with the objects.

\rmk{R200.1}
On the category of algebraic stacks there is the assignment that
takes a stack $X$ to its derived category $\Dqc(X)$. It can
be made into a 2-functor in several ways, in
this remark we recall the following two.
Given a 1-morphism
$f:X\la Z$ there are 1-morphisms $f^*,f^\times:\Dqc(Z)\la\Dqc(X)$, and
given a 2-morphism $f\la g$ there are induced 2-isomorphisms
$f^*\la g^*$ and $f^\times\la g^\times$. Since $(gf)^*$ is canonically isomorphic
to $f^*g^*$ and $(gf)^\times$ is canonically isomorphic to $f^\times g^\times$,
we have 2-functors $(-)^*$ and $(-)^\times$, defined
on the entire 2-category of algebraic stacks. 
We will mostly be interested in the restriction
to the 2-subcategory $\seq$, that is
in the
2-functors $(-)^*:\seq\la\Tri$ and
$(-)^\times:\seq\la\Tri$. Here $\Tri$ is the
2-category whose objects are triangulated categories, the 1-morphisms
are triangulated functors, and the 2-morphisms are natural transformations.
The 2-functors $(-)^*$ and $(-)^\times$  give, for every pair of objects
$X,Z\in\seq$, 
functors
\[
\begin{array}{lcccc}
(-)^*&:&\seq(X,Z)&\la&\Tri\big(\Dqc(Z),\Dqc(X)\big),\\
(-)^\times&:&\seq(X,Z)&\la&\Tri\big(\Dqc(Z),\Dqc(X)\big).
\end{array}\]
In the last section we saw the categories $\nseq(X,Z)$ and
$\gnseq(X,Z)$ and studied some functors between them.
We obtain composites
\[
\CD
\nseq(X,Z) @>\pi>> \gnseq(X,Z)@>F'>>\seq(X,Z) @>(-)^*>> \Tri\big(\Dqc(Z),\Dqc(X)\big)\\
\nseq(X,Z)@>\pi>>\gnseq(X,Z)@>F'>>\seq(X,Z) @>(-)^\times>> \Tri\big(\Dqc(Z),\Dqc(X)\big)\\
\endCD
\]
where $F'$ is an equivalence by Proposition~\ref{P6.2055.5}.
Concretely these composites take an object in $\nseq(X,Z)$, i.e.\ a
2-commutative triangle
\[\xymatrix@C+30pt@R-15pt{
       & Y\ar^p[dr]& \\
X\ar^u[ur]\ar[rr]|f & &  Z \\
}\] 
to $f^*$ (respectively $f^\times$). But on the category $\nseq(X,Y)$ there
are equivalent functors: the definition of an object in $\nseq(X,Y)$
gives a 2-isomorphism $pu\la f$, hence
2-isomorphisms $u^*p^*\la f^*$ and $ u^\times p^\times\la f^\times$.
The reader can convince herself that the map taking an object in
$\nseq(X,Z)$ to $u^*p^*$ (respectively $u^\times p^\times$) extends to a functor
in such a way 
that the comparison isomorphism
$u^*p^*\la f^*$ (respectively $u^\times p^\times\la f^\times$) is a natural
transformation.
\ermk

\con{C200.3}
So far all the categories $\nseq(X,Z)$ seem to have bought us are complicated
functors equivalent to the simple $(-)^*$ and $(-)^\times$. Now we begin
to change this.
\be
\item
  For any object $a\in\lnseq(X,Z)$, that is for a 2-commutative triangle
\[\xymatrix@C+30pt@R-15pt{
       & Y\ar^p[dr]& \\
X\ar^u[ur]\ar[rr]|f & &  Z \\
}\] 
we define $P(a)$ to
be the object in $\Tri\big(\Dqc(Z),\Dqc(X)\big)$ given
by the composite
$P(a)=u^*p^\times$. That is $P(a)$ is the composite 1-morphism
$\Dqc(Z)\stackrel{p^\times}\la\Dqc(Y)\stackrel{u^*}\la\Dqc(X)$ in
$\Tri$.
\item
  For any object $a\in\lnseq(X,Y,Z)$ we define $P(a)$ to be a 2-morphism
  in $\Tri$, more specifically a morphism in $\Tri\big(\Dqc(Z),\Dqc(X)\big)$.
  Even more specifically it will be a 2-morphism $P\big(\pi_{13}^{}(a)\big)\la
  P\big(\pi_{12}^{}(a)\big)P\big(\pi_{23}^{}(a)\big)$.

  After all the preliminaries it's time to give the definition.
  An object $a\in\lnseq(X,Y,Z)$ is a diagram
  \[
  \xymatrix@C+10pt{
    X\ar[dr]_-f \ar[r]^-{u'}  &X' \ar[r]^{u} \ar[d]^{p'}\ar@{}[dr]|-{(\diamondsuit)} &  X''\ar[d]^p\\
 &   Y\ar[dr]_-g \ar[r]^-v & Y'\ar[d]^q\\
 & &  Z
}
  \]
while $\pi_{12}^{}(a)$, $\pi_{13}^{}(a)$ and $\pi_{22}^{}(a)$ are, respectively,
the 2-commutative diagrams
\[\xymatrix{
  X\ar[dr]_-f \ar[r]^-{u'}  &X' \ar[d]^{p} && X\ar[dr]_-{gf} \ar[r]^-{uu'}  &
  X''\ar[d]^{qp}  & &  Y\ar[dr]_-g \ar[r]^-v & Y'\ar[d]^q\\
 &   Y  & & & Z
 & &  &Z
}\]
We have promised the reader a 2-morphism
 $P\big(\pi_{13}^{}(a)\big)\la
P\big(\pi_{12}^{}(a)\big)P\big(\pi_{23}^{}(a)\big)$, meaning a map
$(uu')^*(qp)^\times\la ({u'}^*{p'}^\times)(v^*q^\times)$. And the formula is
that $P(a)$ is the composite
\[
\CD
(uu')^*(qp)^\times @>\sim>>{u'}^*(u^*{p}^\times) q^\times @>{u'}^*\Phi(\diamondsuit)q^\times>>{u'}^*({p'}^\times v^*)q^\times @>\sim>> ({u'}^*{p'}^\times)(v^*q^\times) 
\endCD
\]
where the first map is by the canonical isomorphisms of Remark~\ref{R200.1}
coupled with associativity,
the second map is induced by the base-change $\Phi(\diamondsuit)$
for the 2-cartesian square $(\diamondsuit)$ in the diagram defining $a$,
and the third is associativity of composition. Note that $v$ is
assumed flat
so the base-change map is well defined.
\item
  When $a$ is an element of $\lnseq(W,X,Y,Z)$, it turn out that $P(a)$ is a
  relation. Since something has to be proved we state it as a lemma.
\ee
\econ

\lem{L200.5}
With the notation as in Definition~\ref{D6.2062} and Remark~\ref{R6.2062.5},
any object $a\in\lnseq(W,X,Y,Z)$ yields an equality in
$\Tri\big(\Dqc(Z),\Dqc(W)\big)$ of the morphisms
\[
P\big(\pi_{123}^{}(a)\big)P\big(\pi_{134}^{}(a)\big)\eq
P\big(\pi_{234}^{}(a)\big)P\big(\pi_{124}^{}(a)\big)\ .
\]
\elem

\prf
The object $a\in\lnseq(W,X,Y,Z)$ is a diagram
  \[
  \xymatrix@C+10pt{
    W\ar[dr]_-e \ar[r]^-{t''}  &W' \ar[r]^{t'} \ar[d]^{o''}\ar@{}[dr]|-{(\diamondsuit)} &  W''\ar[d]^{o'}\ar[r]^t \ar@{}[dr]|-{(\heartsuit)}& W'''\ar[d]^o\\
&    X\ar[dr]_-f \ar[r]^-{u'}  &X' \ar[r]^{u} \ar[d]^{p'}\ar@{}[dr]|-{(\spadesuit)} &  X''\ar[d]^p\\
& &   Y\ar[dr]_-g \ar[r]^-v & Y'\ar[d]^q\\
& & &  Z
}
\]
and the objects $\pi_{123}^{}(a)$, $\pi_{134}^{}(a)$, $\pi_{234}^{}(a)$ and
$\pi_{124}^{}(a)$ are, respectively, the diagrams
  \[
  \xymatrix@C+10pt@R+10pt{
    W\ar[dr]_-e \ar[r]^-{t''}  &W' \ar[r]^{t'} \ar[d]^{o''}\ar@{}[dr]|-{(\diamondsuit)} &  W''\ar[d]^{o'} & &
    W\ar[dr]_-{fe} \ar[r]^-{t't''}  &W'' \ar[r]^{t} \ar[d]|-{p'o'}\ar@{}[dr]|-{\sst\left(\!\!\!\!
      \begin{array}{c}\heartsuit\\*[-5pt]
        \spadesuit\end{array}\!\!\!\!\right)} & W'''\ar[d]^-{po}
    \\
    &    X\ar[dr]_-f \ar[r]^-{u'}  &X'  \ar[d]^{p'}& &
   &   Y\ar[dr]_-g \ar[r]^-v & Y'\ar[d]^q    \\
& & Y & & & &  Z
}
  \]
\[
  \xymatrix@C+10pt{
   X\ar[dr]_-f \ar[r]^-{u'}  &X' \ar[r]^{u} \ar[d]^{p'}\ar@{}[dr]|-{(\spadesuit)} &  X''\ar[d]^p  & &
W\ar[dr]_-e \ar[r]^-{t''}  &W' \ar[r]^{tt'} \ar[d]^{o''}\ar@{}[dr]|-{(\diamondsuit\heartsuit)} &   W'''\ar[d]^o
    \\
&   Y\ar[dr]_-g \ar[r]^-v & Y'\ar[d]^q& &      &    X\ar[dr]_-{gf} \ar[r]^-{uu'}  & X''\ar[d]^{qp}\\
& & Z & & & &  Z
}
\]
The maps
$P\big(\pi_{234}^{}(a)\big)P\big(\pi_{124}^{}(a)\big)$
and $P\big(\pi_{123}^{}(a)\big)P\big(\pi_{134}^{}(a)\big)$
are, respectively,
the composites
\[
\xymatrix{
(tt't'')^*(qpo)^\times \ar[r]^-{\sim} & {t''}^*{t'}^*t^*o^\times p^\times q^\times
\ar[r]^-{\Phi(\diamondsuit\heartsuit)} & {t''}^*{o''}^\times {u'}^*u^*p^\times q^\times
\ar[r]^-{\Phi(\spadesuit)} & {t''}^*{o''}^\times {u'}^*{p'}^\times v^* q^\times\\
(tt't'')^*(qpo)^\times \ar[r]^-{\sim} & {t''}^*{t'}^*t^*o^\times p^\times q^\times
\ar[r]^-{\sst\Phi\left(\!\!\!\!\begin{array}{c}\heartsuit\\*[-5pt]
  \spadesuit\end{array}\!\!\!\!\right)} & {t''}^*{t'}^* {o'}^\times {p'}^\times
v^* q^\times
\ar[r]^-{\Phi(\diamondsuit)} & {t''}^*{o''}^\times {u'}^*{p'}^\times v^* q^\times
}
\]
and the equality of these composites comes
from the commutativity of the square in
\[
\xymatrix@C+20pt{
 {t''}^*{t'}^*t^*o^\times p^\times q^\times
 \ar[r]^-{\Phi(\heartsuit)} &
{t''}^*{t'}^*{o'}^\times {u}^*p^\times q^\times
\ar[r]^-{\Phi(\diamondsuit)} \ar[d]_-{\Phi(\spadesuit)}&
 {t''}^*{o''}^\times {u'}^*u^*p^\times q^\times
\ar[d]^-{\Phi(\spadesuit)} \\
  & {t''}^*{t'}^* {o'}^\times {p'}^\times
v^* q^\times
\ar[r]^-{\Phi(\diamondsuit)} & {t''}^*{o''}^\times {u'}^*{p'}^\times v^* q^\times
}
\]
and the fact that base-change maps concatenate, that is
$\Phi(\diamondsuit\heartsuit)=\Phi(\diamondsuit)\Phi(\heartsuit)$ and
$\Phi\left(\!\!\!\!\begin{array}{c}\heartsuit\\*[-5pt]
  \spadesuit\end{array}\!\!\!\!\right)=\Phi(\spadesuit)\Phi(\heartsuit)$.
  \eprf

\con{C200.7}
Now it's time to construct $(-)^!$.
\be
\item
If $a\in\nseq(X,Z)$ is an object we define $a^!=P(a)$, with $P(a)$ as in
Construction~\ref{C200.3}(i). That is $a^!=u^*p^\times$ with the notation
of Construction~\ref{C200.3}(i).
\setcounter{enumiii}{\value{enumi}}
\ee
This deals with the objects of $\nseq(X,Z)$, next we worry about
the morphisms.
Suppose $\ph:a\la b$ \emph{represents} a morphism in $\nseq(X,Z)$. Recall:
the morphisms in $\nseq(X,Z)$ are equivalence classes of 2-commutative
diagrams, and for the construction we choose a representative.
For us right now the relevant part of
the data given by the representative is the 2-commutative
\[\xymatrix@C+30pt@R-15pt{
       & Y\ar^p[dr]\ar[dd]^\alpha & \\
X\ar^u[ur]\ar_{u'}[dr]& &  Z \\
       & Y'\ar_{p'}[ur] & 
}\] 
Out of this we construct the 2-commutative diagram $C$ below
  \[
  \xymatrix@C+20pt{
    X\ar[dr]_-\id \ar[r]^-{\id}  &X \ar[r]^{u}\ar@{}[rd]|-{(\diamondsuit)} \ar[d]^{\id} &  Y\ar[d]^\alpha\\
 &   X\ar[dr]_-f \ar[r]^-{u'} & Y'\ar[d]^{p'}\\
 & &  Z
}
\]
We are given that $p'$ is separated, of finite type
and universally quasi-proper, and Remark~\ref{R6.2052.3} tells
us that so is $\alpha$. Thus the vertical maps are all
of finite type and universally quasi-proper while the horizontal maps
are 
dominant, flat monomorphisms. The diagram is an object
$C\in\nseq(X,X,Z)$.

By Construction~\ref{C200.3}(ii) we have a morphism
 $P(C):P\big(\pi_{13}^{}(C)\big)\la
P\big(\pi_{12}^{}(C)\big)P\big(\pi_{23}^{}(C)\big)$, meaning a map
$P(C):{(u\circ\id)}^*{(p'\alpha)}^\times\la \id^*\id^\times {u'}^* {p'}^\times$.
   Note that
   up to isomorphism the 2-morphism $P(C)\in\seq$
   is given 
  by the base-change map
  $\Phi(\diamondsuit):u^*\alpha^\times
  \la \id^\times {u'}^*$,
  and  Lemma~\ref{L4.7}(ii) applies to the square $(\diamondsuit)$ because
  $\id:X'\la X'$ is of finite Tor-dimension. We deduce that $P(C)$
  is an isomorphism. And now we come to
\be
\setcounter{enumi}{\value{enumiii}}
\item
  With the notation as above---that is with $\ph:a\la b$ a representative
  of a morphism in $\nseq(X,Z)$, and with everything else as constructed
  out of $\ph$ in the preceding paragraphs---we
  define $\ph^!:a^!\la b^!$ to be the composite
  \[
  \xymatrix@C+20pt{
    u^*p^\times \ar[r]^-{\sim} & u^*{(p'\alpha)}^\times \ar[r]^-{P(C)} &
    {u'}^*{p'}^\times
  }\]
  where the first isomorphism is the map induced by the 2-isomorphism
  $p\la p'\alpha$ that is part of the data of $\ph$, and $P(C)$
  is as above. Note that $\ph^!$, being the composite of two
  isomorphisms, must be an isomorphism in $\Tri\big(\Dqc(Z),\Dqc(X)\big)$.
  \ee
There is nothing to stop us from defining $(-)^!$ on all
objects of $\lnseq(X,Z)$ following the recipe above, but then
$(-)^!$ will fail to take the morphisms to isomorphisms. As we will see this
is a key property.
\econ

\lem{L200.9}
If $\ph,\psi$ are two equivalent representatives
of a morphism $a\la b$ in the category $\nseq(X,Z)$
then $\ph^!=\psi^!$, with the notation as in
Construction~\ref{C200.7}(ii).
\elem

\prf
Two representatives $\ph,\psi$ of a morphism
in $\nseq(X,Z)$ and an equivalence between them gives, among
other things, the 2-commutative diagram
\[\xymatrix@C+30pt@R-20pt{
       & Y\ar^p[dr]\ar@/_/[dd]_\alpha\ar@/^/[dd]^{\alpha'} & \\
X\ar^u[ur]\ar_{u'}[dr] & &  Z \\
       & Y'\ar_{p'}[ur] & 
}\]
of Definition~\ref{D6.2051}(ii). By Remark~\ref{R6.2052.3} $\alpha$
and $\alpha'$ are of finite type and universally quasi-proper.
Out of this we cook up the 2-commutative diagram $C$ below
 \[
  \xymatrix@C+10pt{
    X\ar[dr]_-\id \ar[r]^-{\id}  &X \ar[r]^{u} \ar[d]^{\id}\ar@{}[dr]|-{(\diamondsuit)} &  Y\ar[d]^{\alpha}\ar[r]^\id & Y\ar[d]^{\alpha'}\\
&    X\ar[dr]_-f \ar[r]^-{u'}  &Y' \ar[r]^{\id} \ar[d]^{p'} &  Y'\ar[d]^{p'}\\
& &   Z\ar[dr]_-\id \ar[r]^-\id & Z\ar[d]^\id\\
& & &  Z
}
  \]
All the horizontal maps
are dominant, flat monomorphisms, all the vertical maps are of finite type
and universally quasi-proper, and Corollary~\ref{CNag.22} tells us that the
squares are all 2-cartesian. Therefore $C$ is an object
in $\nseq(X,X,Z,Z)$. By Lemma~\ref{L200.5} this gives the relation
$P\big(\pi_{123}^{}(C)\big)P\big(\pi_{134}^{}(C)\big)=
P\big(\pi_{234}^{}(C)\big)P\big(\pi_{124}^{}(C)\big)$.
But $P\big(\pi_{234}^{}(C)\big)=\id$ and $P\big(\pi_{134}^{}(C)\big)$ is
the isomorphism $u^*(p'\alpha')^\times\la u^*(p'\alpha)^\times$
  induced by the isomorphism $\alpha'\la\alpha$, and
  we conclude that in the diagram below the triangle commutes
\[
\xymatrix@C+20pt@R+10pt{
 u^*p^\times \ar[r]^-\sim & u^*(p'\alpha')^\times\ar[dr]_-{P\big(\pi_{124}^{}(C)\big)} \ar[r]^-{\sim} & u^*(p'\alpha)^\times \ar[d]^-{P\big(\pi_{123}^{}(C)\big)}\\
 & &   {u'}^*{p'}^\times
}\]
The two equal composites, from top left to bottom right,
are by definition $\ph^!$ and $\psi^!$.
\eprf

\lem{L200.11}
The assignment $(-)^!$ defines a functor $\nseq(X,Z)\la
\Tri\big(\Dqc(Z),\Dqc(X)\big)$.
\elem

\prf
In Construction~\ref{C200.7}(i) we defined $a^!$ for $a$ an object of
$\nseq(X,Z)$, and in Construction~\ref{C200.7}(ii)
we defined $\ph^!:a^!\la b^!$,
whenever $\ph:a\la b$ is the \emph{representative} of a morphism.
In Lemma~\ref{L200.9} we showed that equivalent representatives map
to equal morphisms in $\Tri\big(\Dqc(Z),\Dqc(X)\big)$, hence we have
an unambiguous construction that sends objects to objects and morphisms
to morphisms. Identities obviously map
to identities, and it remains to show that $(-)^!$ respects composition. 

Let $a\stackrel\ph\la b\stackrel\psi\la c$ be composable morphisms
in $\nseq(X,Z)$. Choose representatives; part of the data this gives
is the 2-commutative diagram
\[\xymatrix@C+50pt@R-2pt{
       & Y\ar^p[dr]\ar[d]_{\alpha}& \\
X\ar^u[ur]\ar_{u''}[dr]\ar[r]|{u'} &
   Y'\ar[r]|{p'}  \ar[d]_{\beta}&  Z \\
       & Y''\ar_{p''}[ur] & 
}\]
Out of this we construct the 2-commutative diagram
 \[
  \xymatrix@C+10pt{
    X\ar[dr]_-\id \ar[r]^-{\id}  &X \ar[r]^{\id} \ar[d]^{\id} &  X\ar[d]^{\id}\ar[r]^{u}\ar@{}[dr]|-{(\heartsuit)} & Y\ar[d]^{\alpha}\\
&    X\ar[dr]_-\id \ar[r]^-{\id}  &X \ar[r]^{u'} \ar[d]^{\id} &  Y'\ar[d]^{\beta}\\
& &   X\ar[dr]_-f \ar[r]^-{u''} & Y\ar[d]^{p''}\\
& & &  Z
}
  \]
We are given that $p''$ is
of finite type and universally quasi-proper, and Remark~\ref{R6.2052.3}
says that 
so are $\alpha$ and $\beta$. Hence the vertical maps
are all of finite type and universally quasi-proper.
The horizontal maps
are given to be dominant, flat monomorphisms, and Corollary~\ref{CNag.22}
informs us that the squares are all 2-cartesian.
Thus the diagram defines an object $C\in\nseq(X,X,X,Z)$. From
Lemma~\ref{L200.5} we learn that 
$P\big(\pi_{123}^{}(C)\big)P\big(\pi_{134}^{}(C)\big)=
P\big(\pi_{234}^{}(C)\big)P\big(\pi_{124}^{}(C)\big)$;
observing that $P\big(\pi_{123}^{}(C)\big)=\id$ this simplifies to
$P\big(\pi_{134}^{}(C)\big)=
P\big(\pi_{234}^{}(C)\big)P\big(\pi_{124}^{}(C)\big)$.
Let $C'\in\nseq(X,X,Z)$ be the object
 \[
  \xymatrix@C+10pt{
    X\ar[dr]_-\id \ar[r]^-{\id}  &  X\ar[d]^{\id}\ar[r]^{u}\ar@{}[dr]|-{(\heartsuit)} & Y\ar[d]^{\alpha}\\
&    X\ar[dr]_-f \ar[r]^-{u'}   &  Y'\ar[d]^{p'}\\
 & &  Z
}
  \]
  which is isomorphic to $\pi_{124}^{}(C)$ via the isomorphism $p'\la p''\beta$,
  and consider the diagram
\[
\xymatrix@C+30pt@R+10pt{
  u^*p^\times \ar[dr]_-\sim \ar[r]^-\sim & u^*(p'\alpha)^\times\ar[d]_-{\sim} \ar[r]^-{P(C')} & {u'}^*(p')^\times \ar[d]^-{\sim}\\
  & u^*(p''\beta\alpha)^\times\ar[dr]_-{P\big(\pi_{134}^{}(C)\big)} \ar[r]^-{P\big(\pi_{124}^{}(C)\big)} & {u'}^*(p''\beta)^\times \ar[d]^-{P\big(\pi_{234}^{}(C)\big)}\\
 & &   {u''}^*{p''}^\times
}\]
The commutativity of the triangle in the bottom right is the identity
$P\big(\pi_{134}^{}(C)\big)=
P\big(\pi_{234}^{}(C)\big)P\big(\pi_{124}^{}(C)\big)$,
and the commutativity of the square is by
definition of the horizontal maps,
both are
induced by the base-change map $\Phi(\heartsuit)$
(see Construction~\ref{C200.3}(ii)). The commutativity
of the perimeter tells us that $\psi^!\ph^!=(\psi\ph)^!$.
\eprf

\rmk{R200.13}
Our lemmas so far tell us
that we have constructed a
functor $(-)^!:\nseq(X,Z)\la\Tri\big(\Dqc(Z),\Dqc(X)\big)$. In
Construction~\ref{C200.7}(ii) we noted that, for any morphism
$\ph:a\la b$ in $\nseq(X,Z)$, the morphism
$\ph^!\in\Tri\big(\Dqc(Z),\Dqc(X)\big)$ is an isomorphism.
Therefore $(-)^!$ factors canonically through the groupoid
completion: there is a canonical factorization as
\[
\xymatrix@C+30pt{
  \nseq(X,Z)\ar[r]^-{\pi} &
   \gnseq(X,Z)\ar[r]^-{(-)^!}& \Tri\big(\Dqc(Z),\Dqc(X)\big)\ .
  }
\]
In Proposition~\ref{P6.2055.5} we showed that the forgetful map
$F':\gnseq(X,Z)\la\seq(X,Z)$ is an equivalence. This means that,
up to canonical equivalence, we have defined a functor
$(-)^!:\seq(X,Z)\la \Tri\big(\Dqc(Z),\Dqc(X)\big)$. More precisely:
the 
functor $F':\gnseq(X,Z)\la\seq(X,Z)$ is not only an
equivalence, $F'$ is surjective on objects and has a right
inverse---constructing one amounts to choosing a preimage for every object.
For every pair of objects $X,Z\in\seq$ choose a right
inverse $R(X,Z):\seq(X,Z)\la\gnseq(X,Z)$ for $F'$.
Make the choice is such a way that,
for every object $X\in\seq$, the object $\id:X\la X$ of $\seq(X,X)$ maps to
the object
\[\xymatrix@C+30pt@R-15pt{
       & Y\ar^\id[dr]& \\
X\ar^\id[ur]\ar[rr]|\id & &  Z \\
}\]
in $\nseq(X,X)$. Now define $(-)^!:\seq(X,Z)\la\Tri\big(\Dqc(Z),\Dqc(X)\big)$
to be the composite
\[\xymatrix@C+30pt{
  \seq(X,Z)\ar[r]^-{R(X,Z)} &
  \gnseq(X,Z) \ar[r]^-{(-)^!} & \Tri\big(\Dqc(Z),\Dqc(X)\big)
}\]
Since all right inverses of $F'$ are canonically isomorphic, the functor
$(-)^!$ is unique up to canonical isomorphism.
And with our choice for what happens to identities we have guaranteed
that $\id^!=\id$.

Note that we are committing the notational crime of writing $(-)^!$
for all of these functors, not distinguishing whether the source
category is $\nseq(X,Z)$, $\gnseq(X,Z)$ or $\seq(X,Z)$.
\ermk

\con{C200.15}
We have three functors $\pi_{12}^{}$, $\pi_{13}^{}$ and $\pi_{23}^{}$ out
of $\nseq(X,Y,Z)$, to $\nseq(X,Y)$, $\nseq(X,Z)$ and $\nseq(Y,Z)$
respectively. We can now compose them with the functors $(-)^!$
of Remark~\ref{R200.13}, and wonder how these three composite
functors might be related. To address this
\be
\item
  Let $a\in\nseq(X,Y,Z)$ be an object. We define the morphism $\rho(a):
  \pi_{13}^{}(a)^!\la \pi_{12}^{}(a)^!\pi_{23}^{}(a)^!$ to be
  $P(a):P\big(\pi_{13}^{}(a)\big)\la P\big(\pi_{12}^{}(a)\big)P\big(\pi_{23}^{}(a)\big)$
  as in Construction~\ref{C200.3}(ii). 
  \ee
\econ

\lem{L200.17}
The assignment $\rho$ of Construction~\ref{C200.15} is a natural
transformation of functors on $\nseq(X,Y,Z)$. 
\elem

\prf
Let $\ph:a\la b$ be a representative of a morphism in $\nseq(X,Y,Z)$,
that is a 2-commutative diagram
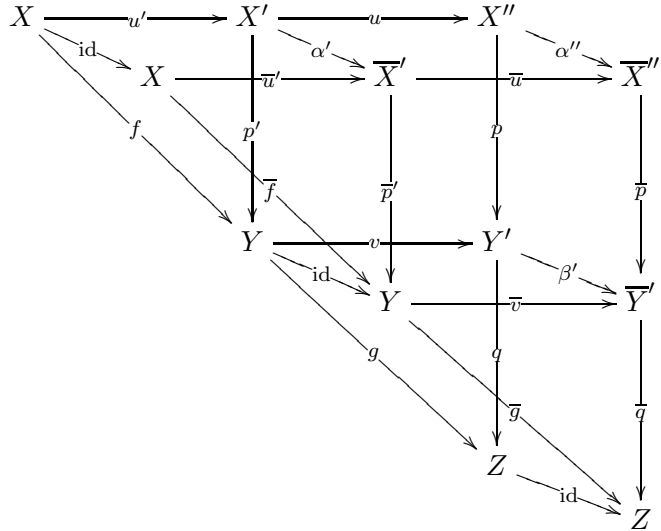
\[
\xymatrix@R-18pt@C-18pt{
X\ar[dddddrrrrr]|-f \ar[drrr]|{\id}\ar[rrrrr]|{u'} &&&&&
 X' \ar[drrr]|{\alpha'}\ar[ddddd]|{p'}  
          \ar[rrrrr]|{u} &&&&& X''
     \ar[drrr]|{\alpha''}\ar[ddddd]|{p}&&&\\
&&& X\ar[dddddrrrrr]|-{\ov f} \ar[rrrrr]|{\ov u'} &&&&& \ov X'\ar[ddddd]|{\ov p'}  
          \ar[rrrrr]|{\ov u} &&&&& \ov X''\ar[ddddd]|{\ov p}\\
&&&&&&&&&&&&&\\
&&&&&&&&&&&&&\\
&&&&&&&&&&&&&\\
&&&&&
Y\ar[dddddrrrrr]|-{g} \ar[drrr]|{\id} \ar[rrrrr]|{v} &&&&& Y' \ar[drrr]|{\beta'}\ar[ddddd]|{q}
&&&\\
& &&&&&&&  Y\ar[dddddrrrrr]|-{\ov g} \ar[rrrrr]|{\ov v} &&&&& \ov Y'\ar[ddddd]|{\ov q}\\
&&&&&&&&&&&&&\\
&&&&&&&&&&&&&\\
&&&&&&&&&&&&&\\
&&&&&&&&&&Z \ar[drrr]|{\id}&&&\\
&&&&&&&&&&&&& Z
}
\]
satisfying all the hypotheses. Out of this we cook up the following
2-commutative diagrams, which we will refer to as $B$ and $C$
respectively
  \[
  \xymatrix{
    X\ar[dr]_-\id \ar[r]^-{\id}  &X \ar[r]^{u'} \ar[d]^{\id}\ar@{}[rd]|-{(\clubsuit)} &  X'\ar[d]^{\alpha'}\ar[r]^u \ar@{}[rd]|-{(\diamondsuit)} & X''\ar[d]^-{\alpha''} & & X\ar[dr]_-f \ar[r]^-{u'}  &X' \ar[r]^{\id} \ar[d]^{p'} &  X'\ar[d]^{p'}\ar[r]^u \ar@{}[rd]|-{(\dagger)} & X''\ar[d]^-{p} \\
    &    X\ar[dr]_-{\ov f} \ar[r]^-{\ov u'}  &\ov X' \ar[r]^{\ov u}\ar@{}[rd]|-{(\heartsuit)} \ar[d]^{\ov p'} &  \ov X''\ar[d]^{\ov p}
    & \ar@{}[d]|-{\ds\text{and}}& &    Y\ar[dr]_-{\id} \ar[r]^-{\id}  &Y \ar[r]^{v} \ar[d]^{\id} \ar@{}[rd]|-{(\spadesuit)}&  Y'\ar[d]^{\beta'}\\
    & &   Y\ar[dr]_-{\ov g} \ar[r]^-{\ov v} & \ov Y'\ar[d]^{\ov q} & &
    & &   Y\ar[dr]_-{\ov g} \ar[r]^-{\ov v} & \ov Y'\ar[d]^{\ov q}\\
& & &  Z
& & & & &  Z
}
\]
Of the vertical maps $p'$, $p$,  $\ov p'$, $\ov p$ and
$\ov q$ are given to
be separated, of finite type and universally quasi-proper, and
Remark~\ref{R6.2057.5} guarantees that so are 
$\alpha'$, $\alpha''$ and $\beta'$.
Thus all the vertical morphisms are of finite type and universally quasi-proper
and the horizontal maps are dominant, flat monomorphisms.
Corollary~\ref{CNag.22} tells us that the squares are 2-cartesian,
hence the diagram $B$ is an object of $\nseq(X,X,Y,Z)$ and the diagram
$C$ is an object
in $\nseq(X,Y,Y,Z)$.

By Lemma~\ref{L200.5} the objects $B\in\nseq(X,X,Y,Z)$ and
$C\in\nseq(X,Y,Y,Z)$ give one relation each. We have the identities
\begin{eqnarray*}
P\big(\pi_{123}^{}(B)\big)P\big(\pi_{134}^{}(B)\big)&=&
P\big(\pi_{234}^{}(B)\big)P\big(\pi_{124}^{}(B)\big)\,,\\
P\big(\pi_{123}^{}(C)\big)P\big(\pi_{134}^{}(C)\big)&=&
P\big(\pi_{234}^{}(C)\big)P\big(\pi_{124}^{}(C)\big)\,.
\end{eqnarray*}
Now recall that $B$ and $C$ were constructed out of
the representative $\ph$ of a morphism $a\la b$ in $\nseq(X,Y,Z)$, and
I assert that the Lemma will be proved by untangling what
each of the entries in the two relations is, in terms of $a$, $b$ and $\ph$. Let
us begin with the relation that comes from $C$.

We have that $\id=P\big(\pi_{123}^{}(C)\big)$, hence the relation given
by $C$ simplifies
to $P\big(\pi_{134}^{}(C)\big)=
P\big(\pi_{234}^{}(C)\big)P\big(\pi_{124}^{}(C)\big)$.
Let $\lambda:q\la\ov q\beta'$ be the isomorphism, and consider the diagram
\[\xymatrix@C+30pt{
  (uu')^*(qp)^\times \ar[r]^-{\lambda^\times}\ar[d]_-{\rho(a)} & (uu')^*(\ov q\beta'p)^\times\ar[d]^-{P\big(\pi_{124}^{}(C)\big)} & \\
  ({u'}^*{p'}^\times)(v^*q^\times)\ar[r]^-{\lambda^\times} & ({u'}^*{p'}^\times)(v^*[\ov q\beta']^\times)\ar[r]^-{P\big(\pi_{234}^{}(C)\big)} &  ({u'}^*{p'}^\times)({v'}^*\ov q^\times)
}\]
The square obviously commutes, and the composite along the bottom row
is by definition the map induced by $\pi_{23}^{}(\ph)^!:v^*q^\times\la
{v'}^*\ov q^\times$. The commutativity, combined the relation given by
$C$, tells us that
$\pi_{23}^{}(\ph)^!\rho(a)=P\big(\pi_{134}^{}(C)\big)\lambda^\times$.

Now we study the relation coming from $B$. If $\mu:qp\la\ov q\ov p\alpha''$
is the isomorphism in the data, then by
definition $\pi_{13}^{}(\ph)^!
=P\big(\pi_{124}^{}(B)\big)\mu^\times$,
while $\rho(b)$ is by definition equal to $P\big(\pi_{234}^{}(B)\big)$.
The relation coming from $B$ rewrites
as $\rho(b)\pi_{13}^{}(\ph)^!=
P\big(\pi_{123}^{}(B)\big)P\big(\pi_{134}^{}(B)\big)\mu^\times$.
Let $\nu:p'\la\ov p'\alpha'$ be the isomorphism,
and we have identities
\begin{eqnarray*}
  \rho(b)\pi_{13}^{}(\ph)^! & = & P\big(\pi_{123}^{}(B)\big)P\big(\pi_{134}^{}(B)\big)\mu^\times \\
  & = &P\big(\pi_{123}^{}(B)\nu^\times P\big(\pi_{134}^{}(C)\big)\lambda^\times\\
  &=& \pi_{12}^{}(\ph)^!\pi_{23}^{}(\ph)^!\rho(a)
\end{eqnarray*}
The second equality is because
$P\big(\pi_{134}^{}(B)\big)\mu^\times=\nu^\times P\big(\pi_{134}^{}(C)\big)\lambda^\times$,
which comes from the isomorphism of $\pi_{134}^{}(B)$ with
$\pi_{134}^{}(C)$. The third equality is because
$\pi_{23}^{}(\ph)^!\rho(a)=
P\big(\pi_{134}^{}(C)\big)\lambda^\times$
by the paragraph above studying the identity
coming from $C$, and because $\pi_{12}^{}(\ph)^!=\pi_{123}^{}(B)\nu^\times$
by definition.
But the equality of the first and last terms
means that the square
\[\xymatrix@C-5pt{
 \pi_{13}^{}(a)^!\ar@{=}[r] & (uu')^*(qp)^\times \ar[rrr]^-{\pi_{13}^{}(\ph)^!}\ar[d]_-{\rho(a)}
 & & & (\ov u\ov u')^*(\ov q\ov p)^\times  \ar[d]^-{\rho(b)}
 & \pi_{13}^{}(b)^!\ar@{=}[l]\\
 \pi_{12}^{}(a)^!\pi_{23}^{}(a)^!\ar@{=}[r] &
 ({u'}^*{p'}^\times)(v^*q^\times) \ar[rrr]^{\pi_{12}^{}(\ph)^!\pi_{23}^{}(\ph)^!}&&& ({\ov u'}^*{\ov p'}^\times)(\ov v^*\ov q^\times)
 & \pi_{12}^{}(b)^!\pi_{23}^{}(b)^!\ar@{=}[l]
}\]
commutes, and hence $\rho$ is a natural transformation.
\eprf

\rmk{R200.19}
The functors $\pi_{12}^{}(-)^!$, $\pi_{23}^{}(-)^!$ and $\pi_{13}^{}(-)^!$
all factor through the map $\pi:\nseq(X,Y,Z)\la \gnseq(X,Y,Z)$,
and hence so does the natural transformation 
$\rho: \pi_{13}^{}(-)^!\la
\pi_{12}^{}(-)^!\pi_{23}^{}(-)^!$.
But  Proposition~\ref{P6.2061} tells us that the
map $F':\gnseq(X,Y,Z)\la \seq(X,Y)\times\seq(Y,Z)$
is an equivalence.
In Remark~\ref{R200.13} we chose, for every pair of objects $X,Z\in\seq$,
a functor $R(X,Z):\seq(X,Z)\la\gnseq(X,Z)$ right inverse to the projection,
and defined the functor $(-)^!$ on $\seq(X,Z)$ using this
right inverse. We leave it to the reader to check that
we can now define an unambiguous natural transformation
which takes the
object $X\stackrel f\la Y\stackrel g\la Z$ in
$\seq(X,Y)\times\seq(Y,Z)$ to the morphism $\rho(f,g):(gf)^!\la f^!g^!$
in $\Tri\big(\Dqc(Z),\Dqc(X)\big)$. The map $\rho(f,g)$ does not
depend on any new choices: it is determined by the $R(X,Y)$,
$R(Y,Z)$ and $R(X,Z)$ which were fixed in Remark~\ref{R200.13}.
\ermk

\thm{T200.19}
We have so far constructed the following map from the 2-category
$\seq$ to the 2-category $\Tri$:
\be
\item
  On objects: we take $X\in\seq$ to $\Dqc(X)\in\Tri$.
\item
  On 1-morphisms and 2-morphisms: For each pair of objects
  $X,Y\in\seq$ we have defined a functor
  $(-)^!:\seq(X,Y)\la\Tri\big(\Dqc(Y),\Dqc(X)\big)$.
\item
  Composition: for any
triple of objects $X,Y,Z\in\seq$ we have defined a natural
transformation $\rho:[(-)\circ(?)]^!\la(?)^!\circ(-)^!$.
\ee
These combine
to make $(-)^!$ an oplax 2-functor of 2-categories.
\ethm

\prf
What needs proof is
the associativity property of $\rho$. We need to show that, for any object
$W\stackrel e\la X\stackrel f\la Y\stackrel g\la Z$ in the
category $\seq(W,X)\times\seq(X,Y)\times\seq(Y,Z)$, the following
square commutes
\[\xymatrix@C+20pt{
(gfe)^! \ar[r]^-{\rho(e,gf)}\ar[d]_{\rho(fe,g)}& e^!(gf)^! \ar[d]^-{e^!\rho(f,g)}\\
 (fe)^!g^! \ar[r]^-{\rho(e,f)g^!} & e^!f^!g^!
}\]
The square involves objects obtained by applying the functor $(-)^!$ to
\[\begin{array}{ccccccccccc}
e&\in&\seq(W,X)&\quad & f&\in&\seq(X,Y) &\quad& g&\in&\seq(Y,Z)\\
fe&\in&\seq(W,Y)& & gf&\in&\seq(X,Z)&& gfe& \in& \seq(W,Z)
\end{array}\]
and the definition of $(-)^!$ suggests that it might be
convenient to first choose preimages for these under
the equivalences of categories $F':\gnseq(-,?)\la\seq(-,?)$.
And we might as well choose preimages in a way that facilitates
the computation.
Lemma~\ref{L6.2063} tells us that every object
$a\in\seq(W,X)\times\seq(X,Y)\times\seq(Y,Z)$
is in
the essential image of $\nseq(W,X,Y,Z)$. If we let
$a$ be the object $W\stackrel e\la X\stackrel f\la Y\stackrel g\la Z$
we may
extend it to an object $A\in\nseq(W,X,Y,Z)$, that is a 2-commutative
diagram
 \[
  \xymatrix{
    W\ar[dr]_-e \ar[r]^-{t''}  &W' \ar[r]^{t'} \ar[d]^{o''} &  W''\ar[d]^{o'}\ar[r]^t & W'''\ar[d]^o\\
&    X\ar[dr]_-f \ar[r]^-{u'}  &X' \ar[r]^{u} \ar[d]^{p'} &  X''\ar[d]^p\\
& &   Y\ar[dr]_-g \ar[r]^-v & Y'\ar[d]^q\\
& & &  Z
}
  \]
Now choose 
\[\begin{array}{ccccccccccc}
\pi_{12}^{}(A)&\in&\gnseq(W,X)&\quad & \pi_{23}^{}(A)&\in&\gnseq(X,Y) &\quad& \pi_{34}^{}(A)&\in&\gnseq(Y,Z)\\
\pi_{13}^{}(A)&\in&\gnseq(W,Y)& & \pi_{24}^{}(A)&\in&\gnseq(X,Z)&& \pi_{14}^{}(A)& \in& \gnseq(W,Z)
\end{array}\]
to be (respectively) the liftings via
$F':\gnseq(-,?)\la\seq(-,?)$ of the objects
\[
e,\ \ f,\ \ g,\ \ fe,\ \ gf\ \text{ and }\ gfe
\]
and then 
\[\begin{array}{ccccccc}
\rho(e,f)&=&P\big(\pi_{123}^{}(A)\big)&\qquad&\rho(fe,g)&=&P\big(\pi_{134}^{}(A)\big)\\
\rho(f,g)&=&P\big(\pi_{234}^{}(A)\big) & & \rho(e,gf)&=&P\big(\pi_{124}^{}(A)\big)
\end{array}\]
and the relation 
$P\big(\pi_{123}^{}(A)\big)P\big(\pi_{134}^{}(A)\big)=
P\big(\pi_{234}^{}(A)\big)P\big(\pi_{124}^{}(A)\big)$
of Lemma~\ref{L200.5} is
exactly the identity $\rho(e,f)\rho(fe,g)=\rho(f,g)\rho(e,gf)$.
\eprf

The oplax 2-functor $(-)^!$ is not a pseudofunctor,
meaning there are in general
pairs of  morphisms for which $\rho(f,g):(gf)^!\la f^!g^!$ need
not be an isomorphism. We end the section with

\pro{P200.23}
Let $X\stackrel f\la Y\stackrel g\la Z$ be composable 1-morphisms
in $\seq$. The map $\rho(f,g):(gf)^!\la f^!g^!$ is an isomorphism if
any one of the following holds
\be
\item
$f$ is of finite Tor-dimension. 
\item
  $g$ is of finite type and universally quasi-proper.
\item
  $gf$ is of finite type and universally quasi-proper.
  \item
We restrict to the
bounded-below derived category. More formally: if
$I:\Dqcpl(Z)\la\Dqc(Z)$ is the inclusion, then
for any composable 1-morphisms
$X\stackrel f\la Y\stackrel g\la Z$ in $\seq$ the map
$\rho(f,g)I:(gf)^!I\la f^!g^!I$ is an isomorphism.
\ee  
\epro

\prf
We begin with (ii): assume $g$ is of finite type and universally
quasi-proper. Choose a Nagata compactification
for $f$, that is write $f$ as the composite
$X\stackrel {u'}\la X'\stackrel{p'}\la Y$, with $u'$
a dominant, flat monomorphism and $p'$ of finite type and
universally quasi-proper. Then the diagram $C$
below
\[
  \xymatrix{
    X\ar[dr]_-f \ar[r]^-{u'}  &X' \ar[r]^{\id}\ar[d]^{p'}\ar@{}[dr]|-{(\diamondsuit)} &  X'\ar[d]^{p'}\\
 &   Y\ar[dr]_-g \ar[r]^-\id & Y\ar[d]^g\\
 & &  Z
}
  \]
is 2-commutative, the horizontal maps are all dominant, flat
monomorphisms and the vertical maps are of finite type and universally
quasi-proper. Hence $C$ is an object of $\nseq(X,Y,Z)$.
But up to isomorphism $\rho(f,g)=P(C)$ is given by the base-change map
$\Phi(\diamondsuit)=\id$.

Now for the proof of (i) and (iv). Choose a lifting via the equivalence
\[F':\gnseq(X,Y,Z)\la\seq(X,Y)\times\seq(Y,Z)\] 
of the object $X\stackrel f\la Y\stackrel g\la Z$,
more concretely choose a diagram $A\in\nseq(X,Y,Z)$ of
the form
\[
  \xymatrix@C+20pt{
    X\ar[dr]_-f \ar[r]^-{u'}  &X' \ar[r]^{u}\ar@{}[dr]|-{(\heartsuit)} \ar[d]^{p'} &  X''\ar[d]^p\\
 &   Y\ar[dr]_-g \ar[r]^-v & Y'\ar[d]^q\\
 & &  Z
}
  \]
In (i)
the map
$f\cong p'u'$ is assumed of finite Tor-dimension while $u'$ is 
given to be flat, hence the image of $u'$ is
contained in the open set on which the finite-type
map $p'$ is of finite Tor-dimension. Theorem~\ref{T4.13}
informs us that ${u'}^*\Phi(\heartsuit):{u'}^*u^*p^\times\la {u'}^*{p'}^\times v^*$
is an isomorphism. Composing on the right with $q^\times$
we deduce that
${u'}^*\Phi(\heartsuit)q^\times:{u'}^*u^*p^\times q^\times\la {u'}^*{p'}^\times v^*q^\times$
is also an isomorphism, but up to isomorphism this is
the morphism $\rho(f,g)=P(A):(gf)^!\la f^! g^!$. This proves (i).

Now
for (iv):
$f$ and $g$ are unrestricted,
but let $E\in\Dqcpl(Z)$ be a bounded-below object.
Then $q^\times E$ is an object of
$\Dqcpl(Y')$, and by Lemma~\ref{L4.7}(ii) the map
$\Phi(\heartsuit)(q^\times E):u^*p^\times q^\times E\la {p'}^\times v^*q^\times E$
is an isomorphism. Applying ${u'}^*$ we deduce that
$\rho(f,g):(gf)^!\la f^!g^!$ evaluates at $E$ to give an isomorphism.

Finally we prove (iii): assume that $gf$ is of finite type
and universally quasi-proper. Choose a Nagata compactification for
$g$, that is write $g$ as the composite $Y\stackrel v\la Y'\stackrel q\la Z$
with $q$ of finite type and universally quasi-proper and $v$ a dominant, flat
monomorphism.
We know that $q$ and $qvf\cong gf$ are of
finite type
and universally quasi-proper and Lemma~\ref{LNag.1} tells us that
so is $vf$. 
By Theorem~\ref{T200.19} we have
the identity $\rho(v,q)\rho(f,qv)=\rho(f,v)
\rho(vf,q)$. We wish to show that $\rho(f,g)\cong\rho(f,qv)$ is an
isomorphism, and by (ii) we know that $\rho(v,q)$ and $\rho(vf,q)$
both are---from the identity it suffices to prove that
$\rho(f,v)$ is an isomorphism. Thus we are reduced to proving (iii)
in the case where $g$ is a dominant, flat monomorphism.

Let us factor $gf$ as $gf\cong ip'$, where $i$ is the (closed) immersion of the
stack-theoretic closure of the image of $gf$.
By Lemma~\ref{LNag.-1} the image
of the finite-type, universally quasi-proper map $gf:X\la Z$ is closed;
as a set the closure of the image is just the image, 
and $p'$ is surjective. We have a 2-commutative square
$(\clubsuit)$ and form the 2-cartesian square $(\diamondsuit)$ below
\[\xymatrix{
  X  \ar[r]^{p'}\ar[d]_f \ar@{}[dr]|-{(\clubsuit)}& Q\ar[d]^i & &
  P  \ar@{}[dr]|-{(\diamondsuit)} \ar[r]^\beta\ar[d]_j & Q\ar[d]^i\\
 Y\ar[r]_g &     Z & &  Y\ar[r]_g &     Z
}\]
The surjective map $p':X\la Q$ factors (up to isomorphism) through
the map $\beta$ which, being the pullback of the flat monomorphism
$g$, must be a surjective, flat monomorphism. Lemma~\ref{LNag.4} tells
us that $\beta$ is an isomorphism. Therefore $f$ is isomorphic to
$j\beta^{-1}p'$. Put $p=\beta^{-1}p'$ and we have $f\cong jp$ where $j$,
being the pullback of the closed immersion $i$, is a closed immersion.

Now Theorem~\ref{T200.19} gives the identity
$\rho(p,j)\rho(jp,g)=\rho(j,g)\rho(p,gj)$. Because $j$ and $gj\cong i$
are closed immersions (ii) tells us that $\rho(p,gj)$ and $\rho(p,j)$
are isomorphisms. We want to show that $\rho(f,g)\cong\rho(jp,g)$ is
an isomorphism, and the identity tells us it suffices to prove
that $\rho(j,g)$ is an isomorphism. That is: we are reduced to proving
that $\rho(f,g)$ is an isomorphism when $g$ is a dominant, flat monomorphism
and $f$ and $gf$ are both closed immersions.
The 2-commutative diagram $C$ below
\[
  \xymatrix@C+10pt{
    X\ar[dr]_-f \ar[r]^-{\id}  &X\ar@{}[rd]|-{(\heartsuit)} \ar[r]^{\id}\ar[d]^{f} &  X\ar[d]^{gf}\\
 &   Y\ar[dr]_-g \ar[r]^-g & Z\ar[d]^\id\\
 & &  Z
}
  \]
has the property that the horizontal maps are all dominant, flat monomorphisms
while the vertical maps are of finite type and universally quasi-proper,
hence $C$ is an object of $\nseq(X,Y,Z)$
lifting $X\stackrel f\la Y\stackrel g\la Z$,
and $P(C)$ computes $\rho(f,g)$. We are reduced to showing
that $P(C)$ is an isomorphism, meaning that the base-change
map $\Phi(\heartsuit):(gf)^\times\la f^\times g^*$ is an isomorphism.
Since $f$ is an affine morphism $f_*$ is conservative, and
it suffices to show
that $f_*\Phi(\heartsuit):f_*(gf)^\times\la f_*f^\times g^*$ is an isomorphism.
The fact that $g$ is a concentrated, flat monomorphism means that the natural
map $g^*g_*\la\id$ is an isomorphism, hence it suffices to
show that the composite
\[
\CD
g^*g_* f_*(gf)^\times @>>> f_*(gf)^\times @>f_*\Phi(\heartsuit)>>
f_*f^\times g^*
\endCD
\]
is an isomorphism. But for the closed immersions $gf$ and $f$
the functors $g_*f_*(gf)^\times$ and $f_*f^\times g^*$ simplify (respectively)
to $\HHom_Z^{}(g_*f_*\co_X^{},-)$ and $\HHom_Y^{}\big(f_*\co_X^{},g^*(-)\big)
\cong\HHom_Y^{}\big(g^*g_*f_*\co_X^{},g^*(-)\big)$.
We are reduced to showing that the natural map
\[
\CD
g^*\HHom_Z^{}(g_*f_*\co_X^{},-)
@>\ph>> \HHom_Y^{}\big(g^*g_*f_*\co_X^{},g^*(-)\big)
\endCD
\]
is an isomorphism. Since $g^*g_*\cong\id$ it certainly suffices to
show that the composite
\begin{equation}
  \label{eqn55}
\xymatrix{
g_*g^*\HHom_Z^{}(g_*f_*\co_X^{},-)
\ar[r]^-{g_*\ph} & g_*\HHom_Y^{}\big(g^*g_*f_*\co_X^{},g^*(-)\big) \ar[d]^-{\wr}\\
 & \HHom_Z^{}\big(g_*f_*\co_X^{},g_*g^*(-)\big)
}
\end{equation}
is an isomorphism.

Define $\cm\subset\Dqc(Z)$ to be the full subcategory
of all objects $M$ for which the natural map $M\la g_*g^*M$ is an
isomorphism. Clearly $\cm$ is 
a
localizing subcategory of $\Dqc(Z)$.
In the 2-cartesian square $(\heartsuit)$ above the vertical
maps are concentrated and the horizontal maps are flat, hence
the base-change map $g^*(gf)_*\la f_*$ is an isomorphism. Therefore
$g_*g^*(gf)_*\cong g_*f_*$ and we deduce that every object in
$(gf)_*\Dqc(X)$ belongs to the category $\cm$. Since $\cm$ is
localizing it follows that the localizing subcategory generated
by $(gf)_*\Dqc(X)$ is contained in $\cm$, that is $\cm$ contains
all of $\D_{\mathbf{qc},X}^{}(Z)$, the subcategory of $\Dqc(Z)$
of objects supported on the closed subset $X$. The
object $g_*f_*\co_X^{}$ is bounded below and pseudo-coherent---the pseudo-coherence is
because $gf$ is a
closed immersion---and moreover $g_*f_*\co_X^{}$ belongs
to $\D_{\mathbf{qc},X}^{}(Z)$. Lemma~\ref{L4.-4}
tells us that
$\HHom_Z^{}(g_*f_*\co_X^{},-)$ belongs to $\D_{\mathbf{qc},X}^{}(Z)$
for every $(-)$. Because
$\D_{\mathbf{qc},X}^{}(Z)\subset\cm$ the horizontal map below is an isomorphism
\begin{equation}
  \label{eqn56}
  \xymatrix{
    \HHom_Z^{}(g_*f_*\co_X^{},-) \ar[r]^-\sim &
g_*g^*\HHom_Z^{}(g_*f_*\co_X^{},-)
\ar[d]^-{(\ref{eqn55})}\\
 & \HHom_Z^{}\big(g_*f_*\co_X^{},g_*g^*(-)\big)
}
\end{equation}
and, to show that (\ref{eqn55}) is an isomorphism, it suffices to prove
the composite (\ref{eqn56}) an isomorphism.

We are assuming that $\Dqc(Z)$ satisfies Thomason's condition, hence
the category $\D_{\mathbf{qc},X}^{}(Z)$ is generated by the compact objects it
contains, and $\cm$ contains all of them. For every compact
object $C\in\D_{\mathbf{qc},X}^{}(Z)$ and any object $E\in\Dqc(Z)$ we have
that the map $C\oo E\la (g_*g^*C)\oo E\cong C\oo (g_*g^*E)$ is an
isomorphism, where the
first map is an isomorphism because $C\in\cm$.
Because $C$ is strongly dualizable we
can rewrite this as saying that $\HHom(C^\vee,E)\la\HHom(C^\vee,g_*g^*E)$
is an isomorphism. Let $\cl\subset\Dqc(Z)$ be the full subcategory of all
objects $L\in\Dqc(Z)$ such that the natural map
$\HHom(L,-)\la\HHom\big(L,g_*g^*(-)\big)$ is an isomorphism. We have
proved that every compact object in $\D_{\mathbf{qc},X}^{}(Z)$ belongs
to $\cl$, but $\cl$ is localizing and the compact objects
in $\D_{\mathbf{qc},X}^{}(Z)$ generate $\D_{\mathbf{qc},X}^{}(Z)$. Hence
$\D_{\mathbf{qc},X}^{}(Z)$ is contained in $\cl$, in particular
$g_*f_*\co_X^{}\in\cl$ and the composite in (\ref{eqn56}) is
an isomorphism.
\eprf

\section{The oplax natural transformation $\psi:(-)^\times\la(-)^!$}
\label{S202}

We want to construct something on $\seq(X,Z)$ but
(as usual) begin with $\nseq(X,Z)$.

\con{C202.1}
Let $a\in\nseq(X,Z)$ be an object, that is a 2-commutative diagram
  \[
  \xymatrix@C+30pt@R-20pt{
    & Y\ar[dr]^-{p} & \\
    X\ar[rr]_-{f}\ar[ur]^-{u}  & & Z 
  }\]
  with $u$ a dominant, flat monomorphism (and $p$ of finite type and
  universally quasi-proper).  
We define $\psi(a):u^\times p^\times\la u^*p^\times$ to be the morphism
$\Phi(\heartsuit)p^\times$, where $\heartsuit$ is the 2-cartesian square
\[
\xymatrix{
  X\ar[r]^\id\ar[d]_\id\ar@{}[dr]|-{(\heartsuit)} & X\ar[d]^u \\
  X\ar[r]^u & Y
}\]
Note that the square is 2-cartesian because $u$ is a monomorphism, and
the base-change map $\Phi(\heartsuit):\id^*u^\times\la\id^\times u^*$
exists because the square is 2-cartesian and $u$ is flat. 
\econ

\lem{L202.3}
On the category $\seq(X,Z)\cong\gnseq(X,Z)$ the
$\psi(-)$ of Construction~\ref{C202.1} gives a natural transformation of
$\psi:(-)^\times\la(-)^!$.
\elem

\prf
We recall that the map $F:\nseq(X,Z)\la\seq(X,Z)$ is
a groupoid completion, hence has
a universal property with respect to functors and natural transformations:
any natural transformation between functors that factor through
$F$ must factor through $F$.
It therefore suffices to show that Construction~\ref{C202.1}
yields a natural transformation between the composite functors
$\xymatrix{\nseq(X,Z)\ar[r]^-{F} &\seq(X,Z)\ar@<0.5ex>[r]^-{(-)^\times}
  \ar@<-0.5ex>[r]_-{(-)^!} &\Tri\big(\Dqc(Z),\Dqc(X)\big).}$
The composite $(-)^!\circ F$
is, by the construction of $(-)^!$,
naturally isomorphic
to the functor taking the object $a\in\nseq(X,Z)$ to
$P(a)=u^*p^\times:\Dqc(Z)\la\Dqc(X)$.
By
Remark~\ref{R200.1} the composite $(-)^\times\circ F$ is naturally
isomorphic to the functor taking the object $a\in\nseq(X,Z)$
to $u^\times p^\times\cong(pu)^\times$ and the reader might find it
convenient to think of this as $P(a')$ where $a'\in\lnseq(X,Z)$ is the
2-commutative diagram
  \[
  \xymatrix@C+30pt@R-20pt{
    & Y\ar[dr]^-{pu} & \\
    X\ar[rr]_-{f}\ar[ur]^-{\id}  & & Z 
  }\]
A representative of a morphism $\ph:a\la b$ in $\nseq(X,Z)$ gives, among other data, a 2-commutative
diagram
  \[\xymatrix@C+30pt@R-20pt{
       & Y\ar^p[dr]\ar_\alpha[dd] & \\
X\ar^u[ur]\ar_{u'}[dr] & &  Z \\
       & Y'\ar_{p'}[ur] & 
}\]
Out of this we cook up the diagram $C$ below
   \[
  \xymatrix@C+10pt{
    X\ar[dr]_-\id \ar[r]^-{\id}  &X \ar[r]^{\id} \ar[d]^{\id} &
    X\ar@{}[dr]|-{(\diamondsuit)}\ar[d]^{\id}\ar[r]^\id & X\ar[d]^{u}\\
    &    X\ar[dr]_-\id \ar[r]^-{\id}  &X
    \ar@{}[dr]|-{(\heartsuit)}\ar[r]^{u} \ar[d]^{\id} &  Y\ar[d]^{\alpha}\\
& &   X\ar[dr]_-{f} \ar[r]^-{u'} & Y'\ar[d]^{p'}\\
& & &  Z
}
  \]
In the diagram $C$
the horizontal maps are all given to be dominant, flat monomorphisms,
and the map $\alpha$ is of finite type and universally quasi-proper
by Remark~\ref{R6.2052.3}.
The square $(\diamondsuit)$ is 2-cartesian because $u$ is a monomorphism, and
the square $(\heartsuit)$ is 2-cartesian by Corollary~\ref{CNag.22}. The
remaining square is trivially 2-cartesian, and the diagram
belongs to $\lnseq(X,X,X,Z)$.

Lemma~\ref{L200.5} gives the relation
$P\big(\pi_{123}^{}(C)\big)P\big(\pi_{134}^{}(C)\big)=
P\big(\pi_{234}^{}(C)\big)P\big(\pi_{124}^{}(C)\big)$,
and since $P\big(\pi_{123}^{}(C) =\id$ this simplifies to
$P\big(\pi_{134}^{}(C)\big)=
P\big(\pi_{234}^{}(C)\big)P\big(\pi_{124}^{}(C)\big)$.
Let $\mu:p\la p'\alpha$ and $\lambda:\alpha u\la u'$ be the given
isomorphisms determined by 
$\ph:a\la b$. Then
\[
P\big(\pi_{124}^{}(C)\big)=\Phi(\diamondsuit)(p'\alpha)^\times \quad
\text{while}\quad \psi(a) =\Phi(\diamondsuit)p^\times,
\]
hence $P\big(\pi_{124}^{}(C)\big)\mu^\times=\mu^\times\psi(a)$.
Therefore
\begin{eqnarray*}
  \ph^!\psi(a) &=& P\big(\pi_{234}^{}(C)\big)\mu^\times\psi(a) \\
  & =& P\big(\pi_{234}^{}(C)P\big(\pi_{124}^{}(C)\big)\mu^\times\\
  &=& P\big(\pi_{134}^{}(C)\big)\mu^\times
\end{eqnarray*}
where the first equality is by the definition of
$\ph^!:a^!\la b^!$, and the second and third equalities
are from the discussion above.
Now in $\lnseq(X,X,Z)$ we have the isomorphic
objects  $\pi_{134}^{}(C)$  and $A$ below
   
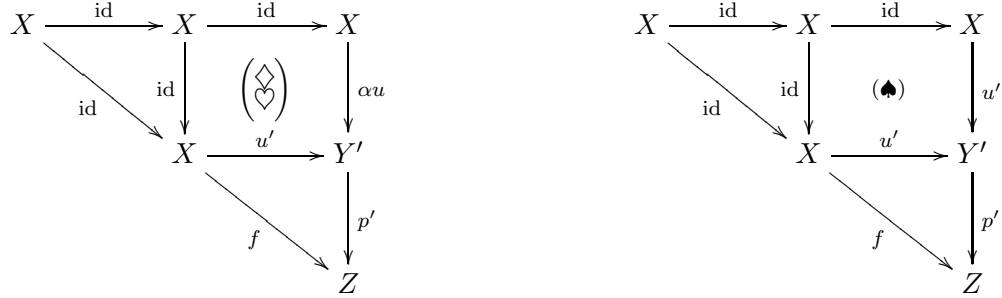
\[
  \xymatrix@C+20pt@R+10pt{
  X \ar[r]^{\id} \ar[dr]_-{\id} &
    X\ar@{}[dr]|-{\left(\!\!\!\!\begin{array}{c}\diamondsuit\\*[-5pt]
      \heartsuit\end{array}\!\!\!\!\right)}\ar[d]_-{\id}\ar[r]^\id & X\ar[d]^{\alpha u}
    & & X \ar[r]^{\id} \ar[dr]_-{\id} &
    X\ar@{}[dr]|-{(\spadesuit)}\ar[d]_-{\id}\ar[r]^\id & X\ar[d]^{u'}\\
    &   X\ar[dr]_-{f} \ar[r]^-{u'} & Y'\ar[d]^{p'} & &     &   X\ar[dr]_-{f} \ar[r]^-{u'} & Y'\ar[d]^{p'}\\
& &  Z
& & & &  Z
}
  \]
where the isomorphism is induced by $\lambda:\alpha u\la u'$. By
Construction~\ref{C202.1} we have
$\psi(b)=P(A)$, and by the isomorphism of $A$ with $\pi_{134}^{}(C)$
we have
$P(A)\lambda^\times=P\big(\pi_{134}^{}(C)\big)$.
Hence
\begin{eqnarray*}
  \ph^!\psi(a)  &=& P\big(\pi_{134}^{}(C)\big)\mu^\times \\
  &=& \psi(b)\lambda^\times\mu^\times
\end{eqnarray*}
in other words the square
\[
\xymatrix@C+30pt{
  a^\times \ar[r]^-{\ph^\times=\lambda^\times\mu^\times} \ar[d]_-{\psi(a)} & b^\times\ar[d]^-{\psi(b)}\\
  a^! \ar[r]^-{\ph^!}  & b^!
}\]
commutes, that is $\psi$ is a natural transformation.
\eprf

\thm{T202.5}
Define $\psi:(-)^\times\la(-)^!$ as follows:
\be
\item
  For an object $X\in\seq$ we define $\psi(X):\Dqc(X)\la\Dqc(X)$ to
  be the identity functor.
\item
  For a morphism $f:X\la Z$ in $\seq$ we define $\psi(f):f^\times\la f^!$
  to be as in Lemma~\ref{L202.3}.
  \ee
These data define an oplax natural transformation of 2-functors.
\ethm

\rmk{R202.7}
We should perhaps expand a little, reminding the reader of the
definitions. We have a pair of oplax 2-functors $(-)^\times$
and $(-)^!$, where it so happens that $(-)^\times$ is actually
a pseudofunctor. To say that the data given define
a natural transformation asserts
\be
\item
  For every pair of objects $X,Z\in\seq$ the data in Theorem~\ref{T202.5}(ii)
  yield a natural transformation of functors $\seq(X,Z)\la\Tri\big(\Dqc(Z),\Dqc(X)\big)$; this much we have already proved, see Lemma~\ref{L202.3}.
\item
  The construction respects identities. This is easy.
\item
  The construction is compatible with composition in $\seq$. Concretely:
for composable 1-morphisms $X\stackrel f\la Y\stackrel g\la Z$ let $\delta$
be the canonical isomorphism $\delta(f,g):(gf)^\times\la f^\times g^\times$, let
$\rho(f,g):(gf)^!\la f^!g^!$ be the natural map of Theorem~\ref{T200.19} and
let $\psi:(-)^\times\la(-)^!$ be the natural transformation of
Lemma~\ref{L202.3}. These are compatible in the sense that the square
below commutes
\[\xymatrix@C+15pt{
{(gf)}^\times \ar[r]^{\delta(f,g)}\ar[d]_{\psi(gf)} &f^\times g^\times \ar[d]^{\psi(f)\psi(g)}\\
{(gf)}^! \ar[r]_{\rho(f,g)} & f^!g^!
}\]
This last part still needs to be proved.
\ee
\ermk

\prf
Our object $X\stackrel f\la Y\stackrel g\la Z$ may
be lifted to $\nseq(X,Y,Z)$; choose
such a lifting, that is a 2-commutative diagram
\[
  \xymatrix{
X\ar[dr]_-f \ar[r]^-{u'}  &X' \ar[r]^{u} \ar[d]^{p'} &  X''\ar[d]^p\\
 &   Y\ar[dr]_-g \ar[r]^-v & Y'\ar[d]^q\\
 & &  Z
}
  \]
satisfying 
a bunch of conditions. Out of this we cook up the
2-commutative
diagram $B$ below
\[
  \xymatrix@C+12pt{
    X\ar[dr]_-\id \ar[r]^-{\id}  &X'\ar@{}[dr]|-{(\clubsuit)} \ar[r]^{\id}
    \ar[d]^{\id} &  X\ar@{}[dr]|-{(\diamondsuit)}\ar[d]^{u'}
    \ar[r]^\id & X\ar[d]^{uu'} \\
    &  X\ar[dr]_-f \ar[r]^-{u'}  &X' \ar@{}[dr]|-{(\heartsuit)}
    \ar[r]^{u} \ar[d]^{p'} &  X''\ar[d]^p \\
    & &   Y\ar[dr]_-g \ar[r]^-v & Y'\ar[d]^q    \\
& & &  Z 
}
\]
The horizontal maps are all given to be dominant, flat monomorphisms. The
square $(\heartsuit)$ is given to be 2-cartesian, the square $(\clubsuit)$
is 2-cartesian by the definition of $u'$ being a monomorphism, and
the square $(\diamondsuit)$ is 2-cartesian by Lemma~\ref{LNag.5} applied
to the composable 1-morphisms $X\stackrel{u'}\la X'\stackrel u\la X''$.
Hence $B$ is an object in $\lnseq(X,X,Y,Z)$ and Lemma~\ref{L200.5} gives
the relation
$P\big(\pi_{123}^{}(B)\big)P\big(\pi_{134}^{}(B)\big)=
P\big(\pi_{234}^{}(B)\big)P\big(\pi_{124}^{}(B)\big)$.
We have $\psi(f)=P\big(\pi_{123}^{}(B)\big)$,
$\psi(gf)=P\big(\pi_{124}^{}(B)\big)$ and
$\rho(f,g)=P\big(\pi_{234}^{}(B)\big)$
and our relation becomes
\begin{equation}
  \label{eq1}
  \psi(f)P\big(\pi_{134}^{}(B)\big)=
  \rho(f,g)\psi(gf).
  \end{equation}
It remains to identify $P\big(\pi_{134}^{}(B)\big)$, and for this it is
handy to consider the two diagrams
\[
  \xymatrix@C+12pt{
  X\ar@{}[dr]|-{(\diamondsuit)}\ar[d]_{u'}
  \ar[r]^\id & X\ar[d]^{uu'}
  & &   X\ar[d]_{p'u'}
    \ar[r]^\id & X\ar[d]^{p'u'}\\
X' \ar@{}[dr]|-{(\heartsuit)}
\ar[r]^{u} \ar[d]_{p'} &  X''\ar[d]^p
 & & Y \ar@{}[dr]|-{(\spadesuit)}
    \ar[r]^{\id} \ar[d]_{\id} &  Y\ar[d]^v \\
 Y \ar[r]^-v & Y' & &  Y \ar[r]^-v & Y'
}
  \]
The diagram on the left
is an extract from the object
$B\in\lnseq(X,X,Y,Z)$, hence the squares $(\diamondsuit)$
and $(\heartsuit)$ are 2-cartesian. The square
$(\spadesuit)$ is 2-cartesian because $v$ is a monomorphism, while
the remaining square is obviously 2-cartesian. And the concatenations
are isomorphic by the isomorphism $\lambda:pu\la vp'$. Hence
$\Phi\left(\!\!\!\!\begin{array}{c}\diamondsuit\\*[-5pt]
\heartsuit\end{array}\!\!\!\!\right)=
(p'u')^\times\Phi(\spadesuit)\lambda^\times$.
This shows that $ \psi(f)P\big(\pi_{134}^{}(B)\big)= \psi(f)\Phi\left(\!\!\!\!\begin{array}{c}\diamondsuit\\*[-5pt]
\heartsuit\end{array}\!\!\!\!\right)q^\times$
is equal to the composite
\[
\xymatrix@C+2pt{
  {u'}^\times u^\times p^\times q^\times \ar[r]^-{\lambda^\times} &
  {u'}^\times {p'}^\times v^\times q^\times \ar[rrrr]^-{(p'u')^\times\Phi(\spadesuit)q^\times  =(p'u')^\times\psi(g)} & & & &
  {u'}^\times {p'}^\times v^* q^\times\ar[r]^-{ \psi(f)} &
  {u'}^* {p'}^\times v^* q^\times
  }
\]
and (\ref{eq1}) yields the commutativity of the square
\[\xymatrix@C+5pt{
  {(gf)}^\times \ar@{=}[r] &{u'}^\times u^\times p^\times q^\times \ar[rrr]^{\delta(f,g)=\lambda^\times}\ar[d]_{\psi(gf)}& & &
 ({u'}^\times {p'}^\times)( v^\times q^\times) \ar[d]^{\psi(f)\psi(g)}\ar@{=}[r] &
  f^\times g^\times \\
{(gf)}^!\ar@{=}[r]  & {u'}^* u^* p^\times q^\times  \ar[rrr]_{\rho(f,g)} && &({u'}^* {p'}^\times) (v^* q^\times) \ar@{=}[r] &f^!g^!
}\]
\eprf

\rmk{R202.9}
In general the natural transformation $\psi$ is not an isomorphism. 
But in passing we note that, if $f:X\la Z$ is of
finite type and universally quasi-proper, then $\psi(f):f^\times\la f^!$
is an isomorphism. The point is that, by Lemma~\ref{L202.3}, we may choose
any preimage of $f$ under the functor $F:\nseq(X,Z)\la\seq(X,Z)$ for
the purpose of computing $\psi(f)$, and the preimage we choose is
  \[
  \xymatrix@C+30pt@R-20pt{
    & X\ar[dr]^-{f} & \\
    X\ar[rr]_-{f}\ar[ur]^-{\id}  & & Z 
  }\]
The recipe of Construction~\ref{C202.1} is to form the 2-cartesian square
\[
\xymatrix{
  X\ar[r]^\id\ar[d]_\id\ar@{}[dr]|-{(\heartsuit)} & X\ar[d]^\id \\
  X\ar[r]^\id & X
}\]
and $\psi(f)=\Phi(\heartsuit)f^\times$ is obviously an isomorphism.
\ermk

\section{Base change}
\label{S8}

We begin by defining the building blocks
for the constructions of this section.

\dfn{D204.1}
The category $\sq$ is defined as follows:
\be
\item
The objects are the 2-cartesian squares in $\seq$
\[\xymatrix{
  W \ar[r]^u\ar[d]_f & X\ar[d]^g\\
  Y\ar[r]^v & Z
  }\]
with $v$ flat.
\item
  The morphisms are the 2-isomorphisms of such squares. That is
  a morphism from the square
 \[\xymatrix{
  W \ar[r]^u\ar[d]_f & X\ar[d]^g&&\ar@{}[d]|-{\ds\text{to the square}} & &W \ar[r]^{u'}\ar[d]_{f'} & X\ar[d]^{g'} \\
  Y\ar[r]^v & Z && & & Y\ar[r]^{v'} & Z
 }\]
 is the data of morphisms $f\la f'$, $g\la g'$, $u\la u'$ and $v\la v'$
 compatible with the isomorphisms $vf\la gu$ and $v'f'\la g'u'$.
\item
  We will also wish to consider two subcategories: the objects
  are the same, but the morphisms in $\sq_H^{}\subset\sq$ are
  the maps where $f\la f'$ and $g\la g'$ are assumed to be identities,
  and the morphisms in
  $\sq_V^{}\subset\sq$ are the ones where
  $u\la u'$ and $v\la v'$ are assumed to be identities.
\ee
\edfn
  
\con{C204.3}
The category $\sq_H^{}$ is a disjoint union of subcategories $\sq_H^{}(f,g)$,
where the vertical maps are fixed to be some $f$ and $g$.
This leads us to the definition of a
2-category $\hseq$:
\be
\item
  The objects of $\hseq$ are the 1-morphisms $f:W\la Y$ in $\seq$.
\item
  Given two objects $f,g\in\hseq$,
  the category $\hseq(f,g)$ is the category $\sq_H^{}(f,g)$ above.
\item
  The composition map $\hseq(f,g)\times\hseq(g,h)\la\hseq(f,h)$
  is the concatenation of 2-cartesian squares.
\setcounter{enumiii}{\value{enumi}}
\ee
There is also an obvious vertical version, the 2-category
$\vseq$ is defined as follows
\be
\setcounter{enumi}{\value{enumiii}}
\item
  The objects of $\vseq$ are the flat 1-morphisms $u:W\la X$ in $\seq$.
\item
  Given objects
  $u,v\in\vseq$, the category
  $\vseq(u,v)$ is the full subcategory $\sq_V^{}$ containing
  the objects (i.e.~the 2-cartesian squares) in which the horizontal
  maps are $u$ and $v$.
\item
  The composition map $\vseq(u,v)\times\vseq(v,w)\la\vseq(u,w)$
  is the concatenation of 2-cartesian squares.
\setcounter{enumiii}{\value{enumi}}
\ee
There are two obvious functors $p_1^{},p_2^{}:\hseq\la\seq$ and two functors
$p_1^{},p_2^{}:\vseq\la\seq$, where $p_1^{}$ takes an object $A\la B$,
of either $\hseq$ or $\vseq$, to the object $A\in\seq$, while
$p_2^{}$ takes $A\la B$ to $B$.
\econ

\rmk{R204.5}
For any object in $\sq$, that is for any 2-cartesian square 
\[\xymatrix{
\ar@{}[d]|-{\ds(\dagger)}&  W \ar[r]^u\ar[d]_f & X\ar[d]^g&\\
 & Y\ar[r]^v & Z &
}\]
with $v$ flat,
there is a base-change map $\Phi:u^*g^\times\la f^\times v^*$.
Furthermore this
base-change map is compatible with the concatenation of squares. In
this Remark we note how this could be rephrased in terms of
2-functors on the categories $\hseq$ and $\vseq$.

Consider first the category $\vseq$. There are the two 2-functors
$p_1^{},p_2^{}:\vseq\la\seq$, and we may compose them with
the 2-functor $(-)^\times:\seq\la\Tri$ to
obtain two 2-functors $\vseq\la\Tri$. Concretely
the functor $(-)^\times\circ p_1^{}$ takes the morphism
$(\dagger)\in\vseq$ to $f^\times:\Dqc(Y)\la\Dqc(W)$,
while $(-)^\times\circ p_2^{}$ takes $(\dagger)$ 
to $g^\times:\Dqc(Z)\la\Dqc(X)$. Now we will define
a natural transformation $\Phi:(-)^\times\circ p_2^{}\la
(-)^\times\circ p_1^{}$.
The definition is
\be
\item
  On objects: the natural transformation $\Phi$ takes
  an object $Y\stackrel v\la Z$ of $\vseq$ to the
  functor $v^*:\Dqc(Z)\la\Dqc(Y)$.
\item
  On 1-morphisms: for the
  1-morphism $(\dagger)\in\vseq$ we need to give a 2-morphism
  in $\Tri$, and we need
  to choose a direction---our natural
  transformation will be oplax.
  This requires us to give 2-morphism $u^*g^\times\la f^\times v^*$. The
  map we choose is the base-change map $\Phi(\dagger)$.
\ee
The assertion that $\Phi$ is an oplax natural transformation says
that the formula in (ii) is functorial and respects composition.

On the category $\hseq$ the construction is similar. We consider the
two composite functors $(-)^*\circ p_1^{}$ and $(-)^*\circ p_2^{}$;
concretely $(-)^*\circ p_1^{}$ takes the 1-morphism $(\dagger)\in\hseq$
to $u^*:\Dqc(X)\la\Dqc(W)$ while $(-)^*\circ p_2^{}$ takes
$(\dagger)$ to $v^*:\Dqc(Z)\la\Dqc(Y)$. 
The natural transformation $\Phi$ takes an object $X\stackrel g\la Z$
of $\hseq$ to $g^\times:\Dqc(Z)\la\Dqc(X)$. And this time the
natural transformation is lax, it takes a 1-morphism $(\dagger)$
to the base-change map $\Phi(\dagger):u^*g^\times\la f^\times v^*$.
\ermk

This section is about mimicking all of this for $(-)^!$
in place of $(-)^\times$. As usual
we begin by constructing a category equivalent to $\sq$.

\con{C204.7}
For objects $W,X,Y,Z$ in $\seq$ we give a category $\nsq(W,X,Y,Z)$ as
follows:
\be
\item
  The objects are 2-commutative diagrams
 
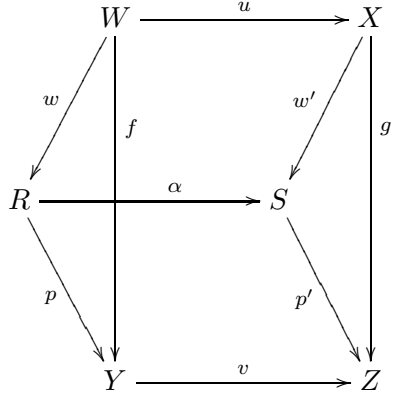
\[\xymatrix@C-5pt@R+30pt{
    & W\ar[dl]_-{w}\ar[rrr]^-{u}\ar[dd]^(0.3){f} & & & X\ar[dl]_-{w'}
    \ar[dd]^(0.3){g} \\
 R\ar[dr]_-{p} \ar[rrr]^(.6){\alpha}& & &S\ar[dr]_-{p'} & \\
& Y\ar[rrr]^-{v}  &&  & Z \\
 }\]
where the squares
  
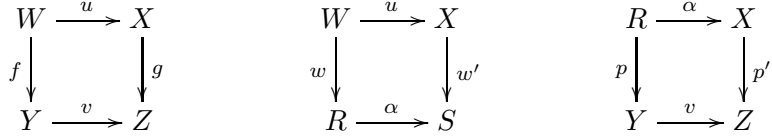
\[\xymatrix{
  W \ar[r]^u\ar[d]_f & X\ar[d]^g& &  W \ar[r]^u\ar[d]_{w} & X\ar[d]^{w'} & &  R \ar[r]^\alpha\ar[d]_{p} & X\ar[d]^{p'} \\
  Y\ar[r]^v & Z & &R\ar[r]^\alpha & S & &Y\ar[r]^v & Z 
  }\]
  are all 2-cartesian, where all the horizontal maps are flat, the maps $p$ and $p'$ are of finite type and universally quasi-proper and the maps $w$ and $w'$
  are dominant, flat monomorphisms.
\item
  A morphism in $\nsq(W,X,Y,Z)$, from object 
 
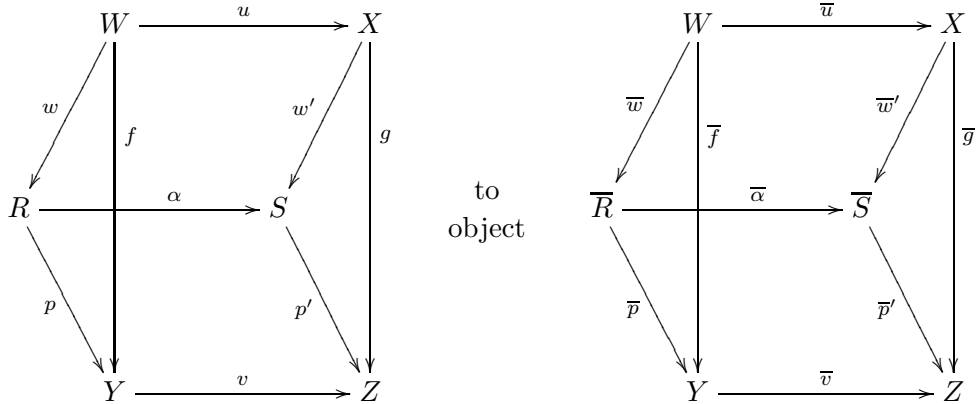
\[\xymatrix@C-5pt@R+30pt{
    & W\ar[dl]_-{w}\ar[rrr]^-{u}\ar[dd]^(0.3){f} & & & X\ar[dl]_-{w'} \ar[dd]^(0.3){g} & &  & & W\ar[dl]_-{\ov w}\ar[rrr]^-{\ov u}\ar[dd]^(0.3){\ov f} & & & X\ar[dl]_-{\ov w'}
    \ar[dd]^(0.3){\ov g} \\
 R\ar[dr]_-{p} \ar[rrr]^(.6){\alpha}& & &S\ar[dr]_-{p'} & 
 & \ar@{}[r]|-{\ds\begin{array}{c}\text{to}\\
 \text{object}\end{array}}& & \ov R\ar[dr]_-{\ov p} \ar[rrr]^(.6){\ov \alpha}& & &\ov S\ar[dr]_-{\ov p'} & \\
& Y\ar[rrr]^-{v}  &&  & Z 
& & && Y\ar[rrr]^-{\ov v}  &&  & Z \\
 }\]
 is an equivalence class of data. The data
 are a pair of 1-morphisms $R\la\ov R$ and $S\la \ov S$, together with a
 bunch of 2-morphisms giving all the 2-commutativity one might expect.
 Two representatives are declared equivalent if they differ from each other
 by perturbing the maps $R\la \ov R$ and $S\la\ov S$ by an isomorphism.
\ee
An object of $\nsq(W,X,Y,Z)$ can be thought of as an object in $\sq$ where
we give compatible Nagata compactifications for the vertical maps. We
also fix the
objects $W,X,Y,Z\in\seq$
for the construction. There is a natural forgetful functor
from $\nsq(W,X,Y,Z)$ to a category we will call $\sq(W,X,Y,Z)\subset\sq$,
the full subcategory of $\sq$ 
whose objects are the 2-cartesian squares
\[\xymatrix{
  W \ar[r]^u\ar[d]_f & X\ar[d]^g\\
  Y\ar[r]^v & Z
}\]
with $W$, $X$, $Y$ and $Z$ fixed.
\econ

\lem{L204.9}
The forgetful map $F:\nsq(W,X,Y,Z)\la\sq(W,X,Y,Z)$
is a groupoid completion: that is $\sq(W,X,Y,Z)$ can be obtained from
$\nsq(W,X,Y,Z)$ by formally inverting all the morphisms.
\elem

\prf
First we factor the functor $F$ as the composite
\[
\CD
\nsq(W,X,Y,Z) @>\rho>>\rsq(W,X,Y,Z) @>\sigma>> \sq(W,X,Y,Z)
\endCD
\]
where $\rho:\nsq(W,X,Y,Z)\la\rsq(W,X,Y,Z)$ is the functor
that remembers only the Nagata compactification on the right.
Thus an object of $\rsq(W,X,Y,Z)$
is a 2-commutative diagram
\[\xymatrix@C-5pt@R+10pt{
    & W\ar[rrr]^-{u}\ar[dd]_-{f} & & & X\ar[dr]^-{w'}
    \ar[dd]_-{g} \\
 &&& & &S\ar[dl]^-{p'}  \\
& Y\ar[rrr]^-{v}  &&  & Z \\
 }\]
where the square is 2-cartesian, the horizontal maps are flat, $w'$ is a
dominant, flat monomorphism and $p'$ is of finite type
and universally quasi-proper. The functor $\rho$ is an equivalence of
categories, a quasi-inverse is given by forming the 2-cartesian square
\[\xymatrix{
  R \ar[r]^\alpha\ar[d]_p & S\ar[d]^{p'}\\
  Y\ar[r]^v & Z
}\]
and letting $w:X\la R$ be the map to the pullback. We leave the details
to the reader.

Since $\rho$ is an equivalence it suffices to prove that $\sigma$ is a
groupoid completion. This is by Proposition~\ref{P6.2055.5}, since
the category $\nsq(W,X,Y,Z)$ is a disjoint union of subcategories
equivalent to
$\seq(W,X)\times\seq(W,Y)\times\seq(Y,Z)\times\nsq(X,Z)$.
\eprf

\con{C204.11}
We consider two composites
\[\xymatrix{
     &\seq(W,X)\times\seq(X,Z)\ar[rd]^-{(-)^*(-)^!} &\\
\sq(W,X,Y,Z)\ar[ru]^-{\pi_1^{}} \ar[rd]_-{\pi_2^{}} & &
\Tri\big(\Dqc(Z),\Dqc(W)\big) \\
 & \seq(W,Y)\times\seq(Y,Z)\ar[ur]_-{(-)^!(-)^*} &
}\]
where $\pi_1^{},\pi_2^{}$ are the obvious projections. Let 
$\fkr,\fks$ be the composites, that is $\fkr,\fks$ are,
respectively, the functors taking the object 
\[\xymatrix{
\ar@{}[d]|-{\ds(\dagger)}&  W \ar[r]^u\ar[d]_f & X\ar[d]^g&\\
 & Y\ar[r]^v & Z &
}\]
of the category $\sq(W,X,Y,Z)$ 
to the objects $\fkr(\dagger)=u^*g^!$ and $\fks(\dagger)=f^!v^*$
in $\Tri\big(\Dqc(Z),\Dqc(W)\big)$. We wish to define a
natural transformation $\theta:\fkr\la\fks$.

Because $F:\nsq(W,X,Y,Z)\la\sq(W,X,Y,Z)$ is
a groupoid completion functor
it suffices to produce a natural transformation $\fkr F\la\fks F$,
it will automatically factor as $\theta F$. Suppose therefore
that we are given an object in $A\in\nsq(W,X,Y,Z)$; the data are
the 2-cartesian squares
\[\xymatrix{
  W \ar[r]^u\ar[d]_f & X\ar[d]^g& &  W \ar[r]^u\ar[d]_{w}\ar@{}[dr]|-{(\diamondsuit)} & X\ar[d]^{w'} & &  R \ar[r]^\alpha\ar[d]_{p}\ar@{}[dr]|-{(\heartsuit)}  & X\ar[d]^{p'} \\
  Y\ar[r]^v & Z & &R\ar[r]^\alpha & S & &Y\ar[r]^v & Z 
  }\]
and 2-isomorphisms $pw\la f$, $p'w'\la g$ rendering 2-commutative the
entire picture. We define $\Theta(A)$ to be the composite
\[
\xymatrix@C+30pt{
u^*{w'}^* {p'}^\times \ar[r]_{(\diamondsuit)^*{p'}^\times } \ar@/^2pc/[rr]^{\Theta(A)} & {w}^*\alpha^*{p'}^\times
\ar[r]_{{w}^*\Phi(\heartsuit)} & {w}^*{p}^\times v^*
}
\]
where $(\diamondsuit)^*:u^*{w'}^* \la {w}^*\alpha^*$ is the isomorphism
induced by the 2-commutative square $(\diamondsuit)$, and
$\Phi(\heartsuit):\alpha^*{p'}^\times
\la {p}^\times v^*$ is the base-change map of the square $(\heartsuit)$.
\econ

\pro{P204.13}
With the notation as in Construction~\ref{C204.11}, the map sending an
object $A\in\nsq(W,X,Y,Z)$ to $\Theta(A)$ defines a natural transformation
$\Theta:\fkr F\la\fks F$ and, as mentioned in
Construction~\ref{C204.11}, the natural transformation $\Theta$ factors
uniquely as $\Theta=\theta F$ with 
$\theta:\fkr\la\fks$ also a natural transformation.
\epro

\prf
What needs proof is the naturality of $\Theta$. A morphism $\ph:A\la B$
in the category $\nsq(W,X,Y,Z)$ has a representative
giving (among other things) the data
of a 2-commutative diagram
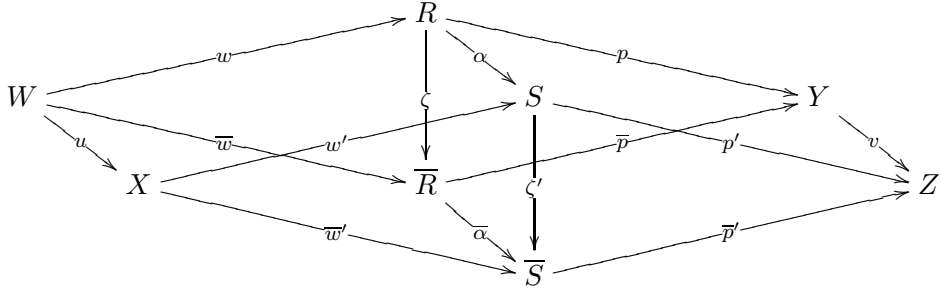
\[\xymatrix@C+2pt@R-25pt{
   &  & && R\ar|{p}[ddddrrrr]\ar|{\zeta}[dddddddd]\ar[ddddr]|{\alpha}  &&& &&&\\
&&& && && &&&\\
&&& && && &&&\\
&&& && && &&&\\
W\ar|{w}[uuuurrrr]\ar|{\ov w}[ddddrrrr]\ar[ddddr]|{u}  &&& & &S\ar|{p'}[ddddrrrr]\ar|{\zeta'}[dddddddd]&& &  Y \ar[ddddr]|{v}  &&\\
&&& && && &&&\\
&&& && && &&&\\
&&& && && &&&\\
     & X\ar|{w'}[uuuurrrr]\ar|{\ov w'}[ddddrrrr] & &  &
    \ov R\ar|{\ov p}[uuuurrrr]\ar[ddddr]|{\ov\alpha} & & &&&Z& \\
&&& && && &&&\\
&&& && && &&&\\
&&& && && &&&\\
&&   & & &  \ov S\ar|{\ov p'}[uuuurrrr] & & & &
}\]
and 2-morphisms $u\la\ov u$ and $v\la\ov v$. The naturality
of $\Theta$ with respect to the 2-morphisms $u\la\ov u$ and $v\la\ov v$ is
easy and left to the reader, in the rest of the proof we will assume
these maps to be identities.
Out of the data we cook up the diagram $C$ below
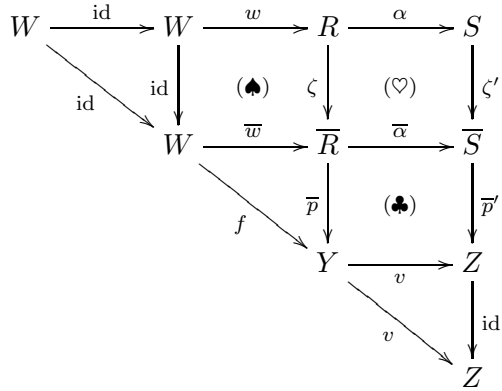
\[
\xymatrix@C+15pt@R+5pt{
W\ar[r]^-\id\ar[dr]_-{\id} & W \ar[d]_\id \ar[r]^{w}\ar@{}[dr]|{(\spadesuit)} & R \ar@{}[dr]|{(\heartsuit)} 
 \ar[d]_{\zeta}\ar[r]^{\alpha} & S \ar[d]^{\zeta'}
 \\
& W\ar[dr]_-{f}  \ar[r]^{\ov w} & \ov R\ar[r]^{\ov\alpha}\ar[d]_{\ov p}\ar@{}[dr]|{(\clubsuit)} &  \ov S\ar[d]^{\ov p'}  \\
& &  Y\ar[r]_{v}\ar[dr]_-{v} & Z\ar[d]^{\id}\\
& & & Z
}
\]
The horizontal maps are all given to be flat.
The map $\zeta$ is of finite type and universally quasi-proper by
Remark~\ref{R6.2052.3}, and Corollary~\ref{CNag.22} applies to
the square $(\spadesuit)$
and tells us that it must be 2-cartesian. The square $(\clubsuit)$ is
given to be 2-cartesian, while the  square
\[
\xymatrix@C+15pt@R+5pt{
  R\ar[r]^\alpha \ar@{}[dr]|-{\left(\!\!\!\!\begin{array}{c}\heartsuit\\*[-5pt]
\clubsuit\end{array}\!\!\!\!\right)}\ar[d]_{\ov p\zeta} & S\ar[d]^-{\ov p'\zeta'} & & &&R\ar[r]^\alpha \ar@{}[dr]|-{(\diamondsuit)}\ar[d]_{p} & S\ar[d]^-{p'} \\
  Y \ar[r]_{v} & Z & &\ar@{}[u]|-{\ds\text{is isomorphic to}} & & Y \ar[r]_{v} & Z 
}  
\]
via the isomorphisms $\lambda:p\la \ov p\zeta$ and $\mu:p'\la\ov p'\zeta'$.
Now $(\diamondsuit)$ is given to be 2-cartesian, hence so
is the isomorphic $\left(\!\!\!\!\begin{array}{c}\heartsuit\\*[-5pt]
\clubsuit\end{array}\!\!\!\!\right)$,
and as $(\clubsuit)$ is also 2-cartesian so is $(\heartsuit)$. Therefore
$C$ is an object of $\lnseq(W,W,Y,Z)$ and Lemma~\ref{L200.5} gives
the relation
$P\big(\pi_{123}^{}(C)\big)P\big(\pi_{134}^{}(C)\big)=
P\big(\pi_{234}^{}(C)\big)P\big(\pi_{124}^{}(C)\big)$.
Note that in the above we also learned that $\Phi\left(\!\!\!\!\begin{array}{c}\heartsuit\\*[-5pt]
\clubsuit\end{array}\!\!\!\!\right)\mu^\times=\lambda^\times\Phi(\diamondsuit)$,
and hence
$P\big(\pi_{134}^{}(C)\big)\mu^\times=w^*\lambda^\times\Phi(\diamondsuit)$.
Therefore
\begin{eqnarray*}
  P\big(\pi_{234}^{}(C)\big)P\big(\pi_{124}^{}(C)\big)\mu^\times& = &
  P\big(\pi_{123}^{}(C)\big)P\big(\pi_{134}^{}(C)\big)\mu^\times \\
  & = & P\big(\pi_{123}^{}(C) w^*\lambda^\times\Phi(\diamondsuit)\\
  & = & \fks(\ph)w^*\Phi(\diamondsuit)
\end{eqnarray*}
where the first equality is the relation coming from $C$,
the second equality is by the discussion above, and the third is
by the definition of $\fks(\ph)$.

Next contemplate the diagrams
\[
\xymatrix@C+15pt@R+2pt{
 W \ar[d]_\id \ar[r]^{w}\ar@{}[dr]|{(\spadesuit)} & R \ar@{}[dr]|{(\heartsuit)} 
 \ar[d]_{\zeta}\ar[r]^{\alpha} & S \ar[d]^{\zeta'}
& & W \ar[d]_\id \ar[r]^{u} & X \ar@{}[dr]|{(\dagger)} 
 \ar[d]_{\id}\ar[r]^{w'} & S \ar[d]^{\zeta'}
 \\
W  \ar[r]_{\ov w} & \ov R\ar[r]_{\ov\alpha} &  \ov S
& & W  \ar[r]_{u} & X\ar[r]_{\ov w'} &  \ov S
}
\]
The diagram on the left is extracted from $C$, hence the squares are
2-cartesian. The morphisms $w'$ and $\ov w'$ are given to be dominant,
flat monomorphisms, while $\zeta'$ is of finite type and universally
quasi-proper by Remark~\ref{R6.2052.3}. Corollary~\ref{CNag.22}
tells us the square $(\dagger)$ is 
2-cartesian. The remaining square is trivially 2-cartesian,
and the concatenations are
isomorphic via the isomorphisms $\gamma:w'u\la \alpha w$ and
$\delta:\ov w'u\la\ov\alpha\ov w$. Therefore
$\Phi(\spadesuit\heartsuit)\gamma^*=\delta^*u^*\Phi(\dagger)$.
Observing that
$P\big(\pi_{124}^{}(C)\big)=\Phi(\spadesuit\heartsuit){\ov p'}^\times$
we deduce that 
\begin{eqnarray*}
  \fks(\ph)\Theta(A)& =& \fks(\ph)w^*\Phi(\diamondsuit)\gamma^*\\
  & =&
  P\big(\pi_{234}^{}(C)\big)P\big(\pi_{124}^{}(C)\big)\mu^\times\gamma^*\\
     & =&
  P\big(\pi_{234}^{}(C)\big)\delta^*u^*\Phi(\dagger){\ov p'}^\times\mu^\times\\
  &=& \Theta(B)\fkr(\ph)
\end{eqnarray*}
The first equality is because $\Theta(A)=w^*\Phi(\diamondsuit)\gamma^*$
by definition, the second equality comes from
the previous paragraph, the third equality expands
$\Phi(\spadesuit\heartsuit)\gamma^*$ as right before the string
of equalities, and the last equality is because
$ \Theta(B)= P\big(\pi_{234}^{}(C)\big)\delta^*$ and
$\fkr(\ph)=\Phi(\dagger){\ov p'}^\times\mu^\times$.
Thus the square
\[\xymatrix@C+20pt{
\fkr(A) \ar[r]^-{\fkr(\ph)}\ar[d]_-{\Theta(A)} & \fkr(B)\ar[d]^{\Theta(B)} \\
\fks(A) \ar[r]^-{\fks(\ph)} & \fks(B) \\
}\]  
commutes, and $\Theta$ is indeed a natural transformation.
\eprf

\thm{T204.15}
Let $\hseq$ be the 2-category of Construction~\ref{C204.3}~(i), (ii) and (iii).
Consider the composite 2-functors
\[\xymatrix@C+30pt{
  \hseq\ar@<0.5ex>[r]^-{p_1^{}} \ar@<-0.5ex>[r]_-{p_2^{}} &
  \seq \ar[r]^-{(-)^*} & \Tri
}\]
There is a lax natural transformation
$\theta:(-)^*\circ p_2^{}\la
(-)^*\circ p_1^{}$
given by the rules
\be
\item
  Given an object $f:W\la Y$ of $\hseq$, then $\theta(f)=f^!:\Dqc(Y)\la\Dqc(W)$.
\item
  Given a morphism in $\hseq$, that is a 2-cartesian square with
  flat horizontal morphisms
\[\xymatrix@C+15pt@R+5pt{
  W \ar[r]^u\ar[d]_f \ar@{}[dr]|-{(\diamondsuit)}& X\ar[d]^g\\
  Y\ar[r]^v & Z
}\]
we need to define a 2-morphism $\theta(\diamondsuit):u^*g^!\la f^!v^*$.
The definition is as in Proposition~\ref{P204.13}.
\ee
\ethm

\prf
For any pair of 1-morphisms $f,g\in\seq$ we have two functors
$\hseq(f,g)\la\Tri\big(\Dqc(Z),\Dqc(X)\big)$, namely
$\big[p_1^{}(-)\big]^*\circ g^!$ and $f^!\circ\big[p_2^{}(-)\big]^*$.
In Proposition~\ref{P204.13} we established that
$\theta$ defines a natural transformation
$\theta:\big[p_1^{}(-)\big]^*\circ g^!\la f^!\circ\big[p_2^{}(-)\big]^*$.
Identities are obviously respected,
hence what remains to be proved is that $\theta$ is compatible with
composition. That is: we have to show that, if
\[\xymatrix@C+15pt@R+5pt{
  \wt W \ar[r]^{\wt u}\ar[d]_{\wt f} \ar@{}[dr]|-{(\clubsuit)}& W \ar[r]^u\ar[d]_f \ar@{}[dr]|-{(\diamondsuit)}& X\ar[d]^g\\
  \wt Y\ar[r]^{\wt v} &Y\ar[r]^v & Z
}\]
are two composable morphisms in $\hseq$, then the diagram 
\[
\xymatrix@C-22pt@R+5pt{ 
& {(u \wt u)}^*g^! \ar[rr]^{\theta(\clubsuit\diamondsuit)}\ar[dl]_{\sim} & & 
          {\wt f}^!{( v\wt v)}^*\ar[dr]^{\sim} & \\
{\wt u}^*{u}^*g^!\ar[rr]_{{\wt u}^*\theta(\diamondsuit)} & &{\wt u}^*f^!v^*\ar[rr]_{\theta(\clubsuit)v^*} & &
 {\wt f}^!{\wt v}^*v^* 
}\]
commutes. To prove this we choose compatible preimages of
$(\clubsuit)\in\sq(\wt W,W,\wt Y,Y)$,
$(\diamondsuit)\in\sq(W,X,Y,Z)$
and $(\clubsuit\diamondsuit)\in\sq(\wt W,W,X,Z)$ via the map
$F:\nsq(-,-,-,-)\la
\sq(-,-,-,-)$.
More precisely: let $X\stackrel{w'}\la R'\stackrel{p'}\la Z$
be a Nagata compactification of $g:X\la Z$. From the pullback diagram
\[\xymatrix{
  \wt R \ar[r]^-{\wt\alpha}\ar[d]_-{\wt p} &R \ar[r]^-{\alpha}\ar[d]^-{p} &
  R' \ar[d]^-{p'}\\
\wt Y \ar[r]^-{\wt v} & Y\ar[r]^-{v} & Z
}\]
and then let $w:W\la R$ and $\wt w:\wt W\la \wt R$ be the maps to the pullback.
We obtain a 2-commutative diagram
\[\xymatrix@C+20pt@R+5pt{
  \wt W \ar@{}[dr]|-{(\dagger)}\ar[r]^-{\wt u}\ar[d]_-{\wt w} &W\ar@{}[dr]|-{(\bullet)} \ar[r]^-{u}\ar[d]^-{w} &
  X \ar[d]^-{w'}\\
\wt R\ar@{}[dr]|-{(\heartsuit)} \ar[r]^-{\wt\alpha}\ar[d]_-{\wt p} &R\ar@{}[dr]|-{(\spadesuit)} \ar[r]^-{\alpha}\ar[d]^-{p} &
  R' \ar[d]^-{p'}\\
\wt Y \ar[r]^-{\wt v} & Y\ar[r]^-{v} & Z
}\]
By construction $(\heartsuit)$ and $(\spadesuit)$ are 2-cartesian, and
since $(\clubsuit)\cong\left(\!\!\!\!\begin{array}{c}\dagger\\*[-5pt]
\heartsuit\end{array}\!\!\!\!\right)$
and
$(\diamondsuit)\cong\left(\!\!\!\!\begin{array}{c}\bullet\\*[-5pt]
\spadesuit\end{array}\!\!\!\!\right)$
are given to be 2-cartesian it follows that $(\dagger)$ and $(\bullet)$
must also be 2-cartesian. The 2-commutative diagrams
\[\xymatrix@C+20pt@R+5pt{
  \wt W \ar@{}[dr]|-{(\dagger)}\ar[r]^-{\wt u}\ar[d]_-{\wt w} &W\ar[d]^-{w} & 
  W\ar@{}[dr]|-{(\bullet)} \ar[r]^-{u}\ar[d]_-{w} &
  X \ar[d]^-{w'} & 
  \wt W \ar@{}[dr]|-{(\dagger\bullet)}\ar[r]^-{u\wt u}\ar[d]_-{\wt w} &
  X \ar[d]^-{w'}\\
\wt R\ar@{}[dr]|-{(\heartsuit)} \ar[r]^-{\wt\alpha}\ar[d]_-{\wt p} &R\ar[d]^-{p} & 
R\ar@{}[dr]|-{(\spadesuit)} \ar[r]^-{\alpha}\ar[d]_-{p} &
  R' \ar[d]^-{p'} & 
\wt R\ar@{}[dr]|-{(\heartsuit\spadesuit)} \ar[r]^-{\alpha\wt\alpha}\ar[d]_-{\wt p}  &
  R' \ar[d]^-{p'}\\
\wt Y \ar[r]^-{\wt v} & Y &  
 Y\ar[r]^-{v} & Z & 
\wt Y \ar[r]^-{v\wt v}  & Z
}\]
together with the 2-morphisms $\wt p\wt w\la\wt f$, $pw\la f$ and $p'w'\la g$,
give a lifting of the objects
\[\xymatrix@C+15pt@R+5pt{
  \wt W \ar[r]^{\wt u}\ar[d]_{\wt f} \ar@{}[dr]|-{(\clubsuit)}& W\ar[d]^f & 
 W \ar[r]^u\ar[d]_f \ar@{}[dr]|-{(\diamondsuit)}& X\ar[d]^g & 
  \wt W \ar[r]^{u\wt u}\ar[d]_{\wt f} \ar@{}[dr]|-{(\clubsuit\diamondsuit)} & X\ar[d]^g\\
  \wt Y\ar[r]^{\wt v} &Y & 
  Y\ar[r]^v & Z &
  \wt Y\ar[r]^{v\wt v} & Z
}\]
via the maps $\nsq(-,-,-,-)\la\sq(-,-,-,-)$.
By the definition of $\theta$ it suffices
to prove the commutativity of
\[
\xymatrix@C-15pt@R+10pt{ 
  & {(u \wt u)}^*({w'}^*{p'}^\times) \ar[rr]^{\Theta\sst
\left(\!\!\!\!\begin{array}{cc}\dagger\!\!&\!\!\bullet\\*[-5pt]
\heartsuit\!\!&\!\!\spadesuit\end{array}\!\!\!\!\right)
}\ar[dl]_{\sim} & & 
          (\wt w^*\wt p^\times){( v\wt v)}^*\ar[dr]^{\sim} & \\
{\wt u}^*{u}^*({w'}^*{p'}^\times) \ar[rr]_{{\wt u}^*\Theta\sst
\left(\!\!\!\!\begin{array}{c}\bullet\\*[-5pt]
\spadesuit\end{array}\!\!\!\!\right)} & &{\wt u}^*({w}^*{p}^\times) v^*\ar[rr]_{\Theta\sst
\left(\!\!\!\!\begin{array}{c}\dagger\\*[-5pt]
\heartsuit\end{array}\!\!\!\!\right)v^*} & &
 (\wt w^*\wt p^\times){\wt v}^*v^* 
}\]
which is immediate from the definitions.
\eprf

\thm{T204.17}
Let $\vseq$ be the 2-category of Construction~\ref{C204.3}~(iv), (v) and (vi).
Consider the composite 2-functors
\[\xymatrix@C+30pt{
  \vseq\ar@<0.5ex>[r]^-{p_1^{}} \ar@<-0.5ex>[r]_-{p_2^{}} &
  \seq \ar[r]^-{(-)^!} & \Tri
}\]
There is an oplax natural transformation
$\theta:(-)^!\circ p_2^{}\la
(-)^!\circ p_1^{}$
given by the rules
\be
\item
  Given an object $u:W\la X$ of $\hseq$, then $\theta(u)=u^*:\Dqc(X)\la\Dqc(W)$.
\item
  Given a morphism in $\vseq$, that is a 2-cartesian square with
  flat horizontal morphisms
\[\xymatrix@C+15pt@R+5pt{
  W \ar[r]^u\ar[d]_f \ar@{}[dr]|-{(\diamondsuit)}& X\ar[d]^g\\
  Y\ar[r]^v & Z
}\]
we need to define a 2-morphism $\theta(\diamondsuit):u^*g^!\la f^!v^*$.
The definition is as in Proposition~\ref{P204.13}.
\ee
\ethm

\prf
For any pair of flat 1-morphisms $u,v\in\seq$ we have two functors
$\vseq(u,v)\la\Tri\big(\Dqc(Z),\Dqc(X)\big)$, namely
$u^*\circ \big[p_2^{}(-)\big]^!$ and $\big[p_1^{}(-)\big]^!\circ v^*$.
In Proposition~\ref{P204.13} we established that
$\theta$ defines a natural transformation
$\theta:u^*\circ \big[p_2^{}(-)\big]^!\la \big[p_1^{}(-)\big]^!\circ v^*$.
Identities are obviously respected,
hence what remains to be proved is that $\theta$ is compatible with
composition. That is: we have to show that, if
\[\xymatrix@C+15pt@R+5pt{
\wt W \ar[r]^{\wt u}\ar[d]_{\wt f} \ar@{}[dr]|-{(\clubsuit)}& \wt X\ar[d]^-{\wt g}\\
  W \ar[r]^u\ar[d]_f \ar@{}[dr]|-{(\diamondsuit)}& X\ar[d]^g\\
  Y\ar[r]^v & Z
}\]
are two composable morphisms in $\vseq$, then the pentagon
\[
\xymatrix@C-22pt@R+5pt{ 
& \wt u^*{(g \wt g)}^! \ar[rr]^{\theta\sst\left(\!\!\!\!\begin{array}{c}\clubsuit\\*[-5pt]
\diamondsuit\end{array}\!\!\!\!\right)}\ar[dl]_{\rho(\wt g,g)} & & 
          (f\wt f)^!v^*\ar[dr]^{\rho(\wt f,f)} & \\
{\wt u}^*\wt g^!g^!\ar[rr]_{\theta(\clubsuit)g^!} & &{\wt f}^!u^*g^!\ar[rr]_{{\wt f}^!\theta(\diamondsuit)} & &
 {\wt f}^!f^!v^*
}\]
commutes. For future reference we name this pentagon {\bf Compent}.
We begin with the easy cases:

\medskip

\nin
{\bf Case 1.}\ \  \emph{The pentagon {\bf Compent} commutes if either $g$ is of
finite type and universally quasi-proper, or $\wt g$ is a dominant, flat
monomorphism.}

\medskip

\nin
\emph{Proof of Case 1.}\ \ Let us treat the case where $g$ is assumed
of finite type and universally quasi-proper, the case where $\wt g$ is
a dominant, flat monomorphism is similar and left to the reader. The point is
that, under the assumption, the square $(\diamondsuit)$ has a
particularly simple lifting to $\nsq(W,X,Y,Z)$ namely the diagram
 
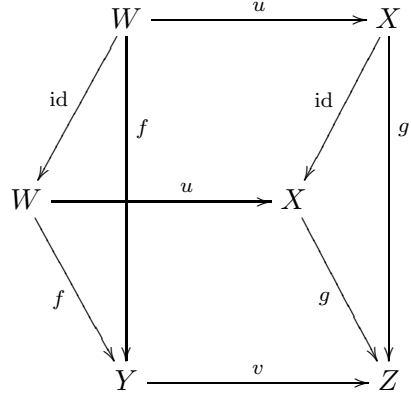
\[\xymatrix@C-5pt@R+30pt{
    & W\ar[dl]_-{\id}\ar[rrr]^-{u}\ar[dd]^(0.3){f} & & & X\ar[dl]_-{\id}
    \ar[dd]^(0.3){g} \\
 W\ar[dr]_-{f} \ar[rrr]^(.6){u}& & &X\ar[dr]_-{g} & \\
& Y\ar[rrr]^-{v}  &&  & Z \\
 }\]
Choose any lifting to $\nsq(\wt W,\wt X,W,X)$ of the square $(\clubsuit)$,
and we have 
that our chosen
 
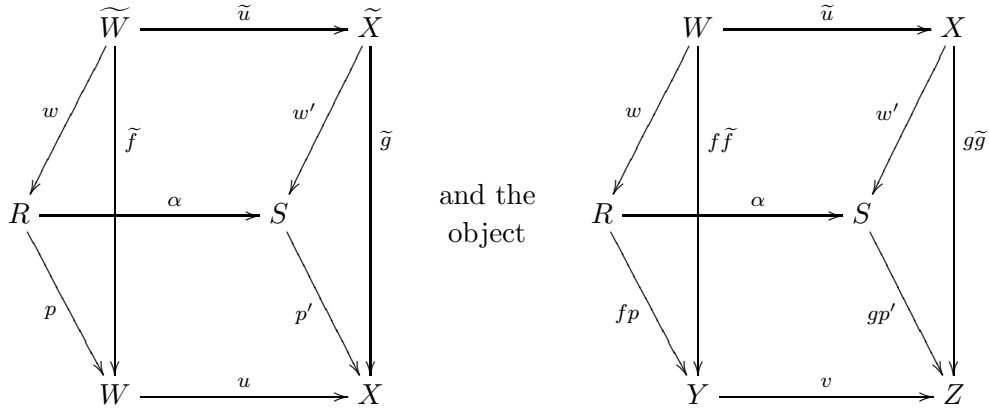
\[\xymatrix@C-5pt@R+30pt{
    &\wt  W\ar[dl]_-{w}\ar[rrr]^-{\wt u}\ar[dd]^(0.3){\wt f} & & &\wt X\ar[dl]_-{w'} \ar[dd]^(0.3){\wt g} & &  & & W\ar[dl]_-{w}\ar[rrr]^-{\wt u}\ar[dd]^(0.3){f\wt f} & & & X\ar[dl]_-{ w'}
    \ar[dd]^(0.3){g\wt g} \\
 R\ar[dr]_-{p} \ar[rrr]^(.6){\alpha}& & &S\ar[dr]_-{p'} & 
 & \ar@{}[r]|-{\ds\begin{array}{c}\text{and the}\\
 \text{object}\end{array}}& & R\ar[dr]_-{f p} \ar[rrr]^(.6){\alpha}& & & S\ar[dr]_-{g p'} & \\
& W\ar[rrr]^-{u}  &&  & X 
& & && Y\ar[rrr]^-{ v}  &&  & Z \\
 }\]
are compatible liftings to $\nsq(\wt W,\wt X,W,X)$ and
$\nsq(\wt W,\wt X,Y,Z)$ of the objects $(\clubsuit)$ and
$\left(\!\!\!\!\begin{array}{c}\clubsuit\\*[-5pt]
\diamondsuit\end{array}\!\!\!\!\right)$. With respect to this choice
the reader will easily check the commutativity of
\[
\xymatrix@C-22pt@R+5pt{ 
& \wt u^*[{w'}^*(gp')^\times] \ar[rr]^{\Theta\sst\left(\!\!\!\!\begin{array}{c}\clubsuit\\*[-5pt]
\diamondsuit\end{array}\!\!\!\!\right)}\ar[dl]_{\rho(\wt g,g)} & & 
          [w^*(fp)^\times]v^*\ar[dr]^{\rho(\wt f,f)} & \\
{\wt u}^*\wt [{w'}^*{p'}^\times]g^\times\ar[rr]_{\Theta(\clubsuit)g^\times} & &(w^*p^\times)u^*g^\times\ar[rr]_{(w^*f^\times)\Theta(\diamondsuit)} & &
 (w^*p^\times)f^\times v^*
}\]
and this completes the proof of Case 1.

\medskip

\nin
{\bf Case 2.} \emph{The composition pentagon commutes if $\wt g$ is of
finite type and universally quasi-proper, and $g$ is a dominant, flat
monomorphism.}

\medskip

\nin
\emph{Proof of Case 2.}\ \ We may choose Nagata compactifications
for the composite $g\wt g:\wt X\la Z$, then pull back
to extend to a diagram
\[\xymatrix@C+15pt@R+5pt{
  \wt W \ar[r]^{\wt u}\ar[d]_{w} \ar@{}[dr]|-{(\heartsuit)}& \wt X\ar[d]^-{w'}
 \\
  R\ar[r]^{\alpha}\ar[d]_{p} \ar@{}[dr]|-{(\spadesuit)}& S\ar[d]^{p'} \\
  Y\ar[r]^v & Z
}\]
This gives us the 2-commutative cube
\[
\xymatrix@R-20pt@C-20pt{
 \wt W \ar[drrr]|{w}\ar[ddddd]|{\wt f}  
          \ar[rrrrr]|{\wt u} &&&&& \wt X
     \ar[drrr]|{w'}\ar[ddddd]|{\wt g}&&&\\
&&& R\ar[ddddd]|{p}  
          \ar[rrrrr]|{\alpha} &&&&& S\ar[ddddd]|{p'}\\
&&&&&&&&\\
&&&&&&&&\\
&&&&&&&&\\
W \ar[drrr]|{f} \ar[rrrrr]|{u} &&&&& X \ar[drrr]|{g}
&&&\\
&&&  Y \ar[rrrrr]|{v} &&&&& Z\\
}
\]
and I assert that all the faces are 2-cartesian. We know this for four of
the faces: for the faces $(\clubsuit)$ and $(\diamondsuit)$ this
is given, and for $(\heartsuit)$ and $(\spadesuit)$ this is by construction.
We have to prove something in the case of the 2-commutative squares
\[\xymatrix@C+15pt@R+5pt{
  \wt W \ar[d]_{\wt f}\ar[r]^{w} \ar@{}[dr]|-{(\dagger)}& R\ar[d]^-{p} & &
  \wt X \ar[d]_{\wt g}\ar[r]^{w'} \ar@{}[dr]|-{(\bullet)}& S\ar[d]^-{p'}
 \\
  W\ar[r]_{f} & Y  & &
  X\ar[r]_g & Z
}\]
We are given that in these squares the horizontal maps are dominant, flat
monomorphisms while the vertical maps are concentrated, of finite
type and universally quasi-proper. Corollary~\ref{CNag.22} therefore
tells us that the squares are 2-cartesian.

To prove the commutativity of the pentagon {\bf Compent} we will consider the
pair of diagrams $B$ and $C$ below 
\[
  \xymatrix{
    \wt W\ar[dr]_-{\wt f} \ar[r]^-{\id}  &\wt W \ar[r]^{w} \ar[d]^{\wt f}\ar@{}[rd]|-{(\dagger)} &  R\ar[d]^{p}\ar[r]^{\alpha} \ar@{}[rd]|-{(\spadesuit)} & S\ar[d]^-{p'} & & \wt W\ar[dr]_-{\wt f} \ar[r]^-{\id}  &\wt W \ar@{}[rd]|-{(\clubsuit)}\ar[r]^{\wt u} \ar[d]^{\wt f} & \wt X\ar[d]^{\wt g}\ar[r]^{w'} \ar@{}[rd]|-{(\bullet)} & S\ar[d]^-{p'} \\
    &    W\ar[dr]_-{f} \ar[r]^-{f}  &Y \ar[r]^{v} \ar[d]^{\id} &  Z\ar[d]^{\id}
    & & &    W\ar[dr]_-{u} \ar[r]^-{u}  &X \ar[r]^{g} \ar[d]^{\id} &  Z\ar[d]^{\id}\\
    & &   Y\ar[dr]_-{v} \ar[r]^-{v} & Z\ar[d]^{\id} & &
    & &   X\ar[dr]_-{g} \ar[r]^-{g} & Z\ar[d]^{\id}\\
& & &  Z
& & & & &  Z
}
\]
The squares are all 2-cartesian and the horizontal maps are all flat, hence
$B$ belongs to $\lnseq(\wt W,W,Y,Z)$ while $C$ belongs to
$\lnseq(\wt W,W,X,Z)$, and Lemma~\ref{L200.5} gives one relation from
each of $B$ and $C$. We have the identities
\begin{eqnarray*}
P\big(\pi_{123}^{}(B)\big)P\big(\pi_{134}^{}(B)\big)&=&
P\big(\pi_{234}^{}(B)\big)P\big(\pi_{124}^{}(B)\big)\,,\\
P\big(\pi_{123}^{}(C)\big)P\big(\pi_{134}^{}(C)\big)&=&
P\big(\pi_{234}^{}(C)\big)P\big(\pi_{124}^{}(C)\big)\,.
\end{eqnarray*}
Note that $P\big(\pi_{234}^{}(B)\big)=\id=P\big(\pi_{234}^{}(C)\big)$,
while $(\heartsuit):w'\wt u\la \alpha w$ and $(\diamondsuit):ug\la vf$ give
isomorphisms of the concatenations of
\[
  \xymatrix{
   \wt W \ar[r]^{w} \ar[d]^{\wt f}\ar@{}[rd]|-{(\dagger)} &  R\ar[d]^{p}\ar[r]^{\alpha} \ar@{}[rd]|-{(\spadesuit)} & S\ar[d]^-{p'} & & \wt W \ar@{}[rd]|-{(\clubsuit)}\ar[r]^{\wt u} \ar[d]^{\wt f} & \wt X\ar[d]^{\wt g}\ar[r]^{w'} \ar@{}[rd]|-{(\bullet)} & S\ar[d]^-{p'} \\
   W  \ar[r]^-{f}  &Y \ar[r]^{v}  &  Z
     & &    W \ar[r]^-{u}  &X \ar[r]^{g}  &  Z
}
\]
from which we deduce that $(\diamondsuit)^*P\big(\pi_{124}^{}(C)=
P\big(\pi_{124}^{}(B)(\heartsuit)^*$. The relations from $B$ and $C$
combine to
\[
P\big(\pi_{123}^{}(B)\big)P\big(\pi_{134}^{}(B)\big)(\heartsuit)^*\eq
(\diamondsuit)^*P\big(\pi_{123}^{}(C)\big)P\big(\pi_{134}^{}(C)\big)
\]
and we deduce
\begin{eqnarray*}
  \rho(\wt f,f)\theta\left(\!\!\!\!\begin{array}{c}\clubsuit\\*[-5pt]
\diamondsuit\end{array}\!\!\!\!\right) &=&
P\big(\pi_{123}^{}(B)\big)P\big(\pi_{134}^{}(B)\big)(\heartsuit)^*\\
&=&
(\diamondsuit)^*P\big(\pi_{123}^{}(C)\big)P\big(\pi_{134}^{}(C)\big)\\
&=& \theta(\diamondsuit)\theta(\clubsuit)\rho(\wt g,g),
\end{eqnarray*}
where the first equality comes from the identities $\rho(\wt f,f)=
P\big(\pi_{123}^{}(B)\big)$ and $\theta\left(\!\!\!\!\begin{array}{c}\clubsuit\\*[-5pt]
\diamondsuit\end{array}\!\!\!\!\right) =
P\big(\pi_{134}^{}(B)\big)(\heartsuit)^*$,
and the third equality is because, for our very simple morphisms
$(\clubsuit)\in\vseq(\wt u,u)$ and $(\diamondsuit)\in\vseq(u,v)$, we have
$\theta(\diamondsuit)=(\diamondsuit)^*$, $\theta(\clubsuit)= P\big(\pi_{123}^{}(C)\big)$ and
$\rho(\wt g,g)=P\big(\pi_{134}^{}(C)\big)$. And this finishes the proof
of Case 2.

Now we move to the general case; the idea is to reduce
to Case 1 and Case 2 using the fact, which we already know from
Theorem~\ref{T200.19}, that the
composition maps $\rho(\alpha,\beta):(\beta\alpha)^!\la\alpha^!\beta^!$
are associative.
We may choose Nagata compactifications
for the morphisms $\wt g:\wt X\la X$ and $g:X\la Z$, and then pull back
to extend to diagrams
\[\xymatrix@C+15pt@R+5pt{
  \wt W \ar[r]^{\wt u}\ar[d]_{w} \ar@{}[dr]|-{(\heartsuit)}& \wt X\ar[d]^-{w'}
  & & 
  W \ar[r]^{u}\ar[d]_{\ov w} \ar@{}[dr]|-{(\bullet)}& X\ar[d]^-{\ov w'}\\
  R\ar[r]^{\alpha}\ar[d]_{p} \ar@{}[dr]|-{(\spadesuit)}& S\ar[d]^{p'} &
  &
\ov  R\ar[r]^{\ov \alpha}\ar[d]_{\ov p} \ar@{}[dr]|-{(\dagger)}& \ov S\ar[d]^{\ov p'} \\
  W \ar[r]^u& X &   &
  Y\ar[r]^v & Z
}\]
In other words: in the 2-category $\vseq$ we factor the
1-morphisms $(\clubsuit)$ and $(\diamondsuit)$, up to 2-isomorphisms,
as composites 
$(\clubsuit)\cong\left(\!\!\!\!\begin{array}{c}\heartsuit\\*[-5pt]
\spadesuit\end{array}\!\!\!\!\right)$
and 
$(\diamondsuit)\cong\left(\!\!\!\!\begin{array}{c}\bullet\\*[-5pt]
\dagger\end{array}\!\!\!\!\right)$
where the vertical maps of $(\heartsuit)$ and $(\bullet)$ are
dominant, flat monomorphisms while the vertical maps of $(\spadesuit)$
and $(\dagger)$ are of finite type and universally quasi-proper.
By the naturality of $\theta$ and $\rho$ we may replace
$(\clubsuit)$ and $(\diamondsuit)$ by the isomorphic
$\left(\!\!\!\!\begin{array}{c}\heartsuit\\*[-5pt]
\spadesuit\end{array}\!\!\!\!\right)$
and 
$\left(\!\!\!\!\begin{array}{c}\bullet\\*[-5pt]
\dagger\end{array}\!\!\!\!\right)$, in other words we may assume
the isomorphisms are equalities.
We wish to prove the commutativity of the pentagon in
\[
\xymatrix@C-22pt@R+5pt{ 
& \wt u^*{(g \wt g)}^! \ar[rr]^{\theta\sst\left(\!\!\!\!\begin{array}{c}\clubsuit\\*[-5pt]
\diamondsuit\end{array}\!\!\!\!\right)}\ar[dl]_{\rho(\wt g,g)} & & 
          (f\wt f)^!v^*\ar[dr]^{\rho(\wt f,f)} & \\
{\wt u}^*\wt g^!g^!\ar[rr]_{\theta(\clubsuit)g^!} & &{\wt f}^!u^*g^!\ar[rr]_{{\wt f}^!\theta(\diamondsuit)} & &
 {\wt f}^!f^!v^* \ar[rrrrrrrr]^{\rho(w,p)\rho(\ov w,\ov p)v^*}&& & &&& & &w^!p^!\ov w^!\ov p^!v^*
}\]
Now $f^!=(\ov p\ov w)^!\stackrel{\rho(\ov w,\ov p)}\la \ov w^!\ov p^!$ and
$\wt f^!=(pw)^!\stackrel{\rho(w,p)}\la w^!p^!$
are isomorphisms by Proposition~\ref{P200.23}~(i) or (ii), and it suffices
to prove that the two composites from top left to bottom right are equal.
But Case 1 gives the commutativity of the pentagons $(1)$ and $(2)$ below
\[
\xymatrix@C-22pt@R+5pt{ 
& & & & & \wt u^*{(g \wt g)}^! \ar[rr]^{\theta\sst\left(\!\!\!\!\begin{array}{c}\clubsuit\\*[-5pt]
\diamondsuit\end{array}\!\!\!\!\right)}\ar[dl]_{\rho(\wt g,g)} & & 
          (f\wt f)^!v^*\ar[dr]^{\rho(\wt f,f)} & \\
 & & &&  {\wt u}^*\wt g^!g^!\ar[rr]^-{\theta(\clubsuit)g^!}
  \ar[dllll]_-{\wt u^*\rho(w',p')g^!}&
  &{\wt f}^!u^*g^!\ar[rr]^-{{\wt f}^!\theta(\diamondsuit)}
  \ar@{}[llld]|-{(1)}\ar@{}[rrrd]|-{(2)}\ar[d]|-{\rho(p,w)u^*g^!}& &
  {\wt f}^!f^!v^* \ar[drrrr]^{\rho(w,p)\rho(\ov w,\ov p)v^*}\\
  {\wt u}^*{w'}^!{p'}^!g^!\ar[rrrr]_-{\theta(\heartsuit){p'}^!g^!} & & && w^!\alpha^*{p'}^!g^!\ar[rr]_-{w^!\theta(\spadesuit)g^!}  && w^!p^!u^*g^!\ar[rr]_-{w^!p^!\theta(\bullet){\ov p'}^!}&  &
  w^!p^!\ov w^!\ov \alpha^*{\ov p'}^!\ar[rrrr]_-{w^!p^!\ov w^!\theta(\dagger)} && & & w^!p^!\ov w^!\ov p^!v^*
}\]
which reduces us to proving the commutativity of the perimeter;
for future reference let us call this perimeter {\bf Perim}. Now consider
the diagram {\bf Bigdia} below
\[
\xymatrix@C-22pt@R+5pt{ 
 & \wt u^*(\ov p'\ov w'p'w')^! \ar[rr]^{\theta\sst\left(\!\!\!\!\begin{array}{c}\heartsuit\\*[-5pt]
\spadesuit\\*[-5pt]
\bullet\\*[-5pt]
\dagger\end{array}\!\!\!\!\right)}\ar[dl]_{{\wt u}^*\rho(w',\ov p'\ov w'p')}\ar@{}[rrd]|-{(1)} & & 
          (\ov p\ov wpw)^!v^*\ar[dr]^{\rho(w,\ov p\ov wp)v^*} & \\
 {\wt u}^*{w'}^!(\ov p'\ov w'p')^! \ar[rr]^-{\theta(\heartsuit)(\ov p'\ov w'p')^!}
  &
 &w^!\alpha^*(\ov p'\ov w'p')^!
\ar[dl]_-{w^!\alpha^*\rho(\ov w'p',\ov p')}
 \ar[rr]|-{{w}^!\theta\sst\left(\!\!\!\!\begin{array}{c}
\spadesuit\\*[-5pt]
\bullet\\*[-5pt]
   \dagger\end{array}\!\!\!\!\right)}
  \ar@{}[rrd]|-{(2)}& &
  w^!(\ov p\ov wp)^!v^* \ar[dr]^{w^!\rho(\ov wp,\ov p)v^*}\\
  & w^!\alpha^*(\ov w'p')^!{\ov p'}^!
\ar@{}[rrd]|-{(3)}
  \ar[dl]_-{ w^!\alpha^*\rho(p',\ov w')^!{\ov p'}^!}\ar[rr]^-{w^!\theta\left(\!\!\!\!\begin{array}{c}
    \spadesuit\\*[-5pt]
    \bullet\end{array}\!\!\!\!\right){\ov p'}^!} & &  w^!(\ov wp)^!\ov\alpha^*{\ov p'}^!\ar[rr]^-{w^!(\ov wp)^!\theta(\dagger)} \ar[rd]^-{w^!\rho(p,\ov w)^!\ov\alpha^*{\ov p'}^!} && w^!(\ov wp)^!{\ov p}^!v^* \\
  w^!\alpha^*{p'}^!{\ov w'}^!{\ov p'}^!\ar[rr]_-{w^!\theta(\spadesuit){\ov w'}^!{\ov p'}^!} && w^!{p}^!u^*{\ov w'}^!{\ov p'}^!\ar[rr]_-{w^!p^!\theta(\bullet){\ov p'}^!} && w^!{p}^!\ov w^!\ov\alpha^*{\ov p'}^!
}\]
In {\bf Bigdia} the pentagons (1) and (2) commute by Case 1 while
pentagon $(3)$
commutes by Case 2.
If we attach to the left and right bottom corners
of {\bf Bigdia} the obviously commutative diagrams
\[
\xymatrix@C-27pt@R+10pt{ 
 & & {\wt u}^*{w'}^!(\ov p'\ov w'p')^!\ar@{}[dr]|-{\ds\text{\bf Rslant}}\ar[dl]_-{{\wt u}^*{w'}^!\rho(\ov w'p',\ov p')} \ar[rr]^-{\theta(\heartsuit)(\ov p'\ov w'p')^!}
  &
 &w^!\alpha^*(\ov p'\ov w'p')^!
\ar[dl]^-{w^!\alpha^*\rho(\ov w'p',\ov p')}
 \\
& {\wt u}^*{w'}^!(\ov w'p')^!{\ov p'}^!
 \ar[dl]_-{{\wt u}^*{w'}^!\rho(p',\ov w')^!{\ov p'}^!}&  & w^!\alpha^*(\ov w'p')^!{\ov p'}^!
 \ar[dl]^-{ w^!\alpha^*\rho(p',\ov w')^!{\ov p'}^!}& \\
 {\wt u}^*{w'}^!{p'}^!{\ov w'}^!{\ov p'}^!\ar[rr]_-{\theta(\heartsuit){p'}^!{\ov w'}^!{\ov p'}^!}&& w^!\alpha^*{p'}^!{\ov w'}^!{\ov p'}^!&& 
}\]
and
\[
\xymatrix@C-22pt@R+5pt{ 
  w^!(\ov wp)^!\ov\alpha^*{\ov p'}^!\ar[rr]^-{w^!(\ov wp)^!\theta(\dagger)} \ar[rd]_-{w^!\rho(p,\ov w)^!\ov\alpha^*{\ov p'}^!} \ar@{}[drrr]|-{\ds\text{\bf Lslant}}&& w^!(\ov wp)^!{\ov p}^!v^*\ar[rd]^-{w^!\rho(p,\ov w)^!{\ov p}^!v^*}  \\
& w^!{p}^!\ov w^!\ov\alpha^*{\ov p'}^!\ar[rr]_-{w^!p^!\ov w^!\theta(\dagger)}  & & w^!{p}^!\ov w^!{\ov p}^!v^* 
}\]
then the perimeter of the union of {\bf Bigdia}, {\bf Rslant} and {\bf Lslant}
agrees with {\bf Perim}, up to the associativity
proved for $\rho$ in Theorem~\ref{T200.19}. The proof is therefore complete.
\eprf

\rmk{R204.19}
In Construction~\ref{C204.11}
we considered the two composites
\[\xymatrix{
     &\seq(W,X)\times\seq(X,Z)\ar[rd]^-{(-)^*(-)^!} &\\
\sq(W,X,Y,Z)\ar[ru]^-{\pi_1^{}} \ar[rd]_-{\pi_2^{}} & &
\Tri\big(\Dqc(Z),\Dqc(W)\big) \\
 & \seq(W,Y)\times\seq(Y,Z)\ar[ur]_-{(-)^!(-)^*} &
}\]
which we called $\fkr$ and $\fks$, and provided a procedure which was
proved in Proposition~\ref{P204.13} to provide a natural transformation
$\theta:\fkr\la\fks$. This was motivated by Remark~\ref{R204.5} where
it was observed that, if we replace $(-)^!$ by $(-)^\times$, then the
natural transformation $\theta$ is analogous to the usual base-change
map $\Phi$. More formally: we could study the
composite functors
\[\xymatrix{
     &\seq(W,X)\times\seq(X,Z)\ar[rd]^-{(-)^*(-)^\times} &\\
\sq(W,X,Y,Z)\ar[ru]^-{\pi_1^{}} \ar[rd]_-{\pi_2^{}} & &
\Tri\big(\Dqc(Z),\Dqc(W)\big) \\
 & \seq(W,Y)\times\seq(Y,Z)\ar[ur]_-{(-)^\times(-)^*} &
}\]
let us name them $\fkr'$ and $\fks'$. To make it concrete:
an object in the category $\sq(W,X,Y,Z)$ is a
2-cartesian square $(\diamondsuit)$ with flat horizontal maps below
\[\xymatrix{
  W\ar[r]^u\ar[d]_f &X\ar[d]^g \\
  Y\ar[r]^v & Z
}\]
and the four functors we consider
take $(\diamondsuit)$ to
$\fkr'(\diamondsuit)=u^*g^\times$, $\fkr(\diamondsuit)=u^*g^!$,
$\fks'(\diamondsuit)=f^\times v^*$ and
$\fks(\diamondsuit)=f^!v^*$. The natural transformations
$\Phi:\fkr'\la\fks'$ and $\theta:\fkr\la\fks$ take the object
$(\diamondsuit)$ (respectively) to the base-change map
$\Phi(\diamondsuit):u^*g^\times\la f^\times v^*$ and the (more complicated)
base-change map $\theta(\diamondsuit):u^*g^!\la f^!v^*$. Of course we also
have the natural transformation $\psi:(-)^\times\la(-)^!$ of Theorem~\ref{T202.5}
and we might wonder how the two are related. We prove
\ermk

\pro{P204.21}
Let the notation be as in Remark~\ref{R204.19}. Then the following square
of natural transformations commutes
\[\xymatrix{
  \fkr' \ar[r]^\Phi \ar[d]_\psi & \fks'\ar[d]^\psi \\
  \fkr\ar[r]^\theta & \fks
}\]  
\epro

\prf
The commutativity may be checked object by object, and we are free to view
our object $(\diamondsuit)\in\sq(W,X,Y,Z)$ as belonging to
$\vseq(u,v)$. In
the category $\vseq$ choose a factorization of $(\diamondsuit)$ as
\[\xymatrix@C+15pt@R+5pt{
  W \ar[r]^{u}\ar[d]_{w} \ar@{}[dr]|-{(\heartsuit)}& X\ar[d]^-{w'}
 \\
  R\ar[r]^{\alpha}\ar[d]_{p} \ar@{}[dr]|-{(\spadesuit)}& S\ar[d]^{p'} \\
  Y\ar[r]^v & Z
}\]
where $w'$ is a dominant, flat monomorphism and $p'$ is of finite type
and universally quasi-proper---in other words choose
a Nagata compactification for the map $g:X\la Z$ and pull back.
We wish to show the commutativity
of the square
\[\xymatrix@C+20pt{
    \fkr'
    {\left(\!\!\!\!\begin{array}{c}
    \heartsuit\\*[-5pt]
    \spadesuit\end{array}\!\!\!\!\right)}
    \ar[r]^{\Phi{\sst\left(\!\!\!\!\begin{array}{c}
    \heartsuit\\*[-5pt]
    \spadesuit\end{array}\!\!\!\!\right)} }\ar[d]_\psi & \fks'{\left(\!\!\!\!\begin{array}{c}
    \heartsuit\\*[-5pt]
    \spadesuit\end{array}\!\!\!\!\right)}\ar[d]^\psi \\
  \fkr{\left(\!\!\!\!\begin{array}{c}
    \heartsuit\\*[-5pt]
    \spadesuit\end{array}\!\!\!\!\right)}\ar[r]^{\theta{\left(\!\!\!\!\begin{array}{c}
    \heartsuit\\*[-5pt]
    \spadesuit\end{array}\!\!\!\!\right)}} & \fks{\left(\!\!\!\!\begin{array}{c}
    \heartsuit\\*[-5pt]
    \spadesuit\end{array}\!\!\!\!\right)}
}\]  
and concretely this is the square $(1)$ in the diagram below
\[\xymatrix@C+30pt{
    u^*(p'w')^\times\ar@{}[rd]|-{(1)}
    \ar[r]^-{\Phi{\sst\left(\!\!\!\!\begin{array}{c}
    \heartsuit\\*[-5pt]
    \spadesuit\end{array}\!\!\!\!\right)} }\ar[d]_{\psi(p'w')} & (pw)^\times v^*\ar[d]^{\psi(pw)} \ar@{}[rd]|-{(2)}\ar[r]^-{\delta(w,p)} & w^\times p^\times v^*\ar[d]^-{\psi(w)\psi(p)}\\
u^*(p'w')^!\ar[r]_-{\theta{\left(\!\!\!\!\begin{array}{c}
    \heartsuit\\*[-5pt]
    \spadesuit\end{array}\!\!\!\!\right)}} & (pw)^!v^* \ar[r]_-{\rho(w,p)} & w^!p^!v^*
}\]
The map $\rho(w,p)$ is an isomorphism by Proposition~\ref{P200.23}~(i) or (ii),
and the square $(2)$ commutes by Theorem~\ref{T202.5}. Hence it suffices
to prove the commutativity of the perimeter. Now
in Theorem~\ref{T204.17} we proved the commutativity of the pentagon
\[
\xymatrix@C-22pt@R+5pt{ 
& \wt u^*{(pw)}^! \ar[rr]^{\theta\sst\left(\!\!\!\!\begin{array}{c}\heartsuit\\*[-5pt]
\spadesuit\end{array}\!\!\!\!\right)}\ar[dl]_{\rho(w,p)} & & 
          (p'w')^!v^*\ar[dr]^{\rho(w'.p')} & \\
{u}^*w^!p^!\ar[rr]_{\theta(\heartsuit)g^!} & &{w}^!u^*{p'}^!\ar[rr]_{{w}^!\theta(\spadesuit)} & &
 {w}^!p^!v^*
}\]
and the commutativity of the pentagon
\[
\xymatrix@C-22pt@R+5pt{ 
& \wt u^*{(pw)}^\times \ar[rr]^{\Phi\sst\left(\!\!\!\!\begin{array}{c}\heartsuit\\*[-5pt]
\spadesuit\end{array}\!\!\!\!\right)}\ar[dl]_{\delta(w,p)} & & 
          (p'w')^\times v^*\ar[dr]^{\delta(w'.p')} & \\
{u}^*w^\times p^\times\ar[rr]_{\Phi(\heartsuit)g^\times} & &{w}^\times u^*{p'}^\times\ar[rr]_{{w}^\times\Phi(\spadesuit)} & &
 {w}^\times p^\times v^*
}\]
is formal, base-change maps concatenate. The perimeter we need to
show commutative is therefore equal to the
perimeter of 
\[\xymatrix@C+30pt@R+20pt{
    u^*(p'w')^\times\ar@{}[rd]|-{(3)}
    \ar[r]^-{\delta(w',p')} \ar[d]_{\psi(p'w')} & u^*{w'}^\times{p'}^\times
\ar@{}[rd]|-{(4)}
\ar[r]^-{\Phi(\heartsuit){p'}^\times} \ar[d]|-{\psi(w')\psi(p')}
 & w^\times \alpha^*{p'}^\times
\ar@{}[rd]|-{(5)}
    \ar[r]^-{w^\times\Phi(\spadesuit)} \ar[d]|-{\psi(w)\psi(p')}
& w^\times p^\times v^*\ar[d]^-{\psi(w)\psi(p)}\\
u^*(p'w')^!\ar[r]_-{\rho(w',p')} & u^*{w'}^!{p'}^!
    \ar[r]_-{\theta(\heartsuit){p'}^!}  &{w}^!\alpha^*{p'}^!
    \ar[r]_-{w^!\theta(\spadesuit)}  & w^!p^!v^*
}\]
and the square $(3)$ commutes by 
Theorem~\ref{T202.5}.
We are reduced to proving that the squares $(4)$ and $(5)$ commute,
in other words we are reduced to the special cases where the vertical maps
are either dominant, flat monomorphisms or of finite type and universally
quasi-proper. We prove these cases below.

\medskip

\nin
{\bf Case 1.}\ \ \emph{The Proposition is true for objects $(\diamondsuit)\in
\sq(W,X,Y,Z)$ where the vertical maps are of finite type and universally
quasi-proper.}

\medskip

\nin
\emph{Proof of Case 1.}\ \
We are given an object $(\diamondsuit)\in\nsq(W,X,Y,Z)$, that is
a 2-cartesian square
\[\xymatrix{
  W\ar[r]^u\ar[d]_f &X\ar[d]^g \\
  Y\ar[r]^v & Z
}\]
and we assume the vertical maps $f$ and $g$ are of finite type and
universally quasi-proper. To evaluate the various maps we lift
the object $(\diamondsuit)$ via the functor
$F:\nsq(W,X,Y,Z)\la
\sq(W,X,Y,Z)$, and the lifting we choose is the 2-commutative
diagram
 \[\xymatrix@C-5pt@R+30pt{
    & W\ar[dl]_-{\id}\ar[rrr]^-{u}\ar[dd]^(0.3){f} & & & X\ar[dl]_-{\id}
    \ar[dd]^(0.3){g} \\
 W\ar[dr]_-{f} \ar[rrr]^(.6){u}& & &X\ar[dr]_-{g} & \\
& Y\ar[rrr]^-{v}  &&  & Z \\
 }\]
For this choice the square we need to prove commutative
reduces to
\[\xymatrix{
  u^*g^\times\ar[r]^{\Phi(\diamondsuit)}\ar[d]_\id &f^\times v^*\ar[d]^\id \\
  u^*g^!\ar[r]^{\Phi(\diamondsuit)} & f^!v^*
}\]
\medskip

\nin
{\bf Case 2.}\ \ \emph{The Proposition is true for objects $(\diamondsuit)\in
\sq(W,X,Y,Z)$ where the vertical maps are dominant, flat monomorphisms.}

\medskip

\nin
\emph{Proof of Case 2.}\ \
In this case the lifting we choose for the object
$(\diamondsuit)$, via the functor
$F:\nsq(W,X,Y,Z)\la
\sq(W,X,Y,Z)$,
is the 2-commutative
diagram
 
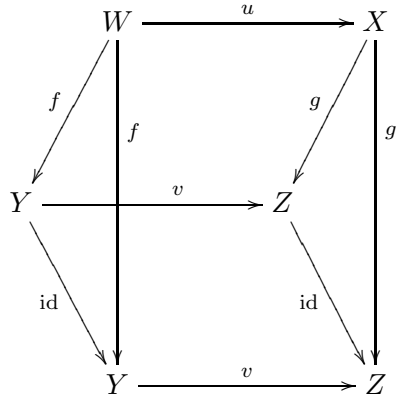
\[\xymatrix@C-5pt@R+30pt{
    & W\ar[dl]_-{f}\ar[rrr]^-{u}\ar[dd]^(0.3){f} & & & X\ar[dl]_-{g}
    \ar[dd]^(0.3){g} \\
 Y\ar[dr]_-{\id} \ar[rrr]^(.6){v}& & &Z\ar[dr]_-{\id} & \\
& Y\ar[rrr]^-{v}  &&  & Z \\
 }\]
Out of this we concoct the 2-commutative
diagram $C$ below
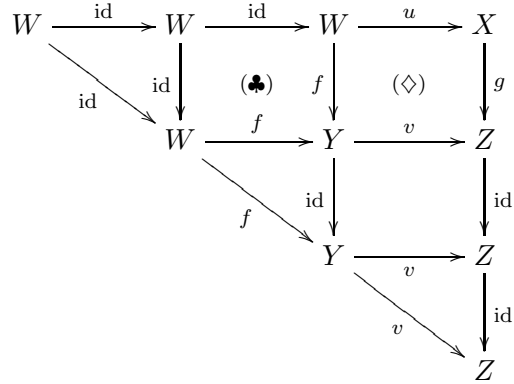
\[
\xymatrix@C+15pt@R+5pt{
W\ar[r]^-\id\ar[dr]_-{\id} & W \ar[d]_\id \ar[r]^{\id}\ar@{}[dr]|{(\clubsuit)} & W \ar@{}[dr]|{(\diamondsuit)} 
 \ar[d]_{f}\ar[r]^{u} & X\ar[d]^{g}
 \\
& W\ar[dr]_-{f}  \ar[r]^{f} & Y\ar[r]^{v}\ar[d]_{\id} &  Z\ar[d]^{\id}  \\
& &  Y\ar[r]_{v}\ar[dr]_-{v} & Z\ar[d]^{\id}\\
& & & Z
}
\]
The square $(\diamondsuit)$ is given to be 2-cartesian, and the square
$(\clubsuit)$ is 2-cartesian because $f$ is a monomorphism. The remaining
square is obviously 2-cartesian. All the horizontal maps are flat, hence
$C$ is an object in $\lnseq(W,W,Y,Z)$ and Lemma~\ref{L200.5} gives
the relation
$P\big(\pi_{123}^{}(C)\big)P\big(\pi_{134}^{}(C)\big)=
P\big(\pi_{234}^{}(C)\big)P\big(\pi_{124}^{}(C)\big)$.

In this relation we have $P\big(\pi_{123}^{}(C)\big)=\psi(f)$,
$P\big(\pi_{134}^{}(C)\big)=\Phi(\diamondsuit)$, and
$P\big(\pi_{234}^{}(C)\big)=\id$. The
relation becomes $\psi(f)\Phi(\diamondsuit)= P\big(\pi_{124}^{}(C)\big)$,
and our problem is to identify $P\big(\pi_{124}^{}(C)\big)$. For this
it helps to consider the two 2-commutative diagrams
\[
\xymatrix@C+15pt@R+5pt{
 W \ar[d]_\id \ar[r]^{\id}\ar@{}[dr]|{(\clubsuit)} & W \ar@{}[dr]|{(\diamondsuit)} 
 \ar[d]_{f}\ar[r]^{u} & X\ar[d]^{g}
 & W \ar[d]_\id \ar[r]^{u}\ar@{}[dr]|{(\heartsuit)} & X \ar@{}[dr]|{(\spadesuit)} 
 \ar[d]_{\id}\ar[r]^{\id} & X\ar[d]^{g}
 \\
W  \ar[r]^{f} & Y\ar[r]^{v} &  Z
& W  \ar[r]^{u} & X\ar[r]^{g} &  Z
}
\]
The diagram on the left is an extract from $C$ and so the squares
$(\clubsuit)$ and $(\diamondsuit)$ are 2-cartesian.
The square $(\spadesuit)$ is 2-cartesian because $g$ is a monomorphism,
and the square $(\heartsuit)$ is trivially 2-cartesian. And the
concatenations are isomorphic via the isomorphism $(\diamondsuit):gu\la vf$.
It follows that
$P\big(\pi_{124}^{}(C)\big)=\Phi(\clubsuit\diamondsuit)=
(\diamondsuit)^*\Phi(\spadesuit)=(\diamondsuit)^*\psi(g)$,
meaning that the square below commutes
\[\xymatrix@C+40pt{
  u^*g^\times\ar[r]^{\Phi(\diamondsuit)}\ar[d]_{\psi(g)} &f^\times v^*\ar[d]^{\psi(f)} \\
  u^*g^*\ar[r]^{\theta(\diamondsuit)=(\diamondsuit)^*} & f^*v^*
}\]
\eprf

\rmk{R204.295}
In Theorem~\ref{T204.15} we gave, on the category $\hseq$,
a natural transformation
$\theta:(-)^*\circ p_2^{}\la (-)^*\circ p_1^{}$ analogous to the base-change
map $\Phi:(-)^*\circ p_2^{}\la (-)^*\circ p_1^{}$. In
Theorem~\ref{T204.17} we gave, on the category $\vseq$,
a natural transformation
$\theta:(-)^!\circ p_2^{}\la (-)^!\circ p_1^{}$ analogous to the base-change
map $\Phi:(-)^\times\circ p_2^{}\la (-)^\times\circ p_1^{}$. On a 1-morphism
$(\diamondsuit)$, in either category, $\theta$ comes down to
the map $\theta(\diamondsuit):u^*g^!\la f^! v^*$ while
$\Phi$ is the map $\Phi(\diamondsuit):u^*g^\times\la f^\times v^*$.

We have produced lots of 2-functors and natural transformations: we can
assemble some of them into the two squares below, where 
the square on the left has the 2-category $\hseq$ for input
while the square on the right
begins with the 2-category $\vseq$
\[\xymatrix{
  (-)^*\circ p_2^{} \ar[d]_-{\id}\ar[r]^-{\Phi} & (-)^*\circ p_1^{}\ar[d]^-{\id} & & (-)^\times\circ p_2^{}\ar[r]^-{\Phi}\ar[d]_-{\psi} & (-)^\times\circ p_1^{}\ar[d]^-{\psi}\\
  (-)^*\circ p_2^{} \ar[r]^-{\theta} & (-)^*\circ p_1^{} & & (-)^!\circ p_2^{}\ar[r]^-{\theta} & (-)^!\circ p_1^{}
}\]
If we combine  Theorems~\ref{T204.15} and \ref{T204.17} 
with Proposition~\ref{P204.21} we learn that these squares
can be filled in by modifications.
Since we are in the world of 2-categories let us elaborate a tiny
bit.

Given an object $f:W\la Y$ of the 2-category $\hseq$, the 2-functor $p_2^{}$
takes it to $Y$ while the 2-functor $p_1^{}$ takes it to $W$. The natural
transformation $\Phi$ takes the object $f$ to
the 1-morphism $f^\times:\Dqc(Y)\la\Dqc(W)$
while the natural transformation $\theta$ takes $f$ to $f^!:\Dqc(Y)\la\Dqc(W)$.
The natural transformation $\id$ is the identity on objects. What is being
asserted is that the assignment taking the object $f\in\hseq$ to
the 2-morphism $\psi(f):f^\times\la f^!$ extends to a
modification of the composites in the square on the left,
i.e.~is compatible with the rest of the structure. This
follows from Proposition~\ref{P204.21}.

Given an object $u:W\la X$ in the 2-category $\vseq$,
the 2-functor $p_2^{}$
takes it to $X$ while the 2-functor $p_1^{}$ takes it to $W$. The natural
transformations $\Phi$ and $\theta$ both take
the object $u$ to the 1-morphism
$u^*:\Dqc(X)\la\Dqc(W)$, and the natural transformation
$\psi$ is the identity on objects. What is being asserted for the 2-category
$\vseq$ is that the square strictly commutes; no modification
is necessary. Once again, the proof amounts to an application of 
Proposition~\ref{P204.21}.
\ermk

We end this section with a result telling us that, under suitable
conditions, the map $\theta(\diamondsuit)$ is an
isomorphism.

\pro{P204.23}
Let $(\diamondsuit)$ be an object of $\sq(W,X,Y,Z)$, given by the
2-cartesian square
\[\xymatrix{
  W\ar[r]^u\ar[d]_f &X\ar[d]^g \\
  Y\ar[r]^v & Z
}\]
The map $\theta(\diamondsuit):u^*g^!\la f^!v^*$ is an
isomorphism if either of the two  conditions below is satisfied
\be
\item
  $f$ is of finite Tor-dimension. More generally: if $u':U\la W$ is
  a map whose image lies in
the subset of $\,W$ on which $f$ is of finite Tor-dimension, then 
${u'}^*\theta(\diamondsuit):{u'}^*u^*g^!\la {u'}^*f^!v^*$
is an isomorphism.
\item
We restrict to the subcategory $\Dqcpl(Z)\subset\Dqc(Z)$.
\ee
\epro

\prf
We begin by choosing a preimage $A$ for $(\diamondsuit)$
under the functor $F:\nsq(W,X,Y,Z)\la \sq(W,X,Y,Z)$, that is
the data of $A$ consists of the 2-cartesian squares
 \[\xymatrix{
   W \ar@{}[dr]|{(\heartsuit)}\ar[r]^u\ar[d]_{w} & X\ar[d]^{w'} & &  R\ar@{}[dr]|{(\spadesuit)} \ar[r]^\alpha\ar[d]_{p} & X\ar[d]^{p'} \\
  R\ar[r]^\alpha & S & &Y\ar[r]^v & Z 
 }\]
as well as 2-isomorphisms $pw\la f$ and $p'w'\la g$ which give
an isomorphism of the concatenation  with $(\diamondsuit)$.
By Proposition~\ref{P204.13} we are free to use any preimage.
The map $\theta(\diamondsuit)$ has the property that
$\theta(\diamondsuit)\cong\Theta(A)$,
where $\Theta(A)$ was defined to be the composite
\[
\xymatrix@C+30pt{
u^*{w'}^* {p'}^\times \ar[r]_{(\heartsuit)^*{p'}^\times } \ar@/^2pc/[rr]^{\Theta(A)} & {w}^*\alpha^*{p'}^\times
\ar[r]_{{w}^*\Phi(\spadesuit)} & {w}^*{p}^\times v^*
}
\]
The map $(\heartsuit)^*$ is an isomorphism unconditionally. Therefore
${u'}^*\theta(\diamondsuit)$ is an isomorphism if and only if
${u'}^*{w}^*\Phi(\spadesuit)$ is an isomorphism.

The map $w$ is flat, and the image of $wu'$ lies in the subset
on which $p$ is of finite Tor-dimension if and
only if the image of $u'$ lies in the subset
on which $pw\cong f$ is
of finite Tor-dimension. The assertions of the Proposition
now follow from Theorem~\ref{T4.13}~(i) and (ii), which
give conditions under which $(wu')^*\Phi(\spadesuit)$
is an isomorphism.
\eprf

\section{The 2-functor $(-)^*$ is a monoid and $(-)_*$ and $(-)^\times$
are modules over it}
\label{S805}

Consider the 2-categories $\seq$ and $\Tri$: we have 
met the 2-functors $(-)^*$, $(-)^\times$ and $(-)^!$, all
of which are 2-functors $\seq\la\Tri$. On objects the three 2-functors
agree---all three
take the object $X\in\seq$ to $\Dqc(X)\in\Tri$. The
category $\Dqc(X)$
is not only triangulated,
it also has a tensor product, and it is natural to
wonder how our three 2-functors behave with respect to the tensor
product.
Of them $(-)^*$ is best: it respects the tensor
product. This section is devoted to formal nonsense
about 2-functors like $(-)^*$ and their adjoints.
We will formulate a language in which we can say that the 2-functor
$(-)^*$ is a monoid and the 2-functor $(-)^\times$ is a module over it.
This turns out to be a compact way of packing a great many
naturality properties.

We begin with

\dfn{D805.-1}
A 2-category $\cb$ is called \emph{premonoidal} if there is a 2-functor
$\fkm:\cb\times\cb\la\cb$ satisfying the obvious associativity conditions.
For this article we will assume that the associativity is strict.
\edfn

\exm{E805.1}
Let $\cb$ be the 2-category whose objects are finite products of
triangulated categories, and whose morphisms are finite  products of
functors $F:\prod_{i=1}^nB_i\la B'$, where $B'$ and each $B_i$ are
triangulated, and where $F$ is triangulated (separately) in
each variable $B_i$. We allow the empty product as an object of $\cb$.
The premonoidal structure takes
two objects of $\cb$, that is $\prod_{i=1}^m B_i$ and $\prod_{j=1}^n B'_j$,
to the object $\prod_{i=1}^m B_i\times\prod_{j=1}^n B'_j$. The
empty product is the unit for this premonoidal structure. 
The 2-morphisms are unrestricted, they
are just the natural transformations.
\eexm

For the next definition it helps to introduce the following notation:
if $\ca$ is any 2-category then $\Delta_n^{}:\ca\la\ca^n$
is the $n$-fold diagonal, the 2-functor
taking an object $X\in\ca$ to $(X,X,\ldots,X)\in\ca^n$. Given 2-functors
$F_1,F_2,\ldots,F_n:\ca\la\cb$, where $\cb$ is premonoidal,
we will be considering composite 2-functors of the form
\[
\xymatrix{
  \ca\ar[r]^-{\Delta_n} & \ca^n \ar[rr]^-{\prod_{i=1}^nF_i} && \cb^n\ar[r]^-{\fkm} &
  \cb
}  
\]
and, in the notation, we will almost always omit the $\fkm$ and sometimes
even omit the $\Delta_n$. Thus $\fkm\circ\big[F\times G\big]\circ\Delta_2$
will usually be abbreviated to $\big[F\times G\big]\circ\Delta_2$ and sometimes
even to $F\times G$.

The next concept is not new, it has occurred in the category-theoretic
literature many times (under different names). We will come to this again
in Remark~\ref{R805.99999}.

\dfn{D805.3}
Let $\ca$ be a 2-category and let $\cb$ be a premonoidal 2-category.
A 2-functor $(-)^*:\ca\la\cb$ is a \emph{premonoid} if
it comes equipped with an associative pseudonatural transformation
 $\mu:[(-)^*\times(-)^*]\circ\Delta_2^{}
  \la(-)^*$.
The associativity means that we are also given a modification isomorphism
between the composites below
\[\xymatrix@C+40pt{
  [(-)^*\times(-)^*\times(-)^*]\circ\Delta_3^{}
  \ar[r]^-{\id\times\mu}\ar[d]_-{\mu\times\id} &
     [(-)^*\times(-)^*]\circ\Delta_2^{}\ar[d]^-{\mu} \\
[(-)^*\times(-)^*]\circ\Delta_2^{}\ar[r]^-{\mu} & (-)^*
}\]
such that the associativity pentagon for the five composites
$[(-)^*\times(-)^*\times(-)^*\times(-)^*]\circ\Delta_4^{}\la(-)^*$
commutes.
\edfn

\exm{E805.5}
With $\cb$ as in Example~\ref{E805.1}, let $\ca$ be the 2-category
of algebraic stacks.
There is a 2-functor $(-)^*:\ca\la\cb$, we have met it before:
\be
\item
  On objects: the object $X\in\ca$ is mapped to $X^*=\Dqc(X)\in\Tri\subset\cb$.
\item
  On 1-morphisms: the 1-morphism $f:X\la Y$ is mapped to the 1-morphism
  $f^*:\Dqc(Y)\la\Dqc(X)$.
\item
  On 2-morphisms: given $\lambda:f\la g\in\ca$, the 2-functor $(-)^*$ maps it to
  $\lambda^*:f^*\la g^*$ in $\Tri\subset\cb$.
\setcounter{enumiii}{\value{enumi}}
\ee
So far we haven't mentioned the tensor product but,
as we have already observed, 
for any
object $X\in\ca$ the object $X^*=\Dqc(X)$ 
is a symmetric premonoid in $\cb$, meaning there is a symmetric tensor
product $\Dqc(X)\times\Dqc(X)\la\Dqc(X)$.
It is the functor taking the
pair of objects
$E,F\in\Dqc(X)$ to $E\oo F\in\Dqc(X)$.
And this does extend to a pseudonatural transformation $\mu$.
\be
\setcounter{enumi}{\value{enumiii}}
\item
  On objects: the object $X\in\ca$ is sent to the
  tensor-product functor
  $\mu_{X}^{}:\Dqc(X)\times \Dqc(X)\la \Dqc(X)$. 
\item
  On morphisms: for any pair of objects $X,Y\in\ca$ we have two
  functors $\ca(X,Y)\la\cb(Y^*\times Y^*,X^*)$, namely
  $\mu_{X}^{}\circ[(-)^*\times(-)^*]$ and $(-)^*\circ\mu_{Y}^{}$.
  We need
  to provide a natural transformation $\mu$ between them
  and, as we are assuming $\mu$ to be \emph{pseudonatural,}
  this map must be an isomorphism---for once
  the direction doesn't matter. Concretely it means
  that for every morphism $f:X\la Y$ and any pair of
  objects $E,F\in\Dqc(Y)$ we give
  an isomorphism $f^*(E\oo F)\cong f^*E\oo f^*F$, and these
  must be  compatible with
  composition and be natural in everything in sight.
\setcounter{enumiii}{\value{enumi}}
\ee
\eexm

Until now we haven't worried much about the variance of our 2-functors,
there was no need to, the 2-functors $(-)^*$, $(-)^\times$ and $(-)^!$
all have the same variance, and they form our main
object of study. This section is the exception, in the next few
results variance plays a role.

The next definition is
a 2-category generalization of an old notion, introduced in Dubuc and
Street~\cite{Dubuc-Street70}. What we  call
natural transformations below goes by the name ``dinatural transformations''
in~\cite{Dubuc-Street70}, 
where the ``di'' stands for diagonal.

\dfn{D805.6}
Let $\ca$ be a 2-category and let $\cb$ be a premonoidal 2-category.
Suppose we are given
functors $E,F,G,H:\ca\la\cb$, with $E,G$ contravariant and $F,H$ covariant.
A \emph{natural transformation} $\mu:E\times F\la G\times H$ is
the following data:
\be
\item
  On objects: for every object $X\in\ca$ we are given a 1-morphism
  $\mu_X^{}:E(X)\times F(X)\la G(X)\times H(X)$ in the category $\cb$.
\item
  Let $X,Y\in\ca$ be objects. We have two functors $\ca(X,Y)\la
  \cb\big[E(Y)\times F(X)\,,\,G(X)\times H(Y)\big]$, namely
  \[
H(-)\mu_X^{}E(-)\quad\text{and}\quad G(-)\mu_Y^{}F(-)\ .
  \]
  We must be given, for every pair of objects $X,Y\in\ca$, a natural
  transformation $\mu(X,Y)$
  between these functors. The direction of this natural
  transformation determines whether $\mu$ is lax or oplax: for us
  the lax direction will be $ G(-)\mu_Y^{}F(-)\la H(-)\mu_X^{}E(-)$.
\item
The $\mu(X,Y)$ must be compatible with composition.
\ee
\edfn

\elb{E805.6.5}
It might be worth elaborating a little. The hypothesis
about the variance of the 2-functors $E,F,G,H$
means that, for any 1-morphism $f:X\la Y$, we have induced
1-morphisms
\[\begin{array}{ccc}
E(f):E(Y)\la E(X),&\quad&
F(f):F(X)\la F(Y), \\
G(f):G(Y)\la G(X),&\quad&
H(f):H(X)\la H(Y),
\end{array}
\]
and any 2-morphism $\lambda:f\la g$ induces 2-morphisms
\[\begin{array}{ccc}
E(\lambda):E(f)\la E(g),&\quad&
F(\lambda):F(g)\la F(f), \\
G(\lambda):G(f)\la G(g),&\quad&
H(\lambda):H(g)\la H(f).
\end{array}
\]
Definition~\ref{D805.6}(ii) says that, for any 1-morphism $f$, we must be given
a 2-morphism between the composites
\[\xymatrix@R+5pt{
  & E(X)\times F(X) \ar[r]^-{\mu_X^{}} & G(X)\times H(X)
            \ar[dr]^-{\id\times H(f)} & \\ 
E(Y)\times F(X)\ar[ur]^-{E(f)\times\id} \ar[dr]_-{\id\times F(f)} & &
& G(X)\times H(Y) \\
  &  E(Y)\times F(Y) \ar[r]^-{\mu_Y^{}} & G(Y)\times H(Y)
            \ar[ur]_-{G(f)\times \id} & 
}\]
The natural transformation is lax if this 2-morphism goes up, oplax if it
goes down.

The meaning of Definition~\ref{D805.6}(iii) is that, if
we are given composable 1-morphisms $X\stackrel f\la Y\stackrel g\la Z$
in the 2-category $\ca$, then the 2-morphism of
the hexagon
\[\xymatrix@R+5pt{
  & E(X)\times F(X) \ar[r]^-{\mu_X^{}} & G(X)\times H(X)
            \ar[dr]^-{\id\times H(gf)} & \\ 
E(Z)\times F(X)\ar[ur]^-{E(gf)\times\id} \ar[dr]_-{\id\times F(gf)} & &
& G(X)\times H(Z) \\
  &  E(Z)\times F(Z) \ar[r]^-{\mu_Z^{}} & G(Z)\times H(Z)
            \ar[ur]_-{G(gf)\times \id} & 
}\]
is compatible with the composite of the 2-morphisms of the hexagons
\[\xymatrix@R+5pt{
  & E(X)\times F(X) \ar[r]^-{\mu_X^{}} & G(X)\times H(X)
            \ar[dr]^-{\id\times H(f)} & \\ 
E(Y)\times F(X)\ar[ur]^-{E(f)\times\id} \ar[dr]_-{\id\times F(f)} & &
& G(X)\times H(Y) \\
  &  E(Y)\times F(Y) \ar[r]^-{\mu_Y^{}} & G(Y)\times H(Y)
\ar[ur]_-{G(f)\times \id}\ar[dr]^-{\id\times H(g)} & \\
E(Z)\times F(Y)\ar[ur]^-{E(g)\times\id} \ar[dr]_-{\id\times F(g)} & &
& G(Y)\times H(Z) \\
  &  E(Z)\times F(Z) \ar[r]^-{\mu_Z^{}} & G(Z)\times H(Z)
            \ar[ur]_-{G(g)\times \id} & 
}\]
once we attach the 2-commutative
\[\xymatrix@C-28pt{
& E(Y)\times F(X) \ar[dr]^-{\id\times F(f)} & & & & & & &
& G(X)\times H(Y)\ar[dr]^-{\id\times H(g)} \\
E(Z)\times F(X)\ar[ur]^-{E(g)\times\id}\ar[dr]_-{\id\times F(f)} & &  E(Y)\times F(Y) & && & & & G(Y)\times H(Y)
\ar[ur]^-{G(f)\times \id}\ar[dr]_-{\id\times H(g)} & & G(X)\times H(Z)\\
& E(Z)\times F(Y)\ar[ur]_-{E(g)\times\id}  & & &&& & &
& G(Y)\times H(Z) \ar[ur]_-{G(f)\times\id}
}\]
and recall the 2-morphisms relating
\[\begin{array}{ccc}
E(Z)\stackrel{E(gf)}\la E(X) & \text{ with the composite } & E(Z)\stackrel{E(g)}\la E(Y) \stackrel{E(f)}\la E(X)\\
G(Z)\stackrel{G(gf)}\la G(X) & \text{ with the composite } & G(Z)\stackrel{G(g)}\la G(Y) \stackrel{G(f)}\la G(X) \\
F(X)\stackrel{F(gf)}\la F(Z) & \text{ with the composite } & F(X)\stackrel{F(f)}\la F(Y) \stackrel{F(g)}\la F(Z) \\
H(X)\stackrel{H(gf)}\la H(Z) & \text{ with the composite } & H(X)\stackrel{H(f)}\la H(Y) \stackrel{H(g)}\la H(Z)
\end{array}\]
where the direction of these 2-morphisms depends on which of
$E,F,G,H$ is lax and which is oplax (we left this vague).
\eelb

\rmd{R805.7}
Let $L:\ca\la\cb$ be a contravariant pseudofunctor of 2-categories. Suppose
that, for every 1-morphism $f:X\la Y$ in $\ca$, the 1-morphism
$L(f):L(Y)\la L(X)$ has a right adjoint. Then there is a canonical
covariant pseudofunctor $R:\ca\la\cb$ so that, for every 1-morphism
$f\in\ca$, we have that $L(f)$ is left adjoint to $R(f)$.
We should perhaps remind the reader how this works:
\be
\item
  On objects: for an object $X\in\ca$ put $R(X)=L(X)$.
\item
  On 1-morphisms: for a 1-morphism $f:X\la Y$ let $R(f):L(X)\la L(Y)$
  be some choice
  of a right adjoint to $L(f):L(Y)\la L(X)$. Recall that
  we assume our 2-functors $L$
  strictly respect identities, meaning $L(\id)=\id$, and we choose
  $R(\id)=\id$.
\item
  On 2-morphisms: given a 2-morphism $\lambda:f\la g$ in $\ca$ we have
  in $\cb$ a 2-morphism $L(\lambda):L(f)\la L(g)$, but the adjunctions
  $L(f)\dashv R(f)$ and $L(g)\dashv R(g)$ give rise
  to further 2-morphism,
  the unit and counit 2-morphisms 
  $\eta_f^{}:\id\la R(f)L(f)$, $\e_f^{}:L(f)R(f)\la\id$,
  $\eta_g^{}:\id\la R(g)L(g)$ and $\e_g^{}:L(g)R(g)\la\id$. We define
  $R(\lambda)$ to be the composite
  \[
  \xymatrix@C+10pt{
    R(g) \ar[r]^-{\eta_f^{}R(g)} &R(f)L(f)R(g)\ar[rr]^-{R(f)L(\lambda)R(g)} &
    &
    R(f)L(g)R(g)\ar[r]^-{R(f)\e_g^{}} & R(f)
  }\]
\item
  Let $X\stackrel f\la Y\stackrel g\la Z$ be composable 1-morphisms
  in $\ca$. Because $L$ is a pseudofunctor we are given an isomorphism
  $\tau:L(gf)\la L(f)L(g)$.
  Now $R(g)R(f)$ is a right adjoint for $L(f)L(g)$
  and $R(gf)$ is a right adjoint for $L(gf)$, and the construction
  in (iii) produces out of $\tau:L(gf)\la L(f)L(g)$ an induced
  isomorphism $\s:R(g)R(f)\la R(gf)$.
\ee
Some checking is required to verify that this recipe delivers a pseudofunctor,
but all of it is straightforward.

Dually with left adjoints: given a pseudofunctor $R:\ca\la\cb$
such that, for every
1-morphism $f\in\ca$ the morphism $R(f)$ has a left adjoint, then
there is a canonical
pseudofunctor $L$ with $L(f)\dashv R(f)$ for all 1-morphisms $f\in\ca$.
\ermd

\lem{L805.9}
As in Definition~\ref{D805.6} let $\ca$ be a 2-category and let 
$\cb$ be a premonoidal 2-category. Suppose we are given 2-functors
$E,F,G,H,L,R,L',R':\ca\la\cb$, with $E,G$ contravariant and $F,H$
covariant. Assume that $L,L',R, R'$ are pseudofunctors,
and that the pairs $L,R$
and $L',R'$ are as in Reminder~\ref{R805.7},
in particular
for every 1-morphism $f\in\ca$ we have $L(f)\dashv R(f)$ and
$L'(f)\dashv R'(f)$.
Then
\be
\item
  If $L,L'$ are contravariant, in which case $R,R'$ must be covariant,
then  there is a
canonical bijection between lax natural transformations
$\lambda:E\times L\times F\la G\times L'\times H$ and lax natural
transformations
$\mu:E\times R\times F\la G\times R'\times H$
which agree on objects.
\item
  If $L,L'$ are covariant, in which case $R,R'$ must be contravariant,
then  there is a
canonical bijection between oplax natural transformations
$\lambda:E\times L\times F\la G\times L'\times H$ and oplax natural
transformations
$\mu:E\times R\times F\la G\times R'\times H$
which agree on objects.
\ee
\elem

\prf
We will prove the case where $L,L'$ are contravariant and leave the
covariant case to the reader.
A natural transformation $\lambda:E\times L\times F\la G\times L'\times H$
gives, for every 1-morphism $f:X\la Y$ in $\ca$, a hexagon
as in Elaboration~\ref{E805.6.5} with a 2-morphism filling it.
If we focus on the part
involving $L$ and $L'$ and compose the composable morphisms in the
hexagon that don't involve $L$ and $L'$, then we reduce to a square
\[\xymatrix{
E(Y)\times L(Y)\times F(X)\ar[d]_-{\id\times L(f)\times\id} \ar[r]^-{p} & H(X)\times L'(Y)\times G(Y)\ar[d]^-{\id\times L'(f)\times\id}\\
E(Y)\times L(X)\times F(X) \ar[r]_-{q} & H(X)\times L'(X)\times G(Y)
}\]
with a 2-cell filling it.
Because we are assuming that $\lambda$ is lax and $L,L'$ are contravariant,
the 2-morphism
goes in the direction $\lambda:L'(f)p\la q L(f)$. Base-change gives
us a 2-morphism $\mu:p R(f)\la R'(f)q$, which delivers the
required 2-morphism in the hexagon for
$\mu:E\times R\times F\la G\times R'\times G$.
Perhaps we should remind the reader:
if $\eta:\id\la R'(f)L'(f)$ is the unit of the adjunction $L'(f)\dashv R'(f)$
and $\e:L(f)R(f)\la\id$ is the counit of the
adjunction $L(f)\dashv R(f)$, then $\mu$ is defined is as
the composite
\[\xymatrix{
  pR(f)\ar[r]^-{\eta} & R'(f)L'(f)pR(f)\ar[r]^-{\lambda} &
  R'(f)qL(f)R(f) \ar[r]^-{\e} & R'(f)q
}\]
and, comparing with Reminder~\ref{R805.7}, we see that this is natural in $f$.
This construction is bijective, we leave to the reader the
formula for the inverse.
In other words we have produced a natural bijection between
\[\xymatrix@R-25pt{
\text{functors}& \ca(X,Y)\ar[r]&  \cb\big[G(-)L'(-)\mu_Y^{}F(-)\,,\,
  H(-)\mu_X^{}L(-)E(-)\big]\\
\text{and functors } &
\ca(X,Y)\ar[r] & \cb\big[G(-)\mu_Y^{}R(-)F(-)\,,\,
  H(-)R'(-)\mu_X^{}E(-)\big]
}\]
The compatibility with composition is clearly respected by this bijection.
\eprf

\dfn{D805.7}
Let $\ca$ be a 2-category, $\cb$ a premonoidal 2-category, and
$(-)^*:\ca\la\cb$ a premonoid as in Definition~\ref{D805.3}. For the
sake of definiteness assume the 2-functor $(-)^*$ is contravariant.

A \emph{module} over $(-)^*$, covariant or contravariant, lax or oplax,
is a (possibly lax or oplax) 2-functor
$F:\ca\la\cb$ so that:
\be
\item
  For all objects $X\in\ca$ we have $F(X)=X^*$.
\item
  We must be given a natural transformation
  $a:\big[(-)^*\times F(-)\big]\circ\Delta_2 \la F(-)$,
  possibly lax or oplax,
  which on objects agrees with $\mu:\big[(-)^*\times (-)^*\big]\circ\Delta_2 \la (-)^*$, and so that the square
  \[\xymatrix@C+40pt{
  [(-)^*\times(-)^*\times F(-)]\circ\Delta_3^{}
  \ar[r]^-{\id\times a}\ar[d]_-{\mu\times\id} &
     [(-)^*\times F(-)]\circ\Delta_2^{}\ar[d]^-{a} \\
[(-)^*\times F(-)]\circ\Delta_2^{}\ar[r]^-{a} & F(-)
}\]
  2-commutes. The modification isomorphism giving the 2-commutativity
  of the square above must
  agree on every object $X\in\ca$ with the modification isomorphism
  of the 2-commutative square below
\[\xymatrix@C+40pt{
  [(-)^*\times(-)^*\times(-)^*]\circ\Delta_3^{}
  \ar[r]^-{\id\times\mu}\ar[d]_-{\mu\times\id} &
     [(-)^*\times(-)^*]\circ\Delta_2^{}\ar[d]^-{\mu} \\
[(-)^*\times(-)^*]\circ\Delta_2^{}\ar[r]^-{\mu} & (-)^*
}\]
\ee
\edfn

Perhaps we should elaborate a little on the last assertion. Let $X$ be
an object of $\ca$: then $F(X)=X^*$ by assumption, and furthermore we
assume that $\mu_X^{}:X^*\times X^*\la X^*$ agrees with
$a_X^{}:X^*\times F(X)\la F(X)$. The modifications in the two
2-commutative squares in Definition~\ref{D805.7}(ii) are both
determined by 2-isomorphisms
between the functors $\mu_X^{}(\mu_X^{}\times\id)=a_X^{}(\mu_X^{}\times\id)$
and $\mu_X^{}(\id\times\mu_X^{})=a_X^{}(\id\times a_X^{})$.
The assumption is that these 2-isomorphisms are identical. Since we
assumed that the associativity pentagon commutes for the five composites
$[(-)^*\times(-)^*\times(-)^*\times(-)^*]\circ\Delta_4^{}\la(-)^*$,
it follows that so does the associativity pentagon
for the five composites
$[(-)^*\times(-)^*\times(-)^*\times F(-)]\circ\Delta_4^{}\la F(-)$.

\rmk{R805.99999}
There is an immense literature about modules over what we
call premonoids, we gave a
tiny sample in the paragraphs immediately following
Theorem~\ref{T0.5}~(iv) and (v).
Given a premonoid in $\Hom(\ca,\cb)$, tensor product
with it is a pseudomonad, and our modules are special cases of
pseudo-algebras or lax algebras over it---most
of the category-theoretic literature treats the more general
situation. But what seems really unusual about our definition is the
restriction that on objects of $\ca$ the monoid and the algebra should be
identical.

This restriction
does spare us from having to worry about what happens to units. As
it happens in our main example $X^*=\Dqc(X)$ is a monoidal category (with
a unit), and each $f^*$ respects the unit. As we assume that
$a_X^{}:X^*\times F(X)\la F(X)$ agrees with $\mu_X^{}$, the action
must be unital. 
\ermk

\exm{E805.6}
The trivial example of a module is $F(-)=(-)^*$ with $a=\mu$, that
is the action map $a$ agrees with the
multiplication map $\mu$. For non-trivial examples we will use
\eexm

\pro{P805.11}
Let $\ca$ be a 2-category, let $\cb$ be a premonoidal 2-category, and let
$(-)^*:\ca\la\cb$ be a premonoid as in Definition~\ref{D805.3}.
Suppose $L,R:\ca\la\cb$ are 2-functors as in Reminder~\ref{R805.7},
in particular $L(f)\dashv R(f)$ for all 1-morphisms $f\in\ca$.

If $L(-)$ is contravariant then $L(-)$ is a a lax
module over $(-)^*$ if and only if $R(-)$ is a lax module. If $L(-)$
is covariant then $L(-)$ is an oplax module over $(-)^*$
if and only if $R(-)$ is.
\epro

\prf
Suppose $L(-)$ is contravariant, and
apply Lemma~\ref{L805.9}(i) with $E=(-)^*$, with $F$, $G$ and $H$ the
trivial 2-functors on the empty product,
and with $L=L'$ and $R=R'$ as given. Then
lax natural transformations $a:\big[(-)^*\times L(-) \big]\circ\Delta_2\la L(-)$
are in bijection with lax natural transformations
$b:\big[(-)^*\times R(-) \big]\circ\Delta_2\la R(-)$ which agree on objects.
If we restrict to the case where $L(X)=R(X)=X^*$ for every object
$X\in\ca$, and where the maps $a_X^{}$ and $b_X^{}$ agree with
$\mu_X^{}:X^*\times X^*\la X^*$, then we have that a lax module
structure on $L$ induces a lax module structure on $R$ and vice
versa. Note that
the 2-commutativity of the square in
Definition~\ref{D805.7}(ii) is an assertion of the equality of
certain pairs of natural transformations whose behavior on objects
is determined, and the fact that the map of Lemma~\ref{L805.9}(i)is a bijection tells us that
equality holds for $L$ if and only if it holds for $R$.

The case where $L$ is covariant is similar,
in the proof simply substitute Lemma~\ref{L805.9}(ii)
for Lemma~\ref{L805.9}(i).
\eprf

\exm{E805.13}
As in Example~\ref{E805.5} let $\ca$ be the 2-category of algebraic stacks,
let $\cb$ be the 2-category of finite products of triangulated
categories, and let $(-)^*:\ca\la\Tri\subset\cb$ be the premonoid we
are familiar with. The natural transformation
$\big[(-)^*\oo(-)^*\big]\circ\Delta_2\la(-)^*$ is pseudonatural, meaning it is both
lax and oplax. In Example~\ref{E805.6} we noted that $(-)^*$
is always a (contravariant) module over $(-)^*$, and from
Proposition~\ref{P805.11} we learn that the 2-functor $(-)_*$,
which always exists and is right adjoint to $(-)^*$, has
a canonical structure of a lax module. If the left adjoint
$(-)_!$ exists it is an oplax module.

Concretely: for a pair of objects $X,Y\in\ca$  the natural transformation
$a:\big[(-)^*\times(-)_*\big]\circ\Delta_2\la (-)_*$ gives a map
$a(X,Y):\ca(X,Y)\la\cb\big[\mu_Y^{}\circ\big((-)^*\times (-)_*\big)\,,\,(-)_*\circ\mu_X^{}\big]$. If $f:X\la Y$ is a 1-morphism in $\ca$, and $E\in Y^*=\Dqc(Y)$
and $F\in X_*=\Dqc(X)$ are objects, this comes down to
a map in $\Dqc(Y)$ of the form $E\oo f_*F\la f_*(f^*E\oo F)$,
and this is nothing other than the projection
natural transformation $p(f,E,F)$
of
Notation~\ref{NThom.73.63}(iii). In general it isn't an isomorphism.

However we can shrink the category $\ca$: if instead of allowing all 1-morphisms
we assume that the 1-morphisms in $\ca$ are concentrated then
the projection formula holds, meaning the
lax natural transformation $a$ becomes pseudonatural, and hence
also oplax. Furthermore for any 1-morphism $f:X\la Y$ in $\ca$ the
1-morphism $f_*$ has a right adjoint $f^\times$, and Reminder~\ref{R805.7}
recalls how to extend this to a 2-functor $(-)^\times$.
Proposition~\ref{P805.11} applies and informs us that $(-)^\times$ is an oplax
contravariant module over $(-)^*$. Concretely: let
$f:X\la Y$ be a morphism in $\ca$ and let $E,F\in Y^*=Y^\times=\Dqc(Y)$
be objects. In Notation~\ref{NThom.73.63}(iv)
we saw a morphism $\chi(f,E,F):f^*E\otimes f^\times F\la f^\times (E\oo F)$,
and what we have now learned is that this is the
evaluation of a natural transformation of 2-functors
$\chi:\big[(-)^*\times(-)^\times \big]\circ\Delta_2\la (-)^\times$ at the 1-morphism $f$ and the pair
of objects $E,F$.

In this section we developed the theory in some generality, but the
situation we care about is where $\cb$ is as in Example~\ref{E805.1},
that is its objects are finite products of triangulated categories, and
$\ca=\seq$ is a 2-subcategory of the 2-category of algebraic stacks
as in Notation~\ref{N0.1}. Because the morphisms in $\seq$ are
all concentrated we are in the situation where $(-)^\times$ is an
oplax module over $(-)^*$.
\eexm

\rmd{R204.923}
Before we end this section we study the compatibility of the module
structure on $(-)^\times$ with base-change, and although
the result generalizes we confine our treatment to
the case  we care about, where $\seq$ is as
in Notation~\ref{N0.1}. We begin by reminding the
reader of Remark~\ref{R204.19}: the category $\sq(W,X,Y,Z)$ has for
its
objects the 2-cartesian squares $(\diamondsuit)$ below
\[\xymatrix{
  W \ar[r]^u\ar[d]_f & X\ar[d]^g\\
  Y\ar[r]^v & Z
  }\]
with $v$ flat. There are four natural projections, the obvious
functors
 \[\xymatrix@R+10pt@C-10pt{
  \seq(W,X)& & & \sq(W,X,Y,Z)\ar[rd]^-{\pi_{\text{bot}}^{}}\ar[rrr]^-{\pi_{\text{left}}^{}}\ar[lll]_-{\pi_{\text{top}}^{}}\ar[ld]_-{\pi_{\text{right}}^{}} & & & \seq(W,Y)\\
 &&\seq(X,Z)&&\seq(Y,Z)&& 
  }\]
If $A$ is any of top, bot (for bottom), left or right let $\pi_A^*$ be the
composite
functor $(-)^*\circ\pi_A^{}$ and let $\pi_A^\times$ be the
composite
functor $(-)^\times\circ\pi_A^{}$. In Remark~\ref{R204.19} we observed
that base-change gives a natural transformation 
$\Phi:\pi_{\text{top}}^*\pi_{\text{right}}^\times\la
\pi_{\text{left}}^\times\pi_{\text{bot}}^*$
of functors on the category $\sq(W,X,Y,Z)$. Even more obviously the
2-commutativity gives a natural isomorphism $\tau:\pi_{\text{top}}^*\pi_{\text{right}}^*\la
\pi_{\text{left}}^*\pi_{\text{bot}}^*$. On the other hand we are
given that $(-)^\times$ is an oplax module over the premonoid $(-)^*$, which
means that we have natural 
transformations
\[\xymatrix@C+5pt@R-15pt{
\mu_W^{}\circ[\pi_{\text{top}}^*\times
\pi_{\text{top}}^*]\ar[r]^-{\mu(W,X)} &\pi_{\text{top}}^*\circ\mu_X^{} ,
&\mu_W^{}\circ[\pi_{\text{top}}^*\times
\pi_{\text{top}}^\times]\ar[r]^-{\chi(W,X)} &
\pi_{\text{right}}^\times\circ\mu_X^{},\\
\mu_X^{}\circ[\pi_{\text{right}}^*\times
\pi_{\text{right}}^*]\ar[r]^-{\mu(X,Z)} &\pi_{\text{right}}^*\circ\mu_Z^{} ,
&\mu_X^{}\circ[\pi_{\text{right}}^*\times
\pi_{\text{right}}^\times]\ar[r]^-{\chi(X,Z)} &
\pi_{\text{right}}^\times\circ\mu_Z^{},\\
\mu_W^{}\circ[\pi_{\text{left}}^*\times
\pi_{\text{left}}^*]\ar[r]^-{\mu(W,Y)} &\pi_{\text{left}}^*\circ\mu_Y^{} ,
&\mu_W^{}\circ[\pi_{\text{left}}^*\times
\pi_{\text{left}}^\times]\ar[r]^-{\chi(W,Y)} &
\pi_{\text{left}}^\times\circ\mu_Y^{},\\
\mu_Y^{}\circ[\pi_{\text{bot}}^*\times
\pi_{\text{bot}}^*]\ar[r]^-{\mu(Y,Z)} &\pi_{\text{bot}}^*\circ\mu_Z^{} ,
&\mu_W^{}\circ[\pi_{\text{bot}}^*\times
\pi_{\text{bot}}^\times]\ar[r]^-{\chi(Y,Z)} & \pi_{\text{bot}}^\times\circ\mu_Z^{}.
}\]
With all these natural transformations we could wonder how they might
be related. We prove
\ermd

\pro{P204.1006}
Consider the diagram below, of natural transformations of functors of
the form $\sq(W,X,Y,Z)\la\cb\big[\Dqc(Z)\times\Dqc(Z)\,,\,\Dqc(W)\big]$,
\[\xymatrix@C+40pt{
\mu_W^{}\circ[\pi_{\text{\rm top}}^*\times \pi_{\text{\rm top}}^*]\circ[\pi_{\text{\rm right}}^*\times
\pi_{\text{\rm right}}^\times]\ar[r]^-{\mu_W^{}\circ[\tau\times\Phi]} \ar[d]_-{\mu(W,X)}&\mu_W^{}\circ[\pi_{\text{\rm left}}^*\times
\pi_{\text{\rm left}}^\times]\circ[\pi_{\text{\rm bot}}^*\times
\pi_{\text{\rm bot}}^*]\ar[d]^-{\chi(W,Y)}\\
\pi_{\text{\rm top}}^*\circ\mu_X^{}\circ[\pi_{\text{\rm right}}^*\times
\pi_{\text{\rm right}}^\times] \ar[d]_-{\chi(X,Z)}& \pi_{\text{\rm left}}^\times\circ \mu_Y^{}\circ[\pi_{\text{\rm bot}}^*\times
\pi_{\text{\rm bot}}^*]\ar[d]^-{\mu(Y,Z)}\\
\pi_{\text{\rm top}}^*\circ \pi_{\text{\rm right}}^\times\circ\mu_Z^{} \ar[r]^-{\Phi\circ\mu_Z^{}}& \pi_{\text{\rm left}}^\times\circ \pi_{\text{\rm bot}}^*\circ\mu_Z^{} 
}\]
where the top arrow is an abbreviation for $\mu_W^{}$ applied to the composite
\[\xymatrix@C+40pt{
[\pi_{\text{\rm top}}^*\times
\pi_{\text{\rm top}}^*]\circ[\pi_{\text{\rm right}}^*\times
\pi_{\text{\rm right}}^\times]\ar@{=}[r] & [\pi_{\text{\rm top}}^*\circ \pi_{\text{\rm right}}^* ]\times
[\pi_{\text{\rm top}}^*\circ\pi_{\text{\rm right}}^\times]\ar[d]^-{\tau\times\Phi}\\
[\pi_{\text{\rm left}}^*\times
\pi_{\text{\rm left}}^\times]\circ[\pi_{\text{\rm bot}}^*\times
\pi_{\text{\rm bot}}^*]\ar@{=}[r] &[\pi_{\text{\rm left}}^*\circ \pi_{\text{\rm bot}}^*]\times[\pi_{\text{\rm left}}^\times\circ\pi_{\text{\rm bot}}^*]
}\]
This diagram commutes.
\epro

\prf
The assertion is that the two composites agree, and these are ordinary natural
transformations of 1-functors on 1-categories. Hence it can be checked
object by object: if $(\diamondsuit)$ below
\[\xymatrix{
  W \ar[r]^u\ar[d]_f & X\ar[d]^g\\
  Y\ar[r]^v & Z
}\]
is an object of the category $\sq(W,X,Y,Z)$ we assert the equality
of the evaluation at $(\diamondsuit)$ of two maps, each of which
produces out of $(\diamondsuit)$
a natural transformation between
functors
$\xymatrix{\Dqc(Z)\times\Dqc(Z)\ar@<0.5ex>[r]\ar@<-0.5ex>[r]&
  \Dqc(W).}$
Each of the natural transformtions sends an object
$(E,F)\in\Dqc(Z)\times\Dqc(Z)$ to a morphism in $\Dqc(W)$, and
it suffices to show that, for each pair of objects $(E,F)\in\Dqc(Z)$,
we obtain the same morphism
$u^*g^* E\otimes u^*g^\times F\la
f^\times v^*(E\oo F)$ in the category $\Dqc(W)$.
But this has already
been proved, way back in Lemma~\ref{L27.9}.
\eprf

\section{The module structure of $(-)^!$ over $(-)^*$}
\label{S703}

In \S\ref{S805} we proved that $(-)^*$ is a premonoid and $(-)^\times$ is a
module over it; this turned out to be formal nonsense. Now we
want to generalize to the 2-functor $(-)^!$. We begin
with

\con{C703.1}
Let $X,Y,Z$ be three objects in $\seq$.
The fact that $(-)^*$ is a premonoid
and $(-)^\times$ is a module over it gives, among other things,
natural transformations
\[\xymatrix@C+10pt{
  \mu_X^{}\circ[(-)_1^*\times(-)_1^*]\ar[r]^-{\mu(X,Y)} & (-)_1^*\circ\mu_Y^{},
  &
\mu_Y^{}\circ[(-)_2^*\times(-)_2^\times]\ar[r]^-{\chi(Y,Z)} & (-)_2^\times\circ\mu_Z^{}.
}\]
where $(-)^*_1$ is the restriction of the 2-functor $(-)^*$ to
$\seq(X,Y)$, while $(-)^*_2$ and $(-)^\times_2$ are the restrictions
of $(-)^*$ and $(-)^\times$ to $\seq(Y,Z)$.
We can combine these to form the composite
\[\xymatrix@C+40pt{
  \mu_X^{}\circ[(-)_1^*\times(-)_1^*]\circ[(-)_2^*\times(-)_2^\times]
  \ar[r]^-{\mu(X,Y)} & (-)_1^*\circ\mu_Y^{}\circ[(-)_2^*\times(-)_2^\times]\ar[d]_-{\chi(Y,Z)} \\
  & (-)_1^*\circ(-)_2^\times\circ\mu_Z^{}.
}\]
Now fix the objects $X,Z\in\seq$ and let $a\in\nseq(X,Z)$ be the object
\[\xymatrix@C+30pt@R-15pt{
       & Y\ar^p[dr]& \\
X\ar^u[ur]\ar[rr]|f & &  Z \\
}\]
that is the 2-commutative diagram
with $u$ a dominant, flat monomorphism and $p$ of finite type and
universally quasi-proper. Applying the above composite
to the object $u\in\seq(X,Y)$ and $p\in\seq(Y,Z)$ we deduce a
2-morphism in $\cb$ of the form
\[\CD\chi(Y,Z)\mu(X,Y):\mu_X^{}(u^*p^*\times u^*p^\times)@>>>
u^*p^\times\mu_Y^{}.\endCD\] 
Define $\s_a$ to be the above, that is
$\chi(Y,Z)\mu(X,Y)$ evaluated at the pair $(u,p)$.
We have a recipe that takes an object $a\in\nseq(X,Z)$ to a
morphism $\s_a^{}:\mu_X^{}(a^*\times a^!)\la a^!\mu_Z^{}$.
We assert:
\econ

\lem{L703.3}
The recipe of Construction~\ref{C703.1} defines a natural transformation
of functors of
the form $\nseq(X,Z)\la\cb\big[\Dqc(Z)\times\Dqc(Z),\Dqc(X)\big]$,
more specifically a natural transformation
\[
\CD
\s(X,Z):\mu_X^{}\circ[(-)^*\times(-)^!]@>>> (-)^!\circ\mu_Z^{}\ .
\endCD
\]
\elem

\prf
We have to prove naturality. Let $\ph:a\la b$ be the representative
of a morphism in
$\nseq(X,Z)$; the data include the 2-commutative diagram 
\[\xymatrix@C+30pt@R-20pt{
       & Y\ar^p[dr]\ar[dd]_\alpha & \\
X\ar^u[ur]\ar_{u'}[dr] & &  Z \\
       & Y'\ar_{p'}[ur] & 
}\]
Let $\lambda:p\la p'\alpha$ be the 2-isomorphism giving
the 2-commutativity of the triangle on the right.
Out of the data we construct the 2-commutative diagram $C$ below
  \[
  \xymatrix@C+20pt{
    X\ar[dr]_-\id \ar[r]^-{\id}  &X \ar[r]^{u}\ar@{}[rd]|-{(\diamondsuit)} \ar[d]^{\id} &  Y \ar[d]^\alpha\\
 &   X\ar[dr]_-f \ar[r]^-{u'} & Y'\ar[d]^{p'}\\
 & &  Z
}
\]
and the morphism $\ph^!$ is defined to be $\ph^!=P(C)\lambda^\times$, see
Construction~\ref{C200.7}(ii).
We wish to prove the commutativity of the square
\[\xymatrix@C+30pt{
  \mu_X^{}(a^*\times a^!)\ar[d]_-{\s_a^{}}\ \ar[r]^-{\mu_X^{}(\ph^*\times\ph^!)}
  \ar@{}[dr]|-{(1)}& \mu_X^{}(b^*\times b^!)\ar[d]^-{\s_b^{}} \\
  a^!\mu_Z^{} \ar[r]_-{\ph^!\mu_Z^{}} & b^!\mu_Z^{}
}\]

Let us begin with the case where $p'=\id$, $p=\alpha$ and
$\lambda:p\la\id\circ\alpha$ is the identity.
Then the square $(1)$ reduces to
\[\xymatrix@C+50pt@R+10pt{
  \mu_X^{}[(u^*\alpha^*)\times (u^*\alpha^\times)] \ar[r]^-{\chi(Y,Y')\mu(X,Y)}
  \ar[d]_-{(\diamondsuit)^*\times\Phi(\diamondsuit)}\ar@{}[dr]|-{(2)}
  & u^*\alpha^\times\mu_{Y'}^{} \ar[d]^-{\Phi(\diamondsuit)}\\
\mu_X^{}[{u'}^*\times{u'}^*]\ar[r]_-{\mu(X,Y')} & 
  {u'}^*\mu_{Y'}^{}
}\]
and it commutes by Proposition~\ref{P204.1006} applied to the cartesian
square $(\diamondsuit)$. The idea of the rest of the proof
is that the general case follows from the commutativity of $(2)$ by
applying very minor adjustments. To spell it out: if we precompose
everything in $(2)$ with $[{p'}^*\times {p'}^\times]$ it gives the
commutative square
\[\xymatrix@C+70pt@R+10pt{
  \mu_X^{}[(u^*\alpha^*)\times (u^*\alpha^\times)][{p'}^*\times {p'}^\times] \ar[r]^-{\chi(Y,Y')\mu(X,Y)}
  \ar[d]_-{(\diamondsuit)^*\times\Phi(\diamondsuit)}\ar@{}[dr]|-{(2')}
  & u^*\alpha^\times\mu_{Y'}^{}[{p'}^*\times {p'}^\times] \ar[d]^-{\Phi(\diamondsuit)}\\
\mu_X^{}[{u'}^*\times{u'}^*][{p'}^*\times {p'}^\times]\ar[r]_-{\mu(X,Y')} & 
  {u'}^*\mu_{Y'}^{}[{p'}^*\times {p'}^\times]
}\]
and if we compose the horizontal maps, on the right, with the map
$\chi(Y',Z):\mu_{Y'}^{}[{p'}^*\times{p'}^\times]\la {p'}^\times\mu_Z^{}$, 
then we obtain the commutative square
\[\xymatrix@C+90pt@R+10pt{
  \mu_X^{}[(u^*\alpha^*)\times (u^*\alpha^\times)][{p'}^*\times {p'}^\times] \ar[r]^-{\chi(Y',Z)\chi(Y,Y')\mu(X,Y)}
  \ar[d]_-{(\diamondsuit)^*\times\Phi(\diamondsuit)}\ar@{}[dr]|-{(3)}
  & u^*\alpha^\times{p'}^\times\mu_{Z}^{}\ar[d]^-{\Phi(\diamondsuit)}\\
\mu_X^{}[{u'}^*\times{u'}^*][{p'}^*\times {p'}^\times]\ar[r]_-{\chi(Y',Z)\mu(X,Y')} & 
  {u'}^*{p'}^\times\mu_{Z}^{}
}\]
But 
because $\chi:\big[(-)^*\times(-)^\times \big]\circ\Delta_2\la(-)^\times$ is an oplax natural
transformation of 2-functors it respects composition, and
the square below commutes
\[\xymatrix@C+50pt@R+10pt{
  \mu_Y^{}[p^*\times p^\times] \ar[r]^-{\chi(Y,Z)}
  \ar[d]_-{\lambda^*\times\lambda^\times}\ar@{}[dr]|-{(4)}
  & p^\times\mu_Z^{} \ar[d]^-{\lambda^\times}\\
\mu_Y^{}[(\alpha^*{p'}^*)\times(\alpha^\times{p'}^\times)]\ar[r]_-{\chi(Y',Z)\chi(Y,Y')} & 
  \alpha^\times{p'}^\times\mu_Z^{}
}\]
Applying the functor $u^*$ gives the commutative square
  \[\xymatrix@C+50pt@R+10pt{
  u^*\mu_Y^{}[p^*\times p^\times] \ar[r]^-{\chi(Y,Z)}
  \ar[d]_-{\lambda^*\times\lambda^\times}\ar@{}[dr]|-{(4')}
  & u^*p^\times\mu_Z^{} \ar[d]^-{\lambda^\times}\\
u^*\mu_Y^{}[(\alpha^*{p'}^*)\times(\alpha^\times{p'}^\times)]\ar[r]_-{\chi(Y',Z)\chi(Y,Y')} & 
  u^*\alpha^\times{p'}^\times\mu_Z^{}
}\]
  Composing the horizontal maps, on the left, with
  the map $\mu(X,Y):\mu_X^{}[u^*\times u^*]\la u^*\mu_Y^{}$
gives
the commutative square
  \[\xymatrix@C+90pt@R+10pt{
  \mu_X^{}[u^*\times u^*][p^*\times p^\times] \ar[r]^-{\chi(Y,Z)\mu(X,Y)}
  \ar[d]_-{\lambda^*\times\lambda^\times}\ar@{}[dr]|-{(5)}
  & u^*p^\times\mu_Z^{} \ar[d]^-{\lambda^\times}\\
\mu_X^{}[u^*\times u^*][(\alpha^*{p'}^*)\times(\alpha^\times{p'}^\times)]\ar[r]_-{\chi(Y',Z)\chi(Y,Y')\mu(X,Y)} & 
  u^*\alpha^\times{p'}^\times\mu_Z^{}
}\]
And the commutative squares (3) and (5) concatenate to (1).
\eprf

\rmk{R703.5}
Because the functors $(-)^*$ and $(-)^!$ both factor through the
groupoid completion functor $F:\nseq(X,Z)\la\seq(X,Z)$, the natural
transformation
$\s(X,Z)$ factors uniquely through $F$. We have defined on $\seq(X,Z)$
a natural transformation
$\s(X,Z):\mu_X^{}\circ\big[(-)^*\times(-)^!\big]\la(-)^!\circ\mu_Z^{}$. We assert:
\ermk

\thm{T703.7}
The following defines an oplax natural transformation
$\s:\big[(-)^*\times(-)^!\big]\circ\Delta_2\la(-)^!$ of 2-functors on $\seq$
\be
\item
On objects: for an object $X\in\seq$ we define $\s_X^{}:X^*\times
X^!\la X$ to be $\mu_X^{}:\Dqc(X)\times\Dqc(X)\la\Dqc(X)$.
\item
On morphisms: for two objects $X,Z\in\seq$ we define the natural
transformation 
$\s(X,Z):\mu_X^{}\circ\big[(-)^*\times(-)^!\big]\la (-)^!\circ\mu_Z^{}$ as in
Lemma~\ref{L703.3} and Remark~\ref{R703.5}.
\ee
The natural transformation $\s$ makes $(-)^!$ into an oplax module over $(-)^*$.
\ethm

\prf
Let us begin by showing that $\s$ is a natural transformation of
2-functors: for this we need to show that it respects composition.
If $X\stackrel f\la Y\stackrel g\la Z$ are composable morphisms in
$\seq$ we need to show the commutativity of the pentagon
\[\xymatrix@C-25pt@R+10pt{
& \mu_X^{}\big[(gf)^*\times (gf)^!\big]
\ar[rr]^-{\s(X,Z)}\ar[ld]_-{\mu_X^{}\circ\big[\tau(f,g)\times\rho(f,g)\big]} 
\ar@{}[rrd]|-{\ds\text{\bf Pent}}& &
(gf)^!\mu_Z^{} \ar[rrd]^-{\rho(f,g)\circ\mu_Z^{}}\\
 \mu_X^{}[f^*g^*\times f^!g^!] \ar[rr]^-{\s(X,Y)} & &
 f^!\mu_Y^{}[g^*\times g^!] \ar[rrr]^-{\s(Y,Z)} & & & f^!g^!\mu_Z^{}
}\]
and we begin with easy cases.

\medskip

\nin
{\bf Case 1.}\ \ 
\emph{The pentagon {\bf Pent} commutes if $f$ is a dominant, flat monomorphism
or if $g$ is of finite type and universally quasi-proper.}

\medskip

\nin
\emph{Proof of Case 1.}\ \ We prove the case where $g$ is of
finite type and universally
quasi-proper, and leave the other case to the reader. To compute the
maps involved it helps to choose a preimage in
$\nseq(X,Y,Z)$ of the object $X\stackrel f\la Y\stackrel g\la Z$
in the category $\seq(X,Y)\times\seq(Y,Z)$, and our choice of
preimage is the 2-commutative diagram $C$ below
\[
  \xymatrix{
X\ar[dr]_-f \ar[r]^-{u}  &X' \ar[r]^{\id} \ar[d]^{p} &  X'\ar[d]^p\\
 &   Y\ar[dr]_-g \ar[r]^-\id & Y\ar[d]^g\\
 & &  Z
}
\]
The map $\rho(f,g)$ was defined to be
$P(C):\pi_{13}^{}(C)^!\la\pi_{12}^{}(C)^!\pi_{23}^{}(C)^!$,
which for our choice of $C$ comes down to
the canonical isomorphism $u^*\delta(p,g):u^*(gp)^\times\la u^*p^\times g^\times$.
Thus {\bf Pent} is the concatenation of the
square
\[\xymatrix@C+20pt{
  \mu_X^{}[u^*\times u^*]\big[(gp)^*\times (gp)^\times\big] \ar[r]^-{\mu(X,X')}
  \ar[d]_-{\mu_X^{}[u^*\times u^*]\big[\tau(p,g)\times\delta(p,g)\big]} \ar@{}[dr]|-{\ds(1)}&
u^*\mu_{X'}^{}\big[(gp)^*\times (gp)^\times\big]\ar[d]^-{u^* \mu_{X'}^{}\big[\tau(p,g)\times\delta(p,g)\big]}\\
\mu_X^{}[u^*\times u^*][p^*g^*\times p^\times g^\times] \ar[r]^-{\mu(X,X')} &
u^*\mu_{X'}^{}[p^*g^*\times p^\times g^\times] 
}\]
and the pentagon
\[\xymatrix@C-25pt@R+10pt{
&u^* \mu_{X'}^{}\big[(gp)^*\times (gp)^\times\big]\ar@{}[drr]|-{\ds(2)}
\ar[rr]^-{\s(X',Z)}\ar[ld]_-{u^*\mu_{X'}^{}\circ\big[\tau(p,g)\times\delta(p,g)\big]} 
& &
u^*(gp)^\times\mu_Z^{} \ar[rrd]^-{u^*\delta(p,g)\circ\mu_Z^{}}\\
u^* \mu_{X'}^{}[p^*g^*\times p^\times g^\times] \ar[rr]^-{\chi(X',Y)} & &
 u^*p^\times\mu_Y^{}[g^*\times g^\times] \ar[rrr]^-{\chi(Y,Z)} & & &
u^* p^\times g^\times\mu_Z^{}
}\]
The square $(1)$ commutes trivially, while $(2)$ is what we get
by applying $u^*$ 
to the commutative pentagon which spells
out the fact that the oplax natural transformation
$\chi:\big[(-)^*\times(-)^\times\big]\circ\Delta_2\la(-)^\times$
respects composition.

\medskip

\nin
{\bf Case 2.}\ \ 
\emph{The pentagon {\bf Pent} commutes if $g$ is a dominant, flat monomorphism
and $f$ is of finite type and universally quasi-proper.}

\medskip

\nin
\emph{Proof of Case 2.}\ \ 
We are given the composable maps $X\stackrel f\la Y\stackrel g\la Z$,
and choose a Nagata compactification for the composite $gf:X\la Z$.
This produces for us the 2-commutative diagram
 $C$ below
\[
  \xymatrix{
X\ar[dr]_-f \ar[r]^-{\id}  &X \ar[r]^{u} \ar[d]^{f}\ar@{}[dr]|-{(\diamondsuit)} & R\ar[d]^p\\
 &   Y\ar[dr]_-g \ar[r]^-g & Z\ar[d]^\id\\
 & &  Z
}
\]
The map $f$ is assumed of finite type and universally quasi-proper, as
is the map $p$ by construction. The map $g$ is assumed to be a dominant, flat
monomorphism, as is $u$ (by construction). Hence
$C$ is an object of the category $\nseq(X,Y,Z)$ lifting
the object $X\stackrel f\la Y\stackrel g\la Z$
via the functor $F:\nseq(X,Y,Z)\la\seq(X,Y)\times\seq(Y,Z)$.
The map $\rho(f,g)$ was defined to be
$P(C):\pi_{13}^{}(C)^!\la\pi_{12}^{}(C)^!\pi_{23}^{}(C)^!$,
and for our $C$ this comes down to
$\Phi(\diamondsuit):u^*p^\times\la f^\times g^*$.
The commutativity of {\bf Pent} now follows by evaluating the
commutative diagram of Proposition~\ref{P204.1006} at the object
$(\diamondsuit)\in\sq(X,R,Y,Z)$.

\medskip

\nin
\emph{Proof of the commutativity of {\bf Pent} in the general case.}\ \ 
Suppose $f$ and $g$ are general, and choose Nagata compactifications
for each of them. That is factor $f$ and $g$ as $f=pu$ and $g=qv$, with
$u:X\la R$, $v:Y\la S$ dominant, flat monomorphisms and $p:R\la Y$,
$q:S\la Z$ of finite type and universally
quasi-proper. Now consider the diagram
\[\xymatrix@C+20pt{
    & & \mu_X^{}[u^*p^*v^*q^*\times u^!p^!v^!q^!]\ar[dd]^{\s(X,R)}\\
  &  \ar@{}[dr]|-{\ds(2)}
  \mu_X^{}[f^*g^*\times f^!g^!]\ar[dd]_-{\s(X,Y)}\ar[ur]^-{\mu_X^{}[\tau\times\rho]}&\\
  \mu_X^{}[(gf)^*\times(gf)^!]\ar[ur]^-{\mu_X^{}[\tau\times\rho]}\ar@{}[ddr]|-{\ds(1)} \ar[dd]_-{\s(X,Z)}&& u^!\mu_R^{}[p^*v^*q^*\times p^!v^!q^!]\ar[d]^{\s(R,Y)} \\
  &f^!\mu_Y^{}[g^*\times g^!]\ar[dd]_{\s(Y,Z)}\ar[r]^{\rho\mu_Y^{}(\tau\times\rho]}&u^!p^!\mu_Y^{}[v^*q^*\times v^!q^!]\ar[d]^{\s(Y,S)}\\
  (gf)^!\mu_Z^{} \ar[dr]_-{\rho\mu_Z^{}} & &u^!p^!v^!\mu_S^{}[q^*\times q^!]\ar[dd]^{\s(S,Z)} \\
& f^!g^!\mu_Z^{}\ar[dr]_{\rho\rho\mu_Z^{}}\ar@{}[ur]|-{\ds(3)}  &\\
 & & u^!p^!v^!q^!\mu_Z^{}
}\]
We wish to prove the commutativity of the pentagon $(1)$. Case 1 gives
the commutativity of $(2)$ and $(3)$, and the maps $\rho(u,p):f^!\la u^!p^!$
and $\rho(v,q):g^!\la v^!q^!$ are both isomorphisms by
Proposition~\ref{P200.23}~(i) or (ii). Thus the map labeled $\rho\rho\mu_Z^{}$
at the bottom of the diagram is an isomorphism, and it suffices to prove
the commutativity of the perimeter; for future reference we call
it {\bf Perim}.

To this end consider the diagram
\[\xymatrix@C-20pt@R+20pt{
  & & & \mu_X^{}[u^*p^*v^*q^*\times u^!p^!v^!q^!]\ar[dd]^{\s(X,R)}\\
    & & &\\
  &  \ar@{}[rr]|-{\ds(7)}
  \mu_X^{}[u^*(gp)^*\times u^!(gp)^!]\ar[dd]_-{\s(X,R)}\ar[uurr]^-{\mu_X^{}[\tau\times\rho]}& & u^!\mu_R^{}[p^*v^*q^*\times p^!v^!q^!]\ar[dd]^-{\s(R,Y)}\\
  \mu_X^{}[(gpu)^*\times(gpu)^!]\ar[ur]^-{\mu_X^{}[\tau\times\rho]}\ar@{}[dddr]|-{\ds(4)} \ar[dd]_-{\s(X,Z)}&& u^!\mu_R^{}[(vp)^*q^*\times (vp)^!q^!]\ar[dd]^{\s(R,S)}\ar[ur]^{u^!\mu_R^{}[\tau\times\rho]} &\\
  &u^!\mu_R^{}[(gp)^*\times (gp)^!]\ar[dd]_{\s(R,Z)}\ar[ur]^{u^!\mu_R^{}(\tau\times\rho]}&& u^!p^!\mu_Y^{}[v^*q^*\times v^!q^!]\ar[dd]^-{\s(Y,S)}\\
  (gpu)^!\mu_Z^{} \ar[dr]_-{\rho\mu_Z^{}} & &u^!(vp)^!\mu_S^{}[q^*\times q^!]\ar[dd]^{\s(S,Z)}\ar[rd]_\rho\ar@{}[ur]|-{\ds(6)} & \\
& u^!(gp)^!\mu_Z^{}\ar[dr]_{\rho\rho\mu_Z^{}}\ar@{}[ur]|-{\ds(5)}  & &u^!p^!v^!\mu_S^{}[q^*\times q^!]\ar[dd]^-{\s(S,Z)}\\
& & u^!(vp)^!q^!\mu_Z^{}\ar[dr]_-{\rho}\ar@{}[ur]|-{\ds(8)}&\\
& & & u^!p^!v^!q^!\mu_Z^{}
}\]
The perimeter agrees with {\bf Perim}, the pentagons $(4)$ and $(5)$ commute by
Case 1, the pentagon $(6)$ commutes by Case 2, and the regions $(7)$ and $(8)$
commute trivially. We have finished the proof that {\bf Pent} commutes, that is
we now know that $\s:\big[(-)^*\times(-)^!\big]\circ\Delta_2\la(-)^!$ is an
oplax
natural transformation.

To finish the proof of the theorem we need to establish the associativity in
Definition~\ref{D805.7}(ii). We should perhaps remind the reader.

The fact that $(-)^*$ is a premonoid gives a pseudonatural transformation
$\mu:\big[(-)^*\times(-)^*\big]\circ\Delta_2\la(-)^*$, and a modification
isomorphism $m:\mu(\mu\times\id)\la\mu(\id\times\mu)$. For every
object $X$ this gives a functor $\mu_X^{}:X^*\times X^*\la X^*$, and a natural
isomorphism $m_X^{}:\mu_X^{}(\mu_X^{}\times\id)\la\mu_X^{}(\id\times\mu_X^{})$.
The 2-functor $(-)^!$
is such that for all objects $X$ we have $X^!=X^*$,
the natural transformation $\s:\big[(-)^*\times(-)^!\big]\circ\Delta_2\la(-)^!$
is such that for every object $X$ the map $\s_X^{}:X^*\times X^!\la X^!$
satisfies $\s_X^{}=\mu_X^{}$, and the definition tells us
that, in order for $(-)^!$ to be a module over $(-)^*$, if we view
$m_X^{}:\mu_X^{}(\mu_X^{}\times\id)\la\mu_X^{}(\id\times\mu_X^{})$
as a natural transformation
$m_X^{}:\s_X^{}(\mu_X^{}\times\id)\la\s_X^{}(\id\times\s_X^{})$,
then it extends to a modification. We have no choice in what
the modification does on objects, what needs checking is that this
is compatible with rest of the structure. Concretely: if $X,Z\in\seq$ are
objects and $f:X\la Z$ is a 1-morphism then we have the following square
\[\xymatrix@C+20pt{
  \mu_X^{}\big[\mu_X^{}(f^*\times f^*)\times f^!\big] \ar[r]^-{\mu(X,Z)}
  \ar[d]_-{m_X^{}}&
  \mu_X^{}\big[f^*\mu_Z^{}\times f^!\big]\ar[r]^-{\s(X,Z)} &
  f^!\mu_Z^{}[\mu_Z^{}\times\id]\ar[d]^-{m_Z^{}}\\
   \mu_X^{}\big[f^*\times\mu_X^{}(f^*\times f^!)\big] \ar[r]^-{\s(X,Z)}&
  \mu_X^{}\big[f^*\times f^!\mu_Z^{}\big]\ar[r]^-{\s(X,Z)} &
  f^!\mu_Z^{}[\id\times\mu_Z^{}]
}\]
and we need to prove it commutative for every 1-morphism $f$. To do this
recall the definition of $\s(X,Z)$. First choose a lifting
$a$ of $f$ via the functor $F:\nseq(X,Z)\la\seq(X,Z)$, that
is choose an object
$a\in\nseq(X,Z)$ with
$F(a)=f$. Concretely $a$ is a Nagata compactification for $f$;
it is a 2-isomorphism $pu\la f$, where $X\stackrel u\la Y\stackrel p\la Z$
are composable 1-morphisms 
with $u$ a dominant, flat monomorphism and $p$ of
finite type and universally quasi-proper. The definition in
Construction~\ref{C703.1} is that $\s_a^{}$ is the composite
\[\xymatrix@C+20pt{
  \mu_X^{}\big[u^*p^*\times u^*p^\times\big] \ar[r]^-{\mu(X,Y)}
  & u^*\mu_Y^{}\big[p^*\times p^\times\big]\ar[r]^-{\chi(Y,Z)}
  & u^*p^\times\mu_Z^{}\ .
}\]
The Theorem now follows from the commutativity of
\[\xymatrix@C+20pt{
  \mu_X^{}\big[\mu_X^{}(u^*\times u^*)\times u^*\big] \ar[r]^-{\mu(X,Y)}
  \ar[d]_-{m_X^{}}&
  \mu_X^{}\big[u^*\mu_Y^{}\times u^*\big]\ar[r]^-{\mu(X,Y)} &
  u^*\mu_Y^{}[\mu_Y^{}\times\id]\ar[d]^-{m_Y^{}}\\
   \mu_X^{}\big[u^*\times\mu_X^{}(u^*\times u^*)\big] \ar[r]^-{\mu(X,Y)}&
  \mu_X^{}\big[u^*\times u^*\mu_Y^{}\big]\ar[r]^-{\mu(X,Y)} &
  u^*\mu_Y^{}[\id\times\mu_Y^{}]
}\]
which is a consequence of the associativity of the multiplication
$\mu:\big[(-)^*\times(-)^*\big]\circ\Delta_2\la(-)^*$, coupled with
the commutativity of
\[\xymatrix@C+20pt{
  \mu_Y^{}\big[\mu_Y^{}(p^*\times p^*)\times p^\times\big] \ar[r]^-{\mu(Y,Z)}
  \ar[d]_-{m_Y^{}}&
  \mu_Y^{}\big[p^*\mu_Z^{}\times p^\times\big]\ar[r]^-{\chi(Y,Z)} &
  p^\times\mu_Z^{}[\mu_Z^{}\times\id]\ar[d]^-{m_Z^{}}\\
   \mu_Y^{}\big[p^*\times\mu_Y^{}(p^*\times p^\times)\big] \ar[r]^-{\chi(Y,Z)}&
  \mu_Y^{}\big[p^*\times p^\times\mu_Z^{}\big]\ar[r]^-{\chi(Y,Z)} &
  p^\times\mu_Z^{}[\id\times\mu_Z^{}]
}\]
which comes because the action
$\chi:\big[(-)^*\times(-)^\times\big]\circ\Delta_2\la(-)^\times$ is
also associative.
\eprf

\pro{P703.9}
The oplax natural transformation $\psi:(-)^\times\la(-)^!$ respects
the action of $(-)^*$. In other words $\psi$ is a homomorphism of
$(-)^*$--modules. Even more precisely: we have a square of oplax
natural transformation of oplax 2-functors
\[\xymatrix@C+40pt{
  \big[(-)^*\times(-)^\times\big]\circ\Delta_2 \ar[r]^-{\chi}
  \ar[d]_-{[\id\times\psi]\circ\Delta_2}\ar@{}[dr]|-{\ds(1)}& (-)^\times\ar[d]^-{\psi} \\
  \big[(-)^*\times(-)^!\big]\circ\Delta_2 \ar[r]^-{\s} & (-)^!
}\]
On objects, both $\chi$ and $\s$ send the object $X\in\seq$ to the functor
$\mu_X^{}:\Dqc(X)\times\Dqc(X)\la\Dqc(X)$, and $\psi$ takes the object $X$ to
$\id:\Dqc(X)\la\Dqc(X)$. The composite natural transformations agree on objects,
both composites take $X\in\seq$ to the functor
$\mu_X^{}:\Dqc(X)\times\Dqc(X)\la\Dqc(X)$.
The assertion is that the assignment taking $X$ to the identity 2-morphism
$\id:\mu_X^{}\la\mu_X^{}$ extends to a modification isomorphism of the
composites. More concretely: the two composites are actually equal,
not just 2-isomorphic.
\epro

\prf
We have studied what happens on objects and it remains to verify the
assertion on 1-morphisms. Any 1-morphism $f:X\la Z$ in $\seq$
can be factored as
$f\cong pu$
with $u$ a dominant, flat monomorphism and $p$ of finite type and
universally quasi-proper, and the map $\rho(u,p):f^!\la u^!p^!$
is an isomorphism by Proposition~\ref{P200.23}~(i) or (ii).
Evaluating the square $(1)$ at $f\in\seq$ produces 
two 2-morphisms
$\xymatrix{
  x\ar@<0.5ex>[r] \ar@<-0.5ex>[r] & f^!\mu_Z^{}
}$
which we wish to show equal, 
 and it suffices to
 prove that the composites
$\xymatrix{
  x\ar@<0.5ex>[r] \ar@<-0.5ex>[r] & f^!\ar[r]^-{\rho(u,p)} & u^!p^!\mu_Z^{}
}$
are equal. But natural transformations respect composition,
allowing us to rewrite this composite, and 
we are reduced to proving the commutativity of $(1)$ when evaluated at $u$ and at $p$.
Thus it suffices to consider two special cases.
\medskip

\nin
{\bf Case 1.}\ \ 
\emph{The two composites in the square of the Proposition agree on 1-morphisms
$f$ which are of finite type and universally quasi-proper.}

\medskip

\nin
\emph{Proof of Case 1.}\ \ We begin by choosing 
the following
object $a\in\nseq(X,Z)$ with $F(a)=f$
\[
  \xymatrix@C+30pt@R-20pt{
    & X\ar[dr]^-{f} & \\
    X\ar[rr]_-{f}\ar[ur]^-{\id}  & & Z 
}\]
For this object $\psi(a):f^\times\la f^!$ is the identity map, and
$\s(X,Z):\mu_X^{}\big[f^*\times f^!\big]\la f^!\mu_Z^{}$ agrees with
$\chi(X,Z):\mu_X^{}\big[f^*\times f^\times\big]\la f^\times\mu_Z^{}$.
Therefore the agreement
of the two composites at $f$ is trivial.

\medskip

\nin
{\bf Case 2.}\ \ 
\emph{The two composites in the square of the Proposition agree on 1-morphisms
$f$ which are dominant, flat monomorphisms.}

\medskip

\nin
\emph{Proof of Case 2.}\ \ Once again we begin by choosing 
an object $a\in\nseq(X,Z)$ with $F(a)=f$, our choice is
\[
  \xymatrix@C+30pt@R-20pt{
    & Z\ar[dr]^-{\id} & \\
    X\ar[rr]_-{f}\ar[ur]^-{f}  & & Z 
}\]
This time $f^!=f^*$, the map
$\s(X,Z):\mu_X^{}\big[f^*\times f^!\big]\la f^!\mu_Z^{}$
agrees with  $\mu(X,Z):\mu_X^{}\big[f^*\times f^*\big]\la f^*\mu_Z^{}$, and
the map $\psi(f):f^\times\la f^*$ identifies, by
Construction~\ref{C202.1},
with the base-change map
$\Phi(\diamondsuit):f^\times\la f^*$,
where $(\diamondsuit)$ is the 2-cartesian square
\[\xymatrix{
  X \ar[r]^\id \ar[d]_\id & X\ar[d]^f \\
  X\ar[r]^f & Z
}\]
We wish to prove the commutativity of
\[\xymatrix@C+40pt{
  \mu_X^{}\big[f^*\times f^\times\big]\ar[r]^-{\chi(X,Z)}
  \ar[d]_-{\mu_X^{}[\id\times\Phi(\diamondsuit)]}& f^\times\mu_Z^{}\ar[d]^-{\Phi(\diamondsuit)} \\
   \mu_X^{}\big[f^*\times f^*\big]\ar[r]^-{\mu(X,Z)}
  & f^*\mu_Z^{}
}\]
and the proof is by evaluating the commutative diagram of
Proposition~\ref{P204.1006} on the object $(\diamondsuit)\in\sq(X,X,X,Z)$.
\eprf

\rmk{R703.10}
We remind the reader of the 2-category $\vseq$: the objects are flat morphisms
$u:W\la X$ in $\seq$, the 1-morphisms $u\la v$ are 2-cartesian squares 
\[\xymatrix{
  W \ar[r]^u\ar[d]_f\ar@{}[dr]|-{(\diamondsuit)} & X\ar[d]^g\\
  Y\ar[r]^v & Z
  }\]
and the 2-morphisms are isomorphisms of 2-cartesian squares which are identities
on $u$ and $v$. The functors $p_1^{},p_2^{}:\vseq\la\seq$ take an object
$u:W\la X$ to $p_1^{}(u)=W$ and $p_2^{}(u)=X$. In Theorem~\ref{T204.17}
we produced an oplax natural transformation
$\theta:(-)^!\circ p_2^{}\la(-)^!\circ p_1^{}$, and in Remark~\ref{R204.295}
we noted the(easier)  oplax natural transformation
$\Phi:(-)^\times\circ p_2^{}\la(-)^\times\circ p_1^{}$. Although we
haven't mentioned it yet, there is also the trivial
pseudonatural transformation
$\tau:(-)^*\circ p_2^{}\la(-)^*\circ p_1^{}$. We should remind the reader.

The three 2-functors $(-)^*\circ  p_1^{}$, $(-)^\times\circ  p_1^{}$ and
$(-)^!\circ  p_1^{}$ all take the object $u:W\la X$ of the category
$\vseq$ to the object $\Dqc(W)\in\Tri$. The three 2-functors
$(-)^*\circ  p_2^{}$, $(-)^\times\circ  p_2^{}$ and
$(-)^!\circ  p_2^{}$ all take the object $u:W\la X$ to $\Dqc(X)\in\Tri$.
When we evaluate at the 1-morphism $(\diamondsuit)$ above the functors
$(-)^*\circ  p_1^{}$, $(-)^\times\circ  p_1^{}$ and
$(-)^!\circ  p_1^{}$ take $(\diamondsuit)$ (respectively) to $f^*$, $f^\times$
and $f^!$, while the functors $(-)^*\circ  p_2^{}$, $(-)^\times\circ  p_2^{}$ and
$(-)^!\circ  p_2^{}$  take $(\diamondsuit)$ (respectively) to $g^*$, $g^\times$
and $g^!$.

So much for the functors. The three natural transformations $\tau$,
$\Phi$ and $\theta$ all take the object $u:W\la X$ to the 1-morphism
$u^*:\Dqc(X)\la\Dqc(W)$ in $\Tri$. Each 1-morphism $(\diamondsuit)$
in $\vseq$ must map to a 2-morphism, and the rule is
$\tau(\diamondsuit):u^*g^*\la f^*v^*$ is the canonical isomorphism,
$\Phi(\diamondsuit):u^*g^\times\la f^\times v^*$ is the base-change map,
and $\theta(\diamondsuit):u^*g^!\la f^!v^*$ is the 2-morphism
of Construction~\ref{C204.11} and Proposition~\ref{P204.13}. Of these
$\Phi$ and $\theta$ have a direction, both are oplax, while $\tau(\diamondsuit)$
is an isomorphism, making $\tau$ a pseudonatural transformation.

These are 2-functors and natural transformations we knew about
back in
\S\ref{S8}, but now that we have the premonoid structure on $(-)^*$
and the fact that $(-)^\times$ and $(-)^!$ are both
$(-)^*$--modules we have many more natural transformations
to play with. We can wonder whether some composites might be equal,
or at least have modifications mapping one to the other. The next Proposition
proves such an assertion, but to state it we introduce a tiny bit of
notation. Until now the map $\Delta_2$ was a strict 2-functor
$\Delta_2:\seq\la\seq\times\seq$, the doubling map. To state the next
Proposition we commit the notational crime of allowing
$\Delta_2:\vseq\la\vseq\times\vseq$ to also be the doubling map. With this
criminal notation we observe that the following squares are strictly
commutative
\[\xymatrix{
\vseq \ar[r]^-{\Delta_2}\ar[d]_-{p_1^{}} & \vseq\times\vseq\ar[d]^-{p_1^{}\times p_1^{}} & & \vseq \ar[r]^-{\Delta_2}\ar[d]_-{p_2^{}} & \vseq\times\vseq\ar[d]^-{p_2^{}\times p_2^{}}\\
\seq \ar[r]^-{\Delta_2} & \seq\times\seq & & \seq \ar[r]^-{\Delta_2} & \seq\times\seq      
}\]
With this notation we prove
\ermk

\pro{P703.10.5}
The following diagrams 2-commute
\[\xymatrix@C+20pt{
  \Big[\big((-)^*\circ p_2^{}\big)\times\big((-)^\times\circ p_2^{}\big)\Big]\circ\Delta_2\ar@{=}[d] \ar[r]^-{\tau\times\Phi}&  \Big[\big((-)^*\circ p_1^{}\big)\times\big((-)^\times\circ p_1^{}\big)\Big]\circ\Delta_2\ar@{=}[d] \\
  \big[(-)^*\times(-)^\times\big]\circ\Delta_2\circ p_2^{}\ar[d]_{\chi}
  \ar@{}[r]|-{\ds(1)}
  & \big[(-)^*\times(-)^\times\big]\circ\Delta_2\circ p_1^{}\ar[d]_{\chi}\\
  (-)^\times\circ p_2^{}\ar[r]^-{\Phi} &
    (-)^\times\circ p_1^{}
}\]
and
\[\xymatrix@C+20pt{
  \Big[\big((-)^*\circ p_2^{}\big)\times\big((-)^!\circ p_2^{}\big)\Big]\circ\Delta_2\ar@{=}[d] \ar[r]^-{\tau\times\theta}&  \Big[\big((-)^*\circ p_1^{}\big)\times\big((-)^!\circ p_1^{}\big)\Big]\circ\Delta_2\ar@{=}[d] \\
  \big[(-)^*\times(-)^!\big]\circ\Delta_2\circ p_2^{}\ar[d]_{\s}
\ar@{}[r]|-{\ds(2)}
  & \big[(-)^*\times(-)^!\big]\circ\Delta_2\circ p_1^{}\ar[d]_{\s}\\
  (-)^!\circ p_2^{}\ar[r]^-{\theta} &
    (-)^!\circ p_1^{}
}\]
More explicitly: we have four natural transformations, two composites
in each diagram. If we trace what happens to an object $u:W\la X$
in the category $\vseq$, the two composites shaped like $\urcorner$
take $u$ to $\mu_W^{}(u^*\times u^*)$, while the two composites
of the shape $\llcorner$ take $u$ to $u^*\mu_X^{}$.
The fact that $(-)^*$ is a premonoid gives us a 2-isomorphism
$\mu(W,X):\mu_W^{}(u^*\times u^*)\la u^*\mu_X^{}$, and we
assert that this extends to a modification that works in both
diagrams.
\epro

\prf
We've analyzed what happens to an object $u\in\vseq$, it remains to
check something on 1-morphisms $(\diamondsuit)\in\vseq$. For the
diagram $(1)$ the proof is by evaluating the
commutative diagram of 
Proposition~\ref{P204.1006} at the object $(\diamondsuit)\in\sq(W,X,Y,Z)$.

For the diagram $(2)$ we actually need to
prove something.
Let $(\diamondsuit)\in\vseq(u,v)$ be a 1-morphism, we need to show
that the diagram commutes when evaluated at $(\diamondsuit)$.
If $(\diamondsuit)$ is the 2-cartesian square 
\[\xymatrix{
  W \ar[r]^u\ar[d]_f & X\ar[d]^g\\
  Y\ar[r]^v & Z
}\]
we may choose a Nagata compactification
$X\stackrel {t'}\la S\stackrel{p'}\la Z$
for the map $g:X\la Z$,
and pull back to form the diagram where
the squares are 2-cartesian 
\[\xymatrix{
  W \ar[r]^u\ar[d]_t\ar@{}[rd]|-{(\heartsuit)} & X\ar[d]^{t'}\\
  R\ar[r]^w\ar[d]_p \ar@{}[rd]|-{(\spadesuit)}& S\ar[d]^{p'}\\
  Y\ar[r]^v & Z
}\]
In other words: $(\diamondsuit)$ is
2-isomorphic in $\vseq(u,v)$ to
the composite
$(\spadesuit)\circ(\heartsuit)$. Replacing $(\diamondsuit)$ by
the 2-isomorphic $(\spadesuit)\circ(\heartsuit)$ we may assume the isomorphism
is an equality.

When we evaluate the diagram $(2)$ at
$(\diamondsuit)=(\spadesuit)\circ(\heartsuit)$ we arrive at two 2-morphisms
 $\xymatrix{
  x\ar@<0.5ex>[r] \ar@<-0.5ex>[r] & f^!v^*\mu_Z^{}, 
}$
which we wish to show equal.
By Proposition~\ref{P200.23}~(i) or (ii) the map $\rho(t,p):f^!\la t^!p^!$
is an isomorphism, hence it suffices to
 prove that the composites
$\xymatrix{
  x\ar@<0.5ex>[r] \ar@<-0.5ex>[r] & f^!v^*\mu_Z^{}\ar[r]^-{\rho(t,p)} & t^!p^!v^*\mu_Z^{}
}$
 are equal. But natural transformations respect composition and
 these composites can be rewritten; it suffices to
 prove the commutativity of $(2)$ when evaluated at $(\spadesuit)$ and at $(\heartsuit)$.
 In other words we are reduced to the two special cases below.

\medskip

\nin
{\bf Case 1.}\ \ 
\emph{The two composites in the diagram $(2)$ agree on 1-morphisms
  $(\diamondsuit)$ where the vertical maps
  are of finite type and universally quasi-proper.}

\medskip

\nin
\emph{Proof of Case 1.}\ \ We begin by choosing
a preimage for $(\diamondsuit)$ under
the functor $F:\nsq(W,X,Y,Z)\la\sq(W,X,Y,Z)$, that is
an object $a\in\nsq(W,X,Y,Z)$ with $F(a)=(\diamondsuit)$.
The
preimage we choose is the diagram
 \[\xymatrix@C-5pt@R+30pt{
    & W\ar[dl]_-{\id}\ar[rrr]^-{u}\ar[dd]^(0.3){f} & & & X\ar[dl]_-{\id}
    \ar[dd]^(0.3){g} \\
 W\ar[dr]_-{f} \ar[rrr]^(.6){u}& & &X\ar[dr]_-{g} & \\
& Y\ar[rrr]^-{v}  &&  & Z \\
}\]
With this 
choice the definitions make the computation explicit, and
diagram $(2)$ evaluated at $a\in\nsq(W,X,Y,Z)$ reduces
to diagram $(1)$ evaluated at $(\diamondsuit)$.
We already know diagram $(1)$ to be commutative.

\medskip

\nin
{\bf Case 2.}\ \ 
\emph{The two composites in the diagram $(2)$ agree on 1-morphisms
  $(\diamondsuit)$ where the vertical maps
  are dominant, flat monomorphisms.}

\medskip

\nin
\emph{Proof of Case 2.}\ \ We begin by choosing
a preimage for $(\diamondsuit)$ under
the functor $F:\nsq(W,X,Y,Z)\la\sq(W,X,Y,Z)$, that is
an object $a\in\nsq(W,X,Y,Z)$ with $F(a)=(\diamondsuit)$.
The
preimage we choose is the diagram
 
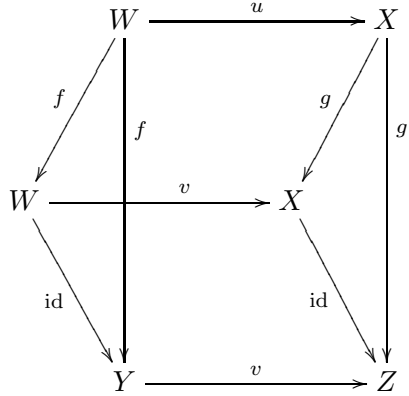
\[\xymatrix@C-5pt@R+30pt{
    & W\ar[dl]_-{f}\ar[rrr]^-{u}\ar[dd]^(0.3){f} & & & X\ar[dl]_-{g}
    \ar[dd]^(0.3){g} \\
 W\ar[dr]_-{\id} \ar[rrr]^(.6){v}& & &X\ar[dr]_-{\id} & \\
& Y\ar[rrr]^-{v}  &&  & Z \\
}\]
With this 
choice the definitions make the computation explicit, and
diagram $(2)$ evaluates at $a\in\nsq(W,X,Y,Z)$ 
to give the following composites
\[\xymatrix@C+40pt@R+2pt{
  \mu_W^{}\big[u^*g^*\times u^*g^*\big] \ar[d]_-{\mu(W,X)}
  \ar[r]^-{(\diamondsuit)^*\times (\diamondsuit)^*}
 &\mu_W^{}\big[f^*v^*\times f^*v^*\big] \ar[d]^-{\mu(W,Y)}\\
u^*\mu_X^{}(g^*\times g^*) \ar[d]_-{\mu(X,Z)} & f^*\mu_Y^{}(v^*\times v^*) \ar[d]^-{\mu(Y,Z)}  \\
u^*g^*\mu_Z^{} \ar[r]^-{(\diamondsuit)^*)}  & f^*v^* \mu_Z^{} 
}\]
and the equality is because of the premonoid structure on $(-)^*$.
\eprf

\rmd{R703.11}
We would like to restate the second part of
Proposition~\ref{P703.10.5} slightly, making it more clearly parallel to
Proposition~\ref{P204.1006}. We briefly remind the reader of the notation
introduced in Reminder~\ref{R204.923}: the category $\sq(W,X,Y,Z)$ has for
its
objects the 2-cartesian squares $(\diamondsuit)$ below
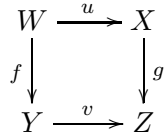
\[\xymatrix{
  W \ar[r]^u\ar[d]_f & X\ar[d]^g\\
  Y\ar[r]^v & Z
  }\]
with $v$ flat. There are four natural projections, the obvious
functors
 \[\xymatrix@R+10pt@C-10pt{
  \seq(W,X)& & & \sq(W,X,Y,Z)\ar[rd]^-{\pi_{\text{bot}}^{}}\ar[rrr]^-{\pi_{\text{left}}^{}}\ar[lll]_-{\pi_{\text{top}}^{}}\ar[ld]_-{\pi_{\text{right}}^{}} & & & \seq(W,Y)\\
 &&\seq(X,Z)&&\seq(Y,Z)&& 
  }\]
If $A$ is any of top, bot (for bottom), left or right let $\pi_A^*$ be the
composite
functor $(-)^*\circ\pi_A^{}$ and let $\pi_A^!$ be the
composite
functor $(-)^!\circ\pi_A^{}$. In
Construction~\ref{C204.11} and Proposition~\ref{P204.13} we produced the
base-change natural transformation 
$\theta:\pi_{\text{top}}^*\pi_{\text{right}}^!\la
\pi_{\text{left}}^!\pi_{\text{bot}}^*$.
With
$\tau:\pi_{\text{top}}^*\pi_{\text{right}}^*\la
\pi_{\text{left}}^*\pi_{\text{bot}}^*$
the canonical isomorphism we have
\ermd

\pro{P703.11}
Consider the diagram below, of natural transformations of functors of
the form $\sq(W,X,Y,Z)\la\cb\big[\Dqc(Z)\times\Dqc(Z)\,,\,\Dqc(W)\big]$,
\[\xymatrix@C+40pt{
\mu_W^{}\circ[\pi_{\text{\rm top}}^*\times \pi_{\text{\rm top}}^*]\circ[\pi_{\text{\rm right}}^*\times
\pi_{\text{\rm right}}^!]\ar[r]^-{\mu_W^{}\circ[\tau\times\theta]} \ar[d]_-{\mu(W,X)}&\mu_W^{}\circ[\pi_{\text{\rm left}}^*\times
\pi_{\text{\rm left}}^!]\circ[\pi_{\text{\rm bot}}^*\times
\pi_{\text{\rm bot}}^*]\ar[d]^-{\s(W,Y)}\\
\pi_{\text{\rm top}}^*\circ\mu_X^{}\circ[\pi_{\text{\rm right}}^*\times
\pi_{\text{\rm right}}^!] \ar[d]_-{\s(X,Z)}& \pi_{\text{\rm left}}^!\circ \mu_Y^{}\circ[\pi_{\text{\rm bot}}^*\times
\pi_{\text{\rm bot}}^*]\ar[d]^-{\mu(Y,Z)}\\
\pi_{\text{\rm top}}^*\circ \pi_{\text{\rm right}}^!\circ\mu_Z^{} \ar[r]^-{\theta\circ\mu_Z^{}}& \pi_{\text{\rm left}}^!\circ \pi_{\text{\rm bot}}^*\circ\mu_Z^{} 
}\]
where the top arrow is an abbreviation for $\mu_W^{}$ applied to the composite
\[\xymatrix@C+40pt{
[\pi_{\text{\rm top}}^*\times
\pi_{\text{\rm top}}^*]\circ[\pi_{\text{\rm right}}^*\times
\pi_{\text{\rm right}}^!]\ar@{=}[r] & [\pi_{\text{\rm top}}^*\circ \pi_{\text{\rm right}}^* ]\times
[\pi_{\text{\rm top}}^*\circ\pi_{\text{\rm right}}^!]\ar[d]^-{\tau\times\theta}\\
[\pi_{\text{\rm left}}^*\times
\pi_{\text{\rm left}}^!]\circ[\pi_{\text{\rm bot}}^*\times
\pi_{\text{\rm bot}}^*]\ar@{=}[r] &[\pi_{\text{\rm left}}^*\circ \pi_{\text{\rm bot}}^*]\times[\pi_{\text{\rm left}}^!\circ\pi_{\text{\rm bot}}^*]
}\]
This diagram commutes.
\epro

\prf
Let $(\diamondsuit)\in\sq(W,X,Y,Z)$ be an object, we need to show
that the diagram commutes when evaluated at $(\diamondsuit)$.
But this is just by evaluating the commutative diagram $(2)$ of
Proposition~\ref{P703.10.5}
at the 1-morphism $(\diamondsuit)\in\vseq$.
\eprf

\rmk{R703.13}
The advantage of Proposition~\ref{P703.11} over
Proposition~\ref{P703.10.5} is that the naturality in $u$ and $v$ is
explicit. One way to explain this is the following: for the purpose of
the 
proof it was convenient to work with natural transformations on the
2-category $\vseq$, it allowed us to factor $(\diamondsuit)$ as 
the composite $(\diamondsuit)=(\spadesuit)\circ (\heartsuit)$ and
reduce our computations to the special cases of $(\spadesuit)$ and
$(\heartsuit)$. The cost was that we fixed $u$ and $v$. But the arrows
in the diagram of Proposition~\ref{P703.11} are all natural, even in
$u$ and $v$, and to
prove commutativity one only needs to check the statement on
objects---and the result of Proposition~\ref{P703.10.5} allows us to evaluate the
diagram at a general object.

As in Remark~\ref{R204.295} the reader can amuse herself by working
out the $\hseq$ analog of what Proposition~\ref{P703.10.5} proved for
$\vseq$. Now that we have proved Proposition~\ref{P703.11}, which is
symmetric
and natural in all arrows, the proofs are straightforward.
\ermk

We end the paper with a result saying that, under certain hypotheses
on $f:X\la Y$, the map $\s$ is an isomorphism.

\pro{P703.15}
Let $X\stackrel f\la Y\stackrel g\la Z$ be composable 1-morphisms in
$\seq$, and assume $f$ is flat and $gf$ is of finite
Tor-dimension. Then $f^*\s(g):f^*\mu_Y^{}(g^*\times g^!)\la f^*g^!\mu_Z^{}$
is an isomorphism.

Also: the 2-morphism $\s(g): \mu_Y^{}(g^*\times g^!)\la g^!\mu_Z^{}$
yields an isomorphism whenever we evaluate it at a pair of objects
$(E,F)\in\Dqc(Z)$ with $E$ a perfect complex.
\epro

\prf
In Proposition~\ref{PThom.73.66} we saw that, for any morphism
$p:\ov Y\la Z$ in $\seq$, any perfect complex $E\in\Dqc(Z)$ and any
object
$F\in\Dqc(Z)$, the
map $\chi(p,E,F):p^*E\otimes p^\times F\la p^\times(E\oo F)$ is an
isomorphism. Now take our map $g:Y\la Z$ and choose a Nagata
compactification $Y\stackrel u\la \ov Y\stackrel p\la Z$. The map
$\s(g)$, evaluated at the pair of objects $(E,F)\in\Dqc(Z)$,
is defined to be the composite
\[
\xymatrix@C+30pt{
u^*p^*E\oo u^*p^\times F\ar[r]^-{\mu(Y,\ov Y)} & u^*(p^*E\oo p^\times
F)\ar[r]^-{\chi(p.E,F)} & u^*p^\times (E\oo F)\ .
}
\]
In this composite the map $\mu(Y,\ov Y)$ is an isomorphism, while the
second map is an isomorphism as long as $E$ is perfect. This proves
the
second part of the Proposition.

Now for the first part: we assume $gf$ of finite Tor-dimension. With the Nagata compactification $Y\stackrel
u\la \ov Y\stackrel p\la Z$
for the map $g:Y\la Z$ as above, consider the diagram
\[
\xymatrix@C+30pt{
X \ar[r]^-{uf}  & \ov Y \ar[r]^-{\id} \ar[d]_-{p} & \ov Y\ar[d]^-{p}\\
                     &   Z \ar[r] ^-{\id} & Z
}
\]
The morphism $uf$ is the composite of the flat maps $u$ and $f$ and is
therefore flat. As $gf=puf$ is of finite Tor-dimension, it follows
that $p$ is of finite Tor-dimension on the image of
$uf$. Theorem~\ref{T4.13} applies and we deduce, among other things,
that the composite $(uf)^*p^\times\id^*=f^*g^!$ respects
coproducts. Therefore the map $f^*\s(g):f^*g^*E\oo f^*g^!F\la f^*g^!(E\oo
F)$ is a natural transformation between functors of $E$ which are
triangulated and respect coproducts. The full subcategory
$\cl\subset\Dqc(Z)$, of all objects 
$E\in\Dqc(Z)$ such that the map $f^*\s(g)$ is an isomorphism for
every $F$, is a
localizing subcategory and contains all the perfect complexes. Hence
$\cl$ is localizing and contains all the
compacts---since $\Dqc(Z)$ is compactly generated we deduce $\cl=\Dqc(Z)$.
\eprf

\appendix

\section{Some counterexamples}
\label{A1}

\exm{E4.7.53}
Let us illustrate the fact that the hypothesis of Lemma~\ref{L4.7},
requiring the morphism $g$ to be 
quasi-proper, cannot easily be weakened---it doesn't suffice
for $g$ to be proper and finitely-presentable. Let $R=k[x,y]$ be the polynomial
ring in two variables, let $M$ be the $R/(y)$--module $M=k(\!(x)\!)/k[[x]]$,
and let $S$ be the Nagata extension $R\oplus M$, with $M$ a square-zero ideal.
Let $g':S\la S'$ be the surjective ring homomorphism with kernel the
ideal $yS$.
We have a pushout square of ring homomorphisms 
\[
\CD
S'[1/x]  @<u'<<  S' \\
@Af'AA        @AAg'A \\
S[1/x] @<v'<< S
\endCD
\]
Applying the functor $\spec-$ we obtain a cartesian square
of affine schemes
\[
\CD
W @>u>> X\\
@VfVV @VVgV \\
Y @>v>> Z
\endCD
\]
with $g$ a finitely-presented closed immersion, 
$f$ of finite Tor-dimension, and $v$ an open
immersion, hence flat. With the exception that $g$ is not
a quasi-proper map, the hypotheses of Lemma~\ref{L4.7} are
all satisfied. We want to show that the conclusion fails, more
precisely we will show that 
$u^*g^\times M$ and
$f^\times v^*M$ are not isomorphic. Note that $M$ belongs to
$\Dqcpl(Z)$,
hence we will have a counterexample to the 
non-pseudo-coherent versions of both Lemma~\ref{L4.7} (i) and (ii).

Tensoring with the map $v':S\la S[1/x]$ takes the $S$--module 
$M= k(\!(x)\!)/k[[x]]$ to zero---in other words $v^* M=0$. Hence
$f^\times v^*M=0$. Next we will show that 
$u^*g^\times M\neq0$; thus the natural map $u^*g^\times\la f^\times v^*$
cannot possibly be an isomorphism.

To this end recall that the ring $S'$, viewed as an $S$--module, fits in
an exact sequence 
\[
\CD
0 @>>> M @>>> S @>y>> S @>>> S' @>>> 0
\endCD
\]
where the map $y:S\la S$ is multiplication by $y$. There is a triangle in
$\Dqc(W)\cong\D(S)$ of the form $\T M\la P\la S'\la$
where $P$ is a perfect complex, more specifically it is the mapping
cone on $y:S\la S$. Now note that in $\D\big(S[1/x]\big)$ we have  the computation
\begin{eqnarray*}
f_*u^*g^\times M &=& v^*g_*g^\times M \\
               &=& v^*\HHom_S^{}\big(S',M\big)\\
               &=&S[1/x]\oo\HHom_S^{}\big(S',M\big),
\end{eqnarray*}
which allows us to compute $f_*u^*g^\times M$
by applying the functor 
$S[1/x]\oo\HHom_S^{}\big(-,M\big)$ to our triangle
$\T M\la P\la S'\la$ in $\D(S)$. As $P$ is perfect we have that
\[
S[1/x]\oo\HHom_S^{}\big(P,M\big)\quad\cong\quad
S[1/x]\oo P^\vee \oo M \eq 0,
\]
 and from the triangle we
learn that 
$f_*u^*g^\times M\cong
\T^{-2} S[1/x]\oo\HHom_S^{}(M,M)$. But then $H^0\HHom_S^{}(M,M)=k[[x]]$ is not
annihilated by the functor $S[1/x]\oo(-)$, hence $H^2(f_*u^*g^\times M)\neq0$.
\eexm

\exm{EA1.1}
Let $k$ be a field, let $R=k[x,y]$ be the polynomial ring, and let the
multiplicative group $\gm$ act by the rule
$t\cdot(x,y)=(t^{-1}x,ty)$. Let $X=\spec R$; then the map
$f:[X/\gm]\la X/\!\!/\gm$, from the stack $[X/\gm]$ to the GIT quotient
$X/\!\!/\gm$, is a quasi-proper map. In this case the invariant subring
$R^\gm$ is the polynomial ring in one variable
$k[xy]\subset k[x,y]$, that is $X/\!\!/\gm=\spec{R^\gm}$ is
just the affine line $\ak^1_k$. The open subset $U=\spec{R[1/x]}\subset\spec R=X$ is
acted on freely by $\gm$ and the quotients $U/\!\!/\gm$ and
$X/\!\!/\gm$ are the same---more precisely, we have an
open immersion $u:[U/\gm]\la [X/\gm]$, and the composite
$fu:[U/\gm]\la X/\!\!/\gm$ is an isomorphism. In this example we will show that
$u^*f^\times\neq\id$; thus, even though for the quasi-proper map $f$
we would expect $f^!\cong f^\times$ and for the open immersion $u$ we
would expect $u^!\cong u^*$, the fact that $fu=\id$ does not imply
that
$\id^!=(fu)^!=u^!f^!=u^*f^\times$. 

Of course $f$ is not separated, hence this is no contradiction to the
theorems of the article. But it does show that something delicate
happens with non-separated morphisms. Note that $\ak^1_k$ is a
finite-dimensional
regular
scheme, so any morphism to it is of finite Tor-dimension---in
particular $f$ is of finite Tor-dimension, while the open immersion $u$
is flat.
If we believed that Theorem~\ref{T0.5} should hold for non-separated
morphisms we might expect \ref{S0.5.3.7} to be true---but for our composite
$fu$ we have that $u$ is flat and $fu=\id$ is proper, meaning that
both \ref{S0.5.3.7}~(1) and (3) would predict that $\rho(u,f):(fu)^!\la u^!f^!$
should be an isomorphism.

If we grade the ring $R=k[x,y]$ by the weight of the representations
of
$\gm$ then $x$ has degree -1 and $y$ has degree $1$. An object in
$\Dqc\big([X/\gm]\big)$ is equivalent to a complex of 
graded $R$--modules, and the functor $f_*$ takes such a complex to its
subcomplex of $\gm$--invariants. The functor $f^\times$ takes a
complex $C$
of $R^\gm$--modules to the (derived) graded Hom--complex
$\Hom_{R^\gm}^{}(R,C)$. We wish to show that $u^*f^\times\neq\id$, and it suffices to produce
an
object $C\in\Dqc(\ak^1_k)$ with $u^*f^\times C=0$.  Our $C$ will be
$C=k$, 
where $k=k[xy]/(xy)$ is the skyscraper sheaf on
$ X/\!\!/\gm$, which we view as a complex concentrated in degree 0. We wish to compute
$u^*f^\times k$. By the above $f^\times k$ is given by
$\Hom_{R^\gm}^{}(R,k)$.

If we write $R=\oplus_{n=-\infty}^\infty R_n$, where $n$ is the
degree, then $R_n$ is a free module of rank 1 over the ring $R^\gm$,
generated by $y^n$ if $n\geq 0$, and $x^n$ if $n\leq0$.
Therefore $\Hom_{R^\gm}^{}(R_n,k)=k$ for every $n$. Thus $f^\times
k=\Hom_{R^\gm}^{}(R,k)$ has cohomology concentrated in degree 0, and
$H^0(f^\times k)$ is a  graded $R$--module, where in each
degree it is $k$. We need to compute
the $R$--module structure on $f^\times k =H^0(f^\times k)$, meaning we must figure
out how
multiplication by $x$ and $y$ acts on $f^\times k$. That is we have
maps $x:(f^\times k)_n\la  (f^\times k)_{n-1}$ and $y:(f^\times k)_n\la
(f^\times k)_{n+1}$ which we need to compute. 

If $n\geq0$ then $R_{-n}$ is the free $R^\gm$--module generated by
$x^n$, so multiplication by $x$ is an isomorphism $R_{-n}\la R_{-n-1}$
and hence induces an isomorphism
\[x:\Hom_{R^\gm}^{}(R_{-n-1},k)\la
\Hom_{R^\gm}^{}(R_{-n},k),\]
which rewrites as $x:(f^\times k)_{n+1}\la
(f^\times k)_{n}$. If $n>0$ multiplication by $y$ takes the generator
$x^n\in R_{-n}$ to $yx^n=[xy]x^{n-1}$, that is to an element of the
submodule
$(xy)R_{-n+1}\subset R_{-n+1}$. Therefore multiplication by $y$ acts as
the
zero map $0:\Hom_{R^\gm}^{}(R_{-n+1},k)\la
\Hom_{R^\gm}^{}(R_{-n},k)$, which rewrites as $0:(f^\times k)_{n-1}\la
(f^\times k)_{n}$. Similar computations, which we leave to the reader,
show that for $n\leq0$ multiplication by $y$ is an isomorphism
$(f^\times k)_{n-1}\la  (f^\times k)_{n}$, whereas for $n<0$
multiplication by $x$ is the zero map $0:(f^\times k)_{n+1}\la  (f^\times
k)_{n}$. This computation shows that $u^*f^\times k=R[1/x]\oo f^\times
k=0$, since multiplication by sufficiently high powers of $x$
annihilates
every cohomology class in $H^0(f^\times k)$.
\eexm

\section{Conflicts of interest}
\label{S2000}

There is no conflict of interest to report.

\def\cprime{$'$}
\providecommand{\bysame}{\leavevmode\hbox to3em{\hrulefill}\thinspace}
\providecommand{\MR}{\relax\ifhmode\unskip\space\fi MR }
\providecommand{\MRhref}[2]{%
  \href{http://www.ams.org/mathscinet-getitem?mr=#1}{#2}
}
\providecommand{\href}[2]{#2}


\begin{thebibliography}{10}

\bibitem{Abramovich-Olsson-Vistoli08}
Dan Abramovich, Martin Olsson, and Angelo Vistoli, \emph{Tame stacks in
  positive characteristic}, Ann. Inst. Fourier (Grenoble) \textbf{58} (2008),
  no.~4, 1057--1091.

\bibitem{Alonso-Jeremias-Lipman11}
Leovigildo Alonso~Tarr{\'\i}o, Ana Jerem{\'\i}as~L{\'o}pez, and Joseph Lipman,
  \emph{Bivariance, {G}rothendieck duality and {H}ochschild homology {I}:
  {C}onstruction of a bivariant theory}, Asian J. Math. \textbf{15} (2011),
  no.~3, 451--497.

\bibitem{Avramov-Iyengar08}
Luchezar~L. Avramov and Srikanth~B. Iyengar, \emph{Gorenstein algebras and
  {H}ochschild cohomology}, Michigan Math. J. \textbf{57} (2008), 17--35,
  Special volume in honor of Melvin Hochster.

\bibitem{Avramov-Iyengar-Lipman-Nayak10}
Luchezar~L. Avramov, Srikanth~B. Iyengar, Joseph Lipman, and Suresh Nayak,
  \emph{Reduction of derived {H}ochschild functors over commutative algebras
  and schemes}, Adv. Math. \textbf{223} (2010), no.~2, 735--772.

\bibitem{Balmer-Dellambrogio-Sanders16}
Paul Balmer, Ivo Dell'Ambrogio, and Beren Sanders, \emph{Grothendieck-{N}eeman
  duality and the {W}irthm\"uller isomorphism}, Compos. Math. \textbf{152}
  (2016), no.~8, 1740--1776.

\bibitem{BenZvi-Francis-Nadler10}
David Ben-Zvi, John Francis, and David Nadler, \emph{Integral transforms and
  {D}rinfeld centers in derived algebraic geometry}, J. Amer. Math. Soc.
  \textbf{23} (2010), no.~4, 909--966.

\bibitem{Bokstedt-Neeman93}
Marcel B\"okstedt and Amnon Neeman, \emph{Homotopy limits in triangulated
  categories}, Compositio Math. \textbf{86} (1993), 209--234.

\bibitem{BondalvandenBergh04}
Alexei~I. Bondal and Michel Van~den Bergh, \emph{Generators and
  representability of functors in commutative and noncommutative geometry},
  Mosc. Math. J. \textbf{3} (2003), no.~1, 1--36, 258.

\bibitem{Chikhladze-Lack-Street10}
Dimitri Chikhladze, Stephen Lack, and Ross Street, \emph{Hopf monoidal
  comonads}, Theory Appl. Categ. \textbf{24} (2010), No. 19, 554--563.

\bibitem{Conrad07}
Brian Conrad, \emph{Deligne's notes on {N}agata compactifications}, J.
  Ramanujan Math. Soc. \textbf{22} (2007), no.~3, 205--257.

\bibitem{Conrad-Lieblich-Olsson12}
Brian Conrad, Max Lieblich, and Martin Olsson, \emph{Nagata compactification
  for algebraic spaces}, J. Inst. Math. Jussieu \textbf{11} (2012), no.~4,
  747--814.

\bibitem{Day-Street97}
Brian Day and Ross Street, \emph{Monoidal bicategories and {H}opf algebroids},
  Adv. Math. \textbf{129} (1997), no.~1, 99--157.

\bibitem{Deligne66}
Pierre Deligne, \emph{Cohomology \`a support propre en construction du foncteur
  $f^!$}, Residues and Duality, Lecture Notes in Mathematics, vol.~20,
  Springer--Verlag, 1966, pp.~404--421.

\bibitem{Dubuc-Street70}
Eduardo Dubuc and Ross Street, \emph{Dinatural transformations}, Reports of the
  {M}idwest {C}ategory {S}eminar, {IV}, Lecture Notes in Mathematics, Vol. 137,
  Springer, Berlin, 1970, pp.~126--137.

\bibitem{Hall-Neeman-Rydh14}
Jack Hall, Amnon Neeman, and David Rydh, \emph{One positive and two negative
  results for derived categories of algebraid stacks}, J. Inst. Math. Jussieu
  \textbf{18} (2019), no.~5, 1087--1111.

\bibitem{Hall-Rydh15}
Jack Hall and David Rydh, \emph{Algebraic groups and compact generation of
  their derived categories of representations}, Indiana Univ. Math. J.
  \textbf{64} (2015), no.~6, 1903--1923.

\bibitem{Hall-Rydh13}
\bysame, \emph{Perfect complexes on algebraic stacks}, Compos. Math.
  \textbf{153} (2017), no.~11, 2318--2367.

\bibitem{Hartshorne66}
Robin Hartshorne, \emph{Residues and duality}, Lecture notes of a seminar on
  the work of A. Grothendieck, given at Harvard 1963/64. With an appendix by P.
  Deligne. Lecture Notes in Mathematics, No. 20, Springer-Verlag, Berlin, 1966.

\bibitem{Illusie71C}
Luc Illusie, \emph{Conditions de finitude}, Th\'eorie des intersections et
  th\'eor\`eme de {R}iemann-{R}och, Springer-Verlag, Berlin, 1971, S\'eminaire
  de G\'eom\'etrie Alg\'ebrique du Bois-Marie 1966--1967 (SGA 6, Expos{\'e}
  III), pp.~222--273. Lecture Notes in Mathematics, Vol. 225.

\bibitem{Illusie71A}
\bysame, \emph{G{\'e}n{\'e}ralit{\'e}s sur les conditions de finitude dans les
  cat{\'e}gories d{\'e}riv{\'e}es}, Th\'eorie des intersections et th\'eor\`eme
  de {R}iemann-{R}och, Springer-Verlag, Berlin, 1971, S\'eminaire de
  G\'eom\'etrie Alg\'ebrique du Bois-Marie 1966--1967 (SGA 6, Expos{\'e} I),
  pp.~78--159. Lecture Notes in Mathematics, Vol. 225.

\bibitem{Iyengar-Lipman-Neeman13}
Srikanth~B. Iyengar, Joseph Lipman, and Amnon Neeman, \emph{Relation between
  two twisted inverse image pseudofunctors in duality theory}, Compos. Math.
  \textbf{151} (2015), no.~4, 735--764.

\bibitem{Kiehl72}
Reinhardt Kiehl, \emph{Ein ``{D}escente''-{L}emma und {G}rothendiecks
  {P}rojektionssatz f\"ur nichtnoethersche {S}chemata}, Math. Ann. \textbf{198}
  (1972), 287--316.

\bibitem{Kresch09}
Andrew Kresch, \emph{On the geometry of {D}eligne-{M}umford stacks}, Algebraic
  geometry---{S}eattle 2005. {P}art 1, Proc. Sympos. Pure Math., vol.~80, Amer.
  Math. Soc., Providence, RI, 2009, pp.~259--271.

\bibitem{Lack00}
Stephen Lack, \emph{A coherent approach to pseudomonads}, Adv. Math.
  \textbf{152} (2000), no.~2, 179--202.

\bibitem{Lipman09}
Joseph Lipman, \emph{Notes on derived functors and {G}rothendieck duality},
  Foundations of {G}rothendieck duality for diagrams of schemes, Lecture Notes
  in Mathematics, vol. 1960, Springer, Berlin, 2009, pp.~1--259.

\bibitem{Lipman-Neeman07}
Joseph Lipman and Amnon Neeman, \emph{Quasi-perfect scheme maps and boundedness
  of the twisted inverse image functor}, Illinois J. Math. \textbf{51} (2007),
  209--236.

\bibitem{MacLane71}
Saunders Mac~Lane, \emph{Categories for the {Working} {Mathematician}}, Lecture
  Notes in Math., vol.~5, Springer-Verlag, New York, 1971.

\bibitem{Nagata62}
Masayoshi Nagata, \emph{Imbedding of an abstract variety in a complete
  variety}, J. Math. Kyoto Univ. \textbf{2} (1962), 1--10.

\bibitem{Nayak09}
Suresh Nayak, \emph{Compactification for essentially finite-type maps}, Adv.
  Math. \textbf{222} (2009), no.~2, 527--546.

\bibitem{Neeman92A}
Amnon Neeman, \emph{The connection between the {K--theory} localisation theorem
  of {Thomason}, {Trobaugh} and {Yao}, and the smashing subcategories of
  {Bousfield} and {Ravenel}}, Ann. Sci. \'Ecole Normale Sup\'erieure
  \textbf{25} (1992), 547--566.

\bibitem{Neeman96}
\bysame, \emph{The {Grothendieck} duality theorem via {Bousfield's} techniques
  and {Brown} representability}, Jour. Amer. Math. Soc. \textbf{9} (1996),
  205--236.

\bibitem{Neeman99}
\bysame, \emph{Triangulated {Categories}}, Annals of Mathematics Studies, vol.
  148, Princeton University Press, Princeton, NJ, 2001.

\bibitem{NeemanTIFR}
\bysame, \emph{The relation between {G}rothendieck duality and {H}ochschild
  homology}, {$K$}-{T}heory---{P}roceedings of the {I}nternational
  {C}olloquium, {M}umbai, 2016, Hindustan Book Agency, New Delhi, 2018,
  pp.~91--126.

\bibitem{Raoult71}
Jean-Claude Raoult, \emph{Compactification des espaces alg\'ebriques normaux},
  C. R. Acad. Sci. Paris S\'er. A-B \textbf{273} (1971), A766--A767.

\bibitem{Raoult74}
\bysame, \emph{Compactification des espaces alg\'ebriques}, C. R. Acad. Sci.
  Paris S\'er. A \textbf{278} (1974), 867--869.

\bibitem{Rydh09B}
David Rydh, \emph{Compactification of stacks and extending stackiness across
  the boundary},  (2009), (Unpublished; available of author's web page).

\bibitem{Rydh09}
\bysame, \emph{Compactification of tame {D}eligne-{M}umford stacks},  (2009),
  (Unpublished; available of author's web page).

\bibitem{Rydh11}
\bysame, \emph{\'{E}tale d\'evissage, descent and pushouts of stacks}, J.
  Algebra \textbf{331} (2011), 194--223.

\bibitem{Rydh15}
\bysame, \emph{Noetherian approximation of algebraic spaces and stacks}, J.
  Algebra \textbf{422} (2015), 105--147.

\bibitem{Shaul13}
Liran Shaul, \emph{Derived hochschild functors over commutative adic algebras},
   (http://arxiv.org/abs/1307.5658).

\bibitem{ThomTro}
Robert~W. Thomason and Thomas~F. Trobaugh, \emph{Higher algebraic {K--theory}
  of schemes and of derived categories}, The Grothendieck Festschrift ( a
  collection of papers to honor Grothendieck's 60'th birthday), vol.~3,
  Birkh{\"a}user, 1990, pp.~247--435.

\bibitem{Verdier68}
Jean-Louis Verdier, \emph{Base change for twisted inverse images of coherent
  sheaves}, vol. Collection: Algebraic Geometry, Tata Inst. Fund. Res., 1968,
  pp.~393--408.

\end{thebibliography}
\end{document}